\documentclass[pdflatex,sn-mathphys-num]{sn-jnl}

\usepackage{graphicx}%
\usepackage{multirow}%
\usepackage{amsmath,amssymb,amsfonts}%
\usepackage{amsthm}%
\usepackage{mathrsfs}%
\usepackage[title]{appendix}%
\usepackage{xcolor}%
\usepackage{textcomp}%
\usepackage{manyfoot}%
\usepackage{booktabs}%
\usepackage{algorithm}%
\usepackage{algorithmicx}%
\usepackage{algpseudocode}%
\usepackage{listings}%
\usepackage{nicematrix}
\usepackage{boites}
\usepackage[bf]{subfigure}
\newcommand{\mbf}[1]{\mathbf{#1}}

\newcommand{\dt}{\varDelta t}
\newcommand{\dx}{\varDelta x}
\newcommand{\mbfu}{\mbf{U}}
\newcommand{\mbfv}{\dot{\mbf{U}}}
\newcommand{\mbfa}{\ddot{\mbf{U}}}
\newcommand{\mbfM}{\mbf{M}}
\newcommand{\mbfC}{\mbf{C}}
\newcommand{\mbfK}{\mbf{K}}
\newcommand{\mbfF}{\mbf{F}}
\newcommand{\T}{^{\mathsf{T}}}

\newcommand{\infrho}{\rho_\infty}

\usepackage{cleveref}
\crefname{figure}{Fig.}{Figs.} 
\crefname{subfigure}{Fig.}{Figs.} 
\crefname{table}{Table}{Tables}
\crefname{equation}{Eq.}{Eqs.}
\usepackage{makecell}
\usepackage{arydshln}
\allowdisplaybreaks[4]

\theoremstyle{thmstyleone}%
\theoremstyle{thmstyletwo}%
\newtheorem{remark}{Remark}%

\theoremstyle{thmstylethree}%
\newtheorem{definition}{Definition}%

\begin{document}

\title[Order elevation of directly self-starting sub-step implicit integrators for transient dynamics]{\bf Order elevation of directly self-starting sub-step implicit integrators for transient dynamics}


\author[1]{\fnm{Jinze} \sur{Li}}\email{pinkie.ljz@hit.edu.cn}

\author[1]{\fnm{Yaokun} \sur{Liu}}\email{liuyk@stu.hit.edu.cn}

\author*[2]{\fnm{Kewei} \sur{Chen}}\email{chenkw@chinatec.com.cn}

\author[3]{\fnm{Hua} \sur{Li}}\email{lihua@ntu.edu.sg}

\author*[1]{\fnm{Kaiping} \sur{Yu}}\email{yukp@hit.edu.cn}


\affil[1]{\small \orgdiv{School of Astronautics}, \orgname{Harbin Institute of Technology}, \orgaddress{\street{No.~92 West Dazhi Street}, \city{Harbin}, \postcode{150001}, \country{China}}}

\affil[2]{\small \orgdiv{Zhongke Ruilong Intelligent Manufacturing (Jiangsu) Co.}, \orgname{Ltd.}, \orgaddress{\street{No.~31 Pujiang Road}, \city{Suzhou}, \postcode{215500}, \country{China}}}

\affil[3]{\small \orgdiv{School of Mechanical and Aerospace Engineering}, \orgname{Nanyang Technological University}, \orgaddress{\street{50 Nanyang Avenue}, \postcode{639798}, \country{Singapore}}}


\abstract{Directly self-starting implicit methods are attractive for transient analysis because they avoid auxiliary starting procedures while retaining the original first- or second-order governing equations. However, most existing formulations usually fix the last sub-step at the end of each time interval, which restricts the attainable order. This study develops a generalized $s$-sub-step implicit framework by releasing this constraint and treating all sub-step locations as design variables. The resulting methods admit a unified Runge--Kutta representation for both first- and second-order transient systems and preserve identical effective matrices over all sub-steps. Accuracy conditions are derived by simultaneously matching the numerical amplification factor and load operator, thereby accounting for both homogeneous and forced responses. For $s=1,~\cdots,~6$, two complementary families are obtained: $s$th-order members with user-controllable high-frequency numerical dissipation and adjustable sub-step locations, and $(s+1)$th-order members obtained by selecting the sub-step locations, with fixed dissipation. The latter reach up to seventh-order accuracy without increasing the number of sub-steps, although some high-order members are $A(\alpha)$-stable with stability angles extremely close to $90^\circ$. Analytical amplitude and phase errors further reveal parity-dependent superconvergence in undamped systems, and appropriate parameter selections can substantially increase either phase or amplitude accuracy beyond the formal order. Numerical benchmarks confirm the predicted convergence orders and the controllable suppression of spurious high-frequency responses. Dispersion analyses of an elastic bar further demonstrate that spatial and temporal errors must be considered simultaneously: optimized mass matrices and CFL numbers can significantly improve the solution accuracy for wave propagations, whereas excessive time-step refinement alone does not necessarily improve the approximation to the original PDE. These results provide a systematic framework for constructing high-order, directly self-starting implicit integrators with flexible accuracy and dissipation trade-offs.}

\keywords{directly self-starting integration, sub-step implicit methods, high-order accuracy, controllable numerical dissipation, dispersion analysis}



\maketitle

\section{Introduction}
The numerical integration of dynamic systems remains one of the fundamental topics in computational mechanics. After spatial discretizations \cite{hughes_FiniteElementMethod_2000}, a broad class of linear dynamic problems can be described by the second-order ordinary differential equation
\begin{subequations}
    \begin{equation}\label{eq:mck}
        \mbfM\mbfa(t)+\mbfC\mbfv(t)+\mbfK\mbfu(t)=\mbfF(t)
    \end{equation}
in linear hyperbolic problems or by the first-order ordinary differential equation
\begin{equation}\label{eq:ck}
\mbfC\mbfv(t)+\mbfK\mbfu(t)=\mbfF(t)
\end{equation}
\end{subequations}
in linear parabolic problems. In \cref{eq:mck} for the structural dynamical problem, $\mbfM$, $\mbfC$, and $\mbfK$ denote the mass, damping, and stiffness matrices, respectively, and $\mbfF(t)$ represents the external load vector, whereas, in \cref{eq:ck} for the heat conduction problem, $ \mbfC $ and $\mbfK$ denote the capacity and conductivity matrices, respectively, and $\mbfF(t)$ is the heat supply vector. The objective of time integration algorithms is to advance the responses ($\mbfu,~\mbfv$ and/or $\mbfa$) from the current time instant to the next while maintaining numerical stability, accuracy, and computational efficiency. Because practical engineering simulations frequently involve large-scale models, nonlinear constitutive behaviors, and long-duration transient responses \cite{li_in_plane_2025}, the development of robust and efficient time integration schemes has remained an active research topic for several decades. Existing time integration algorithms can generally be classified into explicit and implicit methods according to the treatment of equilibrium equations within each time step. Explicit methods \cite{li_IdenticalSecondorderSingle_2021,liu_ExplorationSelfstartingSinglesolve_2026,zhao_SelfstartingDissipativeAlternative_2023,malakiyeh_Explicitv1v2Bathe_2023,li_SecondorderSsubstepExplicit_2023,noh_ExplicitTimeIntegration_2013} determine the structural responses solely from previously computed quantities and therefore avoid solving algebraic equation systems, making them computationally attractive for large-scale wave propagation and highly transient problems. However, their conditional stability imposes a severe restriction on the allowable time step, which is dictated by the highest natural frequency of the discretized system. Consequently, explicit algorithms often become computationally inefficient for stiff structures or quasi-static problems \cite{li_sizeDependent_2026}. In contrast, implicit methods \cite{newmark_MethodComputationStructural_1959,li_SecondorderAccurateThree_2020,yu_NewFamilyGeneralizeda_2008,rezaiee-pajand_EfficientWeightedResidual_2021,li_EnhancedStudiesComposite_2020,shao_ThreeParametersAlgorithm_1988,li_UnifiedFamilyHighorder_2026} require the solution of linear or nonlinear algebraic systems at each time step but possess significantly improved stability. For many structural dynamic analyses, especially those involving low-frequency dominant responses or nonlinear equilibrium iterations, implicit algorithms permit the time increment to be selected according to accuracy considerations rather than stability limitations. As a result, implicit time integration methods continue to dominate practical finite element analyses.

\subsection{Second-order accurate methods}
Over the past several decades, numerous implicit algorithms have been developed by pursuing different combinations of numerical accuracy, stability, algorithmic dissipation, and computational efficiency. Among the most influential methods are the Newmark method \cite{newmark_MethodComputationStructural_1959}, the Hilber--Hughes--Taylor (HHT-$\alpha$) method \cite{hilber_ImprovedNumericalDissipation_1977,hughes_AnalysisTransientAlgorithms_1983}, the Wood--Bossak--Zienkiewicz (WBZ-$\alpha$) method \cite{wood_AlphaModificationNewmarks_1980}, and the TPO/G-$\alpha$ method \cite{shao_ThreeParametersAlgorithm_1988}. These algorithms are all single-step formulations and have become standard integration tools because of their unconditional stability, self-starting capability, and relatively low computational cost. Furthermore, the Newmark-based methods \cite{newmark_MethodComputationStructural_1959,hilber_ImprovedNumericalDissipation_1977,hughes_AnalysisTransientAlgorithms_1983,wood_AlphaModificationNewmarks_1980,shao_ThreeParametersAlgorithm_1988} provide controllable high-frequency numerical dissipation while preserving second-order accuracy, making them particularly effective for suppressing spurious high-frequency oscillations. Nevertheless, the design framework of traditional single-step algorithms inherently limits the available degrees of freedom \cite{Li_JieGouDongLiXue_2025,Li_JieGouDongXiangYing_2025}, and most existing formulations cannot exceed second-order accuracy without introducing additional solution variables \cite{leontiev_ExtensionLMSFormulations_2007,li_DesigningDevelopingSinglestep_2023}, higher-order derivatives \cite{hoff_HigherDerivativeExplicit_1990a}, or substantially increased computational complexity. Moreover, many dissipative single-step methods exhibit undesirable overshoots in displacement and/or velocity responses when strong high-frequency dissipation is imposed, thereby reducing their robustness for certain classes of nonlinear problems.

To overcome the limitations of conventional single-step formulations, the sub-step technique~\cite{li_DevelopmentCompositeSubstep_2021,bathe_CompositeImplicitTime_2005,li_SuiteSecondorderComposite_2023} have attracted increasing attention during the past two decades. Instead of advancing the numerical solution over a single interval, sub-step methods divide each time step into several consecutive sub-steps and employ different integration formulas within individual sub-step. This strategy substantially enlarges the algorithmic design space by introducing intermediate stages while retaining the direct integration framework for second-order equations of motion. A representative example is the two-sub-step method proposed by Bathe et al.~\cite{bathe_CompositeImplicitTime_2005,bathe_ConservingEnergyMomentum_2007}, where the trapezoidal rule is adopted in the first sub-step and the three-point backward difference formula is employed in the second sub-step. The resulting algorithm possesses second-order accuracy, unconditional stability, excellent high-frequency damping, and non-oscillatory responses, demonstrating remarkable robustness for both linear and nonlinear structural dynamics. Subsequent investigations further revealed that the Bathe algorithm is closely related to three-stage diagonally implicit Runge--Kutta (RK) formulations \cite{hosea_AnalysisImplementationTRBDF2_1996} with explicit first stage, providing an important theoretical connection between sub-step methods and implicit RK schemes.

Motivated by the success of the Bathe method \cite{bathe_CompositeImplicitTime_2005,bathe_ConservingEnergyMomentum_2007}, extensive efforts have been devoted to developing generalized sub-step algorithms with improved numerical properties. Various two-~\cite{li_NovelFamilyComposite_2020,kim_ImprovedTimeIntegration_2016} and multi-sub-step \cite{li_NovelFamilyControllably_2019,rezaiee-pajand_MixedMultistepHigherorder_2010} formulations have subsequently been proposed by modifying the time-splitting strategy, introducing controllable dissipation, optimizing spectral properties, or enforcing identical effective matrices throughout all sub-steps. Representative developments include the $\infrho$-Bathe algorithm \cite{noh_BatheTimeIntegration_2019,choi_TimeSplittingRatio_2022}, the $\beta_1/\beta_2$-Bathe algorithm \cite{malakiyeh_BatheTimeIntegration_2019,malakiyeh_NewInsightsv1_2021}, optimal two-sub-step schemes \cite{li_NovelFamilyComposite_2020}, and several generalized multi-sub-step \cite{zhang_OptimizationNsubstepComposite_2020,li_HighorderAccurateMultisubstep_2024} implicit methods capable of preserving identical effective matrices while achieving controllable high-frequency dissipation. These studies have demonstrated that identical effective matrices not only reduce computational cost by allowing repeated matrix factorizations to be avoided but also provide desirable spectral properties, including minimal period elongation and effective suppression of spurious high-frequency responses. Consequently, the sub-step framework has gradually evolved into one of the most promising approaches for constructing robust algorithms with superior stability and dissipation abilities.

Despite these significant advances, almost all existing sub-step algorithms remain fundamentally second-order accurate. Their algorithmic flexibility has primarily been exploited to improve spectral behavior, numerical dissipation, or robustness, whereas relatively little attention has been devoted to systematically increasing the temporal accuracy without sacrificing the attractive computational properties inherited from second-order dynamics. As modern computational mechanics increasingly demands highly accurate long-term transient simulations involving multiscale structural responses, the development of efficient high-order implicit integration schemes has become an important research direction.

\subsection{High-order accurate methods}
To improve temporal accuracy beyond the second order, considerable research has been devoted to the development of high-order implicit integration algorithms. Existing approaches can generally be classified into two categories according to whether the original second-order dynamics are first transformed into an equivalent first-order system. The first category constructs high-order schemes based on collocation methods \cite{kim_EffectiveHigherorderTime_2017}, or higher-order interpolation techniques after rewriting the governing equations as first-order systems. Representative developments include weighted residual methods \cite{golley_WeightedResidualDevelopment_1998}, time finite element formulations \cite{wang_ImprovedTimeIntegration_2021} and Gauss--Lobatto collocation schemes \cite{liu_TwoHighorderEnergypreserving_2022}. These formulations can achieve arbitrarily high orders of accuracy while maintaining excellent stability. However, the resulting algorithms usually require solving coupled algebraic systems with substantially increased dimensions, leading to significantly higher computational costs. In addition, the associated effective matrices generally lose the sparsity and symmetry inherited from the original second-order dynamics, making these methods less attractive for large-scale simulations.

Recently, Song et al. \cite{song_HighorderImplicitTime_2022} proposed a family of high-order implicit algorithms based on Padé approximations, where rational approximations were employed to obtain highly accurate non-dissipative integration schemes directly for second-order dynamics. Subsequently, mixed-order Padé expansions were introduced to incorporate controllable numerical dissipation while preserving high-order accuracy \cite{song_HighorderImplicitTime_2023}. More recently, rational approximation techniques were further generalized to construct sub-step implicit algorithms with enhanced spectral properties and improved computational efficiency \cite{song_HighorderCompositeImplicit_2024}. These developments demonstrate that rational approximation provides an effective mathematical framework for designing high-order implicit methods with desirable stability and accuracy. Nevertheless, the resulting formulations generally involve relatively complicated coefficient constructions and load approximations, and their algorithmic structures remain closely associated with specific rational approximation frameworks.

Another important research direction concerns the development of high-order singly diagonally implicit Runge--Kutta (SDIRK) methods with improved stability, computational efficiency, and controllable numerical dissipation. Compared with fully implicit RK methods, SDIRK formulations require solving only one stage equation at a time and therefore exhibit considerably lower computational complexity while retaining excellent stability. Consequently, SDIRK methods have received renewed attention in recent years. It is worth noting that many multi-sub-step implicit algorithms developed in computational mechanics are, in essence, SDIRK-type formulations directly operating on the original second-order dynamics. Instead of transforming the governing equations into first-order systems, these methods preserve the original equilibrium equations throughout the integration process, thereby retaining the sparsity, symmetry, and positive definiteness of the effective matrices. Consequently, they achieve significantly lower computational costs than conventional RK methods while maintaining desirable numerical stability and accuracy. Within this framework, Zhang and co-workers \cite{zhang_OptimizationNsubstepComposite_2020} developed a family of $n$-sub-step algorithms capable of achieving controllable high-frequency dissipation with acceptable computational cost. However, subsequent studies \cite{li_DirectlySelfstartingHigherorder_2022,li_HighorderAccurateMultisubstep_2024} revealed that these methods may experience order reduction for forced vibrations owing to the approximation of external loads. To overcome this limitation, Li et al.~\cite{li_DirectlySelfstartingHigherorder_2022} proposed a family of directly self-starting implicit algorithms in which the number of sub-steps equals the attainable order of accuracy. These methods successfully achieve unconditional stability, controllable numerical dissipation, identical effective matrices, and identical high-order accuracy for displacement, velocity, and acceleration without requiring auxiliary starting procedures, thereby establishing a unified directly self-starting framework for constructing high-order implicit integration algorithms.

Building upon this framework, recent developments have continued to improve the attainable temporal accuracy and enlarge the algorithm design space. Lee and Noh~\cite{lee_ImplicitSsubstepTime_2025} proposed a family of implicit $s$-sub-step schemes capable of achieving $(s+1)$-th-order accuracy while preserving controllable numerical dissipation, further demonstrating the effectiveness of multi-stage formulations for second-order dynamics. Meanwhile, Souza et al.~\cite{souza_NovelNewmarkFamily_2026} introduced a fourth-order Newmark family employing complex-valued sub-step sizes, showing that relaxing the restrictions imposed on stage locations can further enhance the approximation capability of integration schemes. These studies indicate that increasing the flexibility of stage distributions has become an effective strategy for improving temporal accuracy without substantially increasing computational complexity.

\subsection{Motivation and contributions}
Despite the recent progress in high-order sub-step integration methods, the role of sub-step locations in determining the attainable accuracy has not yet been fully exploited. In most directly self-starting formulations, the last sub-step is constrained to coincide with the end of the current time interval. Although this choice simplifies the algorithm construction, it also removes one potentially useful parameter from the sub-step locations. From the viewpoint of order conditions, such a prescribed stage location may restrict the attainable temporal accuracy within the given $s$-sub-step method.

The present study revisits directly self-starting sub-step implicit integration methods \cite{li_DirectlySelfstartingHigherorder_2022} by releasing the conventional constraint $\gamma_s=1$ and treating the sub-step locations $\gamma_j~(j=1,~\cdots,~s)$ as design parameters. A generalized $s$-sub-step formulation is first established for both first- and second-order transient problems and is uniformly represented in the RK form without transforming the original second-order dynamics into an equivalent first-order system. The proposed construction retains the directly self-starting property and imposes identical effective matrices throughout all sub-steps, thereby preserving the computational advantages of existing sub-step implicit methods.

Within the proposed framework, the algorithmic parameters are determined by jointly considering the consistency and accuracy conditions, the identical effective matrices, stability, and high-frequency numerical dissipation. In particular, the additional freedoms introduced by the sub-step locations lead to two complementary algorithm configurations for a given number of sub-steps. The first configuration achieves $s$th-order accuracy while retaining a user-specified parameter $\infrho$ and adjustable sub-step locations. In this case, the first sub-step size $\gamma_1$ is associated with $\infrho$, whereas the remaining sub-step sizes can be selected without changing any prescribed spectral properties. Alternatively, by imposing one additional order condition, all available sub-step locations can be exploited to construct an $(s+1)$th-order member without increasing the number of sub-steps. The resulting higher-order method, however, possesses a fixed high-frequency dissipation level because $\gamma_1$ is no longer independently adjustable. Therefore, the proposed framework provides a direct trade-off between controllable numerical dissipation and one-order enhancement of temporal accuracy.

The one-, two-, and three-sub-step cases are developed in detail to demonstrate the parameter construction and the underlying relationship among sub-step locations, temporal accuracy, stability, and numerical dissipation. The same procedure is subsequently extended to the four-, five-, and six-sub-step methods, with the lengthy parameter expressions collected in \cref{app:a}. In addition to the formal accuracy, the amplitude and phase errors are analyzed to reveal the spectral accuracy of the proposed methods. The resulting error structures exhibit a clear dependence on the parity of the number of sub-steps, and appropriate selections of $\gamma_1$ can further eliminate dominant amplitude or phase errors in undamped systems. The stability of the higher-order members are also examined carefully; while many members retain unconditional or $A$-stability, some high-order constructions are only $A(\alpha)$-stable with stability angles very close to $90^\circ$.

The remainder of this paper is organized as follows. Section~\ref{sec:dsucin} establishes the generalized directly self-starting $s$-sub-step implicit framework and presents the theoretical conditions associated with identical effective matrices, accuracy, stability, and controllable high-frequency dissipation. Section~\ref{sec:development} develops and analyzes the proposed algorithms for $s=1,~\cdots,~6$, including their higher-order constructions, spectral properties, amplitude and phase errors, and the recovery of additional quantities. Section~\ref{sec:comparisons} summarizes and compares numerical characteristics of the $s$th- and $(s+1)$th-order members, with particular attention to their accuracy, stability, numerical dissipation, and superconvergence properties. Section~\ref{sec:examples} verifies the theoretical results through benchmark problems and further investigates the interaction between spatial and temporal discretization errors through dispersion analyses using linear and quadratic finite elements. Finally, the major findings are summarized in Section~\ref{sec:conclusions}.

\section{The directly self-starting sub-step implicit framework}\label{sec:dsucin}
For solving second-order transient dynamics described by Eq.~\eqref{eq:mck}, a generalized $s$-sub-step implicit integration framework is developed herein. For each time interval $[t_n,~t_{n+1}]$, $s$ sub-step locations are introduced as $t_{n+\gamma_j}=t_n+\gamma_j\dt$, where $\gamma_j~(j=1,~\cdots,~s)$ denote the corresponding stage abscissae and $\dt$ is the time step size. Unlike a conventional partition, these stage locations are not required to be monotonically ordered or restricted to $[0,~1]$ and they are assumed to be distinct from each other. An individual integration scheme is formulated within each sub-step. Without loss of generality, the numerical formulation over the $j$th sub-step is expressed as
\begin{subequations}\label{eq:alg}
\begin{align}
\mbfM\mbfa_{n+\gamma_j}&+\mbfC\mbfv_{n+\gamma_j}+\mbfK\mbfu_{n+\gamma_j}=\mbfF_{n+\gamma_j}\\
\mbfv_{n+\gamma_j}&=\mbfv_n+\dt\sum_{i=1}^{j}\alpha_{ji}\mbfa_{n+\gamma_i}\\
\mbfu_{n+\gamma_j}&=\mbfu_n+\dt\sum_{i=1}^{j}\alpha_{ji}\mbfv_{n+\gamma_i}
\end{align}
where $ \alpha_{ji}~(j=1,~\cdots,~s;~i=1,~\cdots,~j) $ are real-valued algorithmic parameters to be determined. It should be noted that, when solving the $j$th sub-step, the velocity and acceleration vectors obtained from the preceding $(j-1)$ sub-steps are already available and can therefore be treated as known quantities. Consequently, only three unknown vectors, namely the displacement $\mbfu_{n+\gamma_j}$, velocity $\mbfv_{n+\gamma_j}$, and acceleration $\mbfa_{n+\gamma_j}$, are involved in Eqs.~(\ref{eq:alg}a-c). By successively applying Eqs.~(\ref{eq:alg}a-c) for $j=1,~\cdots,~s$, the numerical formulations over $s$ sub-steps are completely established. After completing the calculations at all sub-steps, the displacement and velocity vectors at the end of the current time interval $t_{n+1}$ are obtained by utilizing the previously computed velocities $\mbfv_{n+\gamma_j}~(j=1,~\cdots,~s)$ and accelerations $\mbfa_{n+\gamma_j}~(j=1,~\cdots,~s)$ as
\begin{align}
\mbfv_{n+1}&=\mbfv_n+\dt\sum_{i=1}^{s}\beta_i\mbfa_{n+\gamma_i}\\
\mbfu_{n+1}&=\mbfu_n+\dt\sum_{i=1}^{s}\beta_i\mbfv_{n+\gamma_i}.
\end{align}
\end{subequations}

\begin{remark}
    The proposed $s$-sub-step implicit method is defined by \cref{eq:alg}. It should be noted that the acceleration vector $\mbfa_n$ does not appear in \cref{eq:alg}. Therefore, for given initial displacement $\mbfu_0$ and velocity $\mbfv_0$, the proposed method \eqref{eq:alg} does not require the evaluation of the initial acceleration vector $\mbfa_0$ at the first time step $t\in[t_0,~t_1:=t_0+\dt]$. Consequently, the proposed method \eqref{eq:alg} possesses the directly self-starting property. The most significant advantage of directly self-starting algorithms is that they eliminate the additional solution of an equation in the form of $\mbfM\mbfa_0=\mbf{r}$. In particular, when the mass matrix is non-diagonal, directly self-starting algorithms not only avoid solving an additional linear algebraic system but also eliminate the assembly and factorization of the global mass matrix.
\end{remark}

For solving first-order transient dynamics described by \cref{eq:ck}, the proposed sub-step implicit method \eqref{eq:alg} can be equivalently reformulated as $(j=1,~\cdots,~s)$
\begin{subequations}\label{eq:alg1}
\begin{align}
\mbfC\mbfv_{n+\gamma_j}&+\mbfK\mbfu_{n+\gamma_j}=\mbfF_{n+\gamma_j}\\
\mbfu_{n+\gamma_j}&=\mbfu_n+\dt\sum_{i=1}^{j}\alpha_{ji}\mbfv_{n+\gamma_i}
\end{align}
with the update of $\mbfu_{n+1}$ as 
\begin{equation}
\mbfu_{n+1}=\mbfu_n+\dt\sum_{i=1}^{s}\beta_i\mbfv_{n+\gamma_i}.
\end{equation}
\end{subequations}

\begin{remark}
    For first-order transient dynamics, the proposed $s$-sub-step implicit method is described by \cref{eq:alg1}. Since the derivative $\mbfv_n$ does not appear in \cref{eq:alg1}, the proposed method does not require the evaluation of the initial derivative $\mbfv_0$ at the first time step $t\in[t_0,~t_1:=t_0+\dt]$ for a given initial value $\mbfu_0$. Similar to the case of second-order transient dynamics, the proposed method also possesses the directly self-starting property for first-order transient problems.
\end{remark}

For first-order transient dynamics, the proposed method \eqref{eq:alg1} is essentially an implicit family of RK methods and can thus be represented by the Butcher tableau \cite{butcher_NumericalMethodsOrdinary_2016} as
\begin{equation}\label{eq:but_alg}
    \begin{NiceArray}{c|c}[cell-space-limits=3pt,columns-width=0.5cm]
    \mbf{c} & \mbf{A}\\ \hline & \mbf{b}\T
    \end{NiceArray}=\begin{NiceArray}{c|ccccc}[cell-space-limits=3pt,columns-width=0.7cm]
    \gamma_1 & \alpha_{11} & \\
    \gamma_2 & \alpha_{21} & \alpha_{22} & \\ 
    \gamma_3 & \alpha_{31} & \alpha_{32} & \alpha_{33} \\
    \vdots & \vdots & \vdots & \vdots & \ddots \\
    \gamma_s & \alpha_{s1} & \alpha_{s2} & \alpha_{s3} & \cdots & \alpha_{ss} \\ \hline 
    & \beta_1 & \beta_2 & \beta_3 & \cdots & \beta_s 
    \end{NiceArray}.
\end{equation}
it can be concluded that the proposed sub-step implicit method for both first- and second-order transient dynamics can be uniformly represented by the Butcher tableau given in \cref{eq:but_alg}. In other words, the high-order implicit integrators developed herein for first-order transient problems can be directly extended to second-order transient problems without reformulating the governing equations into an equivalent first-order system. Most previously developed sub-step implicit algorithms \cite{bathe_CompositeImplicitTime_2005,bathe_ConservingEnergyMomentum_2007,li_HighorderAccurateMultisubstep_2024,li_DirectlySelfstartingHigherorder_2022,zhang_OptimizationNsubstepComposite_2020} share the same feature and can therefore be interpreted as particular classes of RK methods.

In the following sections, the RK theory \cite{butcher_NumericalMethodsOrdinary_2016}, together with the principles of computational structural dynamics, will be employed to systematically design, analyze, and develop a family of high-order implicit integration algorithms.

\subsection{Identical effective matrices}

For the proposed $s$-sub-step implicit method, each sub-step requires the solution of an algebraic system involving an effective matrix. If different effective matrices are generated at different sub-steps, repeated matrix assembly and factorization may considerably increase the computational cost, particularly for large-scale linear systems. Therefore, it is desirable to enforce identical effective matrices throughout all sub-steps.

For first-order transient problems described by \cref{eq:ck}, substituting Eqs.~(\ref{eq:alg1}b) into (\ref{eq:alg1}a) at the $j$th sub-step gives
\begin{equation}
\left( \mbfC+\alpha_{jj}\dt\mbfK \right)\mbfv_{n+\gamma_j}=\mbfF_{n+\gamma_j}-\mbfK\left(\mbfu_n+\dt\sum_{i=1}^{j-1}\alpha_{ji}\mbfv_{n+\gamma_i}\right).
\end{equation}
Accordingly, the effective matrix associated with the $j$th sub-step is
$\widetilde{\mbfC}_{j}:=\mbfC+\alpha_{jj}\dt\mbfK$. As the index $j$ varies from $1$ to $s$, different diagonal coefficients $\alpha_{jj}$ generally lead to different effective matrices and therefore require repeated matrix factorizations.

A similar conclusion can be obtained for second-order transient problems described by \cref{eq:mck}. At the $j$th sub-step, substituting Eqs.~(\ref{eq:alg}b-c) into (\ref{eq:alg}a) yields 
\begin{equation}
\left(\mbfM+\alpha_{jj}\dt\mbfC+\alpha_{jj}^{2}\dt^{2}\mbfK\right)\mbfa_{n+\gamma_j}=\mbfF_{n+\gamma_j}-\left(\mbfC+\alpha_{jj}\dt\mbfK\right)\left(\mbfv_n+\dt\sum_{i=1}^{j-1}\alpha_{ji}\mbfa_{n+\gamma_i}\right)-\mbfK\left(\mbfu_n+\dt\sum_{i=1}^{j-1}\alpha_{ji}\mbfv_{n+\gamma_i}\right). 
\end{equation}
Therefore, the effective matrix for second-order transient problems is $\widetilde{\mbfM}_{j}:=\mbfM+\alpha_{jj}\dt\mbfC+\alpha_{jj}^{2}\dt^{2}\mbfK$. 

It follows from the above expressions that, for both first- and second-order transient problems, identical effective matrices throughout all sub-steps can be achieved by imposing
\begin{equation}\label{eq:iem}
\alpha_{11}=\alpha_{22}=\cdots=\alpha_{ss}.
\end{equation}
With the constraint given by \cref{eq:iem}, the proposed implicit method represented by \cref{eq:but_alg} is mathematically reduced to a SDIRK formulation.

It should be emphasized that the identical effective matrix provides an important computational advantage for both first- and second-order transient problems. For linear dynamics with a constant time step, only one effective matrix needs to be assembled and factorized, and the resulting factorization can be repeatedly reused over all sub-steps within the same time step. This feature becomes particularly attractive for large-scale models, where matrix assembly and factorization usually constitute a major portion of the computational cost. In addition to its computational advantage, the identical effective matrix has also been shown to facilitate desirable spectral properties of implicit algorithms \cite{bathe_InsightImplicitTime_2012,li_NovelFamilyComposite_2020,li_NovelImplicitIntegration_2025}. Therefore, the proposed $s$-sub-step implicit methods are developed herein by primarily enforcing and exploiting the identical effective matrix given in \cref{eq:iem}.

\subsection{Accuracy}

As can be observed, \cref{eq:ck} is generally a coupled system, making the direct analysis of numerical properties of an algorithm applied to \cref{eq:ck} rather complicated. It has been rigorously demonstrated \cite{hughes_FiniteElementMethod_2000,li_NovelImplicitIntegration_2025} in linear analyses that the numerical integration of the uncoupled equations is equivalent to that of the original coupled system. Therefore, the following theoretical investigation mainly focuses on the uncoupled single-degree-of-freedom (SDOF) system described as
\begin{subequations}\label{eq:sdof}
\begin{equation}
\dot{u}(t)+\omega u(t)=f(t)
\end{equation}
with the initial condition
\begin{equation}
u(t_0:=0)=u_n. 
\end{equation}
\end{subequations}
It should be noted that the assumption of the initial instant $t_0=0$ does not affect the order conditions; it is introduced solely to simplify and facilitate the accuracy analysis.

Furthermore, an exponential load function $f(t)=\exp(t)$ is considered in \cref{eq:sdof}. This choice is particularly suitable for accuracy analysis \cite{li_HighorderAccurateMultisubstep_2024} because its Taylor expansion contains non-zero terms of all orders, thereby facilitating the systematic derivation of the algorithmic order conditions. The exact solution of \cref{eq:sdof} with $f(t)=\exp(t)$ is given by
\begin{equation}\label{eq:exact_sdof}
u(t)=\exp(-\omega t)u_n+\dfrac{\exp(t)-\exp(-\omega t)}{\omega+1}
\end{equation}
and the exact solution $u(\dt)$ after one integration step can be obtained by substituting $t=\dt$ into \cref{eq:exact_sdof}, yielding
\begin{equation}\label{eq:exact_one_step}
u(\dt)=\exp(-\omega \dt)u_n+\dfrac{\exp(\dt)-\exp(-\omega \dt)}{\omega+1}=D_\mathsf{exa}u_n+L_\mathsf{exa}
\end{equation} 
where the exact amplification factor $D_\mathsf{exa}$ and the exact load operator $L_\mathsf{exa}$ are defined as
\begin{subequations}\label{eq:exact_DL}
\begin{align}
D_\mathsf{exa} &=\exp(-\omega \dt)\\
L_\mathsf{exa}&=\dfrac{\exp(\dt)-\exp(-\omega \dt)}{\omega+1}.
\end{align}
\end{subequations}

On the other hand, applying the proposed sub-step implicit method \eqref{eq:alg1} to the SDOF system \eqref{eq:sdof} with $f(t)=\exp(t)$ leads to the following one-step recursion relation:
\begin{equation}\label{eq:one_step}
u_{\dt}=D_\mathsf{num}u_n+L_\mathsf{num}
\end{equation}
where $D_\mathsf{num}$ and $L_\mathsf{num}$ denote the numerical amplification factor and numerical load operator, respectively. Their explicit expressions are derived below.

For the SDOF system \eqref{eq:sdof} with $f(t)=\exp(t)$, the integration formulation over the $j$th sub-step described by Eqs.~(\ref{eq:alg1}a-b) can be rewritten as
\begin{subequations}\label{eq:alg1_sdof}
\begin{align}
\dot{u}_{n+\gamma_j}+\omega u_{n+\gamma_j}&=f(t_{n+\gamma_j})\\
u_{n+\gamma_j}&=u_n+\dt\sum_{i=1}^{j}\alpha_{ji}\dot{u}_{n+\gamma_i}. 
\end{align}
\end{subequations}
Substituting Eqs.~(\ref{eq:alg1_sdof}b) into (\ref{eq:alg1_sdof}a) gives
\begin{equation}\label{eq:ljz}
\dot{u}_{n+\gamma_j}+\omega \left( u_n+\dt\sum_{i=1}^{j}\alpha_{ji}\dot{u}_{n+\gamma_i} \right) =f(t_{n+\gamma_j}). 
\end{equation}
With the following definitions $\dot{\mbf{u}}=\begin{bmatrix}
\dot{u}_{n+\gamma_1} & \dot{u}_{n+\gamma_2} & \cdots & \dot{u}_{n+\gamma_s}\end{bmatrix}\T $ and $\mbf{f}=\begin{bmatrix}
f(t_{n+\gamma_1}) & f(t_{n+\gamma_2}) & \cdots & f(t_n+\gamma_s)
\end{bmatrix}\T$, Eq.~(\ref{eq:ljz}) with the index $j =1,~\cdots,~s$ can be expressed in the following matrix-vector form as
\begin{equation}\label{eq:dot_u}
\dot{\mbf{u}}+\omega \left( u_n \mbf{1}+\dt \mbf{A}\dot{\mbf{u}} \right) =\mbf{f}\quad\implies\quad \dot{\mbf{u}}= \left( \mbf{I}+\omega\dt\mbf{A} \right) ^{-1} \left( \mbf{f}-\omega u_n\mbf{1} \right) 
\end{equation} 
where $\mbf{I}$ denotes an $s\times s$ identity matrix and $\mbf{1}$ represents a column vector with all entries equal to one. Therefore, the solution at $t=\dt$ can be obtained as
\begin{equation}\label{eq:one_step_2}
    \begin{aligned}
        u_{\dt}&=u_n+\dt\sum_{i=1}^{s}\beta_i \dot{u}_{n+\gamma_i}=u_n+\dt \mbf{b}\T\dot{\mbf{u}}=u_n+\dt\mbf{b}\T \left( \mbf{I}+\omega\dt\mbf{A} \right) ^{-1} \left( \mbf{f}-\omega u_n\mbf{1} \right) \\
        &= \left[ 1-\omega\dt \mbf{b}\T\left( \mbf{I}+\omega\dt\mbf{A} \right) ^{-1} \mbf{1}\right] u_n+\dt\mbf{b}\T \left( \mbf{I}+\omega\dt\mbf{A} \right) ^{-1}\mbf{f}. 
    \end{aligned}
\end{equation}
Note that Eq.~\eqref{eq:dot_u} has been employed above to eliminate $\dot{\mbf{u}}$. By comparing \cref{eq:one_step,eq:one_step_2}, the numerical amplification factor $D_\mathsf{num}$ and numerical load operator $L_\mathsf{num}$ are obtained as
\begin{subequations}\label{eq:num_DL}
\begin{align}
D_\mathsf{num}&=1-\omega\dt \mbf{b}\T\left( \mbf{I}+\omega\dt\mbf{A} \right) ^{-1} \mbf{1}\label{eq:D_num}\\
L_\mathsf{num}&=\dt\mbf{b}\T \left( \mbf{I}+\omega\dt\mbf{A} \right) ^{-1}\mbf{f}. 
\end{align}
\end{subequations}

Assuming that the solution at the initial instant $t_0$ is given as $u_n$ without numerical errors, the proposed method \eqref{eq:alg1} yields the numerical solution $u_{\dt}$ after one time step, as expressed in Eqs.~\eqref{eq:one_step} or \eqref{eq:one_step_2}. Meanwhile, the corresponding exact solution at the same instant is provided by \cref{eq:exact_one_step}. Based on these two solutions, the local truncation error \cite{li_DirectlySelfstartingHigherorder_2022} of the proposed method \eqref{eq:alg1} after one integration step is defined as follows:

\begin{definition}
The proposed method \eqref{eq:alg1} achieves $p$th-order accuracy for solving the standard modal problem \eqref{eq:sdof} if and only if $u_{\dt}-u(\dt)=\mathcal{O}(\dt^{p+1})$ is satisfied, namely 
\begin{subequations}\label{eq:def_acc}
\begin{align}
D_\mathsf{num}-D_\mathsf{exa}=\mathcal{O}(\dt^{p+1})\\
L_\mathsf{num}-L_\mathsf{exa}=\mathcal{O}(\dt^{p+1})
\end{align}
\end{subequations}
where $D_\mathsf{exa},~L_\mathsf{exa}$ and $D_\mathsf{num},~L_\mathsf{num}$ are given by \cref{eq:num_DL,eq:exact_DL}, respectively. 
\end{definition}

\begin{remark}
    If only the local truncation error associated with the homogeneous solution is considered, the resulting accuracy analysis is reduced to the conventional RK analysis, where only Eq.~(\ref{eq:def_acc}a) is taken into account. In contrast, the accuracy analysis presented in this paper considers not only the homogeneous solution but also the non-homogeneous contribution. Consequently, the present accuracy analysis explicitly accounts for the non-homogeneous contribution and avoids the order reduction observed in some existing high-order methods \cite{zhang_OptimizationNsubstepComposite_2020}. This treatment can be regarded as a complementary extension to the classical RK analysis and also represents a minor but meaningful contribution of the present work.
\end{remark}

\begin{remark}
    In addition to the accuracy conditions defined in \cref{eq:def_acc}, the sub-step implicit method should also satisfy the consistency inherited from the RK framework. Specifically, the sub-step locations $\gamma_j$ and the algorithmic coefficients $\alpha_{ji}$ should satisfy
    \begin{equation}\label{eq:consistency}
        \sum_{i=1}^{j}\alpha_{ji}=\gamma_j,\qquad j=1,~2,~\cdots,~s.
    \end{equation}
    This condition ensures that each internal sub-step is consistent with its corresponding temporal location $t_{n+\gamma_j}$. Although the consistency condition is implicitly included in the classical RK theory, it is explicitly emphasized herein because the sub-step locations $\gamma_j$ are treated as independent design variables in the proposed generalized framework. Therefore, \cref{eq:consistency} serves as a fundamental constraint for determining the sub-step distributions and algorithmic parameters in the subsequent development of high-order implicit integration algorithms.
\end{remark}

\begin{remark}
By rewriting the SDOF system in Eq.~(\ref{eq:sdof}a) as the standard test equation in numerical analysis $\dot{u}(t)=\lambda u(t)$ with $\lambda=-\omega$, the numerical amplification factor $D_\mathsf{num}$ can be directly interpreted as the RK stability function $R(z)$. Therefore, the relationship between them can be established as
\begin{equation}\label{eq:R}
R(z:=-\omega\dt)=D_\mathsf{num}=1+z \mbf{b}\T\left( \mbf{I}-z\mbf{A} \right) ^{-1} \mbf{1}.
\end{equation}
\end{remark}

\subsection{Stability}

With the RK stability function $R(z)$ established above, the stability of the proposed method can be directly analyzed within the RK stability framework. Since the corresponding stability criteria have been rigorously established in the literature \cite{butcher_NumericalMethodsOrdinary_2016}, only the essential results are briefly summarized herein without repeating the mathematical proofs.

For a general RK method described by the Butcher tableau $(\mbf{A},~\mbf{b},~\mbf{c})$, the stability function can be alternatively expressed as the ratio of two polynomials:
\begin{equation}
R(z)=\frac{N(z)}{D(z)},
\end{equation}
where $N(z)$ and $D(z)$ denote the numerator and denominator polynomials of the stability function, respectively. The unconditionally stable condition of the RK method is determined by both the pole distribution of $R(z)$ and the corresponding stability boundary. Specifically, a RK method is unconditionally or $A$-stable if and only if the following two conditions are satisfied:
\begin{itemize}
    \item all poles of $R(z)$, namely all zeros of $D(z)$, are located in the right half-plane;
    \item the polynomial
    \begin{subequations}\label{eq:Ey}
    \begin{equation}
    E(y)=D(iy)D(-iy)-N(iy)N(-iy),
    \end{equation}
    where $i$ denote the imaginary unit, satisfies
    \begin{equation}
    E(y)\geq0,\quad \forall y\in\mathbb{R}.
    \end{equation}
    \end{subequations}
\end{itemize}

For the proposed method described in \cref{eq:but_alg}, the stability function is obtained from the numerical amplification factor by \cref{eq:R} with $z=\lambda\dt=-\omega\dt$. Therefore, the above RK stability criteria can be directly employed to determine the stability of the proposed algorithms. In particular, by analyzing the poles of $R(z)$ and evaluating the sign of the corresponding $E(y)$ polynomial, the stability of the proposed implicit methods can be rigorously verified.


\subsection{Dissipation control}

In addition to accuracy and stability, an appropriate control of numerical dissipation is another essential requirement for time integration algorithms in computational dynamics. After spatial discretization \cite{hughes_FiniteElementMethod_2000}, partial differential equations are transformed into a system of ordinary differential equations, during which artificial high-frequency components may be introduced due to the discretization process. These spurious high-frequency modes not only deteriorate the accuracy of the numerical solution but may also compromise the stability of the time integration procedure, particularly in long-term simulations. Therefore, effective time integration algorithms are generally expected to incorporate controllable high-frequency dissipation, which can selectively damp out these non-physical oscillations while maintaining the accuracy of physically meaningful low-frequency responses.

With the numerical amplification factor $D_\mathsf{num}$ given by \cref{eq:D_num}, the proposed implicit method can be designed to control high-frequency numerical dissipation by imposing
\begin{equation}\label{eq:D_num_infty}
D_\mathsf{num}^\infty=\lim_{\omega\to\infty} D_\mathsf{num}=\left|\infrho\right|
\end{equation}
where $\infrho\in[-1,~1]$ denotes a user-specified parameter. Similarly, the condition given by \cref{eq:D_num_infty} can be equivalently expressed in terms of the RK stability function $R(z)$ as
\begin{equation}
\lim\limits_{z\to-\infty}R(z)=\left|\infrho\right|.
\end{equation} 

As will be demonstrated in the following section, \cref{eq:D_num_infty} plays two important roles. First, by solving \cref{eq:D_num_infty}, the variation of the first sub-step size $\gamma_1$ with respect to the parameter $\infrho$ can be established, thereby enabling the direct prescription of high-frequency numerical dissipation of the integration algorithm. Second, once all unknown parameters contained in $D_\mathsf{num}^{\infty}$ have been determined through other approaches, \cref{eq:D_num_infty} can be utilized to evaluate the corresponding value of $\left|\infrho\right|$, thereby providing a quantitative measure of high-frequency numerical dissipation of the integration algorithm.

\section{Developments and analyses}\label{sec:development}

Based on the directly self-starting framework \eqref{eq:but_alg} established in the previous section, a family of high-order implicit integration algorithms is developed and analyzed herein. Unlike conventional sub-step formulations \cite{li_DirectlySelfstartingHigherorder_2022,lee_ImplicitSsubstepTime_2025}, where the last sub-step location is usually fixed at the end of the current time interval, i.e., $\gamma_s=1$, the present framework treats all sub-step locations as independent design variables by relaxing the constraint $\gamma_s=1$. This additional variable provides enhanced flexibility for designing higher-order accuracy without increasing either the number of sub-steps or the dimension of the governing equations. For a given number of sub-steps $s$, the unknown algorithmic parameters are determined by simultaneously enforcing the identical effective matrix given by \cref{eq:iem}, the order conditions described by \cref{eq:def_acc,eq:consistency}, and the prescribed high-frequency dissipation in \cref{eq:D_num_infty}. The cases with a small number of sub-steps are first investigated in detail because they provide fundamental insights into the relationship between sub-step locations, attainable temporal accuracy, and numerical dissipation. Subsequently, the general formulations for more sub-steps are established in \cref{app:a} to construct high-order implicit algorithms.

The proposed framework reveals a flexible trade-off between temporal accuracy and controllable numerical dissipation. Specifically, an $s$-sub-step implicit algorithm can be constructed to achieve $s$th-order accuracy while preserving unconditional stability and controllable high-frequency dissipation. In this case, all sub-step locations $\gamma_j~(j=1,~\cdots,~s)$ are free. Furthermore, the user-specified parameter $\infrho$ is retained to control numerical dissipation, which determines the first sub-step size $\gamma_1$, whereas the remaining sub-step locations $\gamma_j~(j=2,~\cdots,~s)$ are still free. Alternatively, by fully exploiting all sub-step locations and imposing one additional order condition, the accuracy of the $s$-sub-step algorithm can be further improved to $(s+1)$th order. In this higher-order case, all sub-step locations $\gamma_j~(j=1,~\cdots,~s)$ are uniquely determined, leaving no remaining free parameter for dissipation control. Consequently, the resulting algorithms possess fixed high-frequency numerical dissipation. Therefore, the proposed framework enables the construction of different algorithm variants by appropriately allocating sub-step locations between accuracy enhancement and numerical dissipation control.

\subsection{Single sub-step: $s=1$}

The simplest case corresponds to a single sub-step, i.e., $s=1$. In this case, the proposed sub-step formulation reduces to a conventional single-sub-step implicit algorithm, whose Butcher tableau is given by
\begin{equation}\label{eq:s_eq_1}
\begin{NiceArray}{c|c}[cell-space-limits=3pt,columns-width=0.5cm]
\mbf{c} & \mbf{A}\\ \hline & \mbf{b}\T
\end{NiceArray}\quad=\quad\begin{NiceArray}{c|c}[cell-space-limits=3pt,columns-width=0.7cm]
\gamma_1 & \alpha_{11}  \\ \hline
& \beta_1
\end{NiceArray}.
\end{equation}
Based on \cref{eq:s_eq_1}, the corresponding single-sub-step implicit method can be written explicitly as
\begin{subequations}
\begin{align}
\mbfM\mbfa_{n+\gamma_1}&+\mbfC\mbfv_{n+\gamma_1}+\mbfK\mbfu_{n+\gamma_1}=\mbfF_{n+\gamma_1}\\
\mbfv_{n+\gamma_1}&=\mbfv_n+\alpha_{11}\dt\mbfa_{n+\gamma_1}&
\mbfv_{n+1}&=\mbfv_n+\beta_1\dt\mbfa_{n+\gamma_1}\\
\mbfu_{n+\gamma_1}&=\mbfu_n+\alpha_{11}\dt\mbfv_{n+\gamma_1}& \mbfu_{n+1}&=\mbfu_n+\beta_1\dt\mbfv_{n+\gamma_1}
\end{align}
\end{subequations}
for second-order transient dynamics, and as
\begin{subequations}
\begin{align}
\mbfC\mbfv_{n+\gamma_1}&+\mbfK\mbfu_{n+\gamma_1}=\mbfF_{n+\gamma_1}\\
\mbfu_{n+\gamma_1}&=\mbfu_n+\alpha_{11}\dt\mbfv_{n+\gamma_1}& \mbfu_{n+1}&=\mbfu_n+\beta_1\dt\mbfv_{n+\gamma_1}
\end{align}
\end{subequations}
for first-order transient dynamics. Although the resulting single-sub-step method has a relatively simple structure, it provides an important reference for clarifying the relationship between the proposed framework and classical implicit algorithms.

For $s=1$, the identical effective matrix given in \cref{eq:iem} is satisfied automatically, whereas the consistency condition given by \cref{eq:consistency} requires $\alpha_{11}=\gamma_1$. The numerical amplification factor and numerical load operator are obtained from \cref{eq:num_DL} as
\begin{subequations}
\begin{align}
D_\mathsf{num}&=1-\dfrac{\beta_1\omega\dt}{\alpha_{11}\omega\dt+1} \\
L_\mathsf{num}&=\dfrac{\beta_1\dt\exp(\gamma_1\dt)}{\alpha_{11}\omega\dt+1}.
\end{align}
\end{subequations}
Accordingly, the local truncation errors defined in \cref{eq:def_acc} are expressed as
\begin{subequations}
\begin{align}
D_\mathsf{num}-D_\mathsf{exa}&=(1-\beta_1)\omega\dt+ \left( \beta_1\alpha_{11}-\dfrac{1}{2} \right)\omega^{2}\dt ^{2} +\mathcal{O}(\dt^{3})\\
L_\mathsf{num}-L_\mathsf{exa}&=-(1-\beta_1)\dt+ \left[ \beta_1\gamma_1- \dfrac{1}{2}-\left( \beta_1\alpha_{11}-\dfrac{1}{2} \right) \omega  \right] \dt^{2} +\mathcal{O}(\dt^{3}),
\end{align}
\end{subequations}
from which it follows that the present single-sub-step implicit method achieves first-order accuracy if and only if $\beta_1=1$. Controllable numerical dissipation is imposed through $D_\mathsf{num}^\infty=\lim\limits_{\omega\to\infty} D_\mathsf{num}=1-1/\alpha_{11}=\left|\infrho\right|$, which gives $\alpha_{11}=1/(\infrho+1)$
with $\infrho\in(-1,~1]$.

Consequently, all algorithmic parameters in \cref{eq:s_eq_1} are uniquely determined as
\begin{equation}\label{eq:s1_params}
\alpha_{11}=\gamma_1=\dfrac{1}{\infrho+1}\quad\text{and}\quad \beta_1=1,
\end{equation}
where $\infrho\in(-1,~1]$. The stability polynomial $E(y)$ is computed as $y^{2}(1-\infrho^{2})\ge 0$, which always holds for $\infrho\in(-1,~1]$, thus resulting in unconditional stability. For $\infrho\in(-1,~1)$, the method defined by \cref{eq:s_eq_1} is first-order accurate, whereas it becomes second-order accurate for $\infrho=1$ because $\beta_1\alpha_{11}=\beta_1\gamma_1=1/2$. It should be emphasized that the single-sub-step implicit algorithm considered within the present framework does not constitute a novel contribution of this study, since it has already been developed in the previous work \cite{li_DirectlySelfstartingHigherorder_2022}.

The amplitude ($\delta$) and phase ($\epsilon$) errors \cite{li_NovelImplicitIntegration_2025,li_DesigningDevelopingSinglestep_2023} are computed analytically for the single-sub-step implicit method \eqref{eq:s1_params} as 
\begin{subequations}\label{eq:s1_amp_phas}
\begin{align}
\delta=& \dfrac{(2\xi^{2}-1)(\infrho-1)}{2(1+\infrho)}\omega^{2}\dt+ \dfrac{(4\xi^{2}-3)\xi(\infrho^{2}-\infrho+1)}{3(1+\infrho)^{2}}\omega^3\dt^{2}+ \dfrac{(\infrho-1)(\infrho^{2}+1)}{4(1+\infrho)^{3}}\omega^4\dt^{3}+\mathcal{O}(\dt^4)   \\  
\epsilon=& \dfrac{\xi\sqrt{1-\xi^{2}}(\infrho-1)}{1+\infrho}\omega^{2}\dt+ \dfrac{(4\xi^{2}-1)\sqrt{1-\xi^{2}}(\infrho^{2}-\infrho+1)}{3(1+\infrho)^{2}}\omega^{3}\dt^{2}+ \mathcal{O}(\dt^3)  .
\end{align}
\end{subequations}
For the first-order single-sub-step member with $\infrho\neq1$, the leading terms of both the amplitude $\delta$ and phase $\epsilon$ errors are proportional to $(\infrho-1)\dt$. Therefore, in the damped case $(\xi\neq0)$, both errors are generally of order $\mathcal{O}(\dt)$. In the undamped case $(\xi=0)$, the leading term of the amplitude error remains nonzero, so that $\delta=\mathcal{O}(\dt)$, whereas the corresponding leading term in $\epsilon$ vanishes because it is proportional to $\xi$. Consequently, the phase error exhibits one-order superconvergence and becomes $\epsilon=\mathcal{O}(\dt^{2})$. For the second-order member with $\infrho=1$, all $\dt$ terms vanish. In the damped case, the leading terms of both $\delta$ and $\epsilon$ are therefore of order $\mathcal{O}(\dt^{2})$. In the undamped case, the amplitude error in Eq.~(\ref{eq:s1_amp_phas}a) vanishes identically, indicating that the method introduces no numerical amplitude decay (zero dissipation), whereas the phase error generally remains $\mathcal{O}(\dt^{2})$. Thus, the second-order member is non-dissipative for undamped systems and retains second-order phase accuracy.

\begin{figure}[htbp]
	\centering 
	\subfigure[$\infrho=1.0$]{
		\includegraphics[scale=0.4]{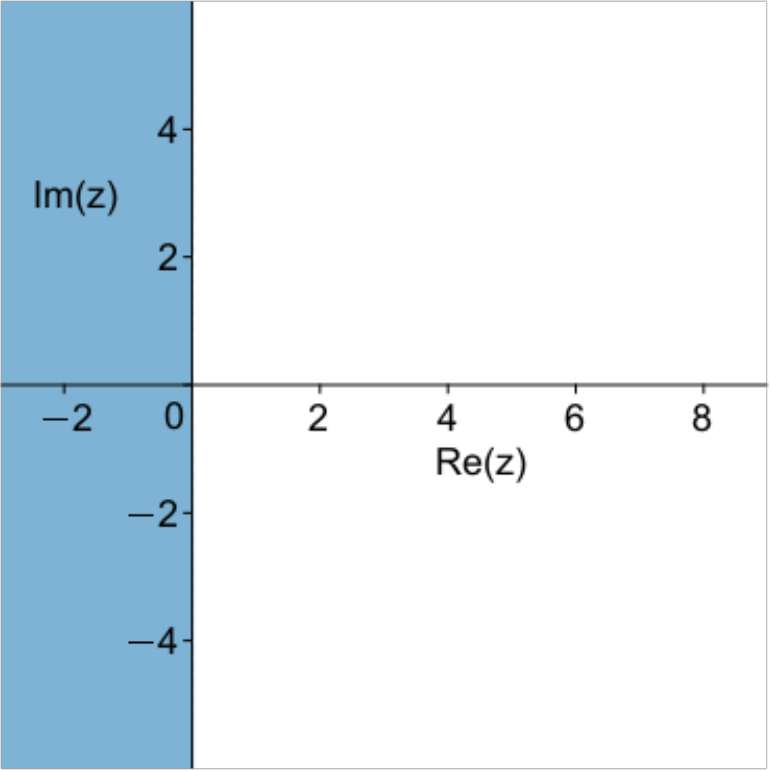}}
	\subfigure[$\infrho=0.6$]{
		\includegraphics[scale=0.4]{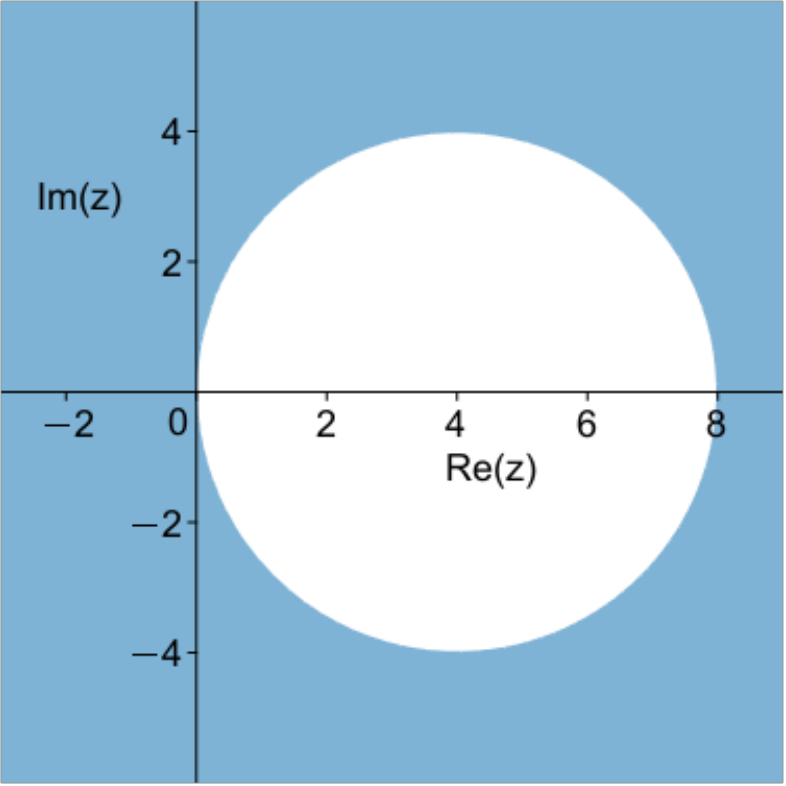}}
	\subfigure[$\infrho=0.3$]{
		\includegraphics[scale=0.4]{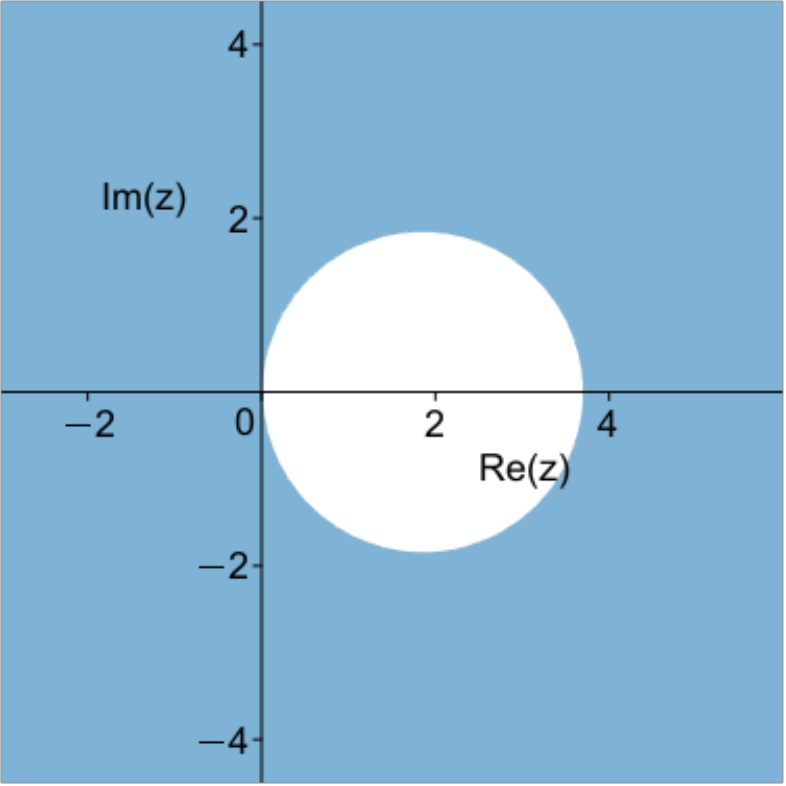}}
	\subfigure[$\infrho=0.0$]{
		\includegraphics[scale=0.4]{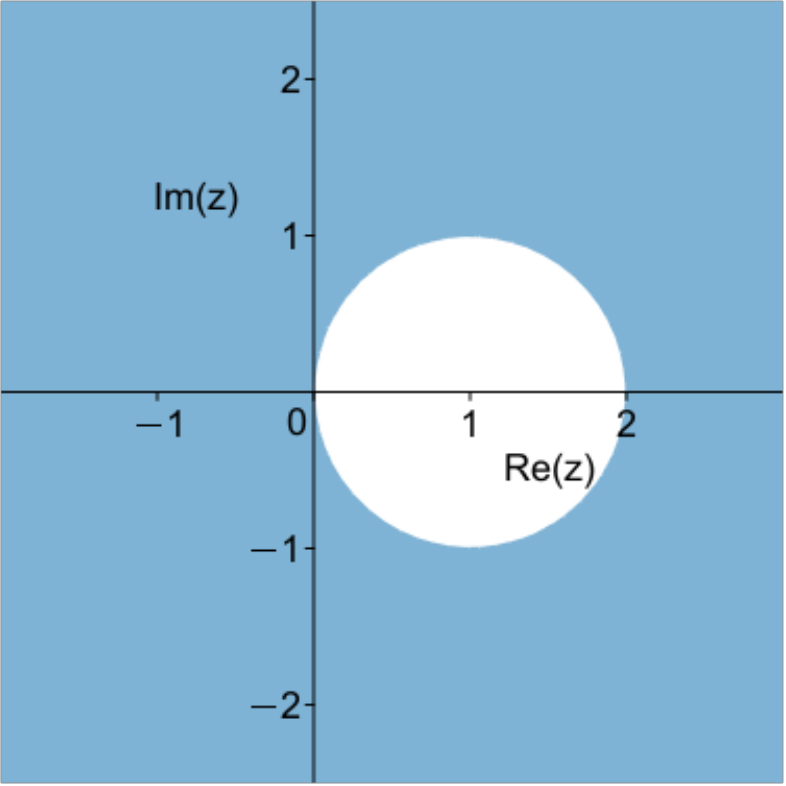}}
	\subfigure[$\infrho=-0.3$]{
		\includegraphics[scale=0.4]{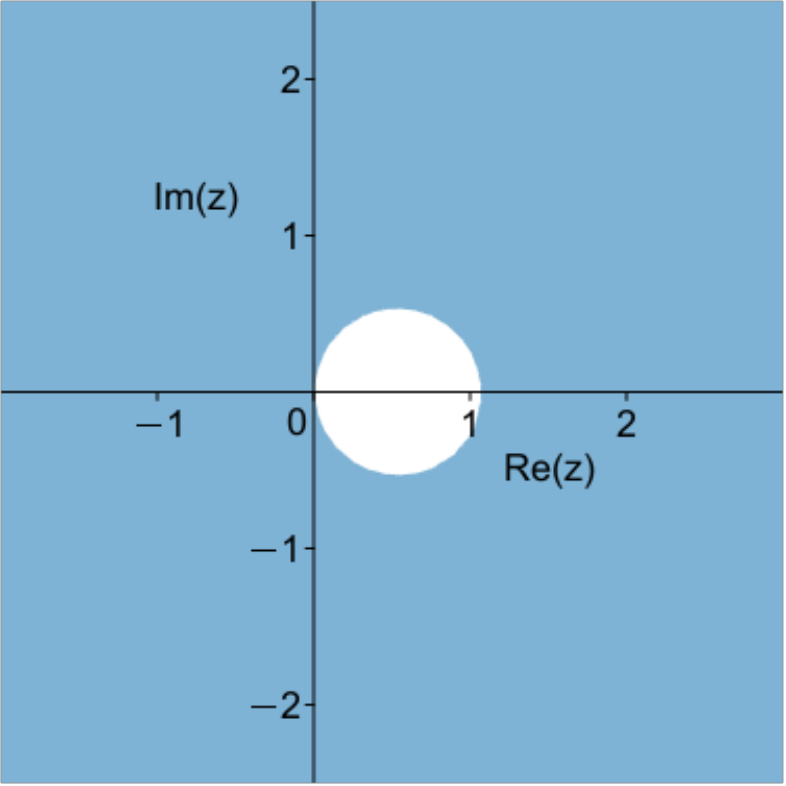}}
	\subfigure[$\infrho=-0.6$]{
		\includegraphics[scale=0.4]{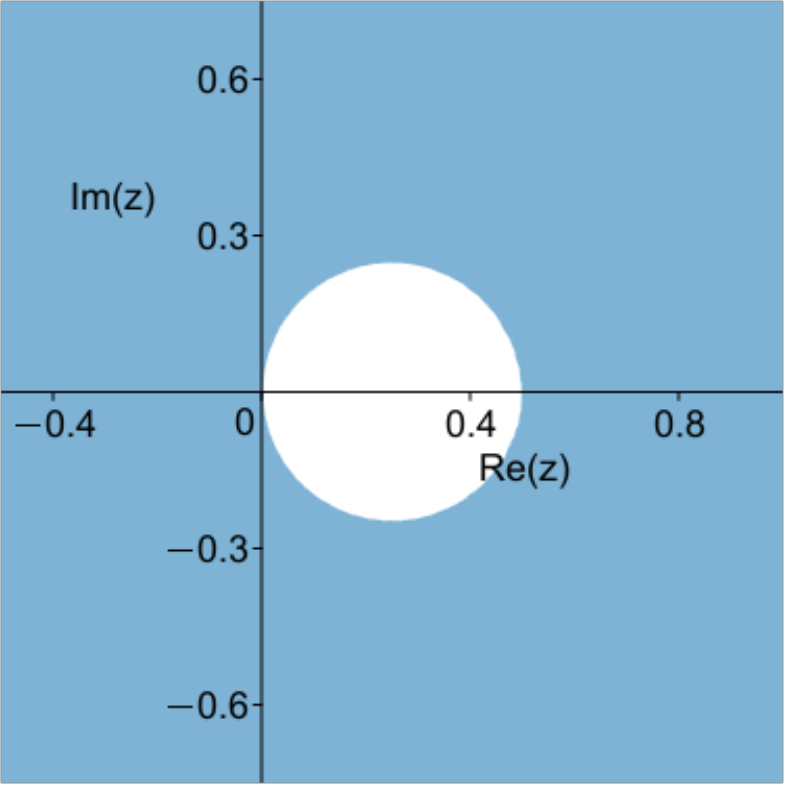}}
	\caption{The stability region for the single-sub-step implicit method defined by \cref{eq:s_eq_1,eq:s1_params}.}
	\label{fig:s1_stability}
\end{figure}

\cref{fig:s1_stability} illustrates the stability regions of the single-sub-step implicit method defined by \cref{eq:s_eq_1,eq:s1_params} for different values of $\infrho$. It can be observed that all cases contain the entire left half-plane, confirming the unconditional stability of the proposed method. When $\infrho=1.0$, the stability function satisfies $\left|R(\infty)\right|=1$, corresponding to the non-dissipative case. As $\infrho$ decreases, the high-frequency spectral radius is reduced, indicating enhanced numerical dissipation for suppressing spurious high-frequency oscillations. Although the shape and size of the stability regions vary with $\infrho$, all cases preserve unconditional stability. 

\subsection{Two sub-steps: $s=2$}
The two-sub-step case represents the simplest configuration where the flexibility of the proposed framework becomes evident. In the case of $ s=2 $, the corresponding Butcher tableau is expressed as
\begin{equation}\label{eq:s_eq_2}
\begin{NiceArray}{c|c}[cell-space-limits=3pt,columns-width=0.5cm]
\mbf{c} & \mbf{A}\\ \hline & \mbf{b}\T
\end{NiceArray}\quad=\quad\begin{NiceArray}{c|cc}[cell-space-limits=3pt,columns-width=0.7cm]
\gamma_1 & \alpha_{11} & 0  \\ 
\gamma_2 & \alpha_{21} & \alpha_{22} \\ \hline
& \beta_1 & \beta_2
\end{NiceArray}.
\end{equation}
Unlike previous directly self-starting two-sub-step methods \cite{li_DirectlySelfstartingHigherorder_2022,kim_ImprovedImplicitMethod_2020}, where the last sub-step is fixed at the end of the time interval, namely $\gamma_2=1$, the present formulation treats the second sub-step location as an independent variable. In the present case, the present two-sub-step implicit method is expressed as 
\begin{subequations}
\begin{align}
\mbfM\mbfa_{n+\gamma_1}&+\mbfC\mbfv_{n+\gamma_1}+\mbfK\mbfu_{n+\gamma_1}=\mbfF_{n+\gamma_1} & \mbfM\mbfa_{n+\gamma_2}&+\mbfC\mbfv_{n+\gamma_2}+\mbfK\mbfu_{n+\gamma_2}=\mbfF_{n+\gamma_2}\\
\mbfv_{n+\gamma_1}&=\mbfv_n+\alpha_{11}\dt\mbfa_{n+\gamma_1}&
\mbfv_{n+\gamma_2}&=\mbfv_n+\dt \left( \alpha_{21}\mbfa_{n+\gamma_1}+\alpha_{22}\mbfa_{n+\gamma_2} \right) \\
\mbfu_{n+\gamma_1}&=\mbfu_n+\alpha_{11}\dt\mbfv_{n+\gamma_1}& \mbfu_{n+\gamma_2}&=\mbfu_n+\dt \left( \alpha_{21}\mbfv_{n+\gamma_1}+\alpha_{22}\mbfv_{n+\gamma_2} \right) \\
\mbfv_{n+1}&=\mbfv_n+\dt \left( \beta_1\mbfa_{n+\gamma_1}+\beta_2\mbfa_{n+\gamma_2} \right) \\
\mbfu_{n+1}&=\mbfu_n+\dt \left( \beta_1\mbfv_{n+\gamma_1}+\beta_2\mbfv_{n+\gamma_2} \right) 
\end{align}
\end{subequations}
for second-order transient dynamics, and as 
\begin{subequations}
\begin{align}
\mbfC\mbfv_{n+\gamma_1}&+\mbfK\mbfu_{n+\gamma_1}=\mbfF_{n+\gamma_1} & \mbfC\mbfv_{n+\gamma_2}&+\mbfK\mbfu_{n+\gamma_2}=\mbfF_{n+\gamma_2}\\
\mbfu_{n+\gamma_1}&=\mbfu_n+\alpha_{11}\dt\mbfv_{n+\gamma_1}& \mbfu_{n+\gamma_2}&=\mbfu_n+\dt \left( \alpha_{21}\mbfv_{n+\gamma_1}+\alpha_{22}\mbfv_{n+\gamma_2} \right) \\
\mbfu_{n+1}&=\mbfu_n+\dt \left( \beta_1\mbfv_{n+\gamma_1}+\beta_2\mbfv_{n+\gamma_2} \right) 
\end{align}
\end{subequations}
for first-order transient dynamics. 

The two-sub-step implicit method \eqref{eq:s_eq_2} achieves the identical effective matrices within two sub-steps, see \cref{eq:iem}, when $\alpha_{11}=\alpha_{22}$, and the consistency in each sub-step given in \cref{eq:consistency} requires $\alpha_{11}=\gamma_1$ and $\alpha_{21}+\alpha_{22}=\gamma_2$. With these conditions in hand, the numerical amplification factor and numerical load operator are simplified as 
\begin{subequations}
\begin{align}
D_\mathsf{num}&=\dfrac{1-(\beta_1+\beta_2-2\gamma_1)\omega\dt-(\beta_1\gamma_1+2\beta_2\gamma_1-\beta_2\gamma_2-\gamma_1^{2})\omega^{2}\dt^{2}}{(1+\gamma_1\omega\dt)^{2}} \\
L_\mathsf{num}&=\dfrac{\beta_1\dt+(\beta_1\gamma_1+\beta_2(\gamma_1-\gamma_2))\omega\dt^{2}}{(1+\gamma_1\omega\dt)^{2}}\exp(\gamma_1\dt)+\dfrac{\beta_2\dt}{1+\gamma_1\omega\dt} \exp(\gamma_2\dt). 
\end{align}
\end{subequations}
Then, the local truncation errors in \cref{eq:def_acc} are computed as 
\begin{subequations}\label{eq:s2_DL_acc}
\begin{align}
D_\mathsf{num}-D_\mathsf{exa}&=-(\beta_1+\beta_2-1)\omega\dt+ \left( \beta_1\gamma_1+\beta_2\gamma_2-\dfrac{1}{2}  \right)\omega^{2}\dt^{2}- \left[ (\beta_1-\beta_2)\gamma_1^{2}+2\beta_2\gamma_1\gamma_2-\dfrac{1}{6}  \right] \omega^{3}\dt^{3}  +\mathcal{O}(\dt^{4})\\
L_\mathsf{num}-L_\mathsf{exa}&= (\beta_1+\beta_2-1)\dt-\left( \beta_1\gamma_1+\beta_2\gamma_2-\dfrac{1}{2}  \right)(\omega-1)\dt^{2}+\bigg[ \left( (\beta_1-\beta_2)\gamma_1^{2}+2\beta_2\gamma_1\gamma_2-\dfrac{1}{6} \right) \omega(\omega-1)\notag\\
&\quad +\dfrac{1}{2} \left( \beta_1\gamma_1^{2}+\beta_2\gamma_2^{2}-\dfrac{1}{3}  \right)   \bigg]\dt^{3} +\mathcal{O}(\dt^{4}). 
\end{align}
\end{subequations}
It follows that the two-sub-step implicit method is second-order accurate when the parameters satisfy 
\begin{equation}
\left.\begin{aligned}
\beta_1+\beta_2-1&=0\\ 
\beta_1\gamma_1+\beta_2\gamma_2-\dfrac{1}{2}&=0
\end{aligned}\right\}\quad \implies\quad \left\{ \begin{aligned}
\beta_1&=\dfrac{2\gamma_2-1}{2(\gamma_2-\gamma_1)}\\
\beta_2&=\dfrac{2\gamma_1-1}{2(\gamma_1-\gamma_2)}. 
\end{aligned}\right.
\end{equation}

By substituting the known algorithmic parameters into the numerical amplification factor, the corresponding stability function $R(z)$ can be simplified accordingly. As a result, the stability polynomial $E(y)$ can be further reduced to $E(y)=y^4(16\gamma_1^{3}-20\gamma_1^{2}+8\gamma_1-1)\ge 0$, and thus the present two-sub-step implicit method achieves unconditional stability when 
\begin{equation}
\dfrac{1}{4}\le\gamma_1\le \infty. 
\end{equation}

\begin{figure}[htbp]
	\centering 
	\includegraphics[scale=1.0]{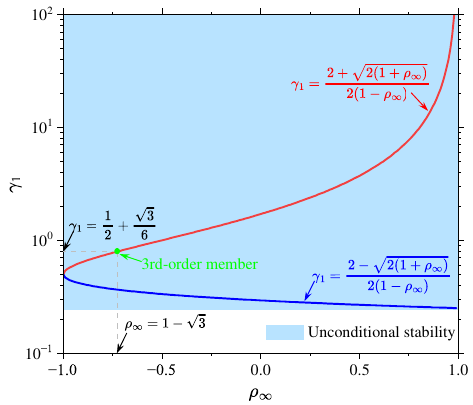}
	\caption{The variations of $\gamma_1$ with $\infrho\in[-1,~1]$ for the two-sub-step implicit method.}
	\label{fig:s2_r1_rho}
\end{figure}

The dissipation control at the high-frequency limit ($\omega\to\infty$) given in \cref{eq:D_num_infty} is simplified as 
\begin{equation}\label{eq:s2_r1}
D_\mathsf{num}^\infty= \dfrac{2\gamma_1^{2}-4\gamma_1+1}{2\gamma_1^{2}}=\left|\infrho\right|\quad\implies\quad \gamma_1= \dfrac{2\pm \sqrt{2(1+\infrho)}}{2(1-\infrho)} .
\end{equation}
\cref{eq:s2_r1} establishes the relationship between the user-specified parameter $\infrho$ and the first sub-step size $\gamma_1$. It can be observed that, for a prescribed value of $\infrho$, two possible solutions of $\gamma_1$ are obtained due to the quadratic form of the dissipation control equation. The variations of these two solutions with respect to $\infrho$ are illustrated in \cref{fig:s2_r1_rho}. The upper branch of $\gamma_1$ increases rapidly as $\infrho$ approaches unity and may even become unbounded, which leads to an excessively large first sub-step size and is therefore unsuitable for practical time integration. In contrast, the lower branch remains within a reasonable range for the entire interval $\infrho\in[-1,~1]$. Spectral analyses further indicate that selecting $\gamma_1\in[1/4,~1/2]$ provides more desirable numerical properties, such as smaller relative period errors. Therefore, the lower branch of the solution is often adopted, namely
\begin{equation}\label{eq:s2_r1_adopted}
\gamma_1=\dfrac{2-\sqrt{2(1+\infrho)}}{2(1-\infrho)}
\quad\text{with}\quad
\infrho\in[-1,~1].
\end{equation}
It should be noted that the selection of $\gamma_1$ based on \cref{eq:s2_r1_adopted} not only guarantees the prescribed high-frequency numerical dissipation but also ensures that the first sub-step size remains within the desirable range. As shown in \cref{fig:s2_r1_rho}, the entire selected branch is located inside the region corresponding to unconditional stability.

In summary, the present two-sub-step implicit method \eqref{eq:s_eq_2} can be represented in the following simplified Butcher tableau:
\begin{equation}\label{eq:s2_2nd_but}
\begin{NiceArray}{c|cc}[cell-space-limits=3pt,columns-width=0.7cm]
\gamma_1 & \gamma_1 & 0  \\ 
\gamma_2 & \gamma_2-\gamma_1 & \gamma_1 \\ \hline
& \dfrac{2\gamma_2-1}{2(\gamma_2-\gamma_1)} & \dfrac{2\gamma_1-1}{2(\gamma_1-\gamma_2)}
\end{NiceArray}
\end{equation}
where the sub-step sizes $\gamma_j~(j=1,~2)$ are free parameters. Of course, $\gamma_1$ can also be  determined by either \cref{eq:s2_r1} or \cref{eq:s2_r1_adopted}, whereas $\gamma_2$ remains the only free parameter within the second-order accurate framework. It is noteworthy that any admissible choice of $\gamma_2$ does not affect spectral properties or stability regions of the resulting algorithm.

\begin{remark}
Different criteria have been adopted in the literature to determine the free parameter $\gamma_2$. Ji et al. \cite{ji_UnconditionallyStableTime_2021} determined $\gamma_2$ by enforcing the BN-stability for nonlinear systems, resulting in $\gamma_2=1-\gamma_1$. Alternatively, Kim \cite{kim_ImprovedTwostageImplicit_2024} derived $\gamma_2$ based on the energy conservation requirement for a conservative nonlinear SDOF system, yielding $\gamma_2=(3\gamma_1-2)/(6\gamma_1-3)$. Li et al. \cite{li_DirectlySelfstartingHigherorder_2022} and Kim \cite{kim_ImprovedImplicitMethod_2020} directly prescribed $\gamma_2$ as unity, i.e., $\gamma_2=1$. Wang et al. \cite{wang_TrulySelfstartingComposite_2023} selected $\gamma_2=1/2$. Similarly, the first two sub-step sizes can also be chosen to be identical, leading to $\gamma_2=2\gamma_1 $. 
\end{remark}

\begin{figure}[htbp]
	\centering 
	\subfigure[$\gamma_1=1/4~(\infrho=1.0)$]{
		\includegraphics[scale=0.407]{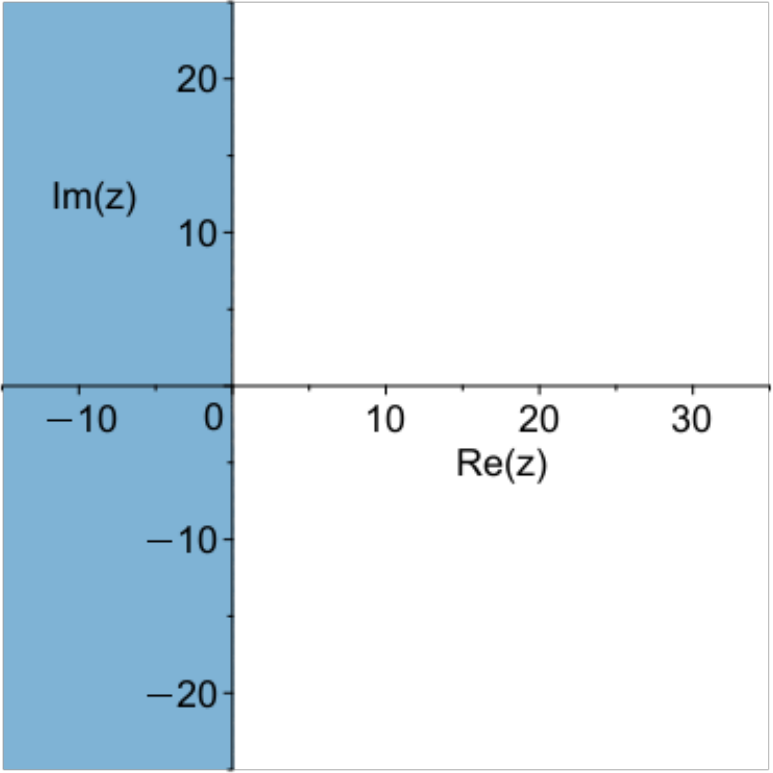}}
	\subfigure[$\gamma_1=0.26794~(\infrho=0.5)$]{
		\includegraphics[scale=0.4]{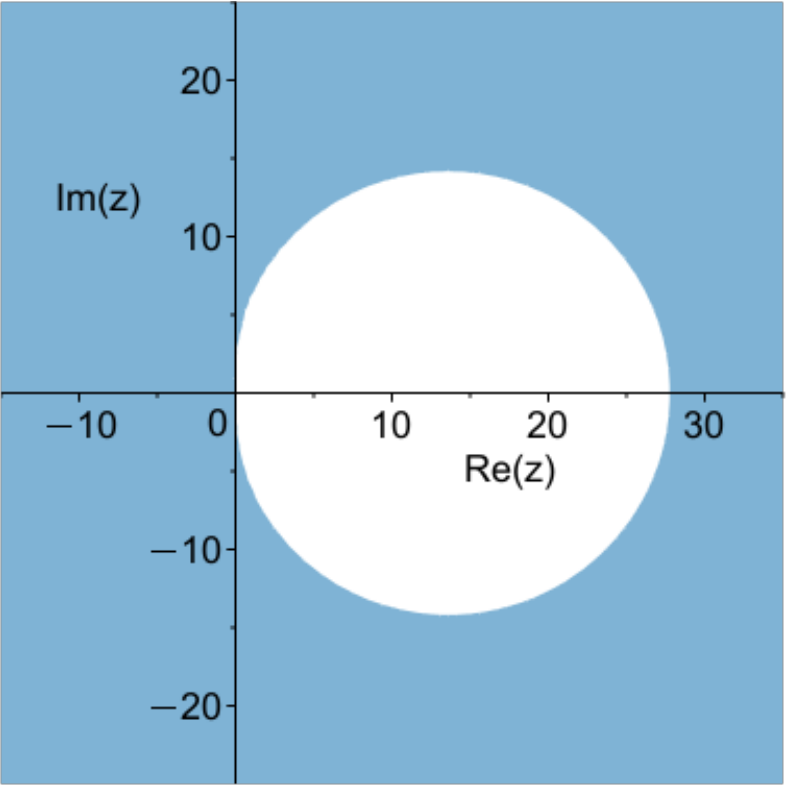}}
	\subfigure[$\gamma_1=0.29289~(\infrho=0.0)$]{
		\includegraphics[scale=0.4]{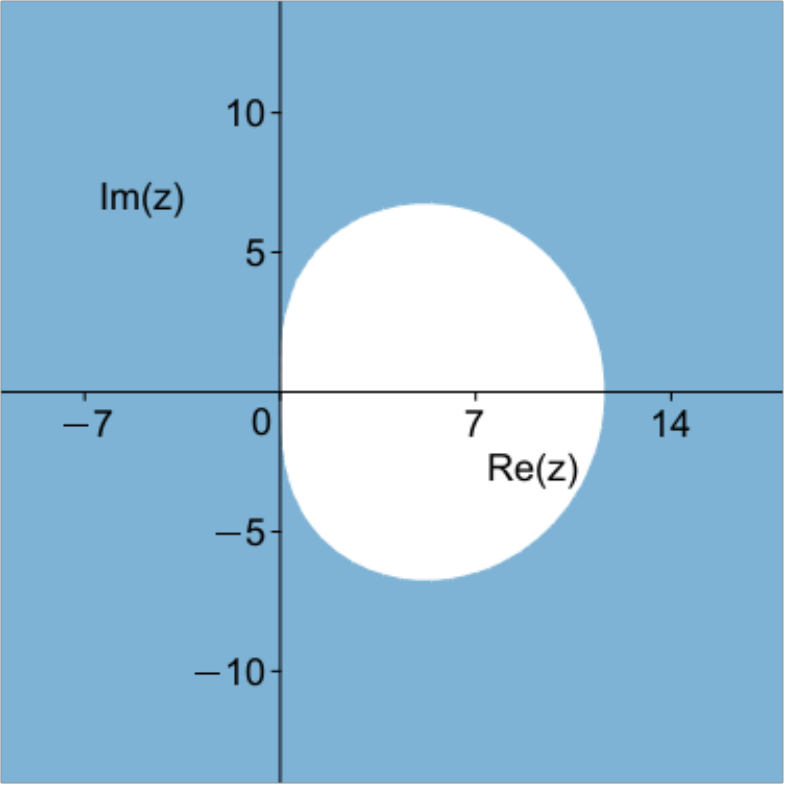}}
	\subfigure[$\gamma_1=1/3~(\infrho=-0.5)$]{
		\includegraphics[scale=0.4]{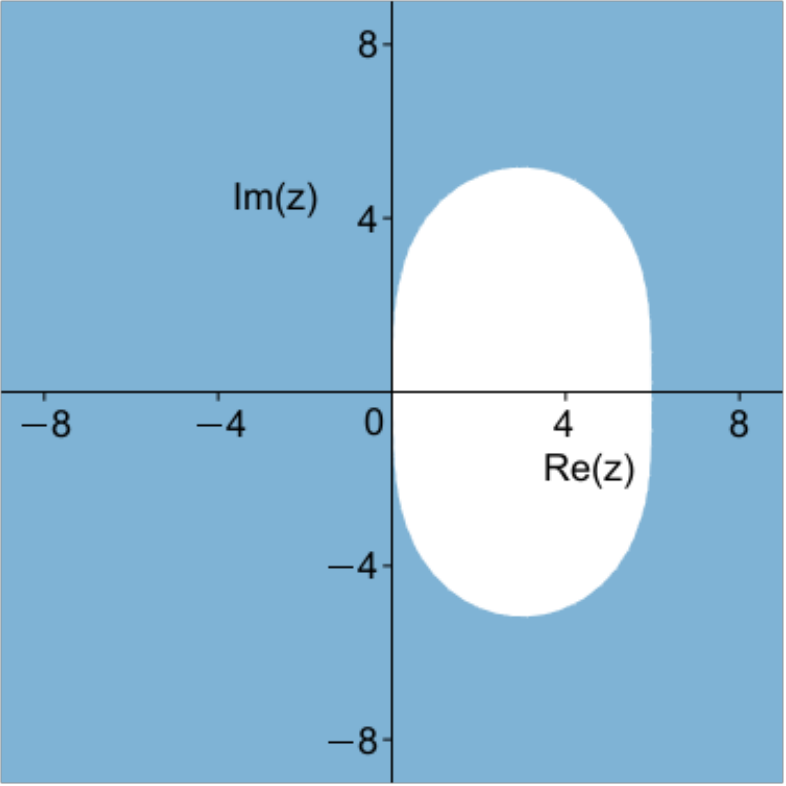}}
	\subfigure[$\gamma_1=1/2~(\infrho=-1.0)$]{
		\includegraphics[scale=0.407]{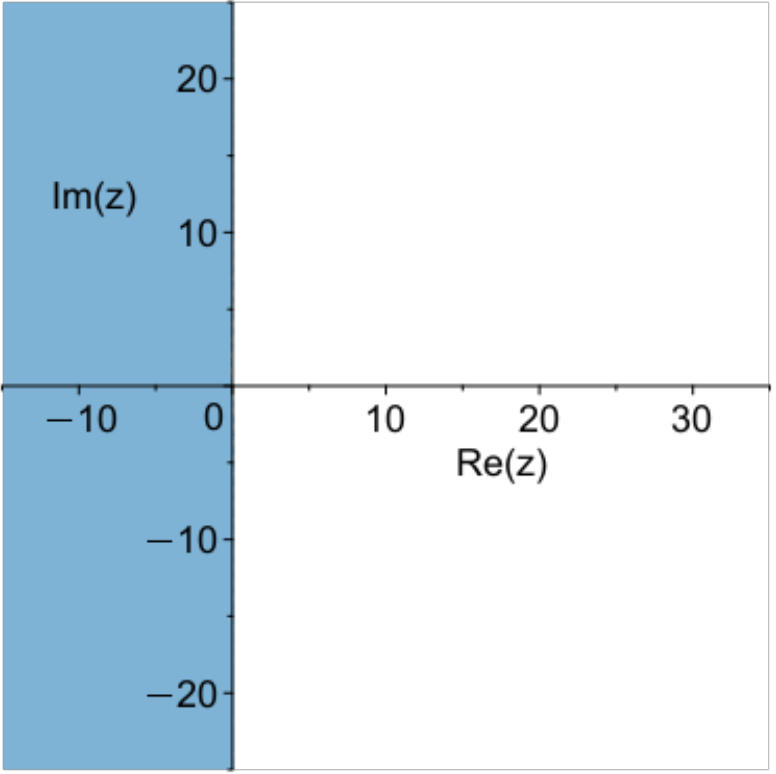}}
	\subfigure[$\gamma_1=1/2+\sqrt{3}/6~(\infrho=1-\sqrt{3})$]{
		\includegraphics[scale=0.4]{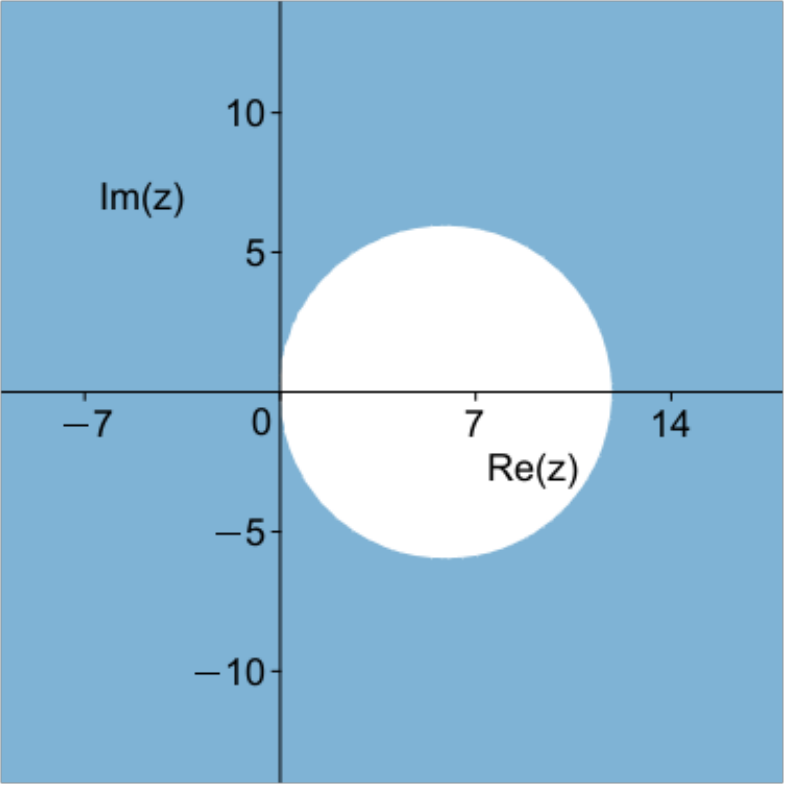}}
    \subfigure[$\gamma_1=1.0~(\infrho=-0.5)$]{
		\includegraphics[scale=0.4]{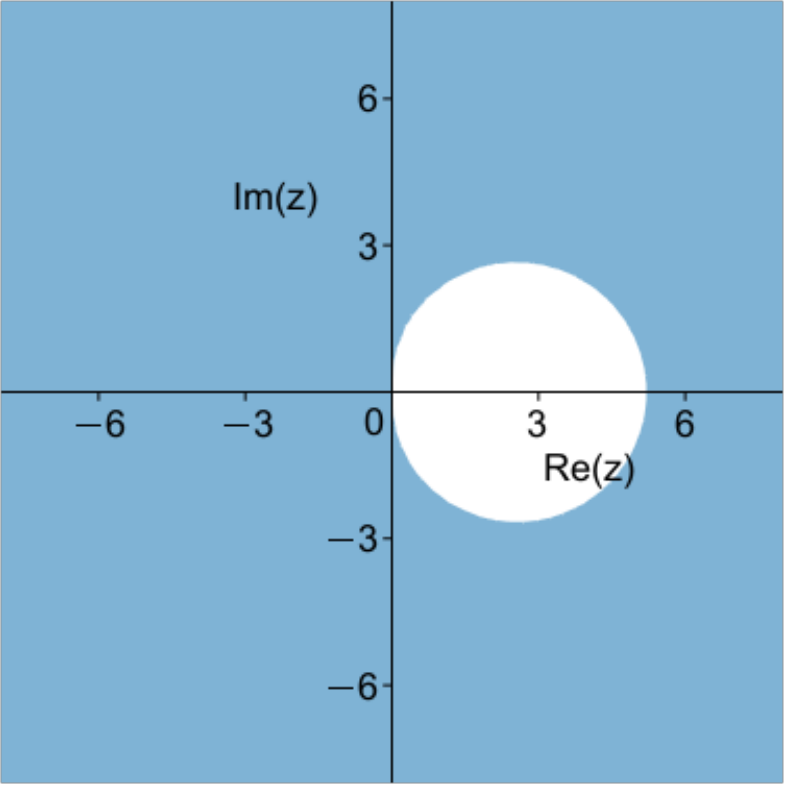}}
    \subfigure[$\gamma_1=1.70711~(\infrho=0.0)$]{
		\includegraphics[scale=0.4]{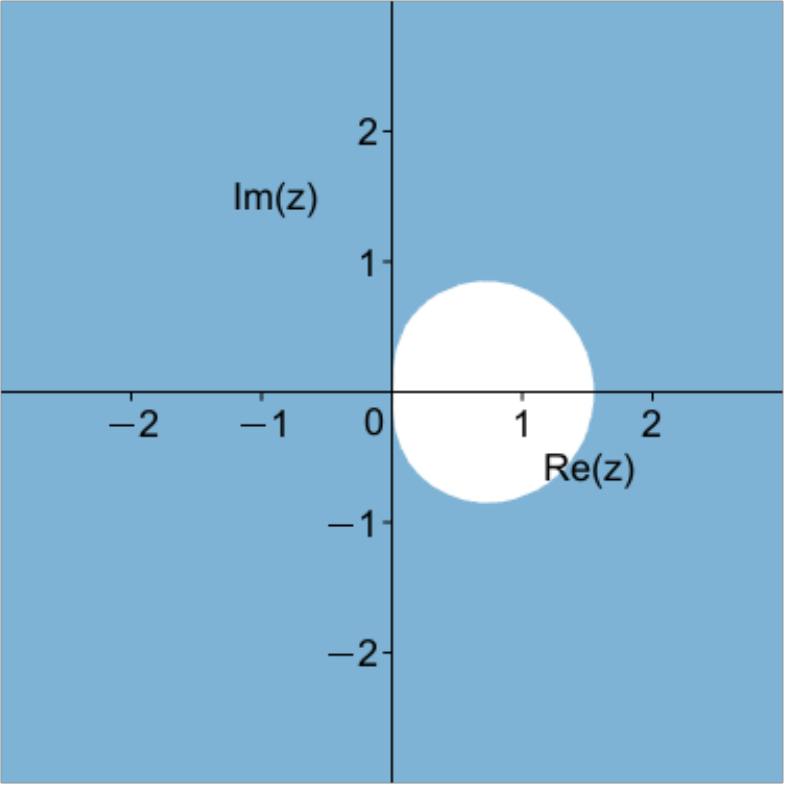}}
    \subfigure[$\gamma_1=3.73205~(\infrho=0.5)$]{
		\includegraphics[scale=0.4]{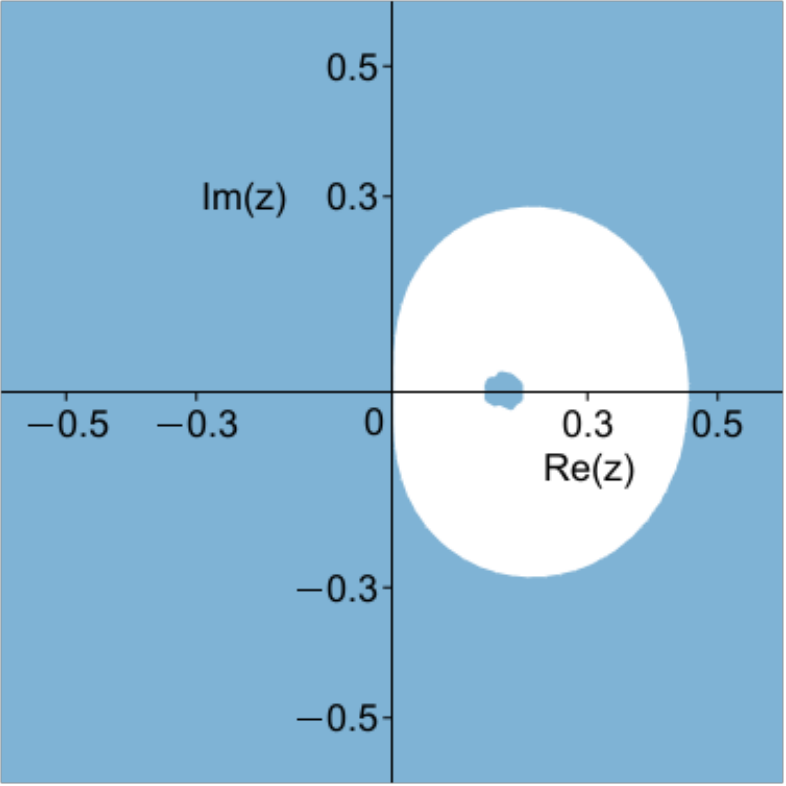}}
	\caption{The stability region for the two-sub-step implicit method defined by \cref{eq:s2_2nd_but}.}
	\label{fig:s2_stability}
\end{figure}

\cref{fig:s2_stability} presents the stability regions of the two-sub-step implicit method defined by \cref{eq:s2_2nd_but} for different selections of the first sub-step size $\gamma_1$ and the corresponding high-frequency spectral radius $\infrho$. The stability regions are obtained based on the RK stability function $R(z)$ introduced previously. It can be observed that all cases preserve unconditional stability, since their stability regions contain the entire left half-plane. Figs.~\ref{fig:s2_stability}(f)-(i) further illustrate the stability regions corresponding to the alternative branch of $\gamma_1>1/2$. Compared with the recommended branch satisfying $\gamma_1\in[1/4,~1/2]$, these choices lead to significantly different stability boundaries but preserve unconditional stability. 


\subsubsection{Higher-order accuracy}

By utilizing the second sub-step size $\gamma_2$, the order conditions can be further satisfied without introducing additional sub-steps. The resulting two-sub-step algorithm is therefore expected to achieve higher-order accuracy than existing directly self-starting schemes \cite{ji_UnconditionallyStableTime_2021,kim_ImprovedTwostageImplicit_2024,li_DirectlySelfstartingHigherorder_2022,kim_ImprovedImplicitMethod_2020,wang_TrulySelfstartingComposite_2023} with the same number of sub-steps. Using the known parameters to simplify \cref{eq:s2_DL_acc} as 
\begin{subequations}
\begin{align}
D_\mathsf{num}-D_\mathsf{exa}&=\dfrac{6\gamma_1^{2}-6\gamma_1+1}{6}\omega^{3}\dt^{3}   +\mathcal{O}(\dt^{4})\\
L_\mathsf{num}-L_\mathsf{exa}&= \left( \dfrac{6\gamma_1^{2}-6\gamma_1+1}{6}\omega(1-\omega)+\dfrac{1}{4}(\gamma_1+\gamma_2)-\dfrac{1}{2}\gamma_1\gamma_2-\dfrac{1}{6}    \right) \dt^{3} +\mathcal{O}(\dt^{4}). 
\end{align}
\end{subequations}
Therefore, the present two-sub-step implicit method \eqref{eq:s2_2nd_but} can achieve third-order accuracy if and only if $6\gamma_1^{2}-6\gamma_1+1=0$ and $\frac{1}{4}(\gamma_1+\gamma_2)-\frac{1}{2}\gamma_1\gamma_2-\frac{1}{6}=0$, yielding (note that $\gamma_1 \ge1/2$ is required to achieve unconditional stability)
\begin{subequations}
\begin{align}
\gamma_1&=\dfrac{1}{2}+ \dfrac{\sqrt{ 3 } }{6}\label{eq:s2_3rd_r1}  \\
\gamma_2&=\dfrac{3\gamma_1-2}{3(2\gamma_1-1)}. \label{eq:s2_3rd_r2}
\end{align}
\end{subequations}
After substituting \cref{eq:s2_3rd_r1} into \cref{eq:s2_r1}, one can compute the amount of high-frequency numerical dissipation described by $\infrho$ as 
\begin{equation}
\infrho=1-\sqrt{ 3 } \quad\implies\quad\left|\infrho\right|=\sqrt{3}-1\approx0.73205081. 
\end{equation}

Furthermore, Fig.~\ref{fig:s2_stability}(f) corresponds to the special case $\gamma_1=1/2+\sqrt{3}/6~(\infrho=1-\sqrt{3})$, where the $\gamma_2$ is utilized to satisfy an additional order condition given in \cref{eq:s2_3rd_r2} and construct the third-order two-sub-step implicit method, whose Butcher tableau is given in \cref{eq:s2_3rd_but}. The stability region of this higher-order member remains consistent with unconditional stability, indicating that the improvement in temporal accuracy does not compromise the fundamental stability of the proposed two-sub-step method.
\begin{equation}\label{eq:s2_3rd_but}
\begin{NiceArray}{c|cc}[cell-space-limits=3pt,columns-width=0.7cm]
\dfrac{1}{2}+ \dfrac{\sqrt{ 3 } }{6} & \dfrac{1}{2}+ \dfrac{\sqrt{ 3 } }{6} & 0  \\ 
\dfrac{1}{2}- \dfrac{\sqrt{ 3 } }{6} & -\dfrac{\sqrt{ 3 } }{3}  & \dfrac{1}{2}+ \dfrac{\sqrt{ 3 } }{6} \\ \hline
& \dfrac{1}{2}  & \dfrac{1}{2}
\end{NiceArray}
\end{equation}

\begin{remark}
For the second-order two-sub-step implicit method with controllable high-frequency dissipation defined by \cref{eq:s2_2nd_but}, the first sub-step size $\gamma_1$ is determined by the user-specified parameter $\infrho$ through \cref{eq:s2_r1} or \cref{eq:s2_r1_adopted}, whereas the second sub-step size $\gamma_2$ remains as a free parameter within the second-order framework. In the present study, the parameter $\gamma_2$ is determined by imposing an additional order condition, resulting in \cref{eq:s2_3rd_r2} and it is naturally adopted in the second-order framework \eqref{eq:s2_2nd_but}, so that the second- and third-order algorithms share the same sub-step splitting strategy. Consequently, the two-sub-step implicit method can be interpreted as a unified family with two possible algorithm configurations: for arbitrary user-specified $\infrho\in[-1,~1]$, the method achieves second-order accuracy with controllable high-frequency dissipation; whereas for the specific selection of $\gamma_1$ (or equivalently $\infrho$) satisfying the additional order condition given in \cref{eq:s2_3rd_r1}, it achieves third-order accuracy with a fixed high-frequency dissipation level $(\left|\infrho\right|=\sqrt{3}-1)$.
\end{remark}

\subsubsection{Spectral properties}
\begin{figure}[htbp]
	\centering 
	\includegraphics[scale=0.6]{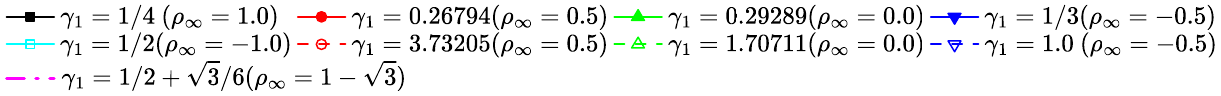}\\
    \subfigure[]{
		\includegraphics[scale=0.75]{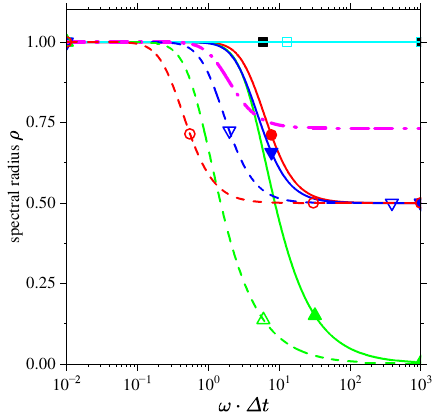}}
    \subfigure[]{
		\includegraphics[scale=0.75]{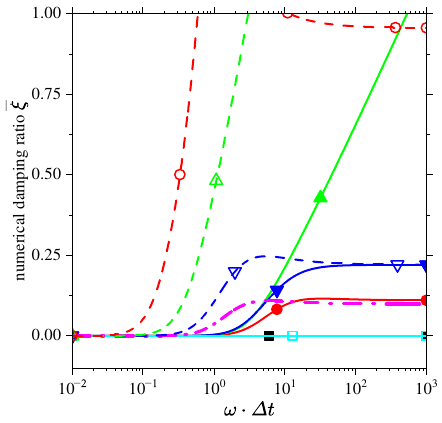}}
    \subfigure[]{
		\includegraphics[scale=0.75]{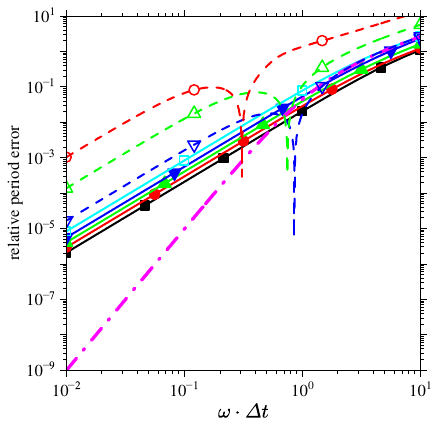}}
	\caption{Spectral properties for the two-sub-step implicit method defined by \cref{eq:s2_2nd_but} in the absence of $\xi$.}
	\label{fig:s2_sp}
\end{figure}

\cref{fig:s2_sp} presents the spectral properties of the two-sub-step implicit method defined by \cref{eq:s2_2nd_but} in the absence of $\xi$. The comparisons are performed for different selections of $\gamma_1$, including both the recommended branch satisfying $\gamma_1\in[1/4,~1/2]$ and the alternative branch discussed previously. The spectral radius, numerical damping ratio, and relative period error are investigated to evaluate the high-frequency dissipation and low-frequency accuracy, and their expressions \cite{hughes_FiniteElementMethod_2000,li_NovelImplicitIntegration_2025,li_DesigningDevelopingSinglestep_2023} are omitted herein for brevity. 

\cref{fig:s2_sp}(a) shows the variation of the spectral radius $\rho$ with respect to the dimensionless frequency $\omega\dt$. It can be observed that all methods preserve $\rho\approx1$ in the low-frequency range, indicating that the amplitude accuracy of the physically meaningful low-frequency responses is well maintained. As the frequency increases, the spectral radius gradually approaches its limiting value $\infrho$, which is consistent with the prescribed high-frequency dissipation parameter. Therefore, the parameter $\infrho$ effectively controls the amount of numerical dissipation introduced into the high-frequency range. In particular, the cases with smaller $\infrho$ exhibit stronger damping of high-frequency components, whereas the non-dissipative case with $\left|\infrho\right|=1$ maintains a unit spectral radius over the entire frequency range. The corresponding numerical damping ratios are presented in \cref{fig:s2_sp}(b). Similar to the spectral radius, the numerical damping ratio remains nearly zero in the low-frequency range, demonstrating that the proposed methods introduce negligible artificial damping for the dominant physical responses. With increasing frequency, numerical damping is gradually activated, and the damping level is determined by the prescribed value of $\infrho$. These results further confirm that the proposed formulation achieves controllable high-frequency dissipation while preserving accurate low-frequency dynamic behavior.

\cref{fig:s2_sp}(c) compares the relative period errors of different algorithmic configurations. It can be observed that the methods associated with the recommended branch of $\gamma_1$ generally provide smaller period errors than those corresponding to the alternative branch. This indicates that the selection of $\gamma_1$ not only determines the high-frequency dissipation but also significantly influences the phase accuracy of the algorithm. Among all cases, the third-order accurate member corresponding to $\infrho=1-\sqrt{3}$ exhibits superior phase accuracy, with the relative period error being substantially reduced in the low-frequency range. This improvement is attributed to the satisfaction of the additional order condition given by \cref{eq:s2_3rd_r2}, which enhances the temporal accuracy without sacrificing the stability demonstrated in \cref{fig:s2_stability}(f).

The amplitude and phase errors \cite{li_NovelImplicitIntegration_2025,li_DesigningDevelopingSinglestep_2023} are computed analytically for the two-sub-step implicit method \eqref{eq:s2_2nd_but} as 
\begin{subequations}\label{eq:s2_amp_phas}
\begin{align}
\delta=& \dfrac{(3-4\xi^{2})\xi\vartheta_1}{6}\omega^{3}\dt^{2}+ \dfrac{(8\xi^4-8\xi^{2}+1)\vartheta_2}{8}\omega^4\dt^{3} +\mathcal{O}(\dt^4) \\  
\epsilon=& \dfrac{(1-4\xi^{2})\sqrt{1-\xi^{2}}\vartheta_1}{6}\omega^{3}\dt^{2} + \dfrac{(2\xi^{2}-1)\xi\sqrt{1-\xi^{2}}\vartheta_2}{2}\omega^4\dt^{3}- \dfrac{\vartheta_3}{20}\omega^5\dt^4+\mathcal{O}(\dt^5)
\end{align}
where $\vartheta_{j}~(j=1,~2,~3) $ are given as 
\begin{align}
\vartheta_1&=6\gamma_1^{2}-6\gamma_1+1\\
\vartheta_2&=(4\gamma_1-1)(2\gamma_1-1)^{2}\\
\vartheta_3&=60\gamma_1^4 - 80\gamma_1^3 + 40\gamma_1^2 - 10\gamma_1 + 1. 
\end{align}
\end{subequations}
For the second-order two-sub-step member, $\vartheta_1\neq0$. In the damped case $(\xi\neq0)$, the leading terms of both the amplitude and phase errors are of order $\mathcal{O}(\dt^{2})$. In the undamped case $(\xi=0)$, the leading $\dt^{2}$ term in $\delta$ vanishes because it is proportional to $\xi$, and hence the amplitude error exhibits one-order superconvergence, i.e., $\delta=\mathcal{O}(\dt^{3})$. By contrast, the corresponding leading term in $\epsilon$ generally remains nonzero at $\xi=0$, so that the phase error remains $\mathcal{O}(\dt^{2})$. For the third-order member, the additional condition $\vartheta_1=0$ eliminates all $\dt^{2}$ terms in the amplitude and phase errors. Consequently, in the damped case, both $\delta$ and $\epsilon$ become $\mathcal{O}(\dt^{3})$. In the undamped case, the amplitude error generally remains $\mathcal{O}(\dt^{3})$, whereas the $\dt^{3}$ term in the phase error vanishes because it is proportional to $\xi$. Therefore, the phase error exhibits one-order superconvergence and becomes $\mathcal{O}(\dt^{4})$. These results indicate that the second-order member possesses enhanced amplitude accuracy for undamped systems, while the third-order member exhibits enhanced phase accuracy under the same condition.

\begin{remark}
    The results given by Eq.~(\ref{eq:s2_amp_phas}a) indicate that the order in amplitude can be improved to third order for the undamped case. Furthermore, when $\gamma_1$ is selected as either $1/2$ or $1/4$, leading to $\vartheta_2=0$, the amplitude error can be further reduced. In fact, these two parameters correspond to the non-dissipative case with $\infrho=1$, meaning that the algorithm introduces no amplitude decay. Consequently, the amplitude error vanishes identically, and its convergence order can be regarded as infinitely high.
\end{remark}
\begin{remark}
    It should be noted that the proposed two-sub-step implicit algorithms, including the third-order member, share the same spectral properties as the existing two-sub-step implicit methods \cite{li_NovelFamilyComposite_2020,ji_UnconditionallyStableTime_2021,kim_ImprovedTwostageImplicit_2024,li_DirectlySelfstartingHigherorder_2022,kim_ImprovedImplicitMethod_2020,wang_TrulySelfstartingComposite_2023,noh_DirectTimeIntegrations_2019,choi_TimeSplittingRatio_2022} with identical effective matrices. Therefore, a comparison of spectral characteristics among these methods is omitted herein, since such a comparison would not provide additional insight. 
\end{remark}

\subsection{Three sub-step: $s=3$}


In this subsection, the three-sub-step algorithm is systematically derived by solving the corresponding algebraic conditions. The influence of the sub-step sizes on the attainable order of accuracy and numerical dissipation is also investigated, providing further insight into the design principles of the proposed family. Moreover, the detailed derivation of the three-sub-step implicit algorithm serves as a representative example for developing higher-order sub-step methods within the proposed framework \eqref{eq:but_alg}. By presenting the procedure for determining algorithmic parameters in detail, the relationship among sub-step sizes, order conditions, and dissipation control can be clearly demonstrated. 

In the case of $s=3$, the proposed framework \eqref{eq:but_alg} reduces to the following three-sub-step formulation:
\begin{equation}\label{eq:s_eq_3}
\begin{NiceArray}{c|c}[cell-space-limits=3pt,columns-width=0.5cm]
\mbf{c} & \mbf{A}\\ \hline & \mbf{b}\T
\end{NiceArray}\quad=\quad\begin{NiceArray}{c|ccc}[cell-space-limits=3pt,columns-width=0.7cm]
\gamma_1 & \alpha_{11}   \\ 
\gamma_2 & \alpha_{21} & \alpha_{22} \\ 
\gamma_3 & \alpha_{31} & \alpha_{32} & \alpha_{33} \\ \hline
& \beta_1 & \beta_2 & \beta_3
\end{NiceArray}.
\end{equation}
The corresponding three-sub-step implicit algorithm can be explicitly written as 
\begin{subequations}
\begin{align}
&\mbfM\mbfa_{n+\gamma_1}+\mbfC\mbfv_{n+\gamma_1}+\mbfK\mbfu_{n+\gamma_1}=\mbfF_{n+\gamma_1} & \mbfM\mbfa_{n+\gamma_2}&+\mbfC\mbfv_{n+\gamma_2}+\mbfK\mbfu_{n+\gamma_2}=\mbfF_{n+\gamma_2}\\
&\quad\mbfv_{n+\gamma_1}=\mbfv_n+\alpha_{11}\dt\mbfa_{n+\gamma_1}&\mbfv_{n+\gamma_2}&=\mbfv_n+\dt \left( \alpha_{21}\mbfa_{n+\gamma_1}+\alpha_{22}\mbfa_{n+\gamma_2} \right) \\
&\quad\mbfu_{n+\gamma_1}=\mbfu_n+\alpha_{11}\dt\mbfv_{n+\gamma_1}& \mbfu_{n+\gamma_2}&=\mbfu_n+\dt \left( \alpha_{21}\mbfv_{n+\gamma_1}+\alpha_{22}\mbfv_{n+\gamma_2} \right) \\
&\mbfM\mbfa_{n+\gamma_3}+\mbfC\mbfv_{n+\gamma_3}+\mbfK\mbfu_{n+\gamma_3}=\mbfF_{n+\gamma_3}\\ 
&\mbfv_{n+\gamma_3}=\mbfv_n+\dt \left( \alpha_{31}\mbfa_{n+\gamma_1}+\alpha_{32}\mbfa_{n+\gamma_2}+\alpha_{33}\mbfa_{n+\gamma_3} \right)&\mbfv_{n+1}&=\mbfv_n+\dt \left( \beta_1\mbfa_{n+\gamma_1}+\beta_2\mbfa_{n+\gamma_2}+\beta_3\mbfa_{n+\gamma_3} \right) \\
&\mbfu_{n+\gamma_3}=\mbfu_n+\dt \left( \alpha_{31}\mbfv_{n+\gamma_1}+\alpha_{32}\mbfv_{n+\gamma_2}+\alpha_{33}\mbfv_{n+\gamma_3} \right)&\mbfu_{n+1}&=\mbfu_n+\dt \left( \beta_1\mbfv_{n+\gamma_1}+\beta_2\mbfv_{n+\gamma_2}+\beta_3\mbfv_{n+\gamma_3} \right) 
\end{align}
\end{subequations}
for second-order transient dynamics, and as 
\begin{subequations}
\begin{align}
\mbfC\mbfv_{n+\gamma_1}&+\mbfK\mbfu_{n+\gamma_1}=\mbfF_{n+\gamma_1} & \mbfC\mbfv_{n+\gamma_2}&+\mbfK\mbfu_{n+\gamma_2}=\mbfF_{n+\gamma_2}\\
\mbfu_{n+\gamma_1}&=\mbfu_n+\alpha_{11}\dt\mbfv_{n+\gamma_1}& \mbfu_{n+\gamma_2}&=\mbfu_n+\dt \left( \alpha_{21}\mbfv_{n+\gamma_1}+\alpha_{22}\mbfv_{n+\gamma_2} \right) \\
\mbfC\mbfv_{n+\gamma_3}&+\mbfK\mbfu_{n+\gamma_3}=\mbfF_{n+\gamma_3} \\
\mbfu_{n+\gamma_3}&=\mbfu_n+\dt \left( \alpha_{31}\mbfv_{n+\gamma_1}+\alpha_{32}\mbfv_{n+\gamma_2}+\alpha_{33}\mbfv_{n+\gamma_3} \right)&
\mbfu_{n+1}&=\mbfu_n+\dt \left( \beta_1\mbfv_{n+\gamma_1}+\beta_2\mbfv_{n+\gamma_2}+\beta_3\mbfv_{n+\gamma_3} \right) 
\end{align}
\end{subequations}
for first-order transient dynamics. 

By enforcing the identical effective matrices and the consistency in each sub-step given by \cref{eq:iem,eq:consistency}, the following algorithmic parameters are obtained:
\begin{subequations}\label{eq:s3_iem_consistency}
\begin{align}
\alpha_{11}&=\alpha_{22}=\alpha_{33}=\gamma_1\\
\alpha_{21}&=\gamma_2-\gamma_1\\
\alpha_{31}&=\gamma_3-\alpha_{32}-\gamma_1.
\end{align}
\end{subequations}
With the above constraints, the numerical amplification factor and numerical load operator are derived from \cref{eq:num_DL} as 
\begin{subequations}
\begin{align}
    D_\mathsf{num}&=\dfrac{\left\{\begin{aligned}
    &1+\bigg[\gamma_1^3-(\beta_1+2\beta_2+2\beta_3)\gamma_1^2+(\beta_3\gamma_3+\beta_3\alpha_{32}+\beta_2\gamma_2)\gamma_1-\beta_3\alpha_{32}\gamma_2\bigg]\omega^3\dt^3\\
    &-\bigg[(2\beta_1+3\beta_2+3\beta_3)\gamma_1-\beta_2\gamma_2-\beta_3\gamma_3-3\gamma_1^2\bigg]\omega^2\dt^2-(\beta_1+\beta_2+\beta_3-3\gamma_1)\omega\dt
    \end{aligned}\right\}}{(1+\gamma_1\omega\dt)^3}\\
    L_\mathsf{num}&=\dfrac{\beta_1+\left(\gamma_{1}^{2}-\beta_{2} \gamma_{1} \gamma_{2} +\beta_{3} \alpha_{32} \gamma_{2} -\beta_{3} \gamma_{1} \gamma_{3}\right) \omega^{2}\dt^{2}+  \left(2 \beta_{1} \gamma_{1} +\beta_{2} \gamma_{1} -\beta_{2} \gamma_{2} -\alpha_{31} \beta_{3}\right)\omega\dt
}{(1+\gamma_1\omega\dt)^{3}}\exp(\gamma_1\dt)\dt\notag\\
&\quad +\dfrac{\big[\beta_2+(\beta_2\gamma_1-\beta_3\alpha_{32})\omega\dt\big]\exp(\gamma_2\dt)+(1+\gamma_1\omega\dt)\beta_3\exp(\gamma_3\dt)}{(1+\gamma_1\omega\dt)^{2}}\dt.
\end{align}
\end{subequations}
Substituting the above expressions into \cref{eq:def_acc} and thus the local truncation errors are obtained as 
\begin{subequations}
\begin{align}
D_\mathsf{num}-D_\mathsf{exa}&=-c_1\omega\dt+c_2\omega^{2}\dt^{2}-c_3\omega^{3}\dt^{3}+\mathcal{O}(\dt^4)\\
L_\mathsf{num}-L_\mathsf{exa}&=c_1\dt-c_2(\omega-1)\dt^{2}+\left[c_3(\omega-1)\omega+\dfrac{1}{2}c_4\right]\dt^{3}+\mathcal{O}(\dt^4)
\end{align}
where the coefficients $c_j~(j=1,~2,~3,~4)$ are defined as 
\begin{align}
c_1 &=\beta_1+\beta_2+\beta_3-1&
c_2 &=\beta_1\gamma_1+\beta_2\gamma_2+\beta_3\gamma_3-\frac12\\
c_3 &= (\beta_1-\beta_2-\beta_3)\gamma_1^{2}+(2\gamma_3-\alpha_{32})\beta_3\gamma_1+2\beta_2\gamma_2\gamma_1+\beta_3\alpha_{32}\gamma_2-\dfrac{1}{6}&
c_4 &= \beta_1\gamma_1^{2}+\beta_2\gamma_2^{2}+\beta_3\gamma_3^{2} -\dfrac{1}{3}.  
\end{align}
\end{subequations}

Therefore, the three-sub-step implicit method \eqref{eq:s_eq_3} achieves third-order accuracy, i.e., $D_\mathsf{num}-D_\mathsf{exa}=L_\mathsf{num}-L_\mathsf{exa}=\mathcal{O}(\dt^4)$, if and only if the conditions $c_1=c_2=c_3=c_4=0$ are simultaneously satisfied. Solving these conditions yields 
\begin{subequations}\label{eq:s3_3rd}
\begin{align}
\alpha_{32}&=\dfrac{(6\gamma_1^{2}-6\gamma_1+1)(\gamma_1-\gamma_3)(\gamma_2-\gamma_3)}{(6\gamma_1\gamma_2-3\gamma_1-3\gamma_2+2)(\gamma_2-\gamma_1)}&
\beta_1&= \dfrac{6\gamma_2\gamma_3-3\gamma_2-3\gamma_3+2}{6(\gamma_1-\gamma_2)(\gamma_1-\gamma_3)} \\
\beta_2&= \dfrac{6\gamma_1\gamma_3-3\gamma_1-3\gamma_3+2}{6(\gamma_2-\gamma_1)(\gamma_2-\gamma_3)} &
\beta_3&= \dfrac{6\gamma_1\gamma_2-3\gamma_1-3\gamma_2+2}{6(\gamma_3-\gamma_1)(\gamma_3-\gamma_2)}.
\end{align}
\end{subequations}
By substituting \cref{eq:s3_iem_consistency,eq:s3_3rd} into the stability polynomial, it can be simplified as $E(y)=y^4\left[\left(216\gamma_1^5-432\gamma_1^4+336\gamma_1^{3}-117\gamma_1^{2}+18\gamma_1-1\right)y^{2}-72\gamma_1^{3}+108\gamma_1^{2}-36\gamma_1+3\right]\ge 0$, which leads to the unconditionally stable condition of the three-sub-step implicit method as 
\begin{equation}\label{eq:s3_uc}
    \dfrac{1}{3}\le\gamma_1\le \dfrac{1}{2}+\dfrac{\sqrt{ 3 } }{3}\cos\left(\dfrac{\pi}{18} \right). 
\end{equation}

Furthermore, the high-frequency dissipation control given by \cref{eq:D_num_infty} can be simplified as 
\begin{equation}\label{eq:s3_D_num_infty}
D_\mathsf{num}^\infty= \dfrac{6\gamma_1^{3}-18\gamma_1^{2}+9\gamma_1-1}{6\gamma_1^{3}}=\left|\infrho\right| 
\end{equation}
and therefore the first sub-step size $\gamma_1$ is determined by 
\begin{equation}\label{eq:s3_r1_set}
\gamma_1\in\left\{\gamma_1~\bigg|~\dfrac{1}{3}\le\gamma_1\le \dfrac{1}{2}+\dfrac{\sqrt{ 3 } }{3}\cos\left(\dfrac{\pi}{18} \right),~6(1-\infrho)\gamma_1^{3}-18\gamma_1^{2}+9\gamma_1-1=0\right\}
\end{equation}
where $\infrho\in[1-\sqrt{3},~1]$. \cref{eq:s3_r1_set} provides the relationship between the first sub-step size $\gamma_1$ and the user-specified parameter $\infrho$. For a prescribed value of $\infrho$, the cubic equation, $6(1-\infrho)\gamma_1^{3}-18\gamma_1^{2}+9\gamma_1-1=0$, yields three candidate solutions for $\gamma_1$. Since not all mathematical solutions necessarily lead to desirable numerical properties, the admissible solution is determined by further imposing unconditional stability given by \cref{eq:s3_uc}. Consequently, a unique relationship between the first sub-step size $\gamma_1$ and the parameter $\infrho$ can be established.

\begin{figure}[htbp]
	\centering 
	\includegraphics[scale=1.0]{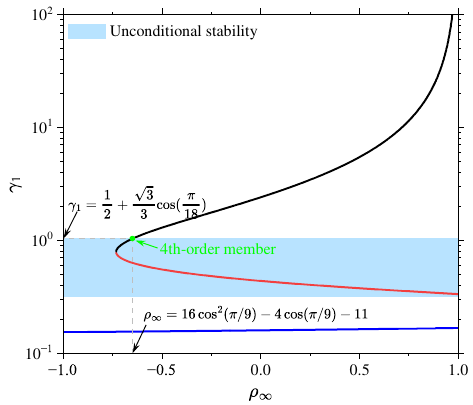}
	\caption{The variations of $\gamma_1$ with $\infrho$ for the three-sub-step implicit method.}
	\label{fig:s3_r1_rho}
\end{figure}

\cref{fig:s3_r1_rho} illustrates the variations of the first sub-step size $\gamma_1$ with respect to the user-specified parameter $\infrho$ for the three-sub-step implicit method. According to \cref{eq:s3_r1_set}, a prescribed value of $\infrho$ may lead to multiple candidate solutions of $\gamma_1$ from the cubic equation. However, only the solutions satisfying unconditional stability given by \cref{eq:s3_uc} are admissible. Therefore, the stable region shown in \cref{fig:s3_r1_rho} provides an effective criterion for selecting the appropriate branch of $\gamma_1$. It can be observed that the admissible solutions of $\gamma_1$ are mainly distributed within the unconditionally stable region, whereas the other mathematical branches may violate the stability requirement. Consequently, the relationship between $\gamma_1$ and $\infrho$ is uniquely determined by selecting the stable solution branch. Similar to the two-sub-step case, this relationship enables the user to prescribe the desired high-frequency dissipation through $\infrho$, while the first sub-step size $\gamma_1$ is subsequently determined from the corresponding algebraic condition.

The proposed three-sub-step implicit method \eqref{eq:s_eq_3} can be further represented by the following simplified Butcher tableau:
\begin{equation}\label{eq:s3_3rd_but}
\begin{NiceArray}{c|ccc}[cell-space-limits=3pt,columns-width=0.7cm]
\gamma_1 & \gamma_1 & 0  \\ 
\gamma_2 & \gamma_2-\gamma_1 & \gamma_1 \\ 
\gamma_3 & \gamma_3-\dfrac{(6\gamma_1^{2}-6\gamma_1+1)(\gamma_1-\gamma_3)(\gamma_2-\gamma_3)}{(6\gamma_1\gamma_2-3\gamma_1-3\gamma_2+2)(\gamma_2-\gamma_1)}-\gamma_1 & \dfrac{(6\gamma_1^{2}-6\gamma_1+1)(\gamma_1-\gamma_3)(\gamma_2-\gamma_3)}{(6\gamma_1\gamma_2-3\gamma_1-3\gamma_2+2)(\gamma_2-\gamma_1)} & \gamma_1 \\ \hline
& \dfrac{6\gamma_2\gamma_3-3\gamma_2-3\gamma_3+2}{6(\gamma_1-\gamma_2)(\gamma_1-\gamma_3)} & \dfrac{6\gamma_1\gamma_3-3\gamma_1-3\gamma_3+2}{6(\gamma_2-\gamma_1)(\gamma_2-\gamma_3)} & \dfrac{6\gamma_1\gamma_2-3\gamma_1-3\gamma_2+2}{6(\gamma_3-\gamma_1)(\gamma_3-\gamma_2)}
\end{NiceArray}
\end{equation}
where the sub-step sizes $\gamma_j~(j=1,~2,~3)$ are free parameters. Of course, the first sub-step size $\gamma_1$ can also be determined by the user-specified parameter $\infrho\in[1-\sqrt{3},~1]$ through \cref{eq:s3_r1_set}, whereas the remaining sub-step sizes $\gamma_2$ and $\gamma_3$ serve as free parameters. These two parameters are independent of the spectral properties and stability of the resulting algorithm.

\begin{figure}[htbp]
	\centering 
	\subfigure[$\gamma_1=1/3~(\infrho=1.0)$]{
		\includegraphics[scale=0.407]{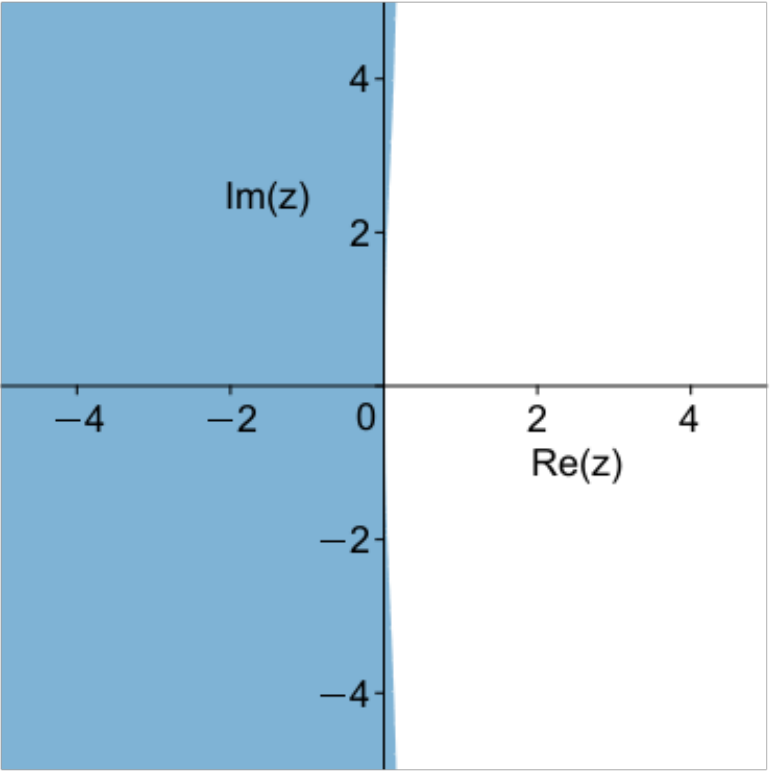}}
	\subfigure[$\gamma_1=0.37560~(\infrho=0.5)$]{
		\includegraphics[scale=0.4]{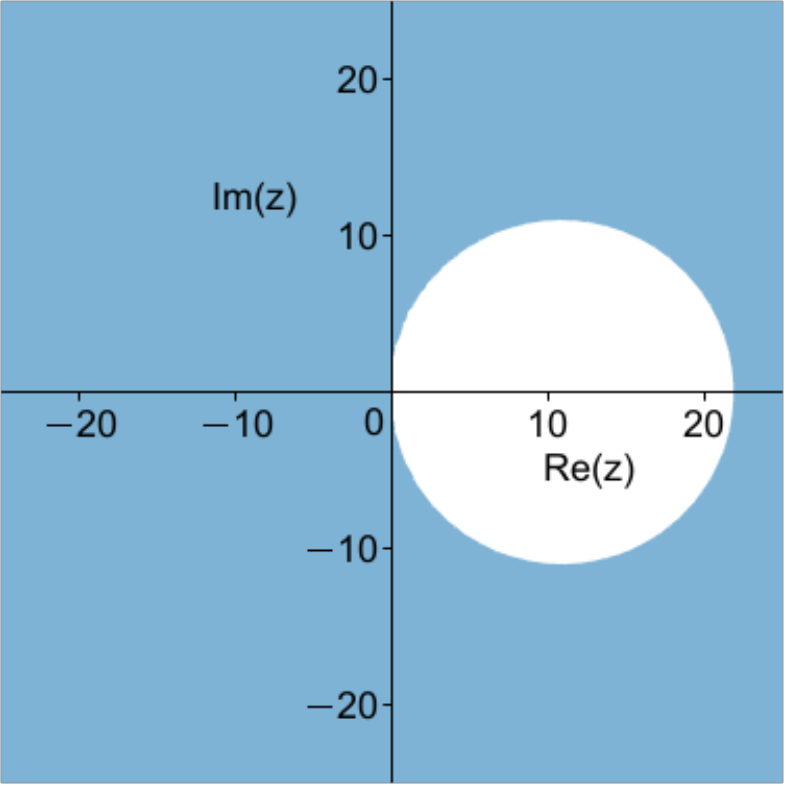}}
	\subfigure[$\gamma_1=0.43587~(\infrho=0.0)$]{
		\includegraphics[scale=0.4]{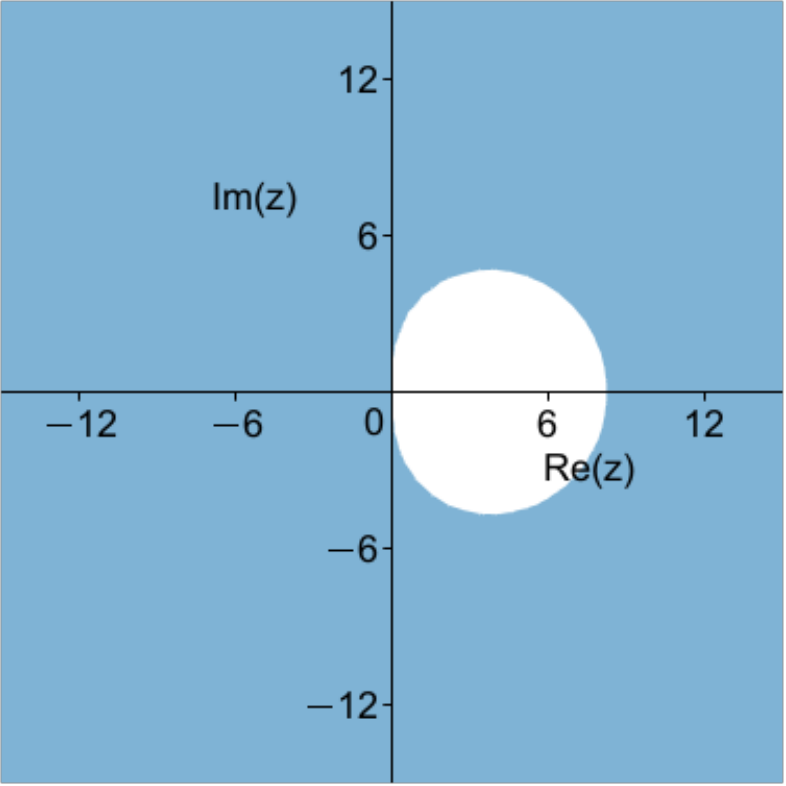}}
	\subfigure[$\gamma_1=0.55090~(\infrho=-0.5)$]{
		\includegraphics[scale=0.4]{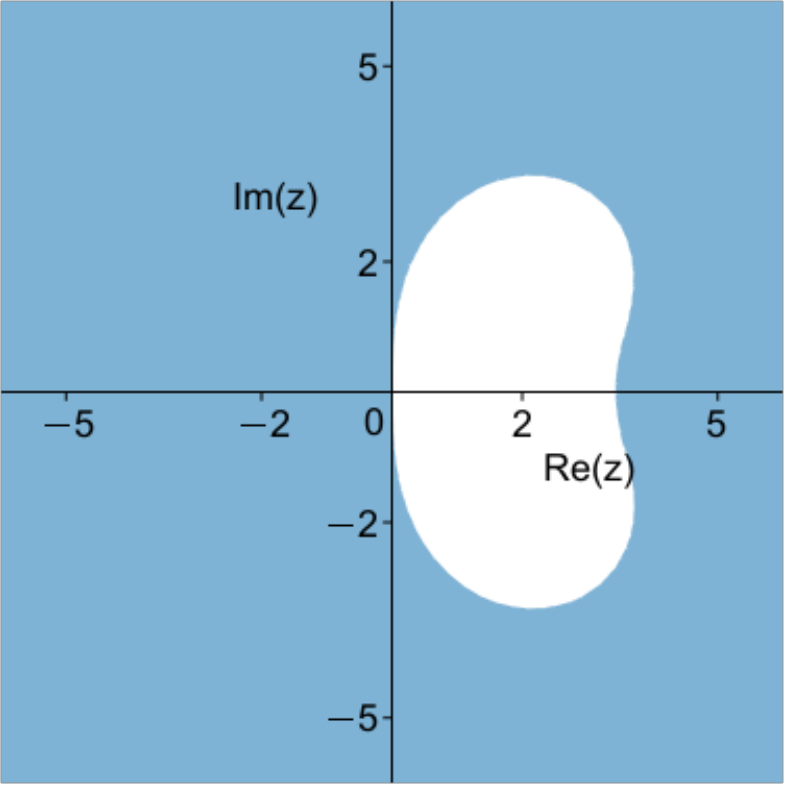}}
	\subfigure[$\gamma_1=1/2+\sqrt{3}/6~(\infrho=1-\sqrt{3})$]{
		\includegraphics[scale=0.4]{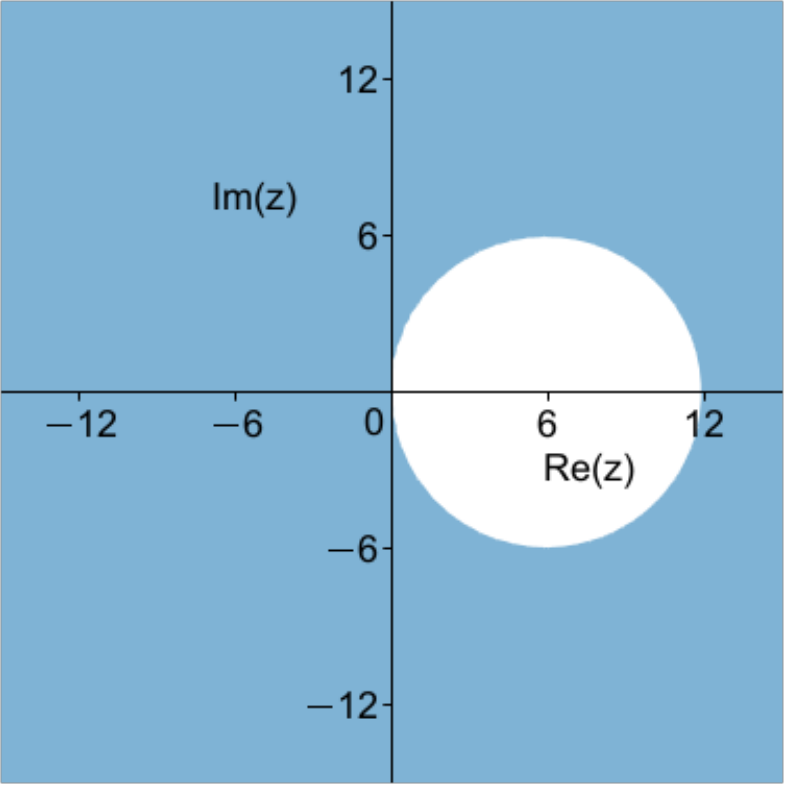}}
	\subfigure[$\gamma_1=1/2+\sqrt{3}\cos(\pi/18)/3~(\infrho\approx-0.63041)$]{
		\includegraphics[scale=0.4]{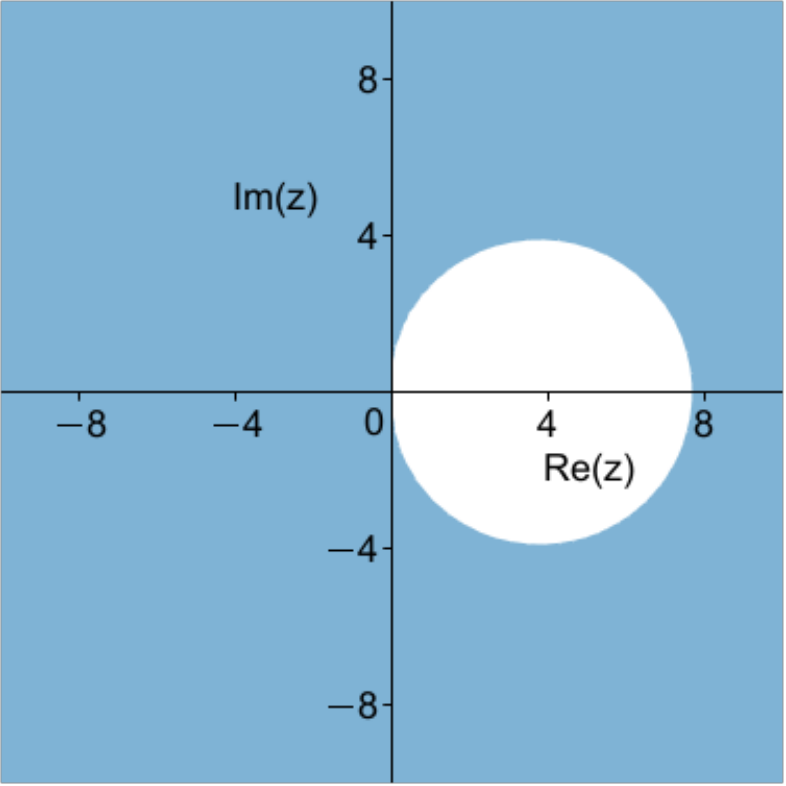}}
	\caption{The stability region for the three-sub-step implicit method defined by \cref{eq:s3_3rd_but}.}
	\label{fig:s3_stability}
\end{figure}

\cref{fig:s3_stability} presents the stability regions of the three-sub-step implicit method defined by \cref{eq:s3_3rd_but} for different admissible values of $\gamma_1$, corresponding to various prescribed high-frequency dissipation $\infrho$. These results are provided to verify that the introduction of the third sub-step and the higher-order accuracy do not compromise unconditional stability established in \cref{eq:s3_uc}. As observed from \cref{fig:s3_stability}, all considered cases include the entire left half-plane of the complex plane within their stability regions, confirming that the proposed three-sub-step implicit algorithms remain unconditionally stable for all admissible values of $\gamma_1$. It should be noted that, unlike the lower-order two-sub-step methods, some stability regions of the three-sub-step algorithms extend across the imaginary axis and partially enter the right half-plane, particularly for the case with $\infrho=1$. Such a phenomenon does not violate the A-stability requirement, since A-stability only requires the complete inclusion of the left half-plane in the stability region. Instead, it reflects the richer stability characteristics introduced by the higher-order construction of the present three-sub-step method.

For the non-dissipative case with $\infrho=1$, the limiting spectral radius satisfies $\left|R(\infty)\right|=1$, indicating that no numerical damping is introduced in the high-frequency range. The corresponding stability region exhibits a large extension along both the real and imaginary directions. Overall, \cref{fig:s3_stability} verifies that the proposed three-sub-step implicit method \eqref{eq:s3_3rd_but} successfully achieves a favorable combination of third-order accuracy, controllable high-frequency dissipation, and unconditional stability.

\subsubsection{Higher-order accuracy}

It is evident that the proposed three-sub-step implicit method \eqref{eq:s3_3rd_but} still contains two free parameters, $\gamma_2$ and $\gamma_3$. Although these parameters do not influence spectral properties or linear stability of the algorithm, their appropriate selections may provide additional flexibility in improving the accuracy.

The local truncation errors associated with $D_\mathsf{num}$ and $L_\mathsf{num}$ are further computed as 
\begin{subequations}
\begin{align}
D_\mathsf{num}-D_\mathsf{exa}&=\dfrac{24\gamma_1^{3}-36\gamma_1^{2}+12\gamma_1-1}{24}\omega^4\dt^4 +\mathcal{O}(\dt^5)\\
L_\mathsf{num}-L_\mathsf{exa}&=\bigg[\dfrac{24\gamma_1^{3}-36\gamma_1^{2}+12\gamma_1-1}{24}(1-\omega)\omega^{2}-\dfrac{1}{2}\left( \gamma_2-\dfrac{1}{2} \right)\left( \gamma_1^{2}-\gamma_1+\dfrac{1}{6} \right)\omega\notag \\
&\qquad + \dfrac{(6\gamma_2\gamma_3-3\gamma_2-3\gamma_3+2)\gamma_1}{36}-\dfrac{(3\gamma_3-2)\gamma_2}{36}+\dfrac{\gamma_3}{18}-\dfrac{1}{24}  \bigg]\dt^4+\mathcal{O}(\dt^5). 
\end{align}
\end{subequations}
Therefore, by setting all coefficients associated with the $\dt^4$ terms to zero, the conditions required for the method to achieve fourth-order accuracy can be obtained as
\begin{subequations}\label{eq:s3_r123}
\begin{align}
\gamma_1&=\dfrac{1}{2}+\dfrac{\sqrt{ 3 } }{3}\cos\left(\dfrac{\pi}{18}\right)\\ 
\gamma_2&=\dfrac{1}{2}\label{eq:s3_r2}\\ 
\gamma_3&=1-\gamma_1\label{eq:s3_r3}. 
\end{align}
\end{subequations}
Substituting \cref{eq:s3_r123} into \cref{eq:s3_3rd_but} gives the Butcher tableau of the present three-sub-step implicit method with fourth-order accuracy as 
\begin{equation}\label{eq:s3_4th_but}
\begin{NiceArray}{c|ccc}[cell-space-limits=3pt,columns-width=3cm]
\dfrac{1}{2}+\dfrac{\sqrt{3}}{3}\cos\left( \dfrac{\pi}{18} \right) & \dfrac{1}{2}+\dfrac{\sqrt{3}}{3}\cos\left( \dfrac{\pi}{18} \right) &                                                                         &                                                                         \\ \rule{0pt}{20pt}
		\dfrac{1}{2}                                                            & -\dfrac{{\sqrt 3 }}{3}\cos \left( {\dfrac{\pi }{{18}}} \right)          & \dfrac{1}{2}+\dfrac{\sqrt{3}}{3}\cos\left( \dfrac{\pi}{18} \right) &                                                                         \\\rule{0pt}{20pt}
		\dfrac{1}{2}-\dfrac{\sqrt{3}}{3}\cos\left( \dfrac{\pi}{18} \right) & 1+\dfrac{2\sqrt{3}}{3}\cos\left( \dfrac{\pi}{18} \right)           & -1-\dfrac{4\sqrt{3}}{3}\cos\left( \dfrac{\pi}{18} \right)          & \dfrac{1}{2}+\dfrac{\sqrt{3}}{3}\cos\left( \dfrac{\pi}{18} \right) \\[6pt]
		\hline \rule{0pt}{15pt}
		& \dfrac{1}{{8{{\cos }^2}\left( {\dfrac{\pi }{{18}}} \right)}}            & 1-\dfrac{1}{{4{{\cos }^2}\left( {\dfrac{\pi }{{18}}} \right)}}          & \dfrac{1}{{8{{\cos }^2}\left( {\dfrac{\pi }{{18}}} \right)}}
\end{NiceArray}.
\end{equation}
Once the value of $\gamma_1$ is given as (\ref{eq:s3_r123}a), the high-frequency spectral radius is computed by \cref{eq:s3_D_num_infty} as 
\begin{equation}\label{eq:s3_4th_rho_infty}
\infrho=16\cos^{2}\left(\dfrac{\pi}{9}\right)-4\cos\left(\dfrac{\pi}{9}\right)-11\approx-0.63041. 
\end{equation}
Furthermore, the special case shown in \cref{fig:s3_stability}(f), corresponding to the fourth-order member, also preserves unconditional stability. This demonstrates that exploiting the remaining algorithmic parameters to achieve higher-order accuracy does not deteriorate the fundamental stability of the three-sub-step method.

\begin{remark}
For the third-order three-sub-step implicit method with controllable numerical dissipation defined by \cref{eq:s3_3rd_but}, the first sub-step size $\gamma_1$ is determined by the user-specified parameter $\infrho$ through \cref{eq:s3_r1_set}, whereas the remaining sub-step sizes $\gamma_2$ and $\gamma_3$ remain as free parameters within the third-order framework. By further imposing an additional order condition, the sub-step sizes $\gamma_1$, $\gamma_2$, and $\gamma_3$ can be completely determined, leading to a fourth-order accurate three-sub-step implicit method described in \cref{eq:s3_4th_but}. However, in this case, the parameter $\gamma_1$ is no longer adjustable and therefore the high-frequency numerical dissipation cannot be freely controlled through $\infrho$. In the present study, the selection rules of $\gamma_2$ and $\gamma_3$ derived from the fourth-order member, namely \cref{eq:s3_r2,eq:s3_r3}, are adopted to determine these two parameters, so that the third- and fourth-order algorithms share the same sub-step splitting strategy. Consequently, the three-sub-step implicit framework provides two possible algorithm configurations: for arbitrary user-specified $\infrho\in[1-\sqrt{3},~1]$, the method achieves third-order accuracy with controllable high-frequency dissipation; whereas for the specific value of $\infrho$ given in \cref{eq:s3_4th_rho_infty}, it achieves fourth-order accuracy with a fixed high-frequency dissipation level ($\left|\infrho\right|\approx0.63041$).
\end{remark}

\subsubsection{Spectral properties}

\begin{figure}[htbp]
	\centering 
	\includegraphics[scale=1.8]{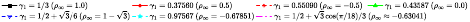}\\
    \subfigure[]{
		\includegraphics[scale=0.72]{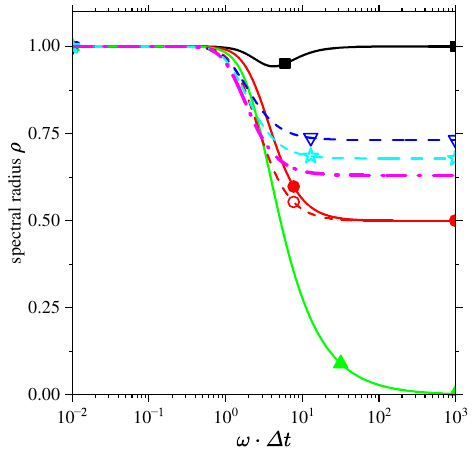}}
    \subfigure[]{
		\includegraphics[scale=0.72]{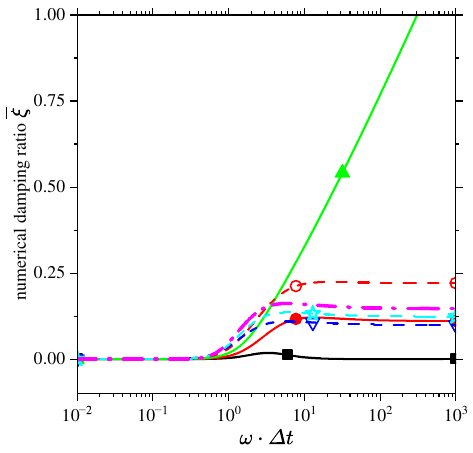}}\\
    \subfigure[]{
		\includegraphics[scale=0.72]{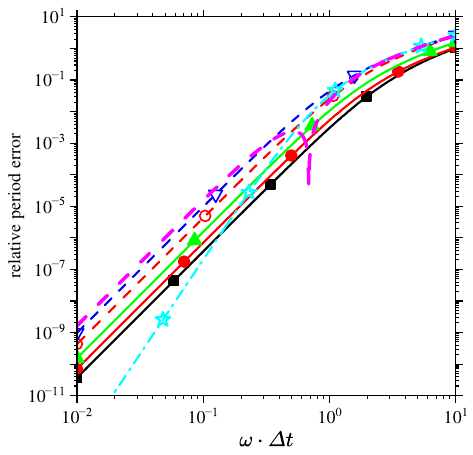}}
    \subfigure[]{
		\includegraphics[scale=0.72]{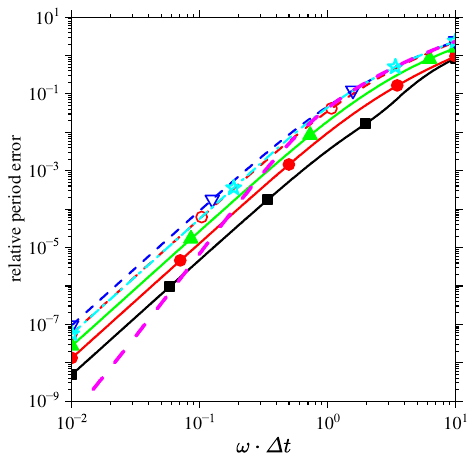}}
	\caption{Spectral properties for the three-sub-step implicit method defined by \cref{eq:s3_3rd_but}: (a-c) $\xi=0.0$ and (d) $\xi=0.35$.}
	\label{fig:s3_sp}
\end{figure}

\cref{fig:s3_sp} further presents the spectral properties of the proposed three-sub-step implicit method defined by \cref{eq:s3_3rd_but}. The selected values of $\gamma_1$ correspond to different prescribed high-frequency spectral radii as determined from \cref{eq:s3_D_num_infty}, together with the fourth-order member satisfying \cref{eq:s3_4th_rho_infty}.

\cref{fig:s3_sp}(a) and (b) illustrate the variation of the spectral radius and the numerical damping ratio with respect to the normalized frequency. As expected, decreasing $\infrho$ generally results in stronger high-frequency dissipation. However, an interesting observation can be made for the case $\infrho=1$. Although this parameter corresponds to the non-dissipative limit at infinite frequency, the spectral radius is not identically equal to unity over the entire frequency range. Instead, a slight reduction of the spectral radius occurs in the mid-frequency range before approaching unity again as $\omega\dt\to\infty$. Consequently, the proposed method still introduces a small amount of numerical damping for mid-frequencies even when $\infrho=1$. This observation is fully consistent with the stability region shown in \cref{fig:s3_stability}(a), where the stability boundary is not exactly coincident with the imaginary axis and therefore the stability region is not the entire left half-plane.

\cref{fig:s3_sp}(c) and (d) present relative period errors for the undamped and damped cases, respectively. In general, as the parameter $\gamma_1$ increases from $1/3$ to $1/2+\sqrt{3}\cos(\pi/18)/3$, the relative period error tends to increase. Nevertheless, \cref{fig:s3_sp}(c) and (d) clearly reveal the existence of two distinct optimal values of $\gamma_1$, each minimizing the relative period error under the corresponding damping condition. Specifically, for the undamped case, the third-order member with $\gamma_1=0.97567$ (corresponding to $\infrho=-0.67851$) produces even smaller relative period errors than the fourth-order member over the low-frequency range. In contrast, the superiority of the fourth-order member becomes evident when the physical damping ratio is present, as shown in \cref{fig:s3_sp}(d). These observations are fully consistent with the analytical amplitude and phase errors derived in \cref{eq:s3_amp_phas}, and can be explained by the leading errors discussed in \cref{rem:s3_1,rem:s3_2}.

The amplitude and phase errors \cite{li_NovelImplicitIntegration_2025,li_DesigningDevelopingSinglestep_2023} are computed analytically for the proposed three-sub-step implicit method \eqref{eq:s3_3rd_but} as 
\begin{subequations}\label{eq:s3_amp_phas}
\begin{align}
\delta=&-\dfrac{(8\xi^4-8\xi^{2}+1)\phi_1}{24}\omega^4\dt^3+ \dfrac{\xi(16\xi^4-20\xi^{2}+5)\phi_2}{30}\omega^5\dt^4 + \dfrac{\phi_3}{72}\omega^6\dt^5+\mathcal{O}(\dt^6)\\  
\epsilon=&- \dfrac{\xi(2\xi^{2}-1)\sqrt{1-\xi^{2}}\phi_1}{6}\omega^4\dt^{3}+ \dfrac{(16\xi^4-12\xi^{2}+1)\sqrt{1-\xi^{2}}\phi_2}{30}\omega^5\dt^4\notag\\
&+ \dfrac{(16\xi^{2}-16\xi^4-3)\xi\sqrt{1-\xi^{2}}\phi_3}{36}\omega^6\dt^5 -\dfrac{\phi_4}{252}\omega^7\dt^6  +\mathcal{O}(\dt^7),
\end{align}
where $\phi_j~(j=1,~2,~3,~4) $ are given as 
\begin{align}
\phi_1&=24\gamma_1^{3}-36\gamma_1^{2}+12\gamma_1-1\\
\phi_2&=90\gamma_1^4-150\gamma_1^{3}+75\gamma_1^{2}-15\gamma_1+1\\
\phi_3&=432\gamma_1^5-756\gamma_1^4+444\gamma_1^{3}-126\gamma_1^{2}+18\gamma_1-1\\
\phi_4&=2520\gamma_1^6 - 4536\gamma_1^5 + 2898\gamma_1^4 - 966\gamma_1^3 + 189\gamma_1^2 - 21\gamma_1 + 1.
\end{align}
\end{subequations}

\begin{remark}\label{rem:s3_1}
For the third-order three-sub-step member, the coefficient $\phi_1$ does not vanish. In the undamped case $(\xi=0)$, the leading term of the amplitude error $\delta$ remains of order $\mathcal{O}(\dt^{3})$, whereas the corresponding $\dt^{3}$ term in the phase error $\epsilon$ vanishes because it is proportional to $\xi$. Consequently, the phase error exhibits a one-order superconvergence and becomes $\mathcal{O}(\dt^{4})$. In contrast, in the damped case $(\xi\neq0)$, the leading terms in both $\delta$ and $\epsilon$ are generally of order $\mathcal{O}(\dt^{3})$. On the other hand, for the fourth-order member, the first sub-step size $\gamma_1$ satisfies $\phi_1=0$ and hence all $\dt^{3}$ terms in \cref{eq:s3_amp_phas} vanish. The leading amplitude and phase errors in the damped case are therefore both of order $\mathcal{O}(\dt^{4})$. In the undamped case, however, the $\dt^{4}$ term in the amplitude error also vanishes because it is proportional to $\xi$, so that $\delta=\mathcal{O}(\dt^{5})$, while the phase error generally remains $\epsilon=\mathcal{O}(\dt^{4})$. Therefore, the fourth-order algorithm exhibits an additional one-order superconvergence in amplitude for undamped systems, whereas its phase accuracy remains fourth order. These analytical observations explain the different relative period error of the damped and undamped cases shown in \cref{fig:s3_sp}(c) and~(d).
\end{remark}
\begin{remark}\label{rem:s3_2}
According to Eq.~(\ref{eq:s3_amp_phas}b), the leading phase error of the third-order implicit method in the undamped case can be eliminated by imposing $\phi_2=0 $. Among the real solutions, the admissible value $\gamma_1=0.9756745886944403$ satisfies the unconditional stability \eqref{eq:s3_uc}. Substituting this value into the high-frequency dissipation relation \eqref{eq:s3_D_num_infty} gives $\infrho\approx-0.67851$. Consequently, the corresponding third-order member exhibits sixth-order phase accuracy $\mathcal{O}(\dt^6)$ for solving undamped problems while retaining unconditional stability and a fixed high-frequency dissipation level. Owing to these favorable properties, $\gamma_1=0.9756745886944403$ is recommended as the default parameter for the third-order three-sub-step implicit method.
\end{remark}

\subsubsection{Comparisons with existing three-sub-step implicit methods}

To ensure a fair comparison, only three-sub-step implicit methods with the same computational complexity as the proposed algorithms are considered in this subsection. Among the available methods, the $\infrho$-MSSBN3 method \cite{ji_AccurateControllablyDissipative_2023} and the Frutos--Serna method \cite{defrutos_EasilyImplementableFourthorder_1992} are selected as representative algorithms. All compared methods employ three implicit stages and therefore require the same number of sub-step equilibrium evaluations per time step. Moreover, the proposed methods and the $\infrho$-MSSBN3 method share the same sub-step splitting strategy, allowing the improvement achieved by the present higher-order construction to be evaluated independently of the time partition. While the $\infrho$-MSSBN3 method is limited to second-order accuracy, the proposed methods systematically construct both third- and fourth-order accurate members. The Frutos--Serna method, which also achieves fourth-order accuracy, is included as a representative high-order three-sub-step implicit scheme for comparison. It should be noted, however, that the Frutos--Serna method does not satisfy the identical effective matrices. Consequently, three different effective matrices must be assembled and factorized when solving linear dynamics, whereas the proposed methods preserve identical effective matrices over all three sub-steps, so that only one matrix assembly and factorization are required.

\begin{figure}[htbp]
	\centering 
	\includegraphics[scale=1.4]{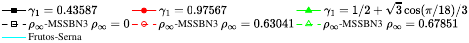}\\
    \subfigure[]{
		\includegraphics[scale=0.72]{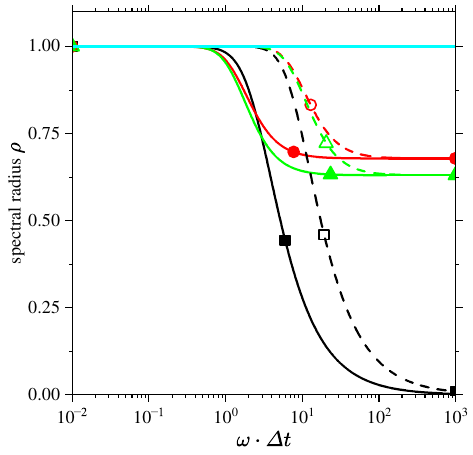}}
    \subfigure[]{
		\includegraphics[scale=0.72]{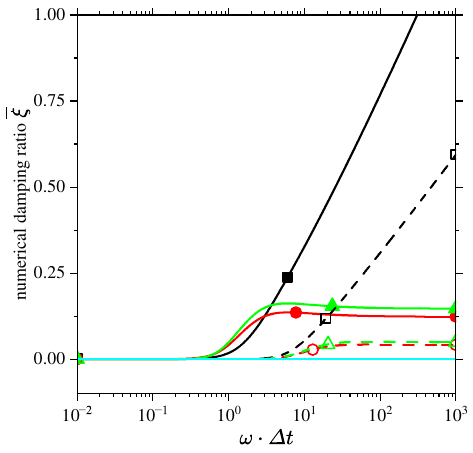}}
    \subfigure[]{
		\includegraphics[scale=0.72]{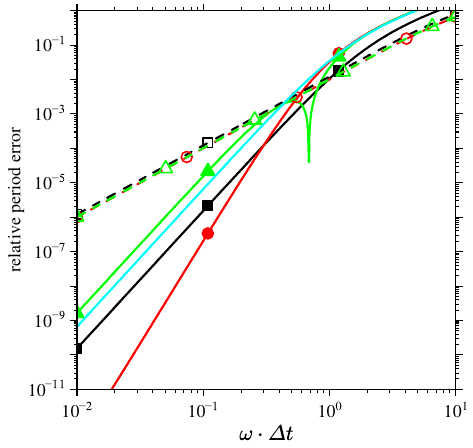}}\\
    \subfigure[]{
		\includegraphics[scale=0.72]{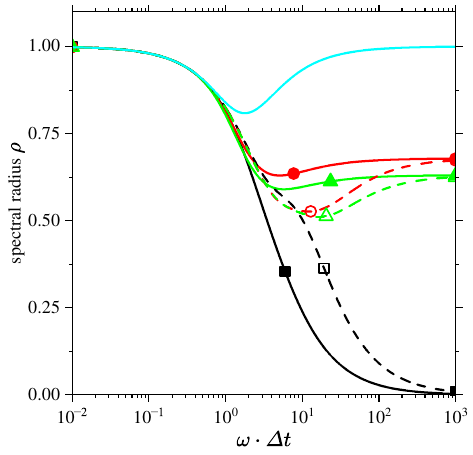}}
    \subfigure[]{
		\includegraphics[scale=0.72]{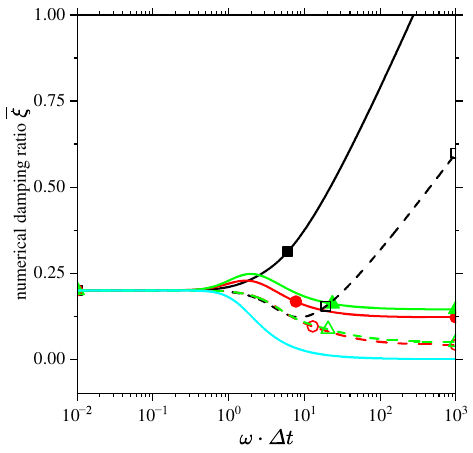}}
    \subfigure[]{
		\includegraphics[scale=0.72]{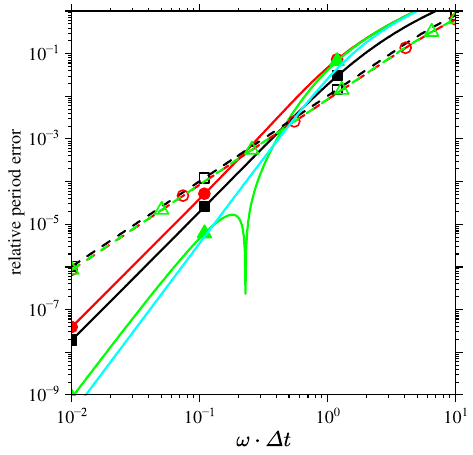}}\\
	\caption{Comparisons of spectral properties for some three-sub-step implicit methods: (a-c) $\xi=0.0$ and (d-f) $\xi=0.2$.}
	\label{fig:s3_com_sp}
\end{figure}

\cref{fig:s3_com_sp} compares the spectral properties of the proposed methods with those of the $\infrho$-MSSBN3 method and the Frutos--Serna method in both undamped ($\xi=0$) and damped ($\xi=0.2$) cases. Three representative members of the proposed family are considered, namely the most-dissipative third-order member ($\gamma_1=0.43587$), the phase-optimized third-order member ($\gamma_1=0.97567$), and the fourth-order member ($\gamma_1=1/2+\sqrt{3}\cos(\pi/18)/3$). \cref{fig:s3_com_sp}(a) and (d) compare the spectral radii. The proposed methods consistently provide stronger dissipation in the mid- and high-frequency ranges. In contrast, the Frutos--Serna method is always non-dissipative over the entire frequency range. The corresponding numerical damping ratios shown in \cref{fig:s3_com_sp}(b) and (e) further confirm these observations. The proposed methods exhibit controllable damping levels determined by the prescribed value of $\infrho$, while the $\infrho$-MSSBN3 method produces similar but generally smaller damping ratios under the same nominal value of $\infrho$.

\cref{fig:s3_com_sp}(c) and (f) compare the relative period errors. As expected, the proposed third- and fourth-order algorithms substantially outperform the second-order $\infrho$-MSSBN3 method throughout the low- and mid-frequency ranges, demonstrating that the present framework significantly improves the phase accuracy without increasing the computational complexity. Among the proposed members, the phase-optimized third-order algorithm ($\gamma_1=0.97567$) exhibits the smallest period error in the undamped case over most of the frequency range considered, even outperforming the fourth-order member. The advantage of the fourth-order member becomes more evident in the damped case, yielding the smallest period errors among all compared methods. The Frutos--Serna method also exhibits excellent phase accuracy owing to its fourth-order accuracy, and its period-error curves are generally comparable to those of the proposed fourth-order member.

\subsection{More sub-steps: $s\ge4$}
Following the development procedure established for the one-, two-, and three-sub-step implicit methods, the proposed framework \eqref{eq:but_alg} can be straightforwardly extended to construct four-, five-, and six-sub-step implicit algorithms with higher-order accuracy. The derivation follows exactly the same sequence of analyses presented in the preceding subsections, including the enforcement of the identical effective matrices \eqref{eq:iem}, the order of accuracy \eqref{eq:def_acc}, the consistency \eqref{eq:consistency}, the unconditional stability \eqref{eq:Ey}, and the controllable high-frequency dissipation \eqref{eq:D_num_infty}. Since the derivation procedure remains unchanged whereas the resulting algebraic expressions become increasingly lengthy and cumbersome as the number of sub-steps increases, the detailed derivations are omitted here for brevity. Instead, the derived algorithmic parameters for the four-, five-, and six-sub-step implicit methods are summarized in \cref{app:a}.

Similar to the cases of $s=1$, $2$, and $3$, the proposed $s$-sub-step method ($s=4,~5,$ and $6$) first constructs $s$th-order accurate implicit members while preserving the identical effective matrices. In this configuration, the sub-step sizes $\gamma_j~(j=1,~\cdots,~s)$ remain undetermined. The first sub-step size $\gamma_1$ can be uniquely associated with the user-specified parameter $\infrho$ through the high-frequency dissipation condition \eqref{eq:D_num_infty}, thereby providing user-specified numerical dissipation. The remaining sub-step sizes $\gamma_j~(j=2,~\cdots,~s)$ are treated as free parameters and do not affect the spectral properties or the unconditional stability of the resulting algorithms.

Furthermore, by imposing one additional order condition, all sub-step sizes are completely consumed, leading to a unique sub-step splitting strategy and consequently an $(s+1)$th-order accurate implicit member. In this case, however, the first sub-step size $\gamma_1$ is no longer adjustable, and the controllable high-frequency numerical dissipation is therefore lost. Motivated by the observations made for the two- and three-sub-step methods, the sub-step sizes $\gamma_j~(j=2,~\cdots,~s)$ obtained from the $(s+1)$th-order construction is also adopted as the default choice for the corresponding $s$th-order members with controllable numerical dissipation. Consequently, for each $s$-sub-step method, two complementary algorithm variants are obtained: for arbitrary admissible user-specified values of $\gamma_1$, equivalently $\infrho$, the resulting method achieves $s$th-order accuracy with controllable high-frequency dissipation; whereas for one specific value of $\gamma_1$, it automatically attains $(s+1)$th-order accuracy with a fixed high-frequency dissipation level.

\subsubsection{Four sub-steps: $s=4$}

\begin{figure}[htbp]
	\centering 
	\subfigure[$\gamma_1=1/4+\sqrt{3}/12~(\infrho=1.0)$]{
		\includegraphics[scale=0.4]{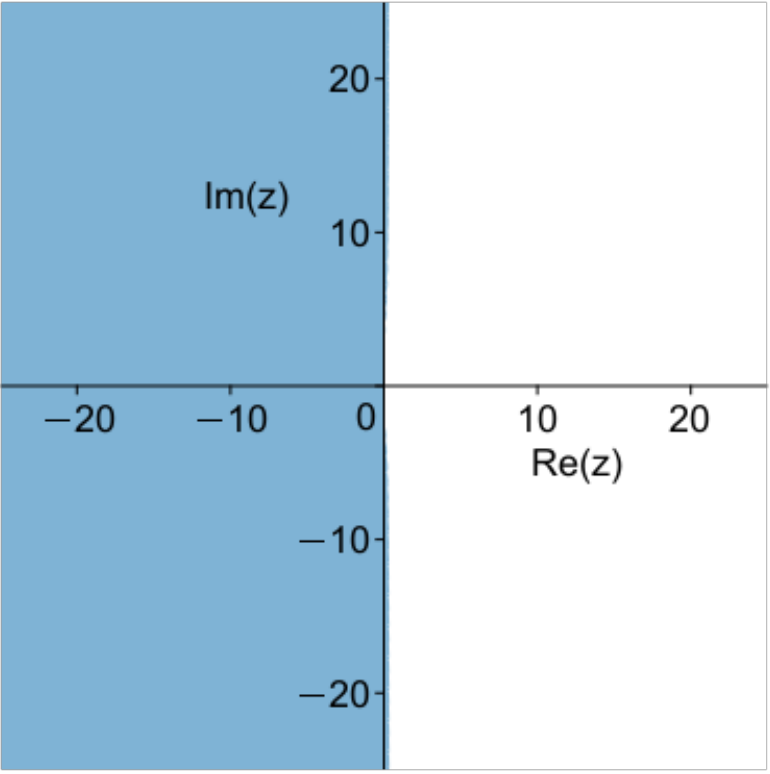}}
	\subfigure[$\gamma_1=0.47048~(\infrho=0.5)$]{
		\includegraphics[scale=0.4]{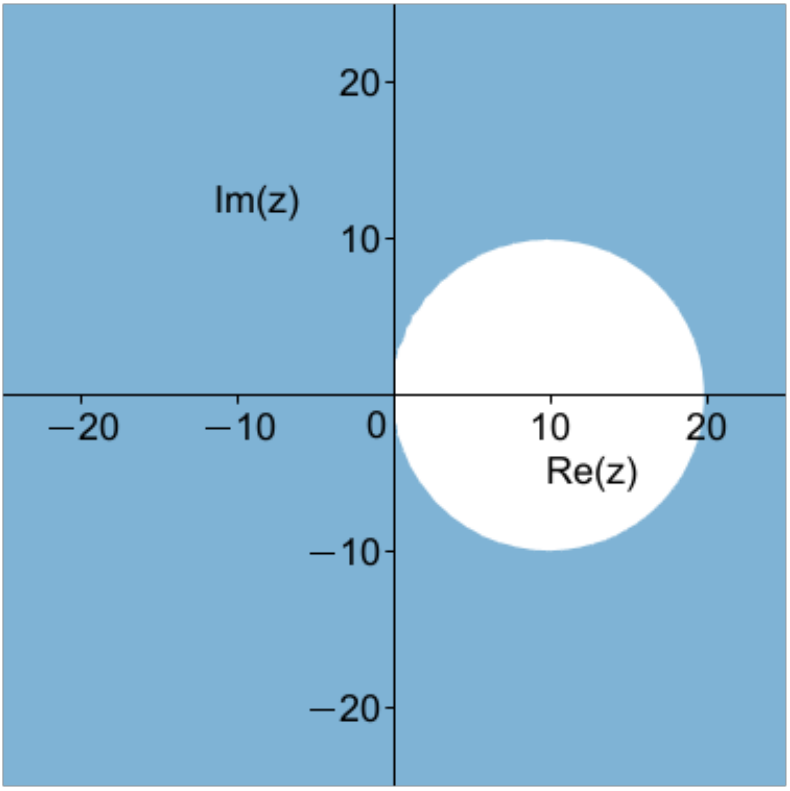}}
	\subfigure[$\gamma_1=0.57282~(\infrho=0.0)$]{
		\includegraphics[scale=0.4]{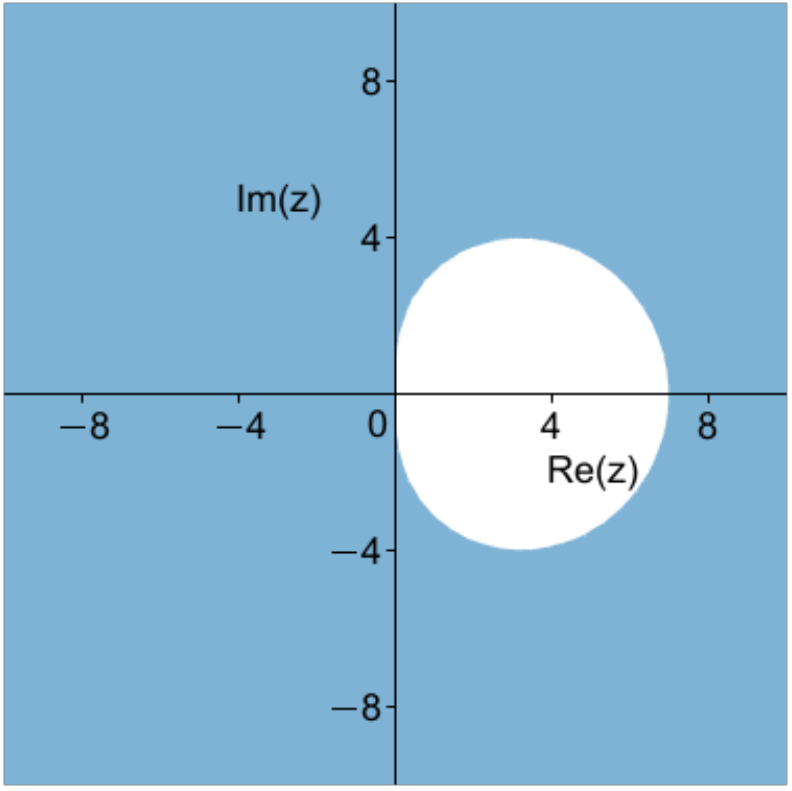}}
	\subfigure[$\gamma_1=0.78868~(\infrho=-0.5)$]{
		\includegraphics[scale=0.4]{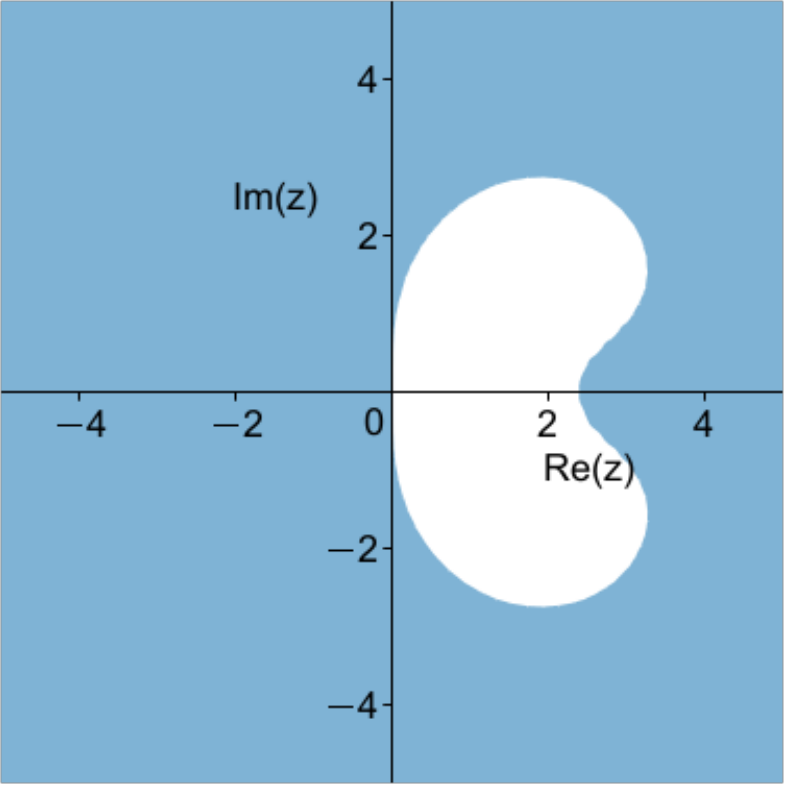}}
	\subfigure[$\gamma_1=1.06858~(\infrho=-0.63041)$]{
		\includegraphics[scale=0.4]{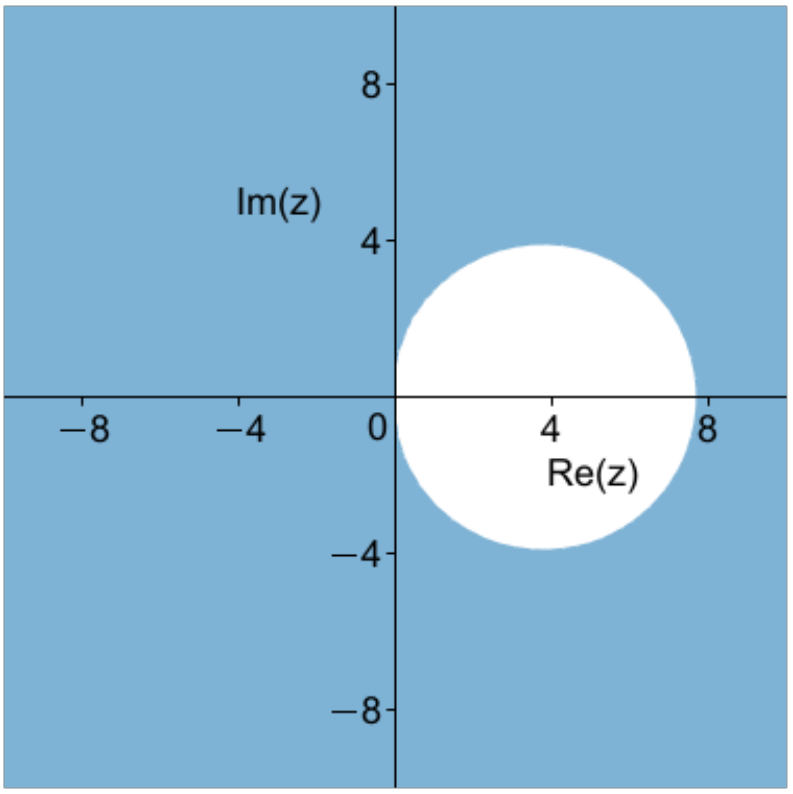}}
	\subfigure[$\gamma_1=1.34537~(\infrho\approx-0.57677)$]{
		\includegraphics[scale=0.4]{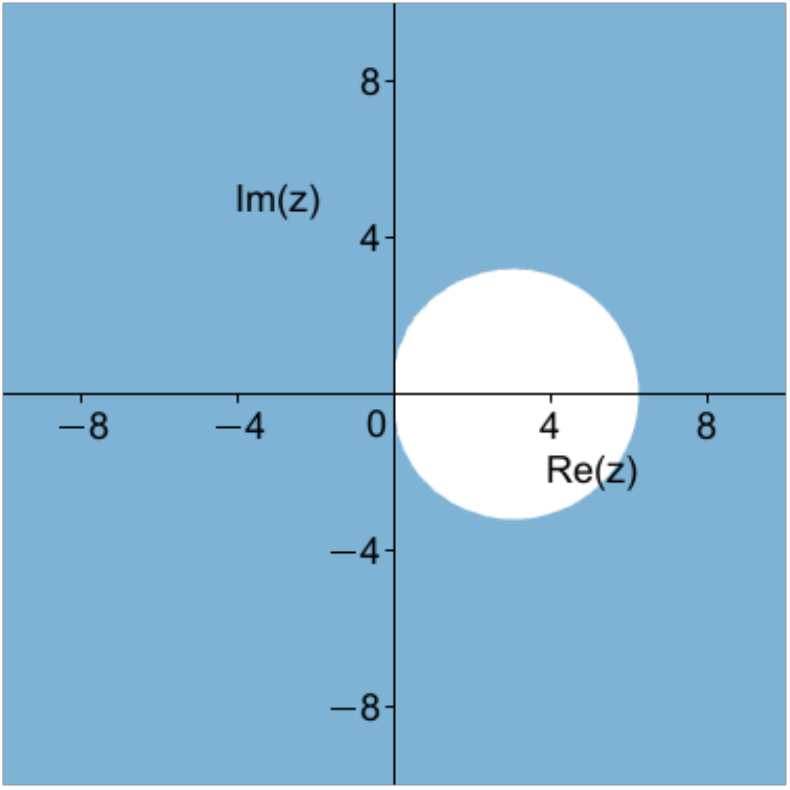}}
	\caption{The stability region for the four-sub-step implicit method defined by \cref{eq:s_eq_4}.}
	\label{fig:s4_stability}
\end{figure}

\begin{figure}[htbp]
	\centering 
	\includegraphics[scale=1.5]{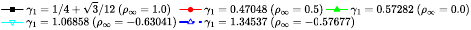}\\
    \subfigure[]{
		\includegraphics[scale=0.72]{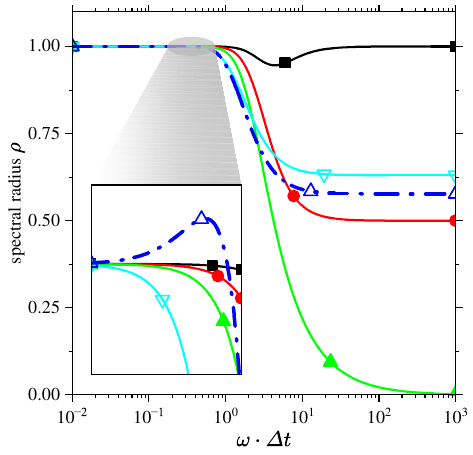}}
    \subfigure[]{
		\includegraphics[scale=0.72]{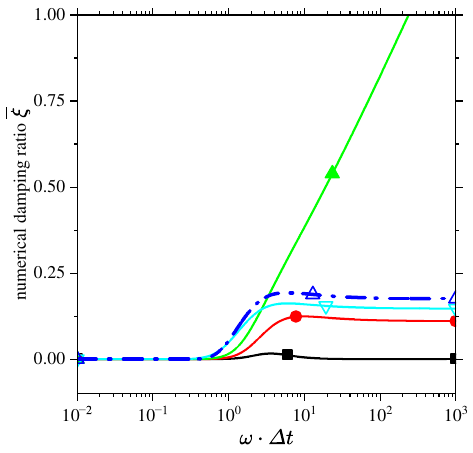}}
    \subfigure[]{
		\includegraphics[scale=0.72]{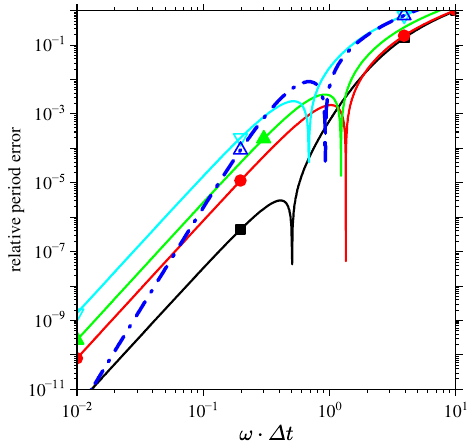}}
	\caption{Spectral properties for the four-sub-step implicit method defined by \cref{eq:s_eq_4} in the absence of $\xi$.}
	\label{fig:s4_sp}
\end{figure}

\cref{fig:s4_stability,fig:s4_sp} present the stability regions and spectral properties, respectively, of the four-sub-step implicit methods defined by \cref{eq:s_eq_4}. Several fourth-order members with different $\gamma_1$, the corresponding $\infrho$, are considered, together with the fifth-order member corresponding to $\gamma_1=1.34537$. The fourth-order members shown in \cref{fig:s4_stability}(a)--(e) are unconditionally stable, as their stability regions contain the entire left half-plane. Although the geometries of the stability boundaries vary with $\gamma_1$, the fundamental unconditional stability is preserved for all these admissible fourth-order configurations. \cref{fig:s4_stability}(f) corresponds to the fifth-order four-sub-step member with $\gamma_1=1.34537$. As discussed in \cref{app:s4}, this value does not satisfy the unconditional stability \eqref{eq:s4_uc}, and the resulting method is only $A(\alpha)$-stable, with $\alpha\approx89.99173663^\circ$. Since this angle is extremely close to $90^\circ$, the unstable sector excluded from the stability domain is very narrow and is barely distinguishable at the scale adopted in \cref{fig:s4_stability}(f). Therefore, although the fifth-order member is not strictly unconditionally stable, its stability region is practically very close to that of an $A$-stable method.

The spectral radius in \cref{fig:s4_sp}(a) further clarify the difference between the fourth- and fifth-order members. For the fourth-order methods, the spectral radius remains no greater than unity over the entire frequency range and approaches the prescribed limiting value $\left|\infrho\right|$ as $\omega\dt\rightarrow\infty$. The parameter $\infrho$ therefore effectively controls the asymptotic high-frequency dissipation. In particular, the case $\infrho=1$ introduces only a small amount of dissipation in the mid-frequency range, whereas the members with smaller limiting spectral radii provide stronger attenuation of high-frequency components. For the fifth-order member with $\gamma_1=1.34537$, however, the spectral radius slightly exceeds unity over a narrow mid-frequency interval, as highlighted in the enlarged view of \cref{fig:s4_sp}(a). This local amplification is a direct consequence of the fact that the method is only $A(\alpha)$-stable rather than strictly $A$-stable. Nevertheless, the maximum spectral radius is only marginally greater than unity, and no pronounced amplification is observed. This spectral behavior is consistent with the extremely small difference between $\alpha$ and $90^\circ$, indicating that the departure from unconditional stability is practically negligible for most applications.

\cref{fig:s4_sp}(b) shows the corresponding numerical damping ratios. All members exhibit negligible damping in the low-frequency range, followed by different levels of numerical damping at mid- and high-frequencies according to the selected value of $\infrho$. The relative period errors shown in \cref{fig:s4_sp}(c) demonstrate the phase accuracy of the four-sub-step methods. All members maintain small period errors in the low-frequency range, while their differences become more evident as $\omega\dt$ increases. The fifth-order member generally exhibits favorable phase accuracy owing to its higher formal order, despite the slight mid-frequency amplification discussed above. Overall, \cref{fig:s4_stability,fig:s4_sp} indicate that the proposed four-sub-step method provides unconditionally stable fourth-order algorithms with controllable high-frequency dissipation, while the fifth-order member offers enhanced temporal accuracy and a stability region that is extremely close to unconditional stability, at the cost of only a very small local amplification in the mid-frequency range.

The amplitude and phase errors \cite{li_NovelImplicitIntegration_2025,li_DesigningDevelopingSinglestep_2023} are computed analytically for the proposed four-sub-step implicit method \eqref{eq:s_eq_4} as 
\begin{subequations}\label{eq:s4_amp_phase}
\begin{align}
\delta=& -\dfrac{(16\xi^4-20\xi^{2}+5)\xi\varpi_1}{120}\omega^5\dt^4+ \dfrac{(32\xi^6-48\xi^4+18\xi^{2}-1)\varpi_2}{144}\omega^6\dt^5 \notag\\
&\quad- \dfrac{(64\xi^6-112\xi^4+56\xi^{2}-7)\xi\varpi_3}{336}\omega^7\dt^6 + \dfrac{\varpi_4}{1152}\omega^8\dt^7 +\mathcal{O}(\dt^8)   \\
\epsilon=& -\dfrac{(16\xi^4-12\xi^{2}+1)\sqrt{1-\xi^{2}}\varpi_1}{120}\omega^5\dt^4+ \dfrac{(16\xi^4-16\xi^{2}+3)\xi\sqrt{1-\xi^{2}}\varpi_2}{72}\omega^6\dt^5 + \dfrac{\varpi_3}{336}\omega^7\dt^6+  \mathcal{O}(\dt^7)
\end{align}
where $\varpi_j~(j=1,~2,~3,~4)$ are given as 
\begin{align}
\varpi_1 &= 120\gamma_1^4 - 240\gamma_1^3 + 120\gamma_1^2 - 20\gamma_1 + 1\\
\varpi_2 &= 576\gamma_1^5 - 1224\gamma_1^4 + 768\gamma_1^3 - 204\gamma_1^2 + 24\gamma_1 - 1\\
\varpi_3 &= 3360\gamma_1^6 - 7392\gamma_1^5 + 5208\gamma_1^4 - 1736\gamma_1^3 + 308\gamma_1^2 - 28\gamma_1 + 1\\
\varpi_4 &= 23040\gamma_1^7 - 51840\gamma_1^6 + 39168\gamma_1^5 - 14832\gamma_1^4 + 3264\gamma_1^3 - 432\gamma_1^2 + 32\gamma_1 - 1.
\end{align}
\end{subequations}
For the fourth-order four-sub-step member, $\varpi_1\neq0$. In the damped case $(\xi\neq0)$, the leading terms of both the amplitude $\delta$ and phase $\epsilon$ errors are of order $\mathcal{O}(\dt^{4})$. In the undamped case $(\xi=0)$, however, the leading $\dt^{4}$ term in $\delta$ vanishes because it is proportional to $\xi$, and hence the amplitude error exhibits one-order superconvergence, i.e., $\delta=\mathcal{O}(\dt^{5})$. By contrast, the leading term in $\epsilon$ generally remains nonzero at $\xi=0$, so that the phase error remains $\mathcal{O}(\dt^{4})$. For the fifth-order member, the additional condition $\varpi_1=0$ eliminates all $\dt^{4}$ terms in the amplitude and phase errors. Consequently, in the damped case, both $\delta$ and $\epsilon$ become $\mathcal{O}(\dt^{5})$. In the undamped case, the amplitude error also remains $\mathcal{O}(\dt^{5})$, whereas the $\dt^{5}$ term in the phase error vanishes because it is proportional to $\xi$. Therefore, the phase error exhibits one-order superconvergence and becomes $\mathcal{O}(\dt^{6})$. These results indicate that the fourth-order member possesses enhanced amplitude accuracy for undamped systems, while the fifth-order member exhibits enhanced phase accuracy under the same condition.

\begin{remark}
According to Eq.~(\ref{eq:s4_amp_phase}a), the amplitude error of the four-sub-step implicit method exhibits an enhanced convergence behavior for undamped systems. When $\xi=0$, all terms containing the physical damping ratio vanish, and the leading term in $\delta$ becomes $\mathcal{O}(\dt^5)$ instead of the formal order $\mathcal{O}(\dt^4)$. Therefore, the four-sub-step implicit method naturally possesses a one-order amplitude superconvergence for undamped problems. Furthermore, if the parameter $\gamma_1$ is selected such that $\varpi_2=0$, the $\dt^5$ term is completely eliminated, and the amplitude error is further improved to $\mathcal{O}(\dt^7)$. Among the real roots of $\varpi_2=0$, the admissible value satisfying the unconditional stability is $\gamma_1=1.2805797612753055$. Consequently, this parameter provides enhanced amplitude accuracy for undamped problems.
\end{remark}

\subsubsection{Five sub-steps: $s=5$}

\begin{figure}[htbp]
	\centering 
	\subfigure[$\gamma_1=0.24650~(\infrho=1.0)$]{
		\includegraphics[scale=0.4]{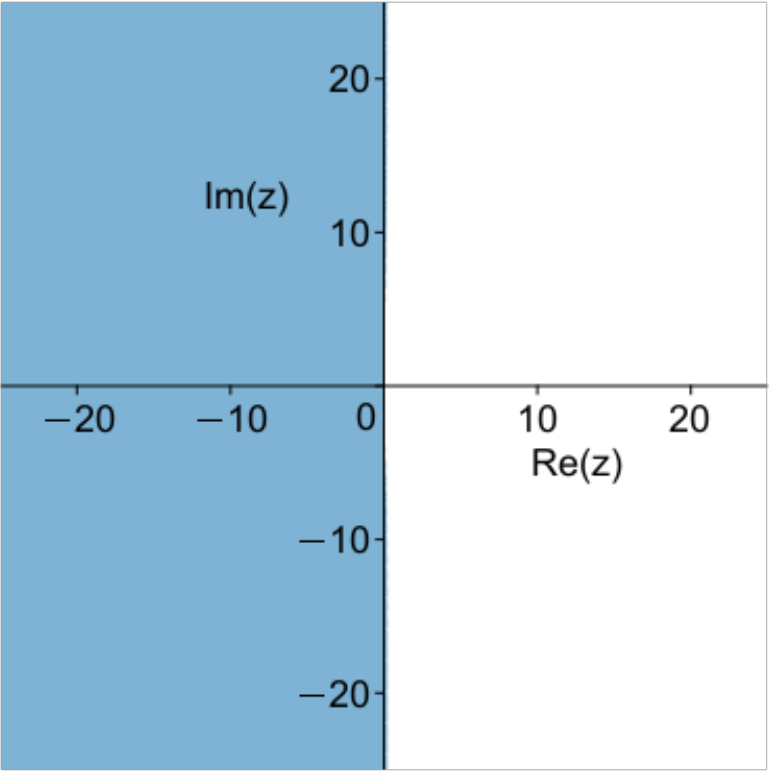}}
	\subfigure[$\gamma_1=0.26052~(\infrho=0.5)$]{
		\includegraphics[scale=0.4]{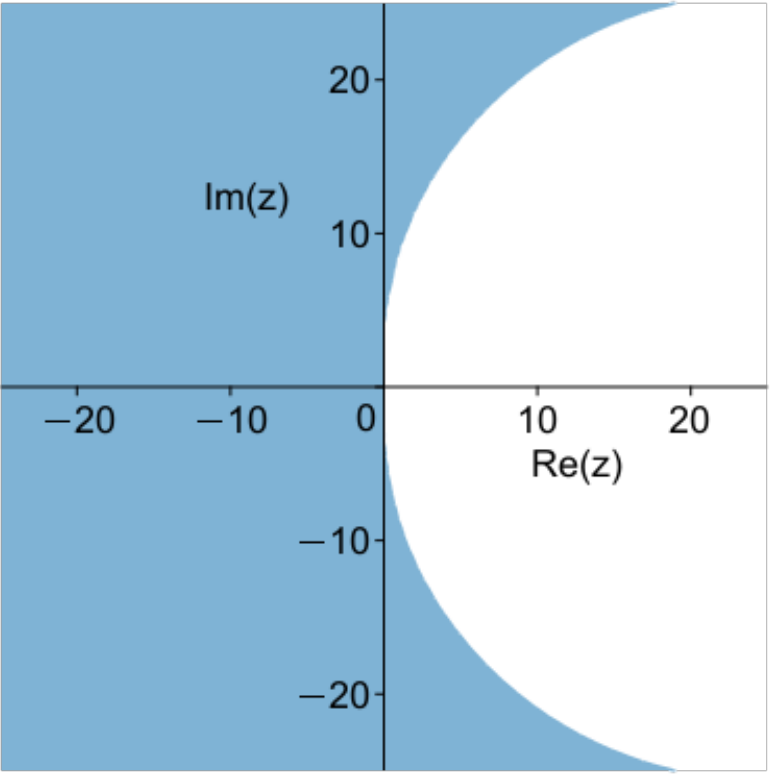}}
	\subfigure[$\gamma_1=0.27805~(\infrho=0.0)$]{
		\includegraphics[scale=0.4]{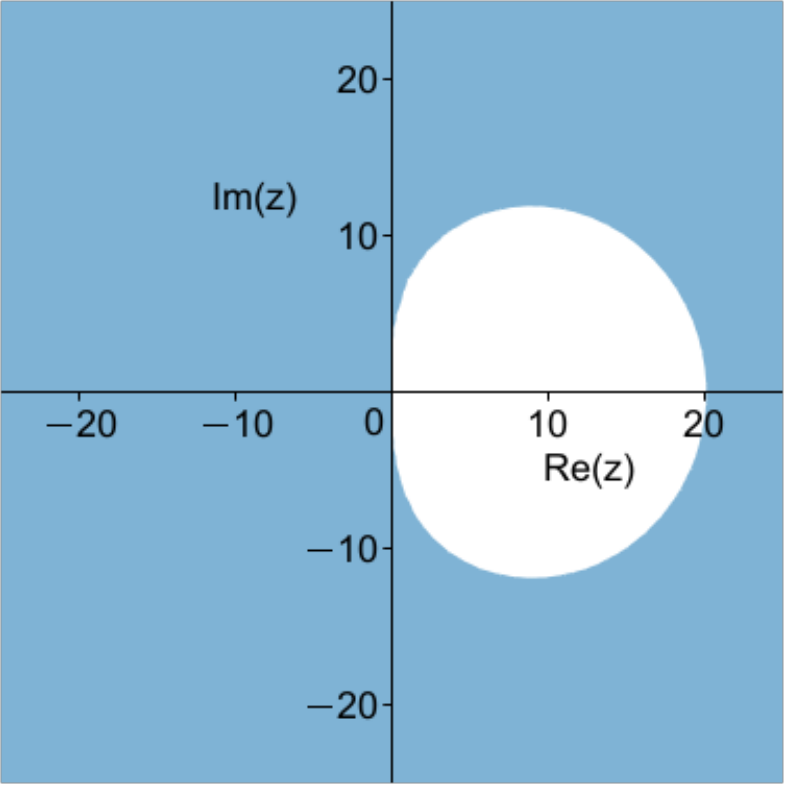}}
	\subfigure[$\gamma_1=0.30300~(\infrho=-0.5)$]{
		\includegraphics[scale=0.4]{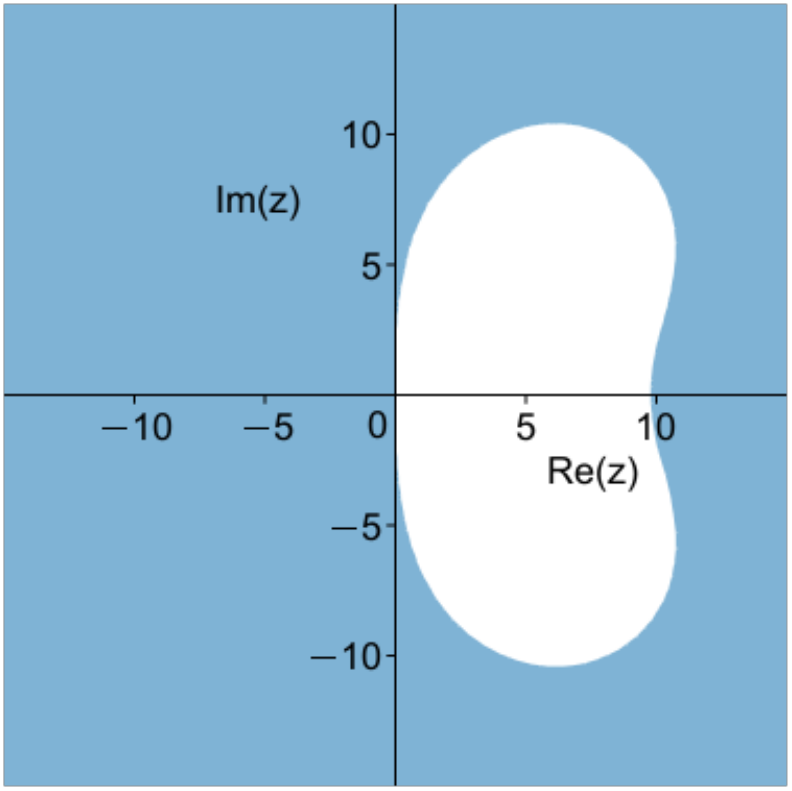}}
	\subfigure[$\gamma_1=0.45155~(\infrho=-0.91419)$]{
		\includegraphics[scale=0.408]{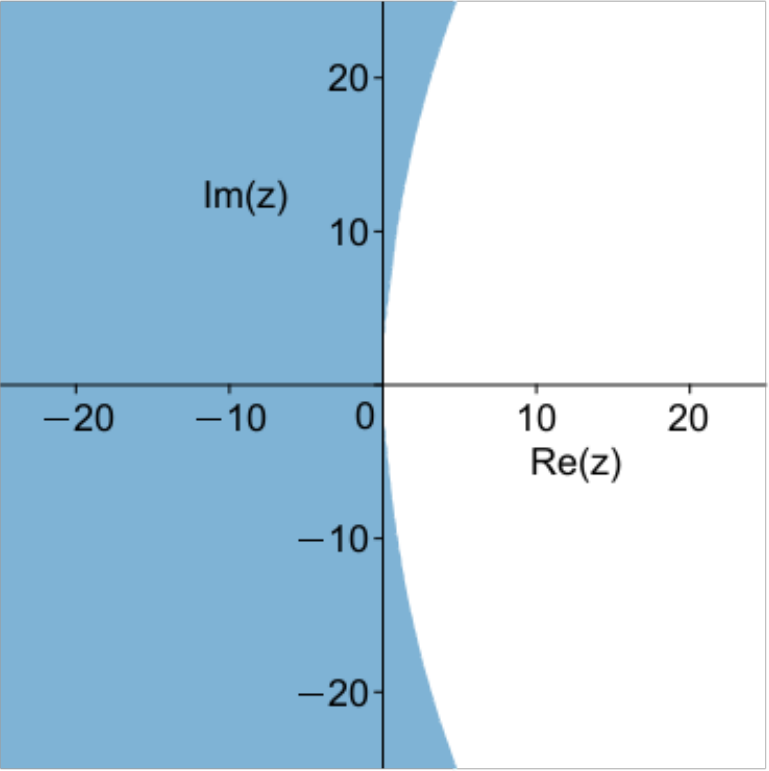}}
	\subfigure[$\gamma_1=0.47327~(\infrho\approx-0.83733)$]{
		\includegraphics[scale=0.408]{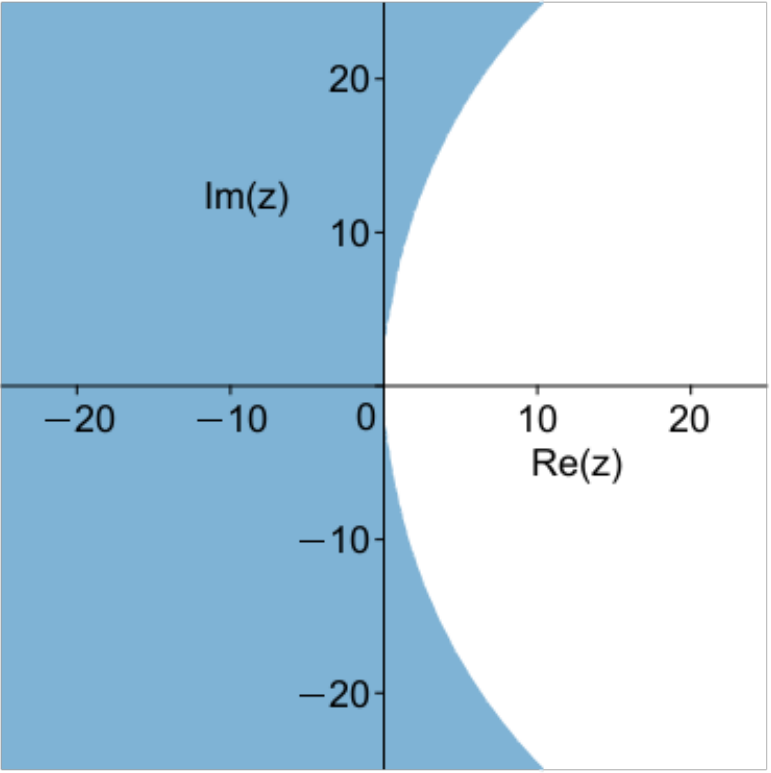}}
	\caption{The stability region for the five-sub-step implicit method defined by \cref{eq:s_eq_5}.}
	\label{fig:s5_stability}
\end{figure}

\cref{fig:s5_stability} illustrates the stability regions of the proposed five-sub-step implicit method defined by \cref{eq:s_eq_5} for different values of the user-specified parameter $\infrho$. As shown in \cref{fig:s5_stability}, the admissible values of $\gamma_1$ determined by the dissipation control always satisfy the stability requirement, demonstrating that the fifth-order five-sub-step implicit method preserves unconditional stability while allowing flexible control of high-frequency dissipation. It is noteworthy that the sixth-order member, corresponding to $\gamma_1=0.47327$ and $\infrho\approx-0.83733$, is also included in \cref{fig:s5_stability}. This member is obtained by imposing the additional order condition \eqref{eq:s5_r12345}, which completely determines the remaining algorithmic parameters. Therefore, unlike the general fifth-order members, the sixth-order method no longer possesses an independently adjustable dissipation parameter.

\begin{figure}[htbp]
	\centering 
	\includegraphics[scale=1.8]{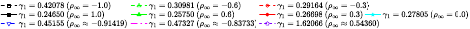}\\
    \subfigure[]{
		\includegraphics[scale=0.72]{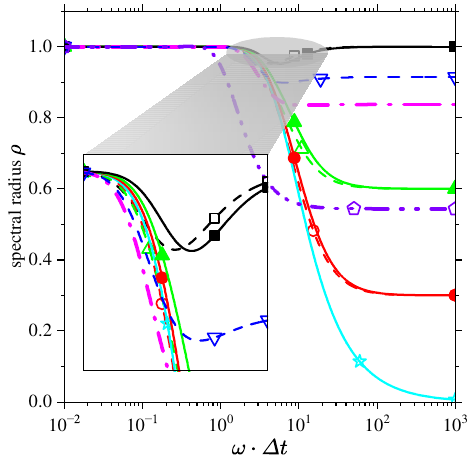}}
    \subfigure[]{
		\includegraphics[scale=0.72]{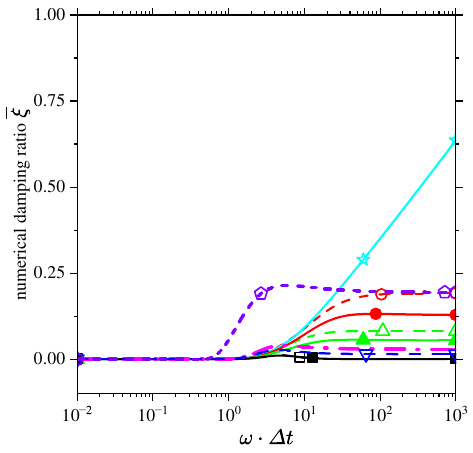}}\\
    \subfigure[]{
		\includegraphics[scale=0.72]{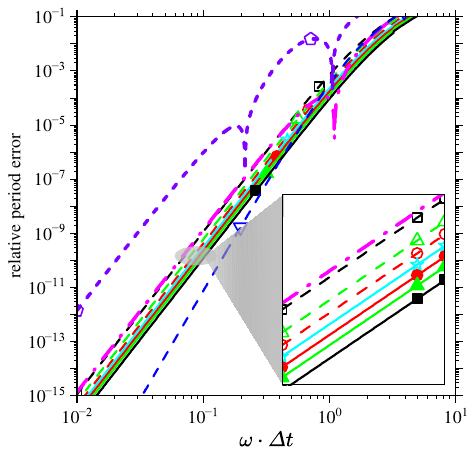}}
    \subfigure[]{
		\includegraphics[scale=0.72]{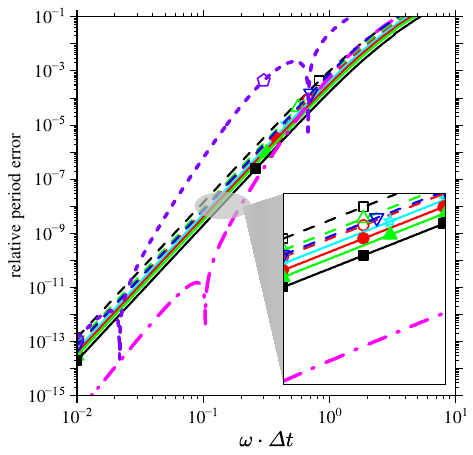}}
	\caption{Spectral properties for the five-sub-step implicit method defined by \cref{eq:s_eq_5}: (a-c) $\xi=0$ and (d) $\xi=0.2$.}
	\label{fig:s5_sp}
\end{figure}

The spectral properties of the five-sub-step implicit method are further examined in \cref{fig:s5_sp}. \cref{fig:s5_sp}(a) presents the spectral radius in the absence of $\xi$. The results demonstrate that different choices of $\infrho$ successfully produce different high-frequency dissipation levels, with the limiting spectral radius approaching the prescribed value $\left|\infrho\right|$. In particular, the two sixth-order members ($\gamma_1=0.47327$ and $\gamma_1=1.62066$) exhibit the fixed asymptotic spectral radius of approximately $0.83733$ and $0.54360$, while the fifth-order members can realize a broad range of high-frequency damping characteristics. This observation further verifies the controllable dissipation property predicted by the theoretical formulation.

\cref{fig:s5_sp}(c) and (d) compare the relative period errors in the undamped and damped cases, respectively. As shown in \cref{eq:s5_amp_phase}, the general fifth-order member exhibits amplitude error of $\mathcal{O}(\dt^5)$, while its phase error is $\mathcal{O}(\dt^6)$ for undamped systems due to the disappearance of the leading phase error proportional to $\xi$. For damped systems, however, the phase error remains $\mathcal{O}(\dt^5)$ because the leading term associated with $\psi_1$ does not vanish. By imposing $\psi_1=0$, the sixth-order members eliminate these leading error terms and achieves improved amplitude and phase accuracy, particularly for damped problems. Furthermore, the $A$-stable sixth-order member ($\gamma_1=0.47327$) exhibits significantly smaller relative period errors than the $A(\alpha)$-stable sixth-order member ($\gamma_1=1.62066$). These theoretical conclusions are consistent with the numerical results shown in \cref{fig:s5_sp}(c) and (d). Furthermore, as discussed in \cref{rem:s5_2}, the fifth-order member with $\gamma_1\approx0.4515512290$ is a special case that eliminates the leading phase error in undamped problems. Consequently, this member exhibits enhanced phase accuracy for undamped problems.

The amplitude and phase errors \cite{li_NovelImplicitIntegration_2025,li_DesigningDevelopingSinglestep_2023} are computed for the proposed five-sub-step implicit method \eqref{eq:s_eq_5} as 
\begin{subequations}\label{eq:s5_amp_phase}
\begin{align}
\delta=&\dfrac{(1-2\xi^{2})(16\xi^4-16\xi^{2}+1)\psi_1}{720}\omega^6\dt^5+ \dfrac{(64\xi^6-112\xi^4+56\xi^{2}-7)\xi\psi_2}{840}\omega^7\dt^6\notag\\
&- \dfrac{(128\xi^8-256\xi^6+160\xi^4-32\xi^{2}+1)\psi_3}{1920}\omega^8\dt^7+\mathcal{O}(\dt^8)   \\
\epsilon=& \dfrac{(4\xi^{2}-3)(1-4\xi^{2})\xi\sqrt{1-\xi^{2}}\psi_1}{360} \omega^6\dt^5+ \dfrac{(64\xi^6-80\xi^4+24\xi^{2}-1)\sqrt{1-\xi^{2}}\psi_2}{840} \omega^7\dt^6\notag\\
&- \dfrac{(16\xi^6-24\xi^4+10\xi^{2}-1)\xi\sqrt{1-\xi^{2}}\psi_3}{240} \omega^8\dt^7+ \dfrac{\psi_4}{6480}\omega^9\dt^8+\mathcal{O}(\dt^9) \label{eq:s5_epsilon}
\end{align}
where $\psi_j~(j=1,~2,~3,~4)$ are given as 
\begin{align}
\psi_1&=720\gamma_1^5 - 1800\gamma_1^4 + 1200\gamma_1^3 - 300\gamma_1^2 + 30\gamma_1 - 1\label{eq:s5_psi_1}\\
\psi_2&=4200\gamma_1^6-10920\gamma_1^5+8400\gamma_1^4-2800\gamma_1^{3}+455\gamma_1^{2}-35\gamma_1+1\label{eq:s5_psi_2}\\
\psi_3&=28800\gamma_1^7-76800\gamma_1^6+64320\gamma_1^5-25200\gamma_1^4+5360\gamma_1^{3}-640\gamma_1^{2}+40\gamma_1-1\\
\psi_4&=226800\gamma_1^8-615600\gamma_1^7+545400\gamma_1^6-236520\gamma_1^5+59130\gamma_1^4-9090\gamma_1^{3}+855\gamma_1^{2}-45\gamma_1+1.
\end{align}
\end{subequations}
For the fifth-order five-sub-step member, the coefficient $\psi_1$ in \cref{eq:s5_psi_1} does not vanish. In this case, the leading amplitude error $\delta$ is generally of order $\mathcal{O}(\dt^{5})$. For the phase error $\epsilon$, the leading $\dt^{5}$ term is proportional to the physical damping ratio $\xi$. Therefore, in the undamped case $(\xi=0)$, the phase error exhibits a one-order superconvergence and becomes $\mathcal{O}(\dt^{6})$, whereas in the damped case $(\xi\neq0)$, the phase error remains $\mathcal{O}(\dt^{5})$. On the other hand, the sixth-order member is obtained by imposing the additional condition $\psi_1=0$, leading to $\gamma_1=0.4732683912582953$. Consequently, all leading amplitude and phase errors associated with $\dt^{5}$ vanish. For the damped case, the leading amplitude and phase errors are therefore improved to $\mathcal{O}(\dt^{6})$. For the undamped case, the amplitude error exhibits a one-order superconvergence, reaching $\mathcal{O}(\dt^7)$. These observations explain the different relative period errors between the fifth-order and sixth-order five-sub-step implicit methods, as shown in \cref{fig:s5_sp}(c) and (d). 
\begin{remark}\label{rem:s5_2}
According to \cref{eq:s5_epsilon}, the leading phase error of the fifth-order five-sub-step implicit method in the undamped case can be further eliminated by imposing $\psi_2=0$. Among the admissible real solutions, the value $\gamma_1\approx 0.4515512289893862$ satisfies the unconditional stability given in \cref{eq:s5_uc}. Substituting this value into the high-frequency dissipation relation gives $\left|\infrho\right|\approx0.914190237$. Consequently, the corresponding fifth-order member achieves eighth-order phase accuracy $\mathcal{O}(\dt^8)$ for undamped problems. Therefore, $\gamma_1\approx 0.4515512289893862$ is recommended as the default parameter for the fifth-order five-sub-step implicit method.
\end{remark}

\subsubsection{Six sub-steps: $s=6$}

\begin{figure}[htbp]
	\centering 
	\subfigure[$\gamma_1=0.28406~(\infrho=1.0)$]{
		\includegraphics[scale=0.4]{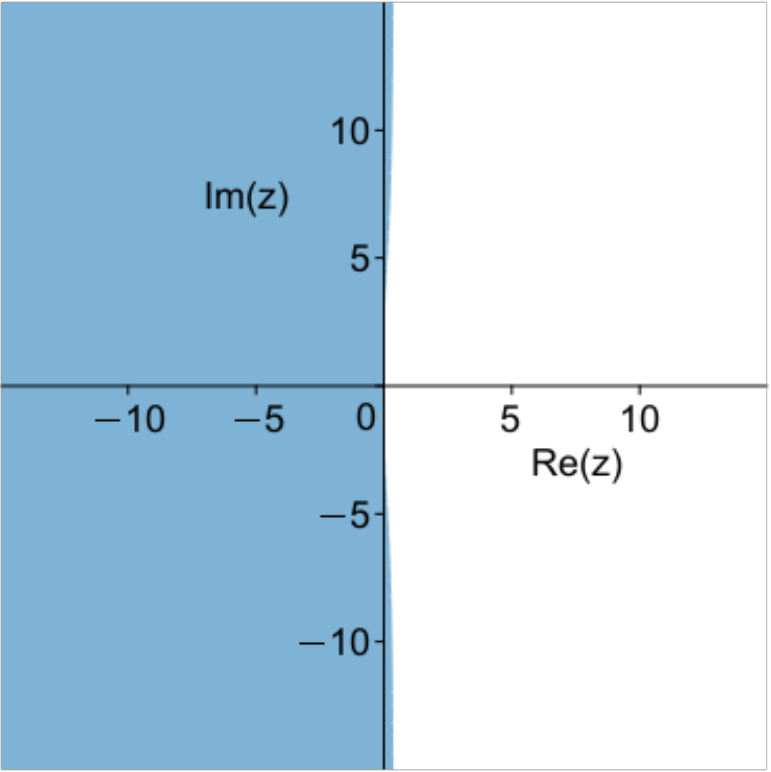}}
	\subfigure[$\gamma_1=0.30633~(\infrho=0.5)$]{
		\includegraphics[scale=0.4]{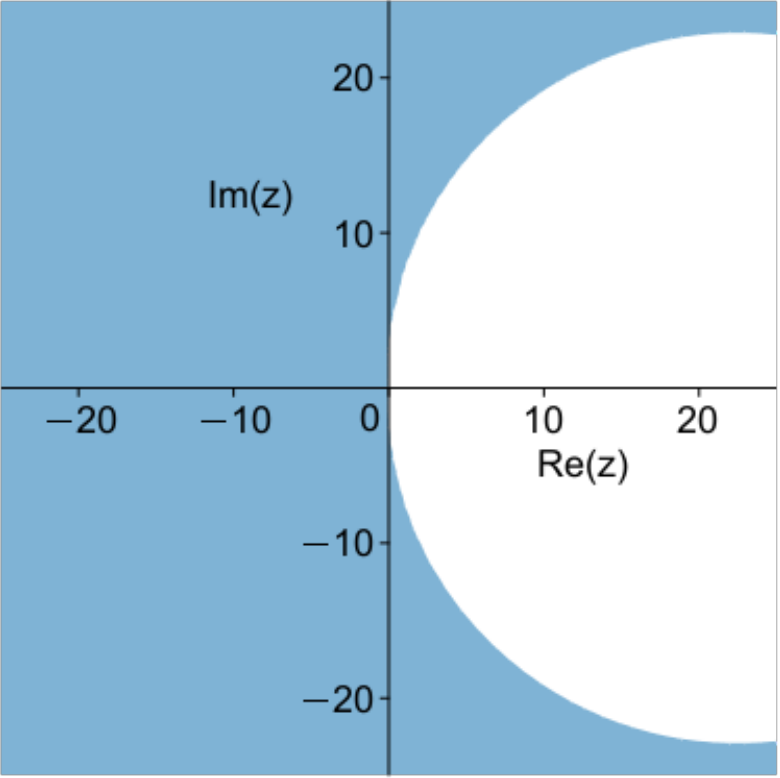}}
	\subfigure[$\gamma_1=0.33414~(\infrho=0.0)$]{
		\includegraphics[scale=0.4]{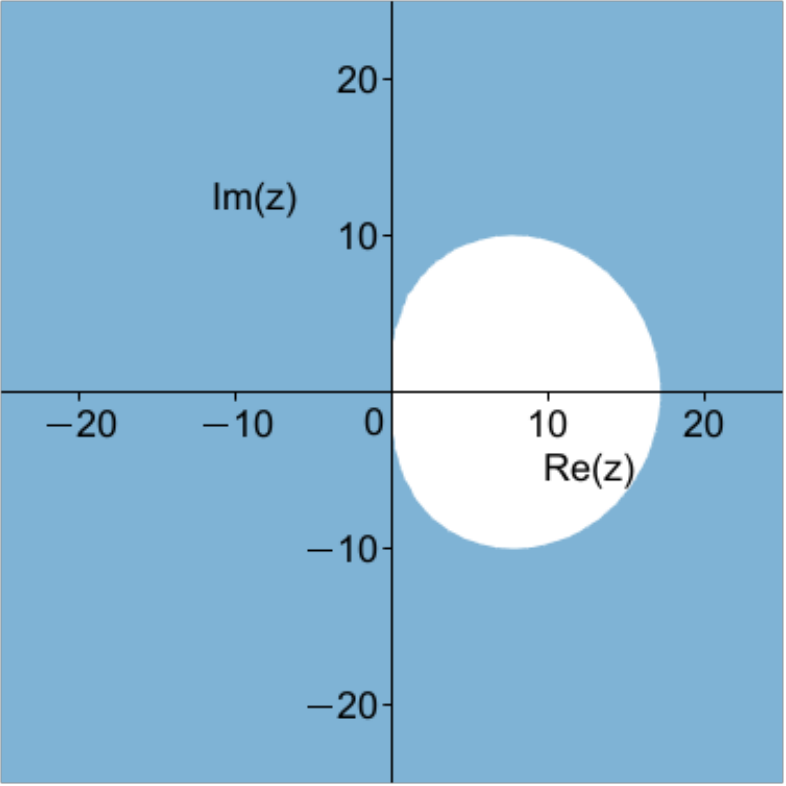}}
	\subfigure[$\gamma_1=0.37670~(\infrho=-0.5)$]{
		\includegraphics[scale=0.4]{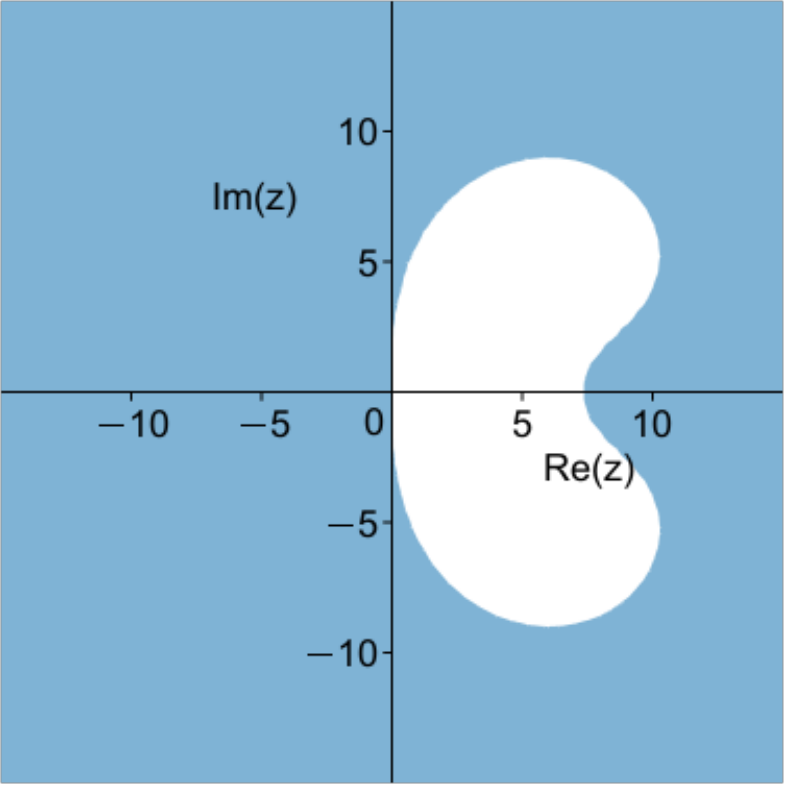}}
	\subfigure[$\gamma_1=0.47327~(\infrho=-0.83733)$]{
		\includegraphics[scale=0.408]{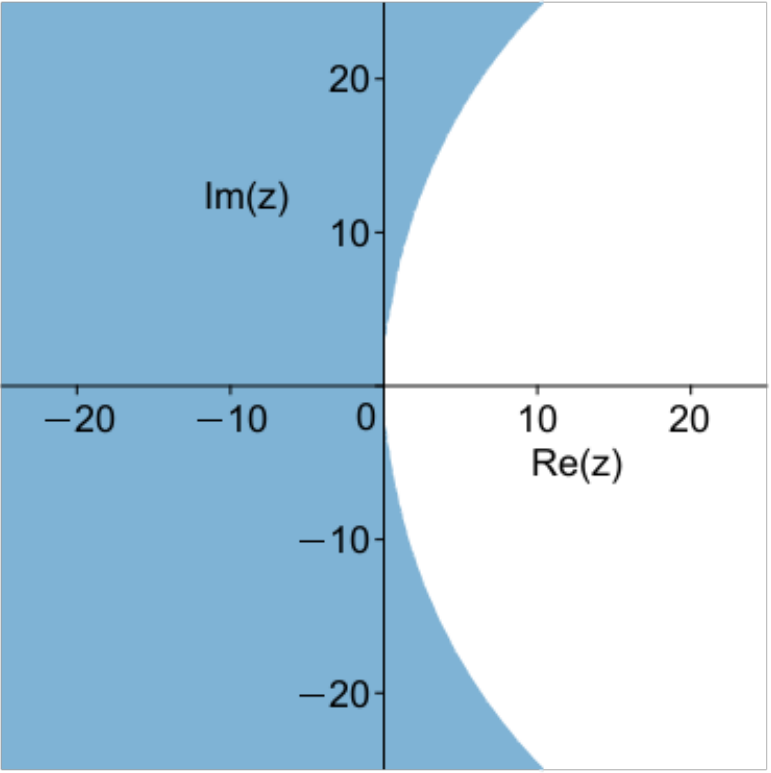}}
	\subfigure[$\gamma_1=0.55670~(\infrho\approx-0.72080)$]{
		\includegraphics[scale=0.408]{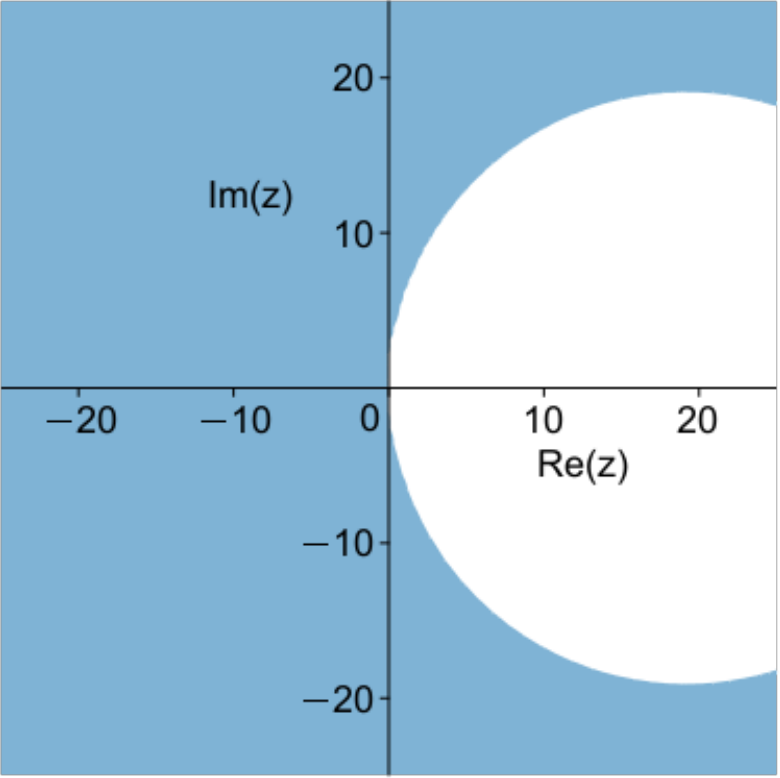}}
	\caption{The stability region for the six-sub-step implicit method defined by \cref{eq:s_eq_6}.}
	\label{fig:s6_stability}
\end{figure}

\cref{fig:s6_stability} presents the stability regions of the proposed six-sub-step implicit method for several representative values of $\gamma_1$, equivalently $\infrho$. Similar to the lower-sub-step members, the stability region evolves continuously as the parameter $\gamma_1~(\infrho)$ varies. For $\infrho=1$, the method exhibits almost no numerical dissipation in the high-frequency limit, while a slight reduction of the spectral radius is still observed in the mid-frequency range, resulting in a stability region that extends slightly into the right half-plane rather than coinciding exactly with the classical $A$-stable region. As $\infrho$ decreases, the high-frequency spectral radius is gradually reduced and the stability region evolves accordingly. In particular, the member with $\infrho=0$ exhibits the strongest high-frequency damping among the sixth-order members. \cref{fig:s6_stability}(e) corresponds to the sixth-order member associated with the lower admissible bound of the user-specified parameter, i.e., $\infrho\approx-0.83733$. \cref{fig:s6_stability}(f) shows the seventh-order member obtained by imposing one additional order condition, corresponding to $\gamma_1\approx0.55670$. Unlike the sixth-order members, this algorithm is only $A(\alpha)$-stable rather than strictly unconditionally stable. Nevertheless, as demonstrated in \cref{fig:s6_r1_rho}, the corresponding stability angle is extremely close to $90^\circ$, indicating that its stability remain very close to unconditional stability. Together with its superior phase accuracy, this observation justifies adopting this parameter set as the recommended seventh-order member within the present framework.

\begin{figure}[htbp]
	\centering 
	\includegraphics[scale=1.8]{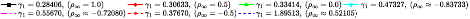}\\
    \subfigure[]{
		\includegraphics[scale=0.72]{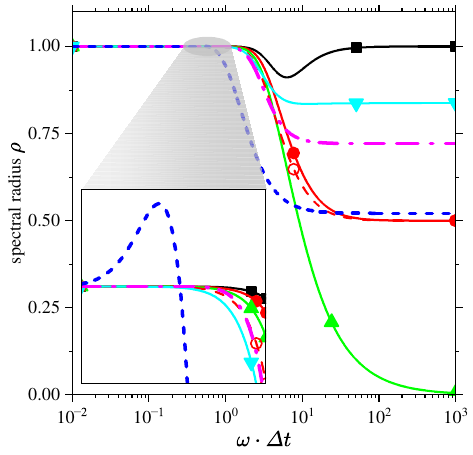}}
    \subfigure[]{
		\includegraphics[scale=0.72]{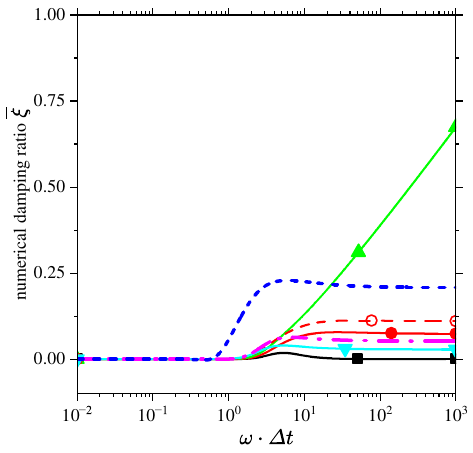}}
    \subfigure[]{
		\includegraphics[scale=0.72]{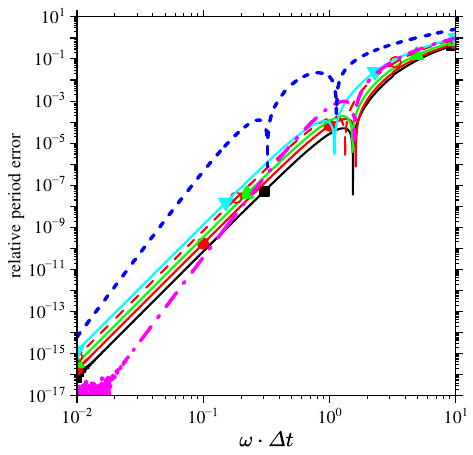}}
	\caption{Spectral properties for the six-sub-step implicit method defined by \cref{eq:s_eq_6} in the absence of $\xi$.}
	\label{fig:s6_sp}
\end{figure}

\cref{fig:s6_sp} compares the spectral properties of the six-sub-step implicit methods in the undamped case. The spectral radius shown in \cref{fig:s6_sp}(a) exhibit the same general trend as observed for the lower-sub-step formulations. As the prescribed value of $\infrho$ decreases, the asymptotic spectral radius approaches the corresponding user-specified high-frequency limit. Meanwhile, the sixth-order member with $\infrho=1$ still introduces a moderate amount of numerical dissipation over the mid-frequency range, explaining why its stability region in \cref{fig:s6_stability}(a) does not exactly coincide with the entire left half-plane.

The preferred seventh-order member, corresponding to $\gamma_1=0.55670$, possesses a fixed high-frequency spectral radius of approximately $0.72080$. Since this member is also $A(\alpha)$-stable rather than strictly unconditionally stable, its spectral radius theoretically exceeds unity over a very narrow mid-frequency interval. However, the overshoot is extremely small and is barely distinguishable even in the enlarged view of \cref{fig:s6_sp}(a). Consequently, its spectral radius and numerical damping remain very close to those of the unconditionally stable sixth-order members. By contrast, the alternative seventh-order member with $\gamma_1\approx1.89513$ exhibits a much more pronounced overshoot of the spectral radius above unity in the mid-frequency range, as clearly shown in the enlarged view of \cref{fig:s6_sp}(a). This behavior is another manifestation of its $A(\alpha)$-stable nature, but indicates a noticeably stronger amplification of mid-frequency components. The corresponding numerical damping ratio shown in \cref{fig:s6_sp}(b) also differs significantly from the other members. Together with its inferior relative period error, these observations further support excluding this parameter choice from practical use.

The relative period errors are compared in \cref{fig:s6_sp}(c). As predicted by the analytical errors given in \cref{eq:s6_amp_phase}, increasing the temporal order from sixth to seventh substantially improves the phase accuracy. The preferred seventh-order member exhibits the smallest relative period error over the low-frequency range, confirming the effectiveness of exploiting the additional degree of freedom to enhance the order of accuracy. In contrast, although the alternative seventh-order member also satisfies the seventh-order conditions, its significantly larger period error, combined with its inferior stability, makes it considerably less attractive in practical computations. These observations are fully consistent with the theoretical analysis presented in the corresponding remarks.

The amplitude and phase errors \cite{li_NovelImplicitIntegration_2025,li_DesigningDevelopingSinglestep_2023} are computed for the proposed six-sub-step implicit method \eqref{eq:s_eq_6} as 
\begin{subequations}\label{eq:s6_amp_phase}
\begin{align}
\delta=& \dfrac{(64\xi^6-112\xi^4+56\xi^{2}-7)\xi\varphi_1}{5040}\omega^7\dt^6+ \dfrac{(128\xi^8-256\xi^6+160\xi^4-32\xi^{2}+1)\varphi_2}{5760}\omega^8\dt^7 \notag\\
& -\dfrac{(4\xi^{2}-3)(64\xi^6-96\xi^4+36\xi^{2}-3)\xi\varphi_3}{12960}\omega^9\dt^8-\dfrac{\varphi_4}{43200}\omega^{10}\dt^9 +\mathcal{O}(\dt^{10})   \\
\epsilon=& -\dfrac{(64\xi^6-80\xi^4+24\xi^{2}-1)\sqrt{1-\xi^{2}}\varphi_1}{5040}\omega^7\dt^6 + \dfrac{(16\xi^6-24\xi^4+10\xi^{2}-1)\xi\sqrt{1-\xi^{2}}\varphi_2}{720}\omega^8\dt^7\notag \\
&- \dfrac{\varphi_3}{12960}\omega^9\dt^8+\mathcal{O}(\dt^9) 
\end{align}
where $\varphi_j~(j=1,~2,~3,~4)$ are determined by 
\begin{align}
\varphi_1 &= 5040\gamma_1^6 - 15120\gamma_1^5 + 12600\gamma_1^4 - 4200\gamma_1^3 + 630\gamma_1^2 - 42\gamma_1 + 1\\
\varphi_2 &=34560\gamma_1^7 - 106560\gamma_1^6 + 97920\gamma_1^5 - 39600\gamma_1^4 + 8160\gamma_1^3 - 888\gamma_1^2 + 48\gamma_1 - 1\\
\varphi_3 &= 272160\gamma_1^8 - 855360\gamma_1^7 + 838080\gamma_1^6 - 382320\gamma_1^5 + 95580\gamma_1^4 - 13968\gamma_1^3 + 1188\gamma_1^2 - 54\gamma_1 + 1\\
\varphi_4 &= 151200\gamma_1^8(16\gamma_1-51)+3600\gamma_1^6(2196\gamma_1-1087)+180\gamma_1^4(6168\gamma_1-1087)+90\gamma_1^{2}(244\gamma_1-17) + 60\gamma_1 - 1. 
\end{align}
\end{subequations}
For the sixth-order six-sub-step member, $\varphi_1\neq0$. In the damped case $(\xi\neq0)$, the leading terms of both the amplitude $\delta$ and phase $\epsilon$ errors are of order $\mathcal{O}(\dt^{6})$. In the undamped case $(\xi=0)$, however, the leading $\dt^{6}$ term in $\delta$ vanishes because it is proportional to $\xi$, so that the amplitude error exhibits one-order superconvergence, i.e., $\delta=\mathcal{O}(\dt^{7})$. By contrast, the leading term in $\epsilon$ does not vanish at $\xi=0$, and the phase error generally remains $\mathcal{O}(\dt^{6})$. For the seventh-order member, the additional condition $\varphi_1=0$ eliminates all $\dt^{6}$ terms in Eqs.~(\ref{eq:s6_amp_phase}a-b). Consequently, in the damped case, both the amplitude and phase errors become $\mathcal{O}(\dt^{7})$. In the undamped case, the amplitude error generally remains $\mathcal{O}(\dt^{7})$, whereas the $\dt^{7}$ term in the phase error further vanishes because it is proportional to $\xi$. Therefore, the phase error exhibits one-order superconvergence and becomes $\mathcal{O}(\dt^{8})$. These results indicate that the sixth-order member possesses enhanced amplitude accuracy for undamped systems, whereas the seventh-order member exhibits enhanced phase accuracy under the same condition.

\begin{remark}
According to Eq.~(\ref{eq:s6_amp_phase}a), the amplitude error of the six-sub-step implicit method exhibits an enhanced convergence behavior for undamped systems. When $\xi=0$, all terms containing the physical damping ratio vanish, and the leading term in $\delta$ becomes $\mathcal{O}(\dt^7)$ instead of the formal order $\mathcal{O}(\dt^6)$. Therefore, the six-sub-step implicit method naturally possesses a one-order amplitude superconvergence for undamped problems. Furthermore, if the parameter $\gamma_1$ is selected such that $\varphi_2=0$, the $\dt^7$ term is eliminated, and the amplitude error is further improved to $\mathcal{O}(\dt^9)$. Among the real roots of $\varphi_2=0$, the admissible value satisfying unconditional stability is $\gamma_1=0.5409068780733081$. Consequently, this parameter provides enhanced amplitude accuracy for undamped problems.
\end{remark}

\subsection{Recovery of additional quantities}

The proposed directly self-starting implicit methods \eqref{eq:but_alg} are developed primarily to eliminate the requirement of initial highest derivatives and avoid additional derivative-related operations within the time-stepping procedure. Therefore, for second-order transient dynamics, the proposed methods directly update the displacement and velocity variables without explicitly providing acceleration responses at each time step. Similarly, for first-order transient dynamics, the proposed methods directly advance the state variables without explicitly evaluating their time derivatives as output quantities. Consequently, additional procedures are required when acceleration responses for second-order systems or derivative responses for first-order systems are desired by users.

It should be emphasized that these additional response quantities are only required for output purposes and do not participate in time integrations. Therefore, the accuracy and construction strategy of the corresponding output procedures do not influence any conclusions established previously, including the identical effective matrices, the accuracy, the stability, spectral properties, and numerical dissipation of the proposed $s$-sub-step implicit methods.

The construction of these additional response quantities is not the focus of the present study. Several existing studies~\cite{li_DirectlySelfstartingHigherorder_2022,lee_ImplicitSsubstepTime_2025,oshea_HighorderImplicitTime_2025} have already investigated this issue and developed efficient approaches for recovering acceleration responses or higher-order derivatives from the computed numerical solutions. For example, the finite difference post-processing techniques and the reconstruction strategies based on governing equations have been successfully applied to integration algorithms. These techniques can be directly incorporated into the proposed methods when additional response quantities are required. Therefore, to avoid unnecessary repetition and maintain the focus on the development of the proposed directly self-starting implicit methods, the detailed construction of these output procedures is not further discussed herein.

\section{Comparisons}\label{sec:comparisons}

\cref{tab:sth,tab:s1th} summarize numerical characteristics of the proposed $s$-sub-step implicit methods with $s$th-order and $(s+1)$th-order accuracy, respectively. As shown in \cref{tab:sth}, the proposed $s$-sub-step method with $s$th-order accuracy provides a unified mechanism for controlling high-frequency dissipation. For each sub-step number, the first sub-step size $\gamma_1$ is associated with the prescribed high-frequency spectral radius $\infrho$, while the remaining sub-step sizes provide additional freedom for improving the accuracy order. The corresponding dissipation parameter can cover the complete range of $\infrho\in[0,~1]$, allowing the method to continuously vary from nearly non-dissipative behavior to strong high-frequency damping. Moreover, some formulations permit negative values of $\infrho$, corresponding to a change in the sign of the high-frequency amplification factor while maintaining the same spectral radius. This additional flexibility provides a wider design space for balancing numerical dissipation and phase accuracy. Some interesting and important observations from \cref{tab:sth} are summarized as follows. 
\begin{itemize}
    \item The first-order single-sub-step method appears to possess the possibility of achieving fourth-order phase accuracy in the undamped case. According to the analytical phase error given in Eq.~(\ref{eq:s1_amp_phas}b), eliminating the corresponding leading term in $\epsilon$ requires $\infrho^{2}-\infrho+1=0$, equivalently $3\gamma_1^{2}-3\gamma_1+1=0$, which yields a complex-valued sub-step parameter. This special case is not considered further in the present study.
    \item The amplitude and phase errors exhibit distinct parity-dependent behaviors with respect to the number of sub-steps. For odd sub-step methods, including the single-, three-, and five-sub-step schemes, the amplitude and phase errors in damped systems remain consistent with the formal accuracy order, i.e., $\delta=\epsilon=\mathcal{O}(\dt^s)$. For undamped systems, however, the leading term in $\epsilon$ vanishes because it is proportional to the physical damping ratio $\xi$, resulting in a one-order phase superconvergence, i.e., $\epsilon=\mathcal{O}(\dt^{s+1})$, while the amplitude error generally remains $\delta=\mathcal{O}(\dt^s)$. In contrast, even sub-step methods, including the two-, four-, and six-sub-step schemes, exhibit a complementary error behavior. For damped systems, both amplitude and phase errors preserve the formal accuracy order, i.e., $\delta=\epsilon=\mathcal{O}(\dt^s)$. For undamped systems, the leading term in $\delta$ vanishes, leading to $\delta=\mathcal{O}(\dt^{s+1})$, whereas the phase error generally remains $\epsilon=\mathcal{O}(\dt^s)$. Therefore, the superconvergence phenomenon alternates between phase and amplitude accuracy depending on the parity of the sub-step number.
    \item For the three- and five-sub-step methods, additional optimization of the parameter $\gamma_1$ can eliminate the dominant error terms in phase in undamped systems. Consequently, the phase accuracy can be improved by three additional orders compared with the formal accuracy. On the other hand, for the four- and six-sub-step methods, the parameter optimization can improve the amplitude accuracy by three additional orders. The two-sub-step method is a special case: for $\gamma_1=1/4$ or $1/2$, it becomes exactly non-dissipative and the amplitude error vanishes in the undamped case.
\end{itemize}

\begin{table*}[htbp]
	\caption{Numerical performance among various $s$-sub-step implicit methods with $s$th-order accuracy.}
	{\small\begin{tabular*}{\textwidth}{@{\extracolsep\fill}cccccccc@{}}\toprule
			\multirow{2}{*}{Integrators} & \multirow{2}{*}{Range of $\gamma_1$} & \multicolumn{2}{c}{Amplitude error $\delta$} & \multicolumn{2}{c}{Phase error $\epsilon$} & \multirow{2}{*}{Stability} &  \multirow{2}{*}{Numerical dissipation $\infrho$} \\ \cmidrule(lr){3-4}\cmidrule(lr){5-6}
& & $\xi\neq0$ & $\xi=0$ & $\xi\neq0$ & $\xi=0$ & \\
			\midrule 
			$s=1$ & $1/2\le\gamma_1\le\infty$ & $\mathcal{O}(\dt)$ & $\mathcal{O}(\dt)$ & $\mathcal{O}(\dt)$ & $\mathcal{O}(\dt^{2})$ & unconditional stability & $\infrho\in(-1,~1]$ \\
			$s=2$ & $1/4\le\gamma_1\le\infty$ & $\mathcal{O}(\dt^{2})$ & $\mathcal{O}(\dt^{3})^\star$ & $\mathcal{O}(\dt^{2})$ & $\mathcal{O}(\dt^{2})$ & unconditional stability & $\infrho\in[-1,~1]$ \\
			$s=3$ & \cref{eq:s3_uc} & $\mathcal{O}(\dt^{3})$ & $\mathcal{O}(\dt^{3})$ & $\mathcal{O}(\dt^{3})$ & $\mathcal{O}(\dt^{4})^*$ & unconditional stability & $\infrho\in[1-\sqrt{3},~1]$ \\
			$s=4$ & \cref{eq:s4_uc} & $\mathcal{O}(\dt^{4})$ & $\mathcal{O}(\dt^{5})^\ddagger$ & $\mathcal{O}(\dt^{4})$ & $\mathcal{O}(\dt^{4})$ & unconditional stability & $\infrho\in[-0.63041,~1]$ \\
			$s=5$ & \cref{eq:s5_uc} & $\mathcal{O}(\dt^{5})$ & $\mathcal{O}(\dt^{5})$ & $\mathcal{O}(\dt^{5})$ & $\mathcal{O}(\dt^{6})^\dagger$ & unconditional stability & $\infrho\in[-1,~1]$ \\
			$s=6$ & \cref{eq:s6_uc} & $\mathcal{O}(\dt^{6})$ & $\mathcal{O}(\dt^{7})^\diamond$ & $\mathcal{O}(\dt^{6})$ & $\mathcal{O}(\dt^{6})$ & unconditional stability & $\infrho\in[-0.83733,~1]$ \\
			\bottomrule
	\end{tabular*}}
    \begin{tablenotes}
\footnotesize
\item[$\star$] The two-sub-step method with $\gamma_1=1/2$ or $1/4$ is non-dissipative ($\rho\equiv1$), so its amplitude error can reach $\mathcal{O}(\dt^\infty)$ in the undamped case. 
\item[$*$] The three-sub-step method with $\gamma_1\approx0.9756745886944403~(\infrho\approx-0.67851)$ can achieve sixth-order phase accuracy $\mathcal{O}(\dt^6)$ in the undamped case. 
\item[$\ddagger$] The four-sub-step method with $\gamma_1\approx 1.2805797612753055~(\infrho\approx -0.59615)$ can achieve seventh-order amplitude accuracy $\mathcal{O}(\dt^7)$ in the undamped case. 
\item[$\dagger$] The five-sub-step method with $\gamma_1\approx 0.4515512289893862~(\infrho\approx-0.91419)$ can achieve eighth-order phase accuracy $\mathcal{O}(\dt^8)$ in the undamped case. 
\item[$\diamond$] The six-sub-step method with $\gamma_1\approx 0.5409068780733081~(\infrho\approx -0.75569)$ can achieve ninth-order amplitude accuracy $\mathcal{O}(\dt^9)$ in the undamped case.
\end{tablenotes}
	\label{tab:sth}
\end{table*}

\begin{figure}[htbp]
	\centering 
	\includegraphics[scale=1.8]{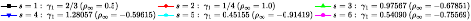}\\
    \subfigure[]{
		\includegraphics[scale=0.71]{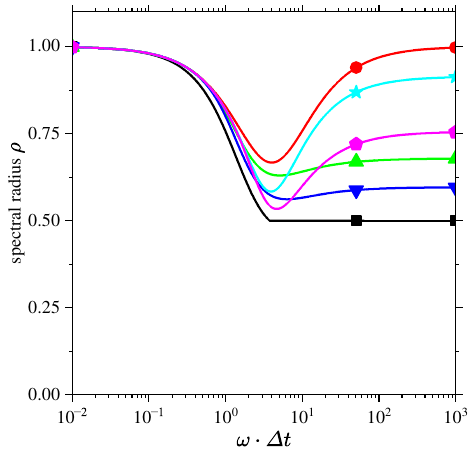}}
    \subfigure[]{
		\includegraphics[scale=0.72]{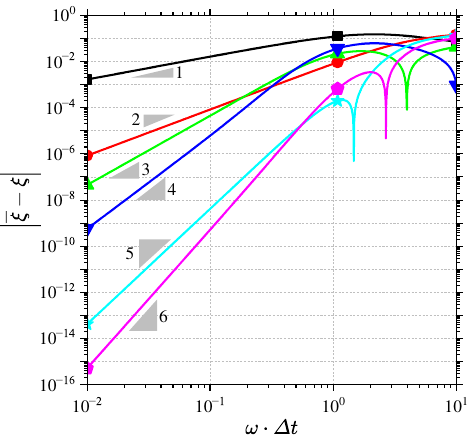}}
    \subfigure[]{
		\includegraphics[scale=0.71]{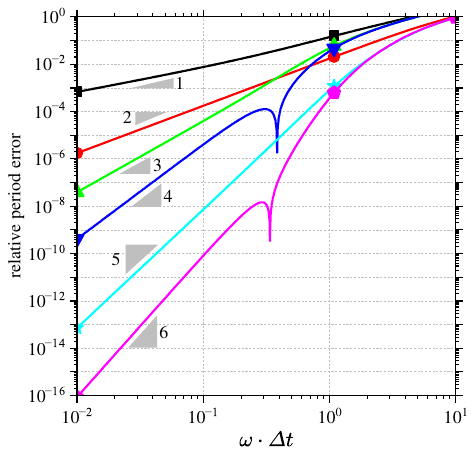}}
    \subfigure[]{
		\includegraphics[scale=0.71]{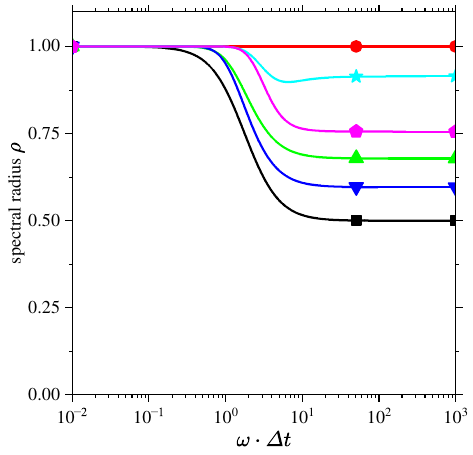}}
    \subfigure[]{
		\includegraphics[scale=0.72]{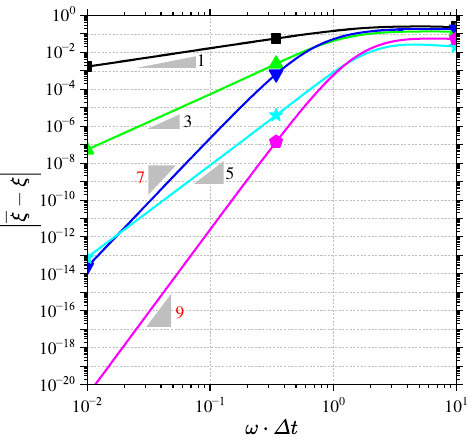}}
    \subfigure[]{
		\includegraphics[scale=0.71]{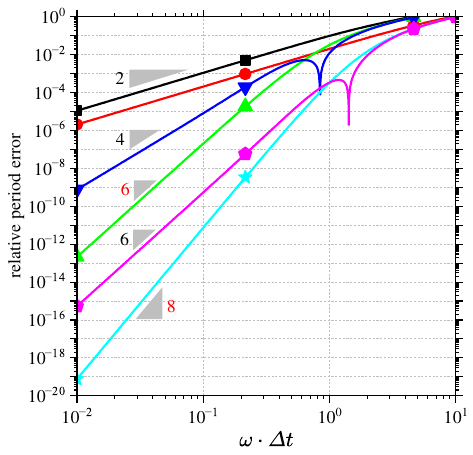}}
	\caption{Comparisons of spectral properties among the proposed $s$th-order $s$-sub-step implicit methods: (a-c) $\xi=0.2$ and (d-f) $\xi=0.0$.}
	\label{fig:com_sp}
\end{figure}

\cref{fig:com_sp} compares the spectral properties of the proposed $s$th-order $s$-sub-step implicit methods with different numbers of sub-steps. The corresponding parameters are selected according to the optimized cases listed in \cref{tab:sth}, where each method possesses enhanced accuracy in the undamped case. \cref{fig:com_sp}(a) and (d) show the spectral radius for damped and undamped cases, respectively, and all considered members maintain $\rho\le1$ over the entire frequency range, confirming their unconditional stability. It is also observed that the optimized parameters do not necessarily correspond to the strongest high-frequency dissipation; instead, they provide a favorable balance between numerical damping and low-frequency accuracy. \cref{fig:com_sp}(b) and (c) present the amplitude and phase errors in the presence of $\xi$, respectively. Since the physical damping ratio is nonzero, the error orders directly correspond to the formal accuracy of each method. Specifically, the amplitude and phase errors decrease with slopes consistent with the theoretical orders listed in \cref{tab:sth}. For undamped systems, \cref{fig:com_sp}(e) and (f) further demonstrate the superconvergence behavior predicted by the analytical error analysis. A clear distinction between odd and even numbers of sub-steps can be observed. For the three- and five-sub-step methods, the specially selected $\gamma_1$ eliminates the dominant phase errors, resulting in three additional orders of phase accuracy compared with the formal order. Consequently, the slopes in \cref{fig:com_sp}(f) become significantly steeper than those corresponding to the nominal accuracy. For example, the three- and five-sub-step methods achieve sixth- and eighth-order phase accuracy, respectively, which agrees well with the theoretical predictions in \cref{tab:sth}. In contrast, the two-, four-, and six-sub-step methods exhibit improved amplitude accuracy in undamped systems. By appropriately selecting $\gamma_1$, the leading amplitude errors can be eliminated, leading to three additional orders of amplitude accuracy compared with the formal order. This phenomenon is clearly observed in \cref{fig:com_sp}(e), where the optimized even-sub-step methods exhibit much faster decay of amplitude errors. In particular, the six-sub-step method achieves ninth-order amplitude accuracy, confirming the effectiveness of exploiting the parameter $\gamma_1$ to enhance the accuracy beyond the formal order.

These results in \cref{fig:com_sp} verify the theoretical errors summarized in \cref{tab:sth}. With increasing sub-step numbers, the proposed methods not only systematically increases the formal accuracy but also preserves the possibility of further improving either phase or amplitude accuracy through appropriate parameter selections. The alternating superconvergence behavior between odd- and even-sub-step methods originates from the different structures of the leading errors and provides an effective strategy for designing high-order implicit integration algorithms with improved spectral accuracy.

\begin{table*}[htbp]
	\caption{Numerical performance among various $s$-sub-step implicit methods with $(s+1)$th-order accuracy.}
	{\small\begin{tabular*}{\textwidth}{@{\extracolsep\fill}cccccccc@{}}\toprule
			\multirow{2}{*}{Integrators} & \multirow{2}{*}{$\gamma_1$} & \multicolumn{2}{c}{Amplitude error $\delta$} & \multicolumn{2}{c}{Phase error $\epsilon$} & \multirow{2}{*}{Stability} &  \multirow{2}{*}{Numerical dissipation $\infrho$} \\ \cmidrule(lr){3-4}\cmidrule(lr){5-6}
& &  $\xi\neq0$ & $\xi=0$ & $\xi\neq0$ & $\xi=0$ & \\
			\midrule 
			$s=1$ & $1/2$ & $\mathcal{O}(\dt^{2})$ & $\mathcal{O}(\dt^{3})$ & $\mathcal{O}(\dt^{2})$ & $\mathcal{O}(\dt^{2})$ & unconditional stability & $\infrho=1$ \\
			$s=2$ & $1/2+\sqrt{3}/6$ & $\mathcal{O}(\dt^{3})$ & $\mathcal{O}(\dt^{3})$ & $\mathcal{O}(\dt^{3})$ & $\mathcal{O}(\dt^{4})$ & unconditional stability & $\infrho=1-\sqrt{3}$ \\
			$s=3$ & $1/2+\sqrt{3}\cos(\pi/18)/3$ & $\mathcal{O}(\dt^{4})$ & $\mathcal{O}(\dt^{5})$ & $\mathcal{O}(\dt^{4})$ & $\mathcal{O}(\dt^{4})$ & unconditional stability & $\infrho\approx-0.63041$ \\ \hdashline
			$s=4$ & $1.3453664197803336$ & $\mathcal{O}(\dt^{5})$ & $\mathcal{O}(\dt^{5})$ & $\mathcal{O}(\dt^{5})$ & $\mathcal{O}(\dt^{6})$ & \makecell[c]{$A(\alpha)$-stability \\ $\alpha=89.99173663^\circ$} & $\infrho\approx-0.57677$ \\ \hdashline
			\multirow{3}{*}{$s=5$} & $0.4732683912582953$ & \multirow{3}{*}{$\mathcal{O}(\dt^{6})$} & \multirow{3}{*}{$\mathcal{O}(\dt^{7})$} & \multirow{3}{*}{$\mathcal{O}(\dt^{6})$} & \multirow{3}{*}{$\mathcal{O}(\dt^{6})$} & unconditional stability & $\infrho\approx-0.83733$ \\
             & $1.6206645011137183$ &  &  &  &  & \makecell[c]{$A(\alpha)$-stability \\ $\alpha=89.94077997^\circ$} & $\infrho\approx0.54360$ \\ \hdashline
			\multirow{3}{*}{$s=6$} & $0.5566999420873480$ & \multirow{3}{*}{$\mathcal{O}(\dt^{7})$} & \multirow{3}{*}{$\mathcal{O}(\dt^{7})$} & \multirow{3}{*}{$\mathcal{O}(\dt^{7})$} & \multirow{3}{*}{$\mathcal{O}(\dt^{8})$} & \makecell[c]{$A(\alpha)$-stability \\ $\alpha=89.99999995^\circ$} & $\infrho\approx-0.72080$ \\
             & $1.8951305922313202$ &  &  &  &  & \makecell[c]{$A(\alpha)$-stability \\ $\alpha=89.85071945^\circ$} & $\infrho\approx0.52105$ \\
			\bottomrule
	\end{tabular*}}\\
	\label{tab:s1th}
\end{table*}

\cref{tab:s1th} further demonstrates the capability of the proposed methods to enhance the accuracy by exploiting the sub-step sizes $\gamma_j~(j=1,~\cdots,~s)$. For each $s$-sub-step formulation, imposing one additional order condition uniquely determines the sub-step locations and increases the accuracy order from $s$ to $(s+1)$. However, this improvement is obtained at the expense of losing the independent control of high-frequency dissipation, because $\gamma_1$ is no longer freely related to the prescribed parameter $\infrho$. Therefore, the proposed $s$-sub-step method naturally provides two different algorithm configurations: one achieves $s$th-order accuracy with user-specified numerical dissipation, whereas the other achieves $(s+1)$th-order accuracy with a fixed dissipation level. It is also observed that increasing the number of sub-steps systematically improves both the attainable accuracy and the flexibility of the algorithm design. Specifically, the $s$-sub-step method can achieve $s$th-order accuracy while preserving the controllability of high-frequency dissipation, and the additional freedoms introduced by the sub-step splitting strategy enable the construction of $(s+1)$th-order members. Meanwhile, the stability gradually becomes more restrictive for the higher-order members, changing from unconditional stability for low-order formulations to $A(\alpha)$-stability for some high-order members. Nevertheless, the stability angles reported in \cref{tab:s1th} remain extremely close to $90^\circ$, indicating that these methods retain stability properties very close to $A$-stability.

Compared with the corresponding $s$th-order methods listed in \cref{tab:sth}, these higher-order members consume the sub-step sizes to further increase the formal accuracy order by one. As a consequence, the first sub-step size $\gamma_1$ is uniquely determined and the high-frequency dissipation can no longer be independently controlled. The amplitude and phase errors of these higher-order members follow the similar parity-dependent tendency, but with one additional order improvement. Specifically, for odd sub-step methods, the amplitude and phase errors in the damped case become $\mathcal{O}(\dt^{s+1})$, while the amplitude error in the undamped case exhibits an additional one-order superconvergence, resulting in $\mathcal{O}(\dt^{s+2})$. For even sub-step methods, the amplitude and phase errors are also improved to $\mathcal{O}(\dt^{s+1})$, whereas the phase error in the undamped case benefits from one-order superconvergence and reaches $\mathcal{O}(\dt^{s+2})$.

Therefore, the two groups of algorithms generated by the proposed framework exhibit a clear trade-off: the $s$th-order members provide controllable high-frequency dissipation through the user-specified parameter $\infrho$, while the $(s+1)$th-order members achieve enhanced temporal accuracy by sacrificing the dissipation control capability. 

\section{Numerical examples}\label{sec:examples}

Although the proposed framework is applicable to both first- and second-order transient problems, the numerical investigations in this section focus on second-order structural dynamics for brevity. Since both formulations share the same Butcher tableau and are constructed from the same accuracy, stability, and dissipation conditions, the theoretical results established previously are equally applicable to first-order transient problems. The second-order dynamical problems considered herein additionally enable the numerical dissipation and spatial--temporal dispersion properties of the proposed methods to be investigated.
Since the proposed framework naturally produces six classes of implicit methods with different numbers of sub-steps ($s=1,~\cdots,~6$), a comprehensive comparison with all existing methods would introduce excessive numerical data and obscure the main observations. Therefore, the numerical examples in this section mainly focus on revealing the performance variations among different $s$-sub-step members within the proposed framework.

All numerical simulations are performed on a 64-bit computer equipped with an Intel(R) Core(TM) Ultra 9 185H CPU, using the Julia programming language \cite{bezanson_JuliaFreshApproach_2017} (version 1.12.1). The three benchmark problems considered in this section have been widely adopted in the literature \cite{bathe_InsightImplicitTime_2012,li_NovelFamilyComposite_2020,li_NovelFamilyControllably_2019} for evaluating time integration algorithms. They possess representative characteristics of structural dynamics, including damping, external excitation, stiffness disparity, and practical applications. Therefore, these examples provide a reliable and reproducible platform for assessing the accuracy, stability, and numerical performance of the proposed methods.

\subsection{The damped and forced oscillator}

The damped and forced oscillator
\begin{equation}
\ddot{u}(t)+0.2\dot{u}(t)+u(t)=\exp(t)
\end{equation}
with $u(0)=\dot{u}(0)=1$ is first considered to examine the numerical accuracy of the proposed algorithms. In particular, this example aims to verify whether the theoretical order of accuracy can be achieved in practical computations, since the theoretical order and observed numerical order may not always coincide for some time integration algorithms \cite{zhang_OptimizationNsubstepComposite_2020}. The global error of the computed solutions is evaluated as
\begin{equation}\label{eq:globalError}
	\text{Error} = \left[\sum_{j=1}^{N}\left(x(t_j)-x_j \right)^2/\sum_{j=1}^{N}\left(x(t_j)\right)^2 \right]^{1/2},
\end{equation}
where $N$ denotes the total number of time steps in the analysis, and $x(t_j)$ and $x_j$ represent the exact and numerical solutions at time $t_j$, respectively.

\begin{figure}[htbp]
	\centering 
	\includegraphics[scale=1]{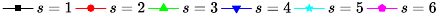}\\
    \subfigure[$\infrho=0.0$]{
		\includegraphics[scale=0.71]{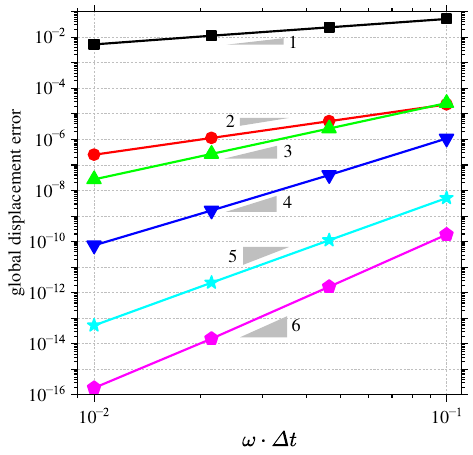}}
    \subfigure[$\infrho=0.5$]{
		\includegraphics[scale=0.71]{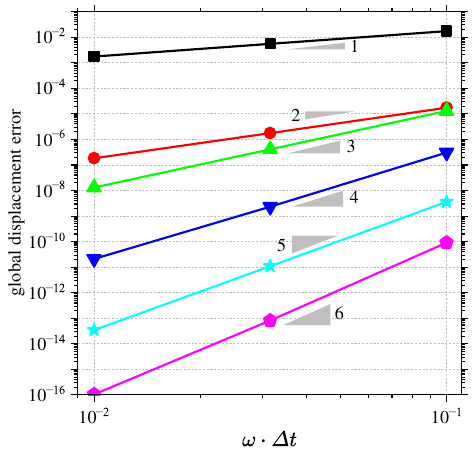}}
    \subfigure[$\infrho=1.0$]{
		\includegraphics[scale=0.71]{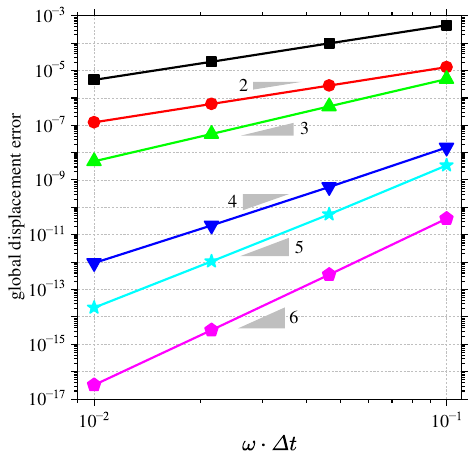}}
    \subfigure[$\infrho=0.0$]{
		\includegraphics[scale=0.71]{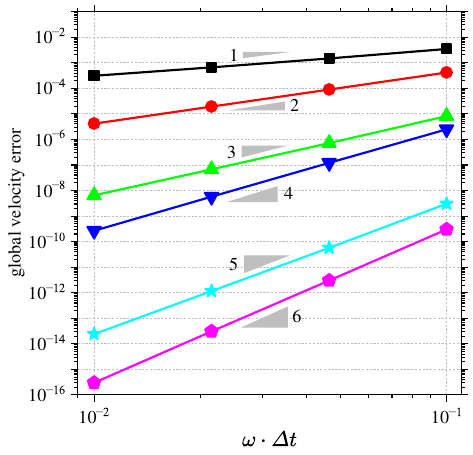}}
    \subfigure[$\infrho=0.5$]{
		\includegraphics[scale=0.71]{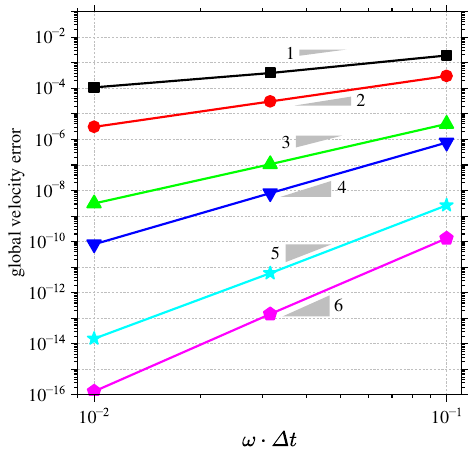}}
    \subfigure[$\infrho=1.0$]{
		\includegraphics[scale=0.71]{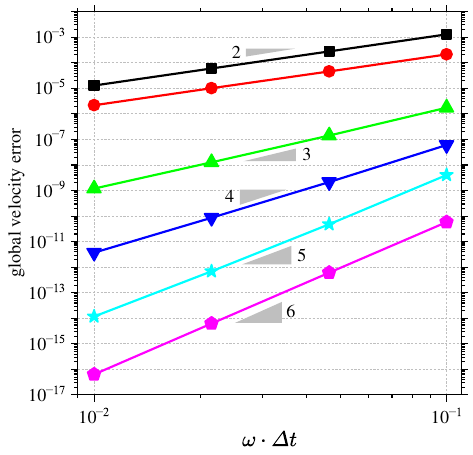}}
	\caption{Convergence ratios of the proposed $s$-sub-step implicit methods with $s$th-order accuracy: (a-c) displacement and (d-f) velocities.}
	\label{fig:con_algs}
\end{figure}

\cref{fig:con_algs} presents the convergence behaviors of the proposed $s$-sub-step implicit methods with $s$th-order accuracy. The displacement and velocity errors are evaluated under different $\infrho\in\{0,~0.5,~1\}$, including the dissipative cases ($\infrho=0$ and $0.5$) and the non-dissipative case ($\infrho=1$). As shown in \cref{fig:con_algs}(a)-(c), the displacement errors exhibit straight lines with slopes approximately equal to $s$ for the methods with $s=1,~\cdots,~6$. Similar convergence features can also be observed for the velocity errors in \cref{fig:con_algs}(d)-(f). These results confirm that the proposed $s$-sub-step implicit framework can achieve the theoretically derived $s$th-order accuracy for both displacement and velocity responses. It is also observed that the convergence rates are almost unaffected by the user-specified parameter $\infrho$. Although different values of $\infrho$ lead to different numerical dissipation, they only modify the high-frequency spectral behavior and do not change the formal order of accuracy except for $s=1$. Furthermore, the convergence behaviors remain consistent for different sub-step numbers. With the increase of $s$, the accuracy order increases correspondingly, demonstrating the systematic high-order constructions.

\begin{figure}[htbp]
	\centering 
	\includegraphics[scale=1.4]{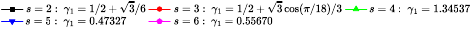}\\
    \subfigure[]{
		\includegraphics[scale=0.71]{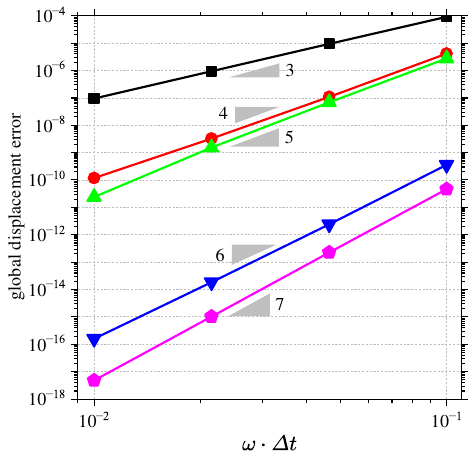}}
    \subfigure[]{
		\includegraphics[scale=0.71]{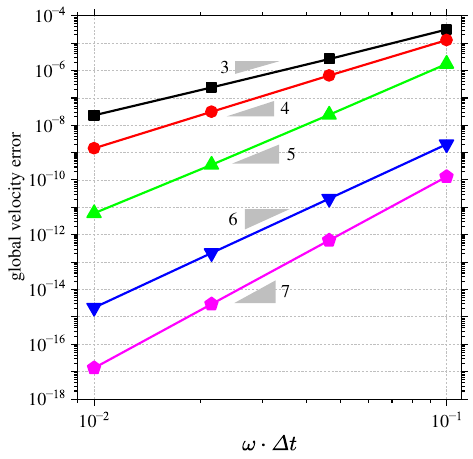}}
	\caption{Convergence ratios of the proposed $s$-sub-step implicit methods with $(s+1)$th-order accuracy.}
	\label{fig:con_algs_s1}
\end{figure}

\cref{fig:con_algs_s1} further investigates the convergence properties of the specially constructed $(s+1)$th-order members. Unlike the general $s$th-order members shown in \cref{fig:con_algs}, the sub-step parameters of these algorithms are completely determined and no longer provide the freely adjustable numerical dissipation. As shown in \cref{fig:con_algs_s1}(a), the displacement errors of all selected $(s+1)$th-order members converge with slopes approximately equal to $s+1$, confirming that the additional order condition successfully increases the temporal accuracy by one order. The same conclusion can be obtained from the velocity errors in \cref{fig:con_algs_s1}(b). Specifically, the two-, three-, four-, five-, and six-sub-step methods achieve third-, fourth-, fifth-, sixth-, and seventh-order accuracy, respectively. It is worth noting that the convergence curve of the fourth-order three-sub-step member in \cref{fig:con_algs_s1}(a) exhibits an apparent slope close to five within the considered time step range. This behavior, however, should not be interpreted as a general fifth-order convergence property. Instead, it is mainly attributed to the small error coefficient for the specific damped and forced oscillator considered in this example. The asymptotic accuracy of this method is still fourth order, as theoretically derived from the order conditions.

These results demonstrate that the additional parameters introduced by the sub-step sizes $\gamma_j~(j=1,~\cdots,~s)$ can be effectively exploited to construct higher-order members within the same $s$-sub-step method. Compared with the general $s$th-order algorithms, the $(s+1)$th-order members provide improved accuracy without increasing the number of sub-steps. However, because all sub-step sizes are fixed by the higher-order conditions, the high-frequency spectral radius $\infrho$ is no longer user-controllable, resulting in a fixed numerical dissipation level for each specific higher-order member.

\subsection{The model problem with stiff and flexible components}

\begin{figure}[htpb]
	\centering
	\includegraphics[scale=0.4]{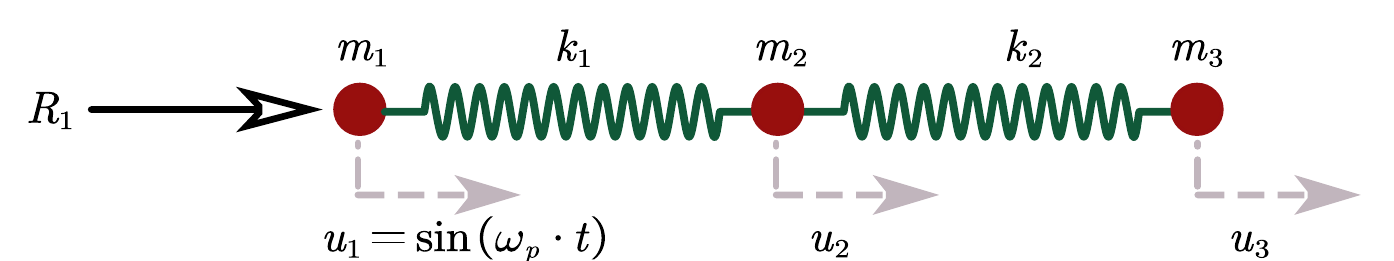}
	\caption{A classical mass-spring system \cite{li_HighorderAccurateMultisubstep_2024} with $k_1=10^7$N/m, $k_2=1$N/m, $m_1=0$kg, $m_2=m_3=1$kg, and $ \omega_\mathsf{p}=1.2$rad/s.}
	\label{fig:modalproblem}
\end{figure}

To evaluate the capability of the proposed algorithms in controlling high-frequency numerical oscillations, the classical stiff-flexible modal problem \cite{bathe_InsightImplicitTime_2012,li_HighorderAccurateMultisubstep_2024} shown in \cref{fig:modalproblem} is considered. This problem is a representative simplified model for structural systems containing components with significantly different stiffness characteristics, where the coexistence of low-frequency physical responses and high-frequency spurious modes provides a suitable benchmark for assessing numerical dissipation properties. The governing equation of the mass-spring system is given by
\begin{equation}\label{eq:modeproblem}
	\begin{bmatrix}
		m_1 & 0 & 0 \\ 0 & m_2 & 0\\ 0 & 0 & m_3\\
	\end{bmatrix}\begin{bmatrix}
		\ddot{u}_1 \\ \ddot{u}_2 \\ \ddot{u}_3\\
	\end{bmatrix}+\begin{bmatrix}
		k_1 & -k_1 & 0 \\ -k_1 & k_1+k_2 & -k_2\\ 0 & -k_2 & k_2\\
	\end{bmatrix}\begin{bmatrix}
		u_1 \\ u_2\\ u_3\\
	\end{bmatrix}=\begin{bmatrix}
		R_1 \\ 0 \\ 0\\
	\end{bmatrix}
\end{equation}
where the prescribed displacement at node 1 is defined as $u_1=\sin(\omega_\mathsf{p}t)=\sin(1.2t)$m, and $R_1$ denotes the corresponding reaction force. The initial conditions are assumed to be zero. By prescribing the displacement history at node 1, the original system can be reduced to
\begin{equation}
	\begin{bmatrix}
		m_2 & 0 \\ 0 & m_3\\
	\end{bmatrix}\begin{bmatrix}
		\ddot{u}_2 \\ \ddot{u}_3\\
	\end{bmatrix}+\begin{bmatrix}
		k_1+k_2 & -k_2 \\ -k_2 & k_2\\
	\end{bmatrix}\begin{bmatrix}
		u_2 \\ u_3\\
	\end{bmatrix}=\begin{bmatrix}
		k_1u_1 \\ 0\\
	\end{bmatrix}\label{eq:standardproblem}
\end{equation}
where $m_2=m_3=1$kg, $k_1=10^7$N/m, and $k_2=1$N/m are selected. For these parameters, the system contains two distinct vibration modes, with the frequencies approximately given by $\omega_1\approx1.0$rad/s and $\omega_2\approx3162.3$rad/s. The first mode represents the dominant low-frequency physical response, whereas the second mode corresponds to a high-frequency component caused by the stiff spring and is typically regarded as a spurious mode in numerical simulations.

Therefore, this problem is particularly suitable for examining the ability of the proposed methods to suppress undesired high-frequency responses while preserving the accuracy of low-frequency components. When both frequency modes are physically meaningful, non-dissipative algorithms ($\infrho=1$) together with sufficiently small time steps should be adopted to avoid excessive attenuation of high-frequency responses. Conversely, for practical structural dynamic simulations where high-frequency modes are mainly numerical artifacts, algorithms with appropriate high-frequency dissipation can effectively improve computational robustness.

\begin{figure}[htbp]
	\centering 
    \subfigure[]{
		\includegraphics[scale=1.03]{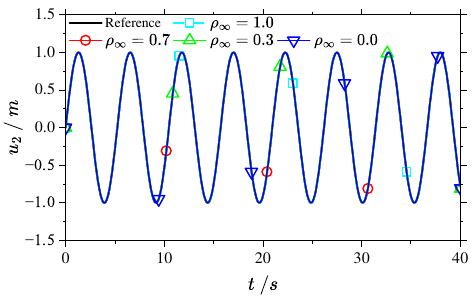}}
    \subfigure[]{
		\includegraphics[scale=1.0]{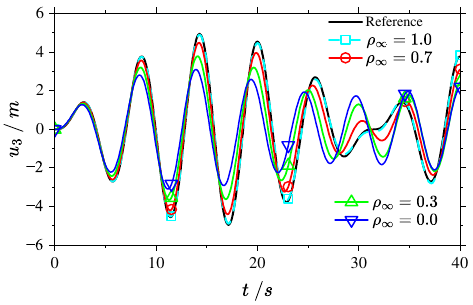}}\\
    \subfigure[]{
		\includegraphics[scale=1.0]{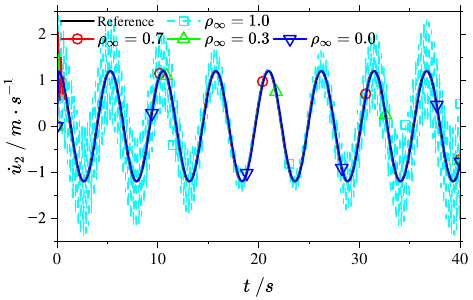}}
    \subfigure[]{
		\includegraphics[scale=1.0]{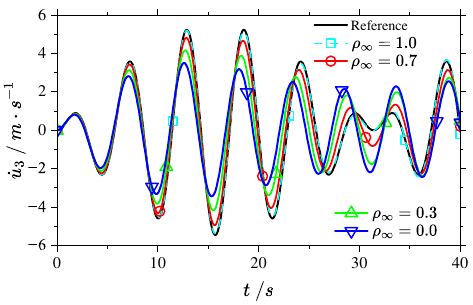}}
	\caption{Numerical response predicted by the single-sub-step implicit method ($s=1$) with $\dt=0.1309$s for solving \cref{eq:standardproblem}.}
	\label{fig:ex2_s1}
\end{figure}

\cref{fig:ex2_s1} presents the numerical responses obtained by the single-sub-step implicit method with different values of $\infrho$, where the influence of the high-frequency dissipation on the stiff-flexible system is investigated. Since the influence of $\infrho$ on the suppression of spurious high-frequency responses is essentially consistent for different $s$-sub-step methods, only the single-sub-step case is presented herein to avoid redundant numerical results. As shown in \cref{fig:ex2_s1}, the value of $\infrho$ plays a crucial role in controlling the high-frequency component associated with the stiff mode. When $\infrho=0$, the high-frequency oscillation is effectively eliminated, while the low-frequency physical response remains accurately preserved. With increasing $\infrho$, the numerical dissipation becomes weaker, and more high-frequency information is retained. In particular, the case of $\infrho=1$ corresponds to the non-dissipative algorithm, for which the spectral radius approaches unity in the high-frequency region. Consequently, the spurious high-frequency oscillation cannot be suppressed and remains visible in the numerical response, as shown in \cref{fig:ex2_s1}(c). It is also worth noting that the single-sub-step method with $\infrho=1$ becomes a second-order accurate method, as predicted by the theoretical analysis, and therefore exhibits improved accuracy compared with its general first-order formulation. These results demonstrate that the parameter $\infrho$ provides an effective mechanism for balancing high-frequency dissipation and solution accuracy.

\begin{figure}[htbp]
	\centering 
    \subfigure[]{
		\includegraphics[scale=1.03]{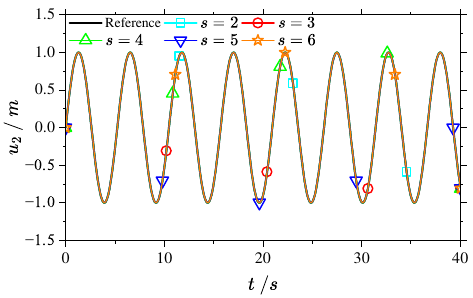}}
    \subfigure[]{
		\includegraphics[scale=1.0]{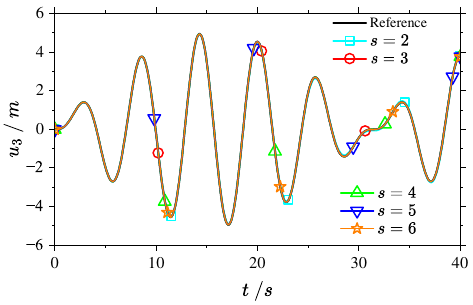}}\\
    \subfigure[]{
		\includegraphics[scale=1.0]{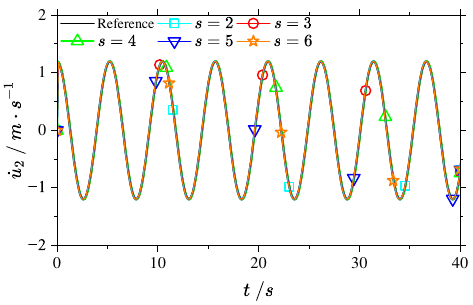}}
    \subfigure[]{
		\includegraphics[scale=1.0]{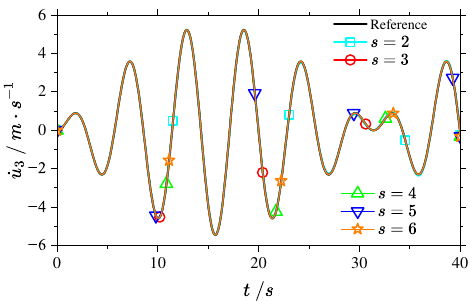}}\\
    \subfigure[]{
		\includegraphics[scale=1.0]{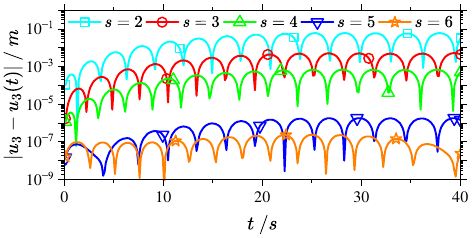}}
    \subfigure[]{
		\includegraphics[scale=1.0]{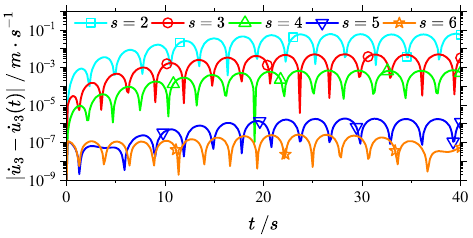}}
	\caption{Numerical response predicted by the $s$th-order $s$-sub-step implicit method with $\infrho=0$ and $\dt=0.1309$s for solving \cref{eq:standardproblem}.}
	\label{fig:ex2_s}
\end{figure}

The numerical responses of the proposed $s$-sub-step implicit methods with $s$th-order accuracy are compared in \cref{fig:ex2_s}. In this comparison, all methods are assigned $\infrho=0$ to provide the same high-frequency dissipation level and to highlight the influence of the sub-step number on the numerical performance. As shown in \cref{fig:ex2_s}(a-d), all proposed methods accurately capture the low-frequency physical responses of the system, while the high-frequency spurious oscillations are effectively removed. With increasing $s$, the numerical solutions become increasingly close to the reference solution, especially for the responses associated with the flexible component. To better distinguish the differences among various algorithms, the absolute errors of the displacement and velocity responses at node 3 are further presented in \cref{fig:ex2_s}(e-f). It can be observed that increasing the number of sub-steps continuously reduces numerical errors. This improvement is attributed to the higher-order temporal accuracy achieved by the proposed $s$-sub-step framework. In particular, the higher-order methods not only preserve the low-frequency response with greater accuracy but also maintain effective suppression of the undesired high-frequency components. These results confirm the advantage of increasing the sub-step number in simultaneously improving accuracy and controlling numerical oscillations.

\begin{figure}[htbp]
	\centering 
    \subfigure[]{
		\includegraphics[scale=1.03]{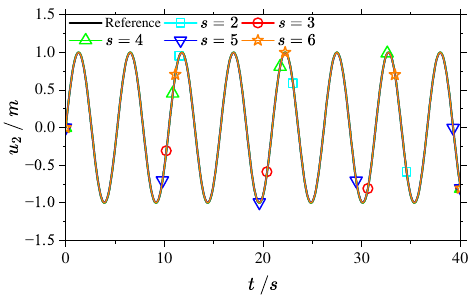}}
    \subfigure[]{
		\includegraphics[scale=1.0]{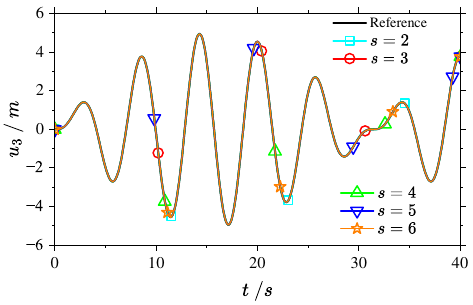}}\\
    \subfigure[]{
		\includegraphics[scale=1.0]{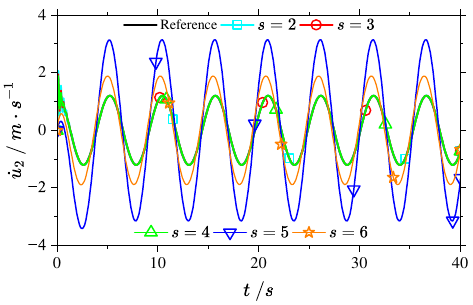}}
    \subfigure[]{
		\includegraphics[scale=1.0]{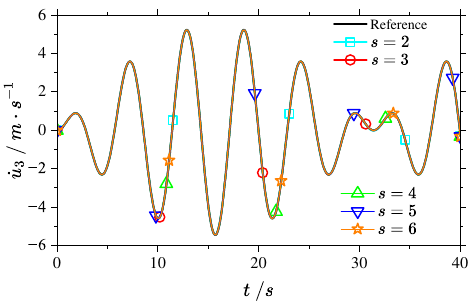}}\\
    \subfigure[]{
		\includegraphics[scale=1.0]{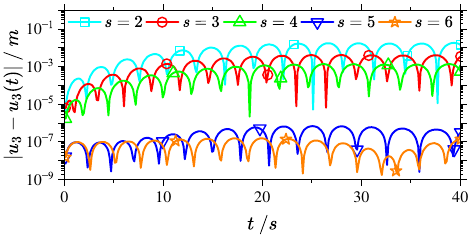}}
    \subfigure[]{
		\includegraphics[scale=1.0]{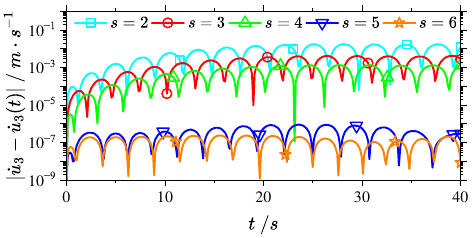}}
	\caption{Numerical response predicted by the $(s+1)$th-order $s$-sub-step implicit method with $\dt=0.1309$s for solving \cref{eq:standardproblem}.}
	\label{fig:ex2_ss1}
\end{figure}

\cref{fig:ex2_ss1} presents the numerical responses of the proposed $s$-sub-step implicit methods with $(s+1)$th-order accuracy for solving \cref{eq:standardproblem}. Compared with the corresponding $s$th-order methods in \cref{fig:ex2_s}, the additional order of accuracy is achieved by appropriately selecting the sub-step parameters $\gamma_j~(j=1,\cdots,s)$. As shown in \cref{fig:ex2_ss1}(a-d), almost all $(s+1)$th-order methods accurately capture the low-frequency physical response while effectively eliminating the spurious high-frequency oscillations. It is worth noting that the velocity response at node 2 predicted by the five- and six-sub-step methods exhibits relatively larger amplitude deviations in \cref{fig:ex2_ss1}(c). This behavior is mainly associated with the selection of the sub-step sizes. For the $(s+1)$th-order methods, the sub-step parameters are determined by satisfying additional accuracy constraints, which reduces the available parameters for simultaneously optimizing other spectral properties, such as high-frequency dissipation. If the additional order condition is relaxed, the five- and six-sub-step methods can adopt alternative sub-step splitting strategies to predict more accurate results for $\dot{u}_2$. 

Similar to the previous comparisons, the absolute errors of the displacement and velocity responses at node 3 are plotted in \cref{fig:ex2_ss1}(e-f) to distinguish the numerical performance among different algorithms. The results indicate that the $(s+1)$th-order methods generally provide smaller errors than their corresponding $s$th-order counterparts, which agrees well with the theoretical analysis. In particular, the six-sub-step method achieves the smallest overall errors among the investigated algorithms, confirming the advantage of increasing the number of sub-steps together with additional order conditions.

It should be noted that the improvement in numerical accuracy observed with increasing the number of sub-steps is obtained under the condition of using the same time step size. As demonstrated in \cref{fig:ex2_s,fig:ex2_ss1}, increasing $s$ leads to higher-order accuracy and consequently provides more accurate numerical predictions. However, this improvement is achieved at the expense of increased computational cost, since a larger number of sub-steps requires more intermediate evaluations within each time step. Therefore, the selection of the number of sub-steps should consider the balance between the accuracy and computational efficiency. The proposed $s$-sub-step methods provides a flexible way to achieve higher accuracy, while allowing users to select an appropriate sub-step number according to the specific characteristics of the dynamical problem. 

\begin{figure}[htbp]
	\centering 
    \includegraphics[scale=1.8]{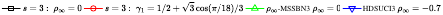}\\
    \subfigure[]{
		\includegraphics[scale=1.0]{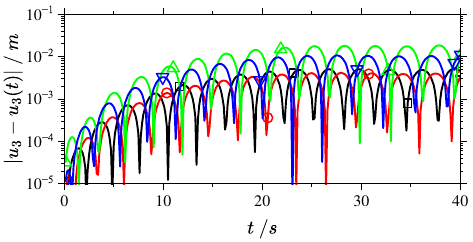}}
    \subfigure[]{
		\includegraphics[scale=1.0]{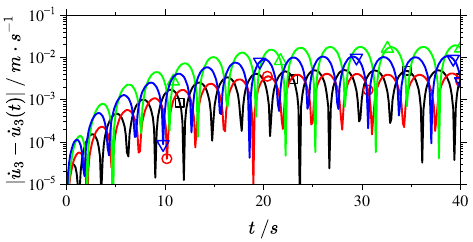}}
	\caption{Numerical response predicted by some three-sub-step implicit methods with $\dt=0.1309$s for solving \cref{eq:standardproblem}.}
	\label{fig:ex2_com}
\end{figure}
Finally, to further demonstrate the advantages of the proposed methods, the three-sub-step implicit method is selected as an example and compared with two existing three-sub-step algorithms \cite{ji_AccurateControllablyDissipative_2023,lee_ImplicitSsubstepTime_2025} in \cref{fig:ex2_com}. The absolute errors of the displacement and velocity responses at node 3 are employed for quantitative comparison. As shown in \cref{fig:ex2_com}, the proposed three-sub-step methods exhibit smaller numerical errors than the compared algorithms while maintaining effective attenuation of the high-frequency spurious mode.


\subsection{An elastic bar with axial tip loading}

\begin{figure}[htbp] 
	\centering 
	\includegraphics[scale=0.5]{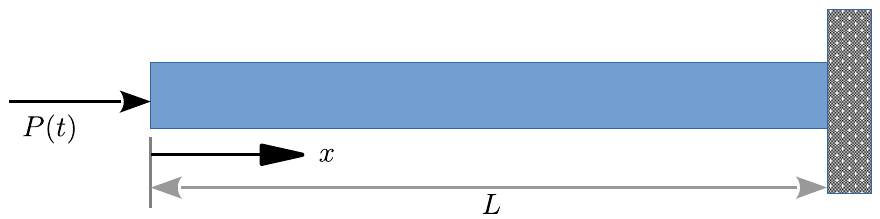}
	\caption{An elastic bar with axial tip loading.}
	\label{fig:2dof}
\end{figure}

The previous examples mainly focus on mathematically dynamic systems. However, in practical dynamical analyses, the governing partial differential equations (PDEs) are usually transformed into semi-discrete ordinary differential equations (ODEs) through spatial discretization \cite{hughes_FiniteElementMethod_2000}. Such discretization procedures inevitably introduce additional high-frequency numerical modes, which may affect the accuracy and stability of time integration algorithms. Therefore, an elastic bar \cite{li_HighorderAccurateMultisubstep_2024} subjected to an axial tip load $P(t)$, as illustrated in \cref{fig:2dof}, is considered to further evaluate the performance of the proposed methods for spatially discretized problems. For continuum dynamic problems, the final numerical responses are affected by both spatial and temporal discretization errors. A reliable evaluation of time integration algorithms for spatially discretized systems should consider the combined influence of these two error sources, as discussed in previous studies \cite{li_DirectlySelfstartingHigherorder_2022,li_UnifiedDispersionError_2026,li_FurtherAssessmentThree_2021}. 

The longitudinal vibration of the bar is governed by
\begin{subequations}
	\begin{align} 
		\rho A\dfrac{\partial^2u(x,~t)}{\partial t^2}-EA\dfrac{\partial^2u(x,~t)}{\partial x^2}&=0 \qquad\ \text{for}\quad 0<x<L,\quad t>0\\
		-EA\dfrac{\partial u(x,~t)}{\partial x}&=P(t)\quad \text{at}\quad x=0\\
		u(x,~t)&=0\qquad\ \text{at}\quad x=L\\
		u(x,~t)=\dfrac{\partial u(x,~t)}{\partial t}&=0\qquad\ \text{at}\quad t=0
	\end{align}
\end{subequations}
where $ A=1.0 $in$ ^2 $, $ E=3.0\times 10^7 $psi, $ \rho=7.4\times 10^{-4} $lb-s$ ^2 $/in$^4$, $ L=20 $in, and $ P(t) =100$lb. The analytical solution of the bar response is available and is adopted as the reference solution for evaluating the numerical accuracy of the proposed integration methods. 

\subsubsection{Spatial discretization with two-node linear elements} 

For the spatial discretization, 200 linear finite elements with equal length are adopted. The corresponding elemental mass and stiffness matrices are given, respectively, as
\begin{subequations}
	\begin{equation}\label{eq:me}
		\mbf{m}^e=(1-r)\cdot\frac{\rho A\dx}{2}\begin{bmatrix}
			1 & 0\\ 0 & 1
		\end{bmatrix}+r\cdot\dfrac{\rho A\dx}{6}\begin{bmatrix}
		2 & 1\\ 1 & 2
		\end{bmatrix}=\dfrac{\rho A\dx}{6}\begin{bmatrix}
		3-r & r\\ r & 3-r 
		\end{bmatrix}
	\end{equation}
	and 
	\begin{equation}
		\mbf{k}^e=\frac{EA}{\dx}\begin{bmatrix}
			1 & -1 \\ -1 & 1
		\end{bmatrix}
	\end{equation}
\end{subequations} 
where $\dx$ denotes the length of two-node finite elements. The parameter $r$ in \cref{eq:me} controls the distribution of the elemental mass matrix between lumped and consistent forms. Specifically, $r=0$, $r=1/2$, and $r=1$ correspond to the lumped, high-order, and consistent mass matrices, respectively. Following the dispersion analysis developed in the literature \cite{li_UnifiedDispersionError_2026,li_FurtherAssessmentThree_2021}, the relationship between the numerical frequency $\overline{\omega}$ and the exact frequency $\omega$ can be obtained for the proposed implicit methods. Since the derivation procedure follows the standard dispersion analysis\cite{li_UnifiedDispersionError_2026,li_FurtherAssessmentThree_2021}, only the resulting analytical expressions are presented herein without repeating the detailed derivations.
\begin{subequations}
\begin{align}
\overline{\omega}_\mathrm{s=1}=&\omega+\left[ \dfrac{(2r-1)\dx^{2}}{24c_0^{2}}-\left( \gamma_1^{2}-\gamma_1+ \dfrac{1}{3}  \right) \dt^{2}   \right]\omega^{3}+\mathcal{O}(\omega^5)\label{eq:omega_2node_1}\\
\overline{\omega}_\mathrm{s=2}=&\omega+\left[ \dfrac{(2r-1)\dx^{2}}{24c_0^{2}}+ \dfrac{\vartheta_1}{6}  \dt^{2}   \right]\omega^{3}+\left[  \dfrac{\vartheta_1(2r-1)\dt^{2}\dx^{2}}{48c_0^{2}} +\dfrac{20r^{2}-20r+1}{1920c_0^4}\dx^4- \dfrac{\vartheta_3}{20}\dt^4    \right]\omega^5+\mathcal{O}(\omega^7)\label{eq:omega_2node_2}\\
\overline{\omega}_\mathrm{s=3}=&\omega+ \dfrac{(2r-1)\dx^{2}}{24c_0^{2}}\omega^3+\left[ \dfrac{20r^{2}-20r+1}{1920c_0^4}\dx^4+ \dfrac{\phi_2}{30}\dt^4  \right]\omega^5\notag\\
&+ \bigg[ \dfrac{1400r^{3}-2100r^{2}+546r-3}{967680c_0^6}\dx^6+ \dfrac{(2r-1)\phi_2\dx^{2}\dt^{4}}{144c_0^{2}}- \dfrac{\phi_4}{252}\dt^6 \bigg]\omega^7+\mathcal{O}(\omega^9)\label{eq:omega_2node_3}\\
\overline{\omega}_\mathrm{s=4}=&\omega+ \dfrac{(2r-1)\dx^{2}}{24c_0^{2}}\omega^3+\left[ \dfrac{20r^{2}-20r+1}{1920c_0^4}\dx^4- \dfrac{\varpi_1}{120}\dt^4  \right]\omega^5\notag\\
&+ \bigg[ \dfrac{1400r^{3}-2100r^{2}+546r-3}{967680c_0^6}\dx^6- \dfrac{(2r-1)\varpi_1\dx^{2}\dt^{4}}{576c_0^{2}}+\dfrac{\varpi_3}{336}\dt^6 \bigg]\omega^7+\mathcal{O}(\omega^9)\label{eq:omega_2node_4}\\
\overline{\omega}_\mathrm{s=5,6}=&\omega+ \dfrac{(2r-1)\dx^{2}}{24c_0^{2}}\omega^3+ \dfrac{20r^{2}-20r+1}{1920c_0^4}\dx^4\omega^5+\mathcal{O}(\omega^7)\label{eq:omega_2node_56}
\end{align}
\end{subequations}
where $\vartheta_{1,3}$, $\phi_{2,4}$ and $\varpi_{1,3}$ are given by \cref{eq:s2_amp_phas,eq:s3_amp_phas,eq:s4_amp_phase}, respectively. 

\paragraph{The single-sub-step implicit method}
As shown in \cref{eq:omega_2node_1}, the leading dispersion error for $s=1$ consists of two independent contributions. The first term, ${(2r-1)\dx^{2}}/{(24c_0^2)}$, is introduced by the spatial discretization, where the parameter $r$ controls the form of the elemental mass matrix. The second term, $-(\gamma_1^2-\gamma_1+1/3)\dt^2$, represents the contribution from the time integration scheme with $s=1$. Therefore, the numerical frequency error is governed by the interaction between the spatial and temporal parameters. It can be observed that, for a fixed spatial discretization, the parameter $\gamma_1$ directly affects the leading dispersion error. \cref{eq:omega_2node_1} indicates that the spatial and temporal discretization errors can compensate each other. Specifically, by properly selecting the time step size, the leading $\omega^3$ term can be eliminated, resulting in a significantly improved dispersion property. Introducing the CFL number $c_0\dt~/\dx$, the optimal CFL number is computed for $s=1$ as 
\begin{equation}\label{eq:s1_cfl_opt}
\mathsf{CFL}_\mathsf{Opt}= \dfrac{1}{4}\sqrt{\dfrac{2(2r-1)}{3\gamma_1^{2}-3\gamma_1+1} } 
\end{equation}
where $\gamma_1=1/(\infrho+1)$. 

\begin{figure}[htbp] 
	\centering 
	\includegraphics[scale=2.0]{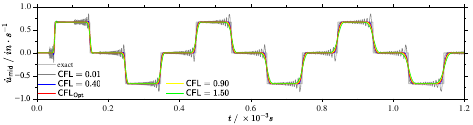}
	\caption{Velocities at the mid-point of the bar predicted by the single-sub-step implicit method ($s=1$) with $\infrho=0.8$.}
	\label{fig:bar_s1_CFLs}
\end{figure}

Herein, the case with $\infrho=0.8$ and $r=1$ is first considered, and the corresponding optimal CFL number predicted by \cref{eq:s1_cfl_opt} is adopted as a reference. As shown in \cref{fig:bar_s1_CFLs}, different CFL numbers lead to noticeably different numerical responses. It is interesting to observe that reducing the CFL number does not necessarily improve the accuracy of the numerical solution. For example, when a very small CFL number is adopted, the time step size becomes smaller, but stronger high-frequency oscillations are observed in the velocity response. This phenomenon is fundamentally different from the conventional understanding that reducing the time step size always improves the accuracy of time integration. For spatially discretized models, the numerical solution represents the response of the semi-discrete finite element model rather than the original PDE. Therefore, reducing the time step size only decreases the temporal discretization error, while the spurious high-frequency components introduced by the spatial discretization may become more apparent. Consequently, for spatially discretized models, an excessively small time step does not necessarily lead to a more accurate approximation of the original PDE solution.

\begin{figure}[htbp]
	\centering 
    \includegraphics[scale=1.0]{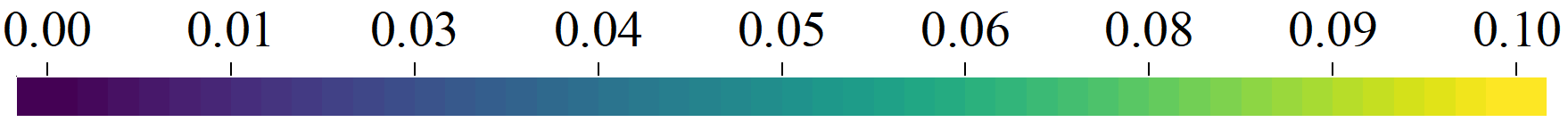}\\
    \subfigure[$\infrho=0.8$]{
		\includegraphics[scale=1.0]{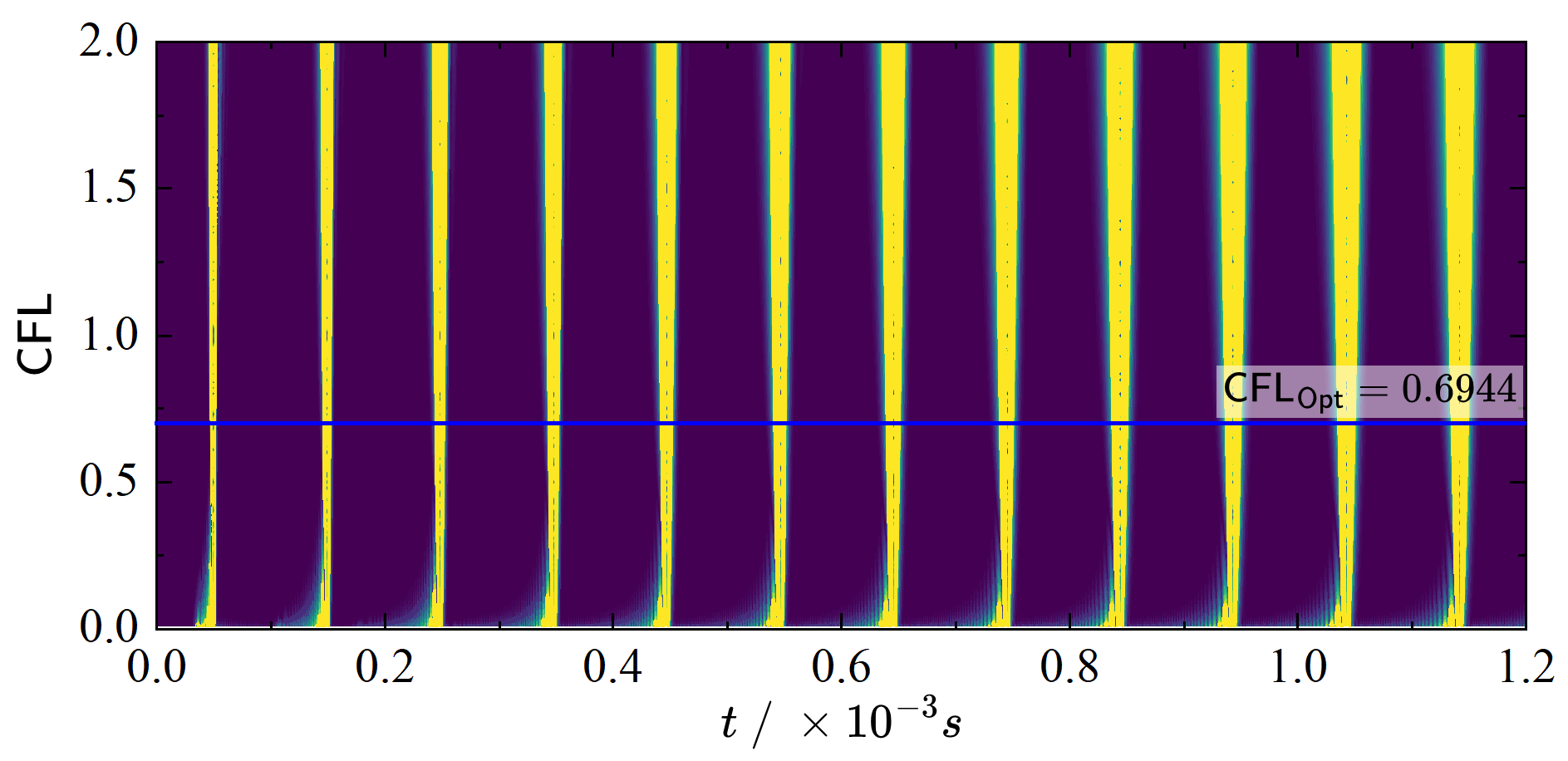}\label{fig:bar_s1_rho_a}}
    \subfigure[$\infrho=1.0$]{
		\includegraphics[scale=1.0]{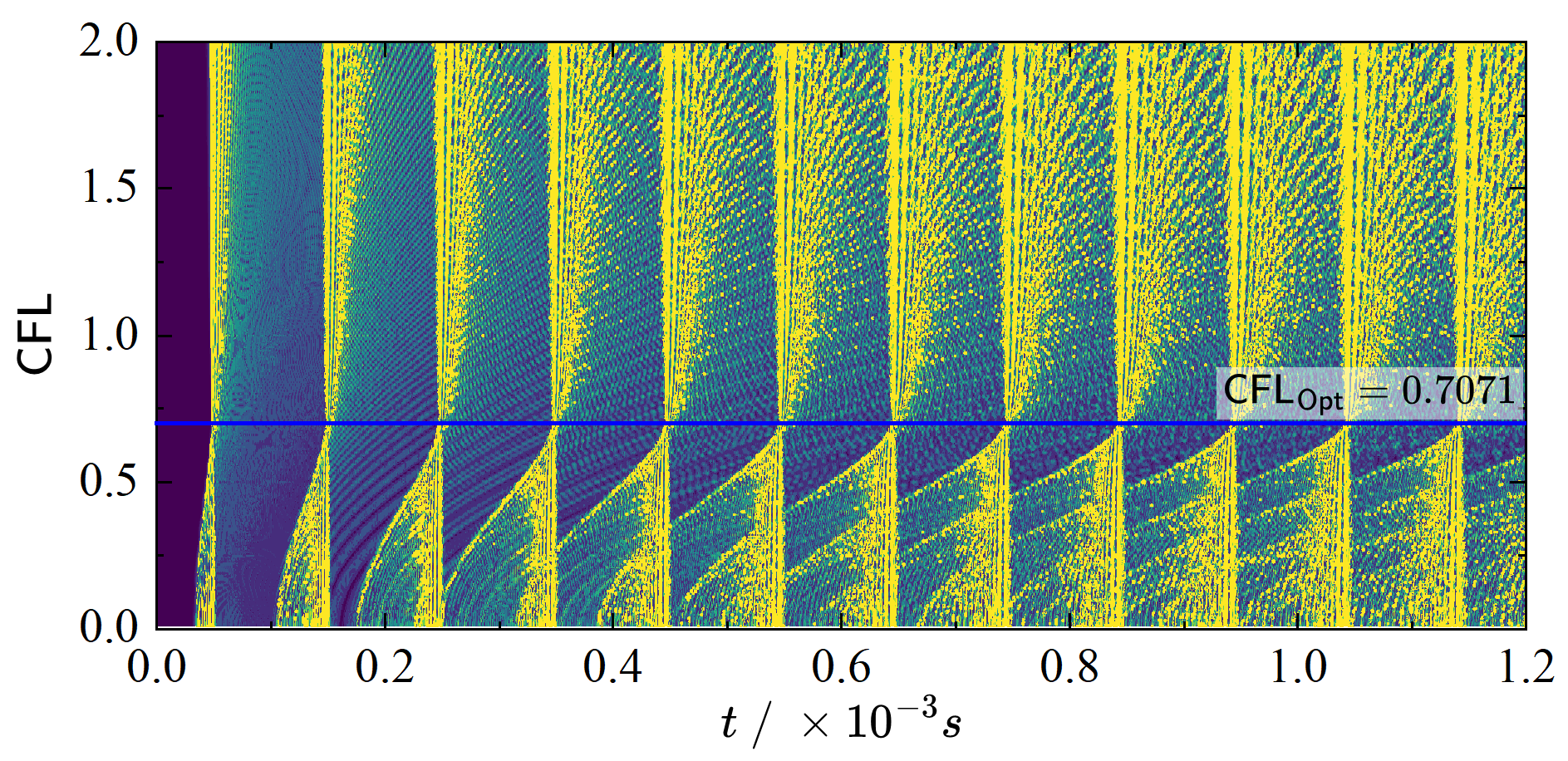}\label{fig:bar_s1_rho_b}}
	\caption{Errors of the mid-point velocity of the bar predicted by the single-sub-step implicit method ($s=1$).}
\end{figure}

For the case with $\infrho=0.8$, the numerical dissipation of the single-sub-step implicit method ($s=1$) is relatively strong due to its first-order accuracy. As a result, the advantage of the theoretically optimal CFL number is not clearly reflected in \cref{fig:bar_s1_CFLs,fig:bar_s1_rho_a}. In fact, the solution obtained with $\mathrm{CFL}=0.4$ appears to be slightly more accurate than that obtained using the theoretically optimal value. This is because the strong algorithmic dissipation suppresses not only the spurious high-frequency components but also partially attenuates the physical wave components, thereby masking the dispersion improvement achieved through the optimal CFL selection.

To further demonstrate the effect of the optimal CFL condition, the non-dissipative case with $\infrho=1.0$ is examined. As shown in \cref{fig:bar_s1_rho_b}, when numerical dissipation is absent, the influence of dispersion error becomes dominant. The error distribution clearly exhibits a minimum around the optimal CFL value, confirming that the leading error can be effectively eliminated through the appropriate coupling of spatial and temporal discretization parameters. This result verifies the theoretical analysis and demonstrates that the optimal CFL number provides significant benefits for reducing dispersion errors.

\begin{figure}[htbp]
	\centering 
        \includegraphics[scale=1.0]{bar_s1_leg}\\
    \subfigure[$r=1~\&~\gamma_1\approx0.25658$]{
		\includegraphics[scale=1.0]{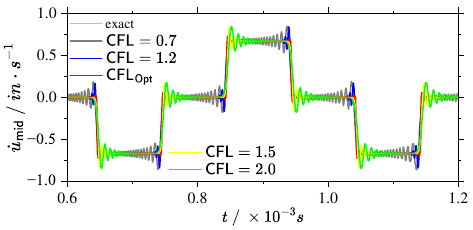}}
    \subfigure[$r=1~\&~\gamma_1\approx0.25658$]{
        \includegraphics[scale=1.0]{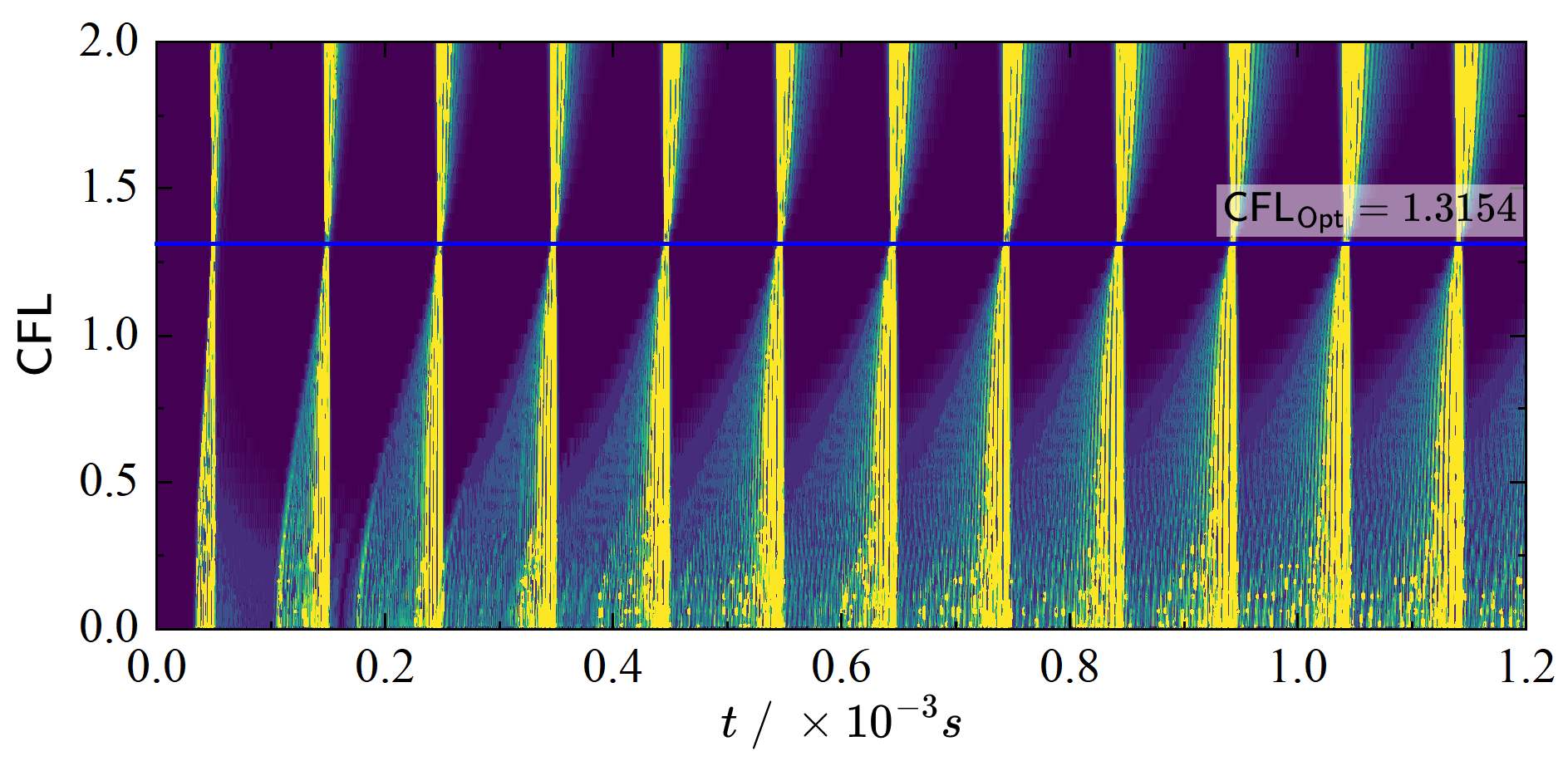}}
    \subfigure[$r=0~\&~\gamma_1\approx9.74342$]{
        \includegraphics[scale=1.0]{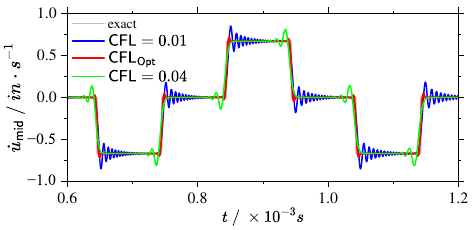}}
    \subfigure[$r=0~\&~\gamma_1\approx9.74342$]{
		\includegraphics[scale=1.0]{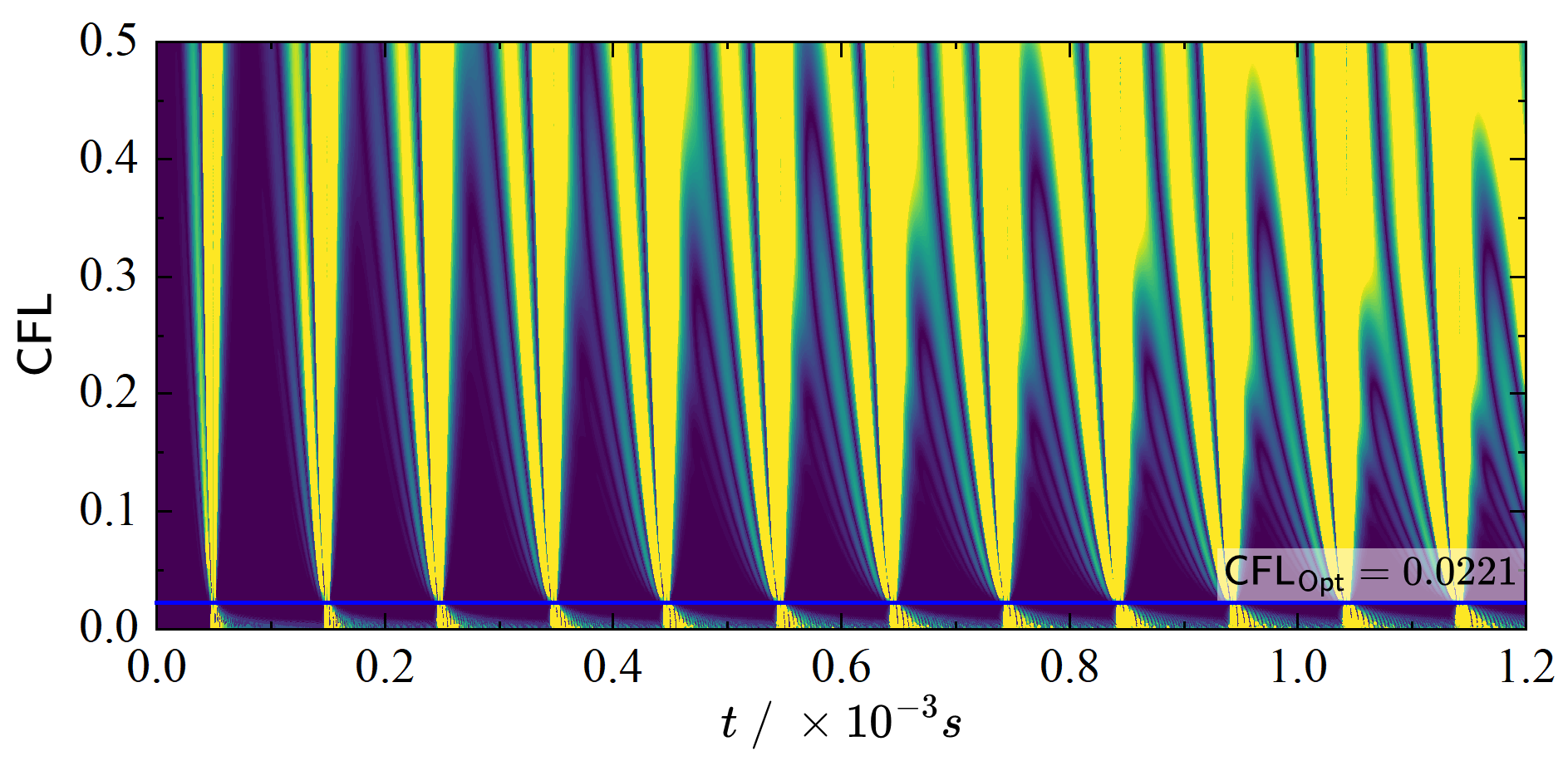}}
	\caption{The mid-point velocity of the bar predicted by the second-order two-sub-step implicit method ($s=2$) with $\infrho=0.8$.}
    \label{fig:bar_s2_CFLs}
\end{figure}

\paragraph{The two-sub-step implicit method}
For the second-order member, $\vartheta_1=6\gamma_1^2-6\gamma_1+1$ is generally nonzero. By eliminating the leading $\omega^3$ term in \cref{eq:omega_2node_2}, the optimal CFL number can be obtained as
\begin{equation}
\mathsf{CFL}_{\mathsf{Opt}}=\dfrac{1}{2}\sqrt{\dfrac{1-2r}{6\gamma_1^2-6\gamma_1+1} } 
\end{equation}
where the values of $r$ and $\gamma_1$ should ensure that $\mathsf{CFL}_{\mathsf{Opt}}$ is real and positive. For example, when $r>1/2$, the admissible range of $\gamma_1$ is $\gamma_1\in\left[1/4,~1/2+\sqrt{3}/6\right)$, whereas $\gamma_1>1/2+\sqrt{3}/6$ is required for $r<1/2$. Under these conditions, the spatial and temporal dispersion errors have opposite signs and can compensate each other, leading to an improved numerical accuracy. The numerical results of the second-order two-sub-step implicit method are presented in \cref{fig:bar_s2_CFLs} to verify the optimal CFL analysis. The effectiveness of $\mathsf{CFL}_{\mathsf{Opt}}$ is expected to depend on the adopted spatial discretization parameters, particularly the mass weighting parameter $r$.

For the consistent mass matrix case ($r=1$), as shown in \cref{fig:bar_s2_CFLs}(a)-(b), the optimal value $\mathsf{CFL}_{\mathsf{Opt}}\approx1.3154$ provides the most accurate numerical response among the tested CFL numbers. Compared with $\mathsf{CFL}=0.7$, $1.2$, $1.5$, and $2.0$, the solution obtained with the optimal CFL number exhibits significantly reduced oscillations around the discontinuities and better agreement with the exact solution. This confirms that the spatial and temporal dispersion errors can effectively compensate each other when the CFL number is appropriately selected. The error distribution shown in \cref{fig:bar_s2_CFLs}(b) further demonstrates that the minimum error occurs around the optimal CFL value.

For the lumped mass matrix case ($r=0$), the same theoretical analysis leads to a much smaller optimal CFL number, namely $\mathsf{CFL}_{\mathsf{Opt}}\approx0.0221$, as shown in \cref{fig:bar_s2_CFLs}(c)-(d). The numerical results indicate that this optimal value also achieves the best agreement with the exact solution. In contrast, when a larger CFL number ($\mathsf{CFL}=0.04$) is adopted, additional dispersion oscillations appear near the wave fronts, whereas a very small CFL number ($\mathsf{CFL}=0.01$) does not provide further accuracy improvement. The error map in \cref{fig:bar_s2_CFLs}(d) clearly identifies the optimal region predicted by the theoretical analysis.

For the third-order two-sub-step method, the parameter is selected as $\vartheta_1=0$, corresponding to $\gamma_1=1/2+\sqrt{3}/6$. In this case, the leading temporal contribution to the dispersion error vanishes, and \cref{eq:omega_2node_2} reduces to
\begin{equation}
\overline{\omega}_{\mathrm{s=2}}=\omega+ \dfrac{(2r-1)\dx^{2}}{24c_0^{2}}   \omega^{3}+\left[  \dfrac{20r^{2}-20r+1}{1920c_0^4}\dx^4- \left(\dfrac{1}{20}+\dfrac{\sqrt{3}}{36}\right)\dt^4    \right]\omega^5+\mathcal{O}(\omega^7).
\end{equation}
Obviously, when the third-order method is combined with either the standard consistent ($r=1$) or lumped ($r=0$) mass matrix, the numerical frequency error remains dominated by the spatial discretization error and is maintained at $\mathcal{O}(\omega^3)$. By employing the high-order mass matrix ($r=1/2$), the $\omega^3$ spatial dispersion error can be eliminated, and the dispersion accuracy is consequently improved to $\mathcal{O}(\omega^5)$. However, the coefficient of the $\omega^5$ term cannot be completely eliminated by selecting a real CFL number. Therefore, no optimal CFL number exists in the real domain for further removing the leading dispersion error of the third-order method.

\begin{figure}[htbp]
	\centering 
        \includegraphics[scale=1.0]{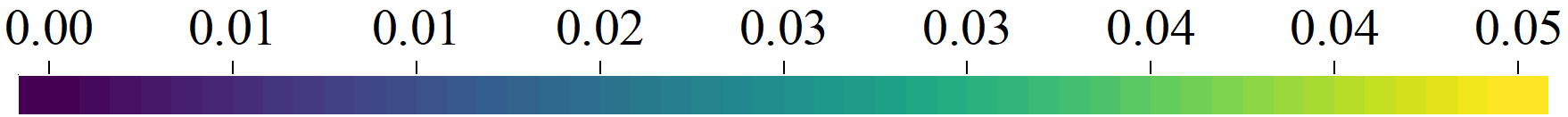}\\
    \subfigure[ ]{
		\includegraphics[scale=1.0]{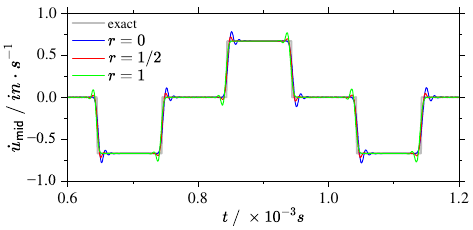}}
    \subfigure[ ]{
        \includegraphics[scale=1.0]{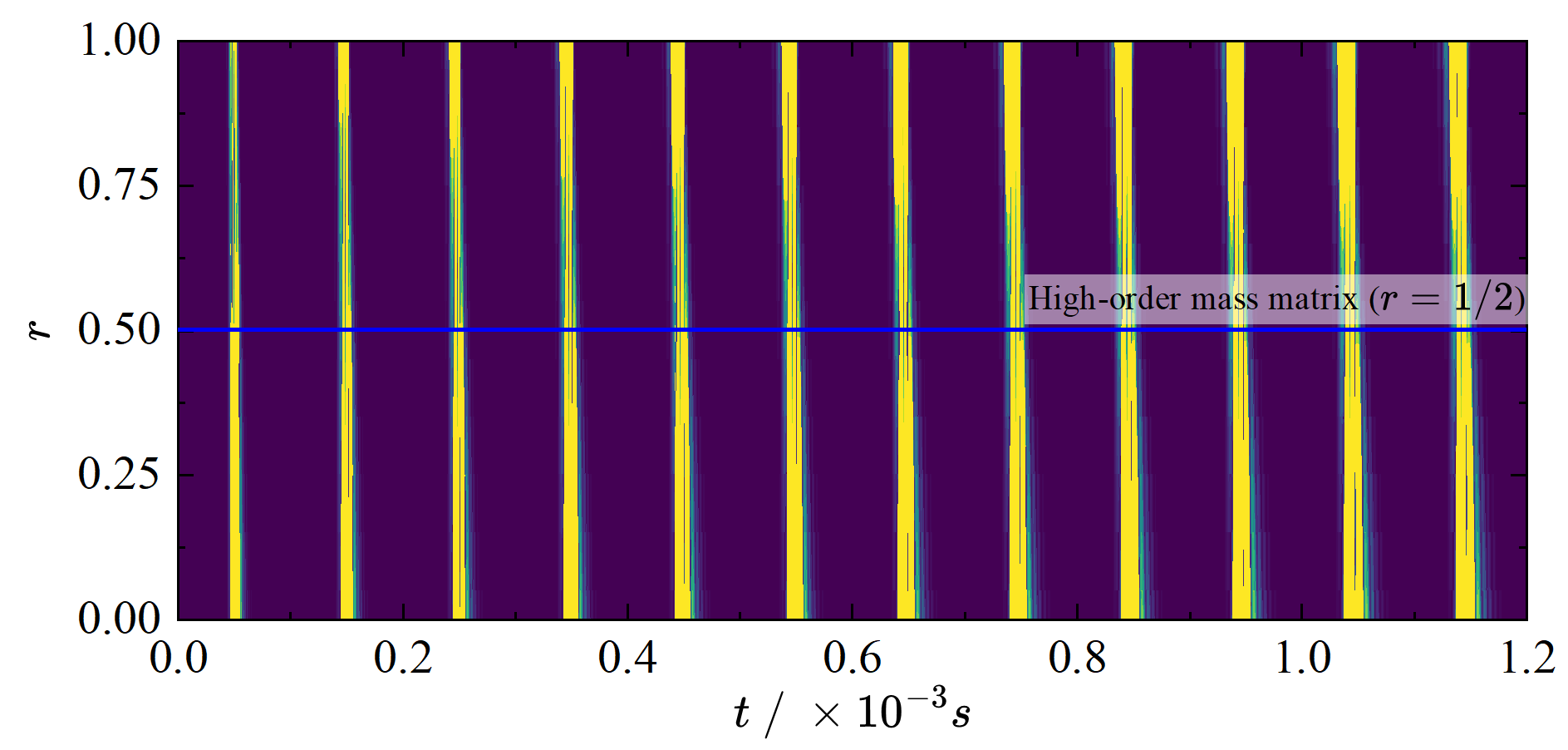}}
	\caption{The mid-point velocity of the bar predicted by the third-order two-sub-step implicit method ($s=2$) with $\mathsf{CFL}=1$.}
    \label{fig:bar_s2_1}
\end{figure}

The numerical results of the third-order two-sub-step implicit method are presented in \cref{fig:bar_s2_1} to further verify the theoretical analysis. As shown in \cref{fig:bar_s2_1}(a), the third-order method with $\mathsf{CFL}=1$ is evaluated using different mass weighting parameters. The results indicate that the numerical solutions obtained with the standard lumped mass matrix ($r=0$) and consistent mass matrix ($r=1$) still exhibit noticeable oscillations near the wave fronts. This behavior is consistent with the theoretical prediction that the dispersion error remains at $\mathcal{O}(\omega^3)$ for these two cases. In contrast, when the high-order mass matrix ($r=1/2$) is employed, the numerical response almost coincides with the exact solution. This improvement results from the fact that the $\omega^3$ spatial dispersion term is completely removed for $r=1/2$, and the dominant frequency error is further reduced to $\mathcal{O}(\omega^5)$. Therefore, the combination of the third-order implicit method and the high-order mass matrix provides a significantly improved accuracy. The error distribution in \cref{fig:bar_s2_1}(b) further confirms this observation. The minimum error region is located around $r=1/2$, where the theoretical cancellation of the leading dispersion term occurs. In contrast, for $r=0$ and $r=1$, the error remains relatively large.

These results indicate that the second- and third-order two-sub-step methods exhibit different mechanisms for improving numerical accuracy. The second-order method benefits from the compensation between spatial and temporal dispersion errors through an optimal CFL number, whereas the third-order method relies primarily on improving the spatial discretization accuracy. Consequently, increasing the temporal order alone does not necessarily guarantee superior dispersion performance unless the spatial discretization error is simultaneously controlled.

\paragraph{The three-sub-step implicit method}
For the third-order three-sub-step implicit method, the dispersion relation is given by \cref{eq:omega_2node_3}. To further improve the dispersion accuracy, the high-order mass matrix with $r=1/2$ must be adopted, which eliminates the $\omega^3$ term in spatial errors. In this case, the leading frequency error becomes $\left[\phi_2\dt^4/30-\dx^4/(480c_0^4)\right]\omega^5$ and thus a real optimal CFL number exists only when $\phi_2>0$, which is given by
\begin{equation}\label{eq:s3_CFL_Opt}
\mathsf{CFL}_{\mathsf{Opt}}=\dfrac{1}{2\sqrt[4]{\phi_2}}= \dfrac{1}{2\sqrt[4]{90\gamma_1^4-150\gamma_1^{3}+75\gamma_1^{2}-15\gamma_1+1}}
\end{equation}
where $\gamma_1\in\left(0.9756745886944403,~1/2+\sqrt{3}\cos(\pi/18)/3\right]$. It indicates that, unlike the second-order two-sub-step method, the third-order three-sub-step method requires the high-order mass matrix ($r=1/2$) to achieve further dispersion improvement. Under this condition, the spatial and temporal dispersion errors of order $\mathcal{O}(\omega^5)$ can be balanced through an optimal CFL number given by \cref{eq:s3_CFL_Opt}, leading to a higher-order numerical frequency approximation $\mathcal{O}(\omega^7)$. 

Since \cref{eq:omega_2node_3} does not involve $\phi_1$, the optimal CFL number derived for the fourth-order ($\phi_1=0$) three-sub-step implicit method with $\gamma_1=1/2+\sqrt{3}\cos(\pi/18)/3$ is also applicable to \cref{eq:s3_CFL_Opt}. Therefore, for spatially discretized models using linear finite elements, increasing the temporal accuracy from third to fourth order through the three-sub-step method does not bring a significant improvement in the dispersion error. The reason is that both methods are constrained by the spatial discretization, and the attainable frequency accuracy remains at the same order. Consequently, further improvement of dispersion accuracy requires not only increasing the temporal order but also employing higher-order spatial discretization schemes or specially designed mass matrices.

\begin{figure}[htbp]
	\centering 
        \includegraphics[scale=1.0]{bar_s1_leg_down}\\
    \subfigure[ ]{
		\includegraphics[scale=1.0]{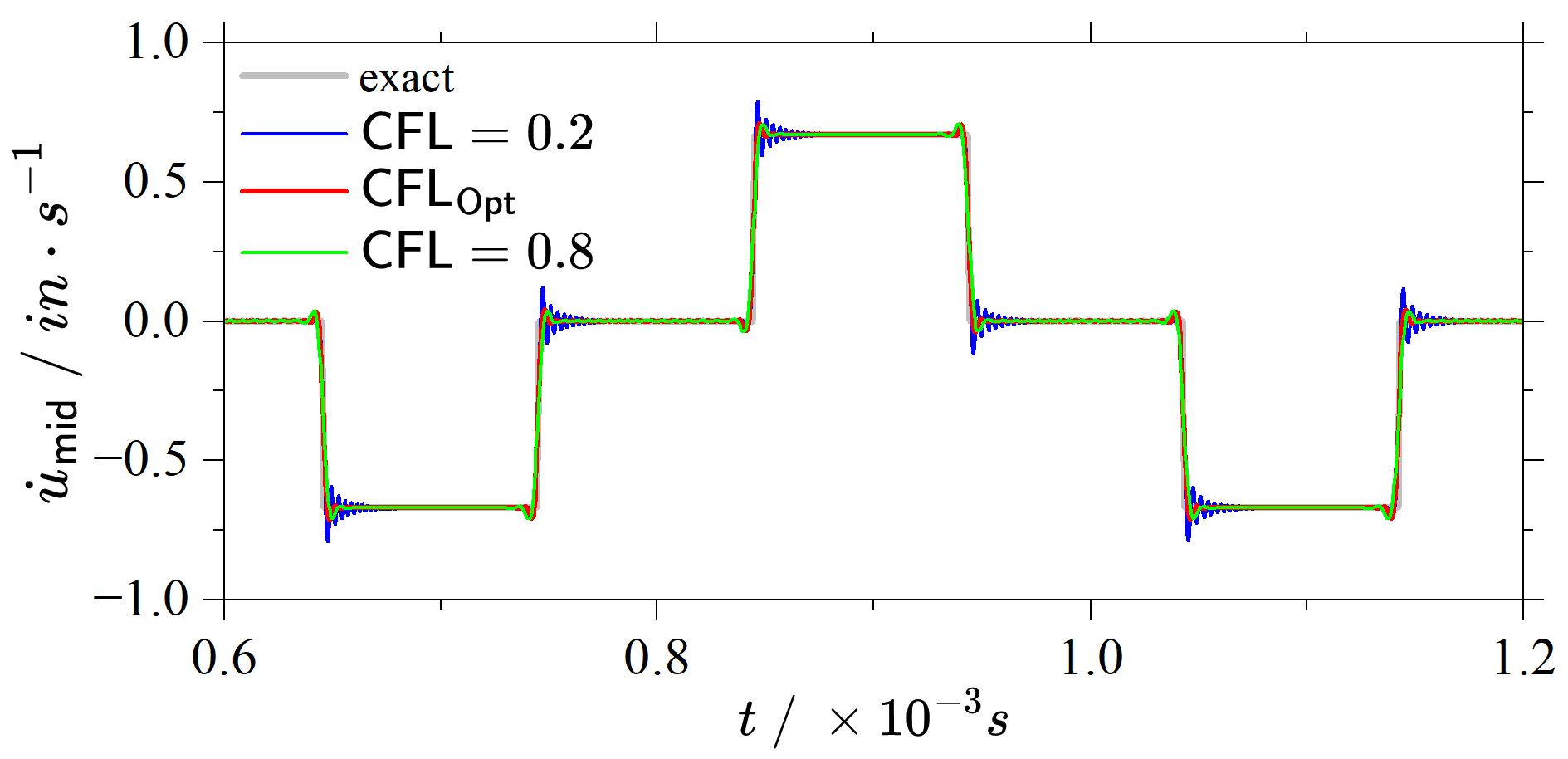}}
    \subfigure[ ]{
        \includegraphics[scale=1.0]{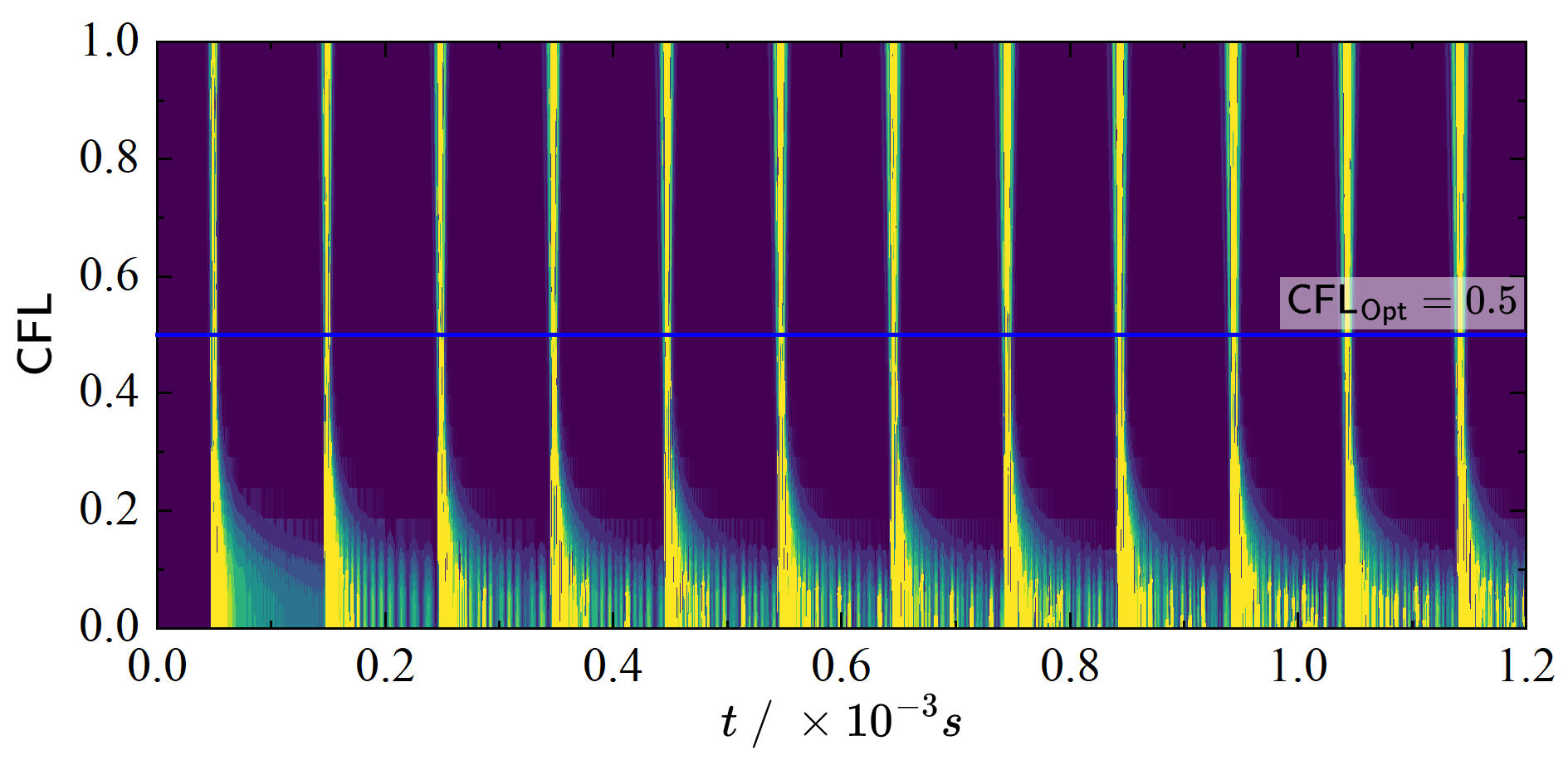}}
    \subfigure[ ]{
		\includegraphics[scale=1.0]{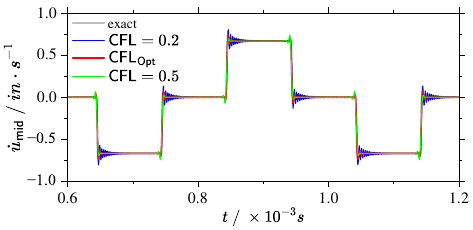}}
    \subfigure[ ]{
        \includegraphics[scale=1.0]{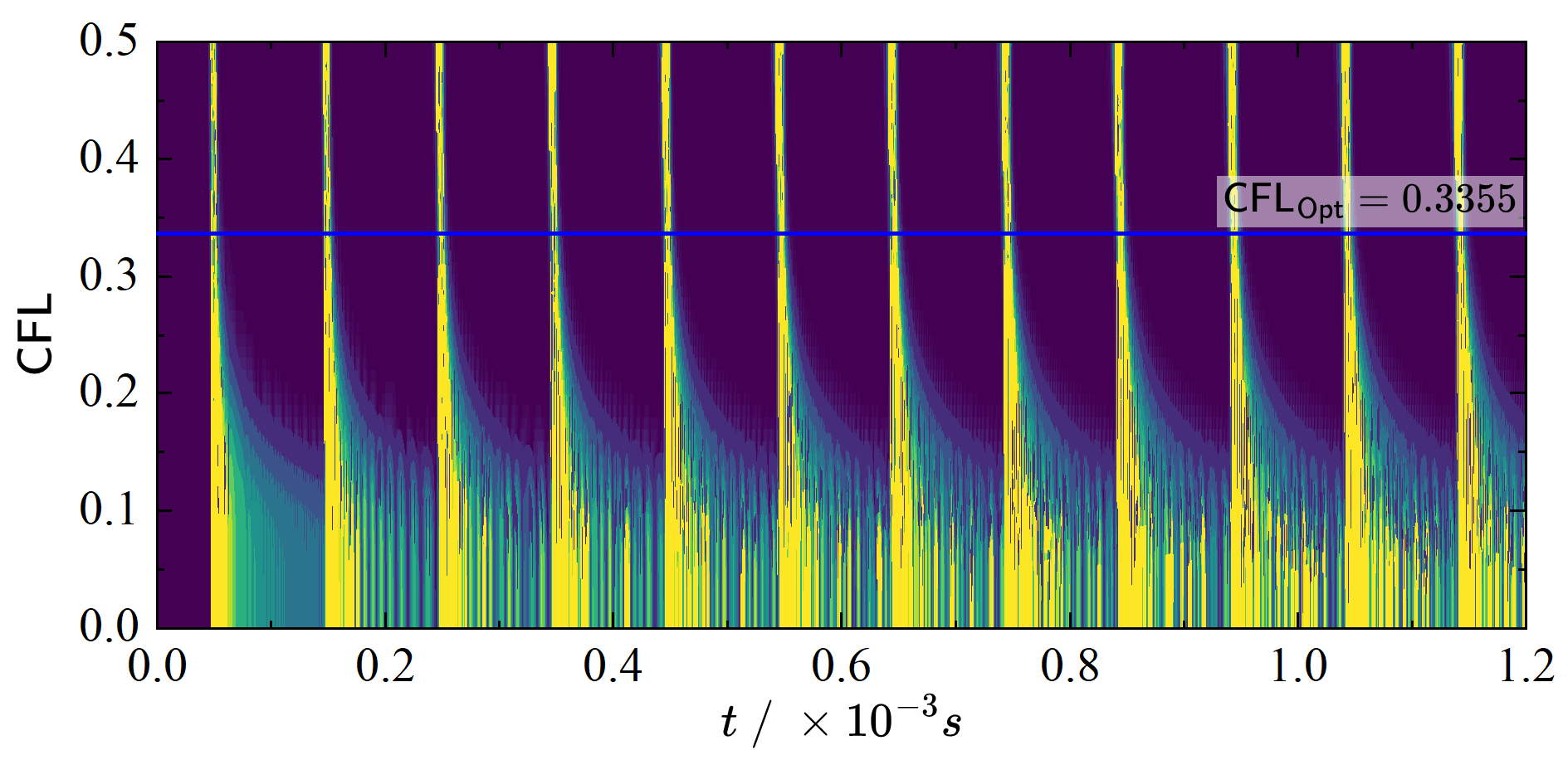}}
	\caption{The mid-point velocity of the bar predicted by the three-sub-step implicit method ($s=3$): (a-b) third-order accuracy with $\gamma_1=1.0$ and (c-d) fourth-order accuracy.}
    \label{fig:bar_s3_1}
\end{figure}

\cref{fig:bar_s3_1} further investigates the influence of the optimal CFL number predicted by the dispersion analysis for the three-sub-step implicit method. Two cases with different accuracy are considered. As shown in \cref{fig:bar_s3_1}(a-b), the third-order three-sub-step implicit method with $\gamma_1=1.0$ is employed. The numerical results indicate that the solution obtained using the optimal CFL number agrees better with the exact solution than those using other CFL values. In particular, compared with smaller or larger CFL numbers, the oscillations near the discontinuity regions are effectively reduced, demonstrating the validity of the optimal CFL selection. For the fourth-order three-sub-step implicit method shown in \cref{fig:bar_s3_1}(c-d), the parameter $\gamma_1=1/2+\sqrt{3}\cos(\pi/18)/3$ is adopted such that $\phi_1=0$. The optimal CFL number for the fourth-order method is smaller than that of the third-order method. The numerical results demonstrate that the solution obtained with the optimal CFL value provides the best agreement with the exact solution among the tested CFL numbers.

It is worth noting that, although the fourth-order method achieves a higher temporal accuracy than the third-order method, both methods are ultimately limited by the spatial discretization error when linear finite elements are employed. Therefore, the improvement in temporal order alone cannot further enhance the attainable dispersion accuracy unless it is combined with an appropriate spatial discretization scheme or a specially designed mass matrix. The results in \cref{fig:bar_s3_1} confirm that the theoretically derived optimal CFL numbers can effectively minimize the numerical dispersion error and provide guidance for selecting appropriate time step sizes in practical simulations.

\paragraph{The four-sub-step implicit method}
For the fourth-order four-sub-step method ($\varpi_1\neq0$), the $\omega^3$ term in \cref{eq:omega_2node_4} represents the dispersion error introduced by the spatial discretization. It is independent of time integration algorithms and remains at $\mathcal{O}(\omega^3)$ when the standard lumped mass matrix ($r=0$) or consistent mass matrix ($r=1$) is adopted. Therefore, increasing the temporal order alone cannot improve the frequency accuracy for these commonly used mass matrices. To eliminate this dominant spatial dispersion error, the high-order mass matrix with $r=1/2$ must be employed. With $r=1/2$, \cref{eq:omega_2node_4} is simplified as
\begin{equation}
\overline{\omega}_\mathrm{s=4}=\omega- \left[ \dfrac{1}{480c_0^4}\dx^4+\dfrac{\varpi_1}{120}\dt^4  \right]\omega^5+ \bigg[ \dfrac{\varpi_3}{336}\dt^6-\dfrac{1}{12096c_0^6}\dx^6 \bigg]\omega^7+\mathcal{O}(\omega^9)
\end{equation}
At this stage, the spatial and temporal dispersion errors are of the same order and can be balanced through an optimal CFL number, which is computed as
\begin{equation}
\mathsf{CFL}_{\mathsf{Opt}}= \dfrac{\sqrt{2}}{2\sqrt[4]{-\varpi_1}}=\dfrac{\sqrt{2}}{2\sqrt[4]{-(120\gamma_1^4 - 240\gamma_1^3 + 120\gamma_1^2 - 20\gamma_1 + 1)}}
\end{equation}
where $\gamma_1$ should satisfy \cref{eq:s4_uc}. Therefore, the fourth-order four-sub-step method combined with the high-order mass matrix ($r=1/2$) can achieve a further improvement of dispersion accuracy from $\mathcal{O}(\omega^5)$ to $\mathcal{O}(\omega^7)$. 

\begin{figure}[htbp]
	\centering 
        \includegraphics[scale=1.0]{bar_s1_leg_down}\\
    \subfigure[ ]{
		\includegraphics[scale=1.0]{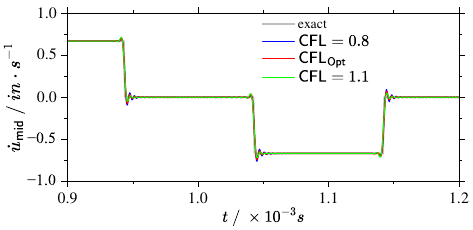}}
    \subfigure[ ]{
        \includegraphics[scale=1.0]{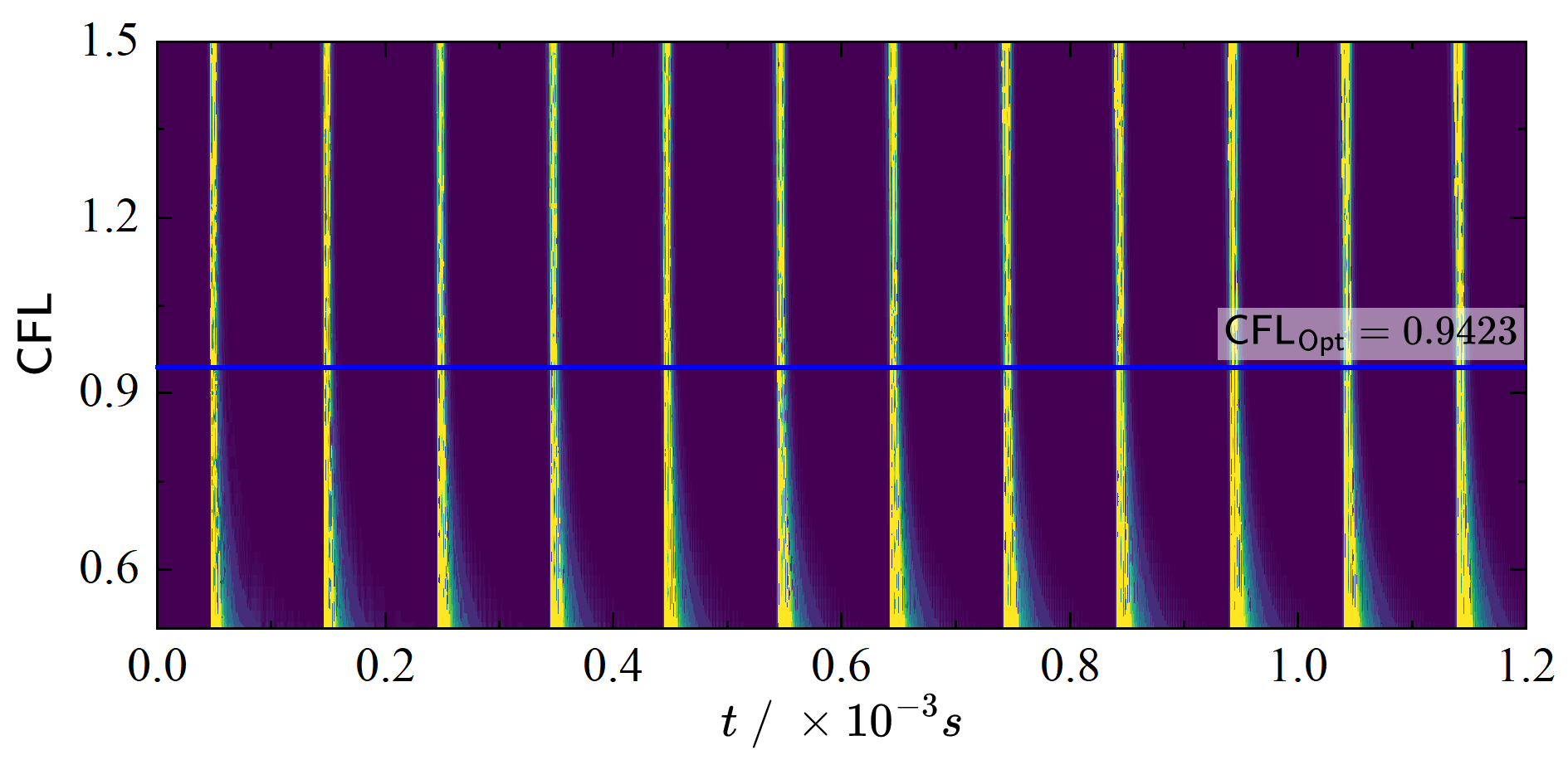}}
    \subfigure[ ]{
		\includegraphics[scale=1.0]{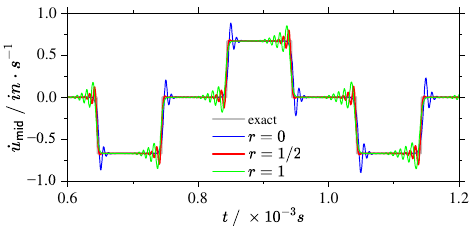}}
    \subfigure[ ]{
        \includegraphics[scale=1.0]{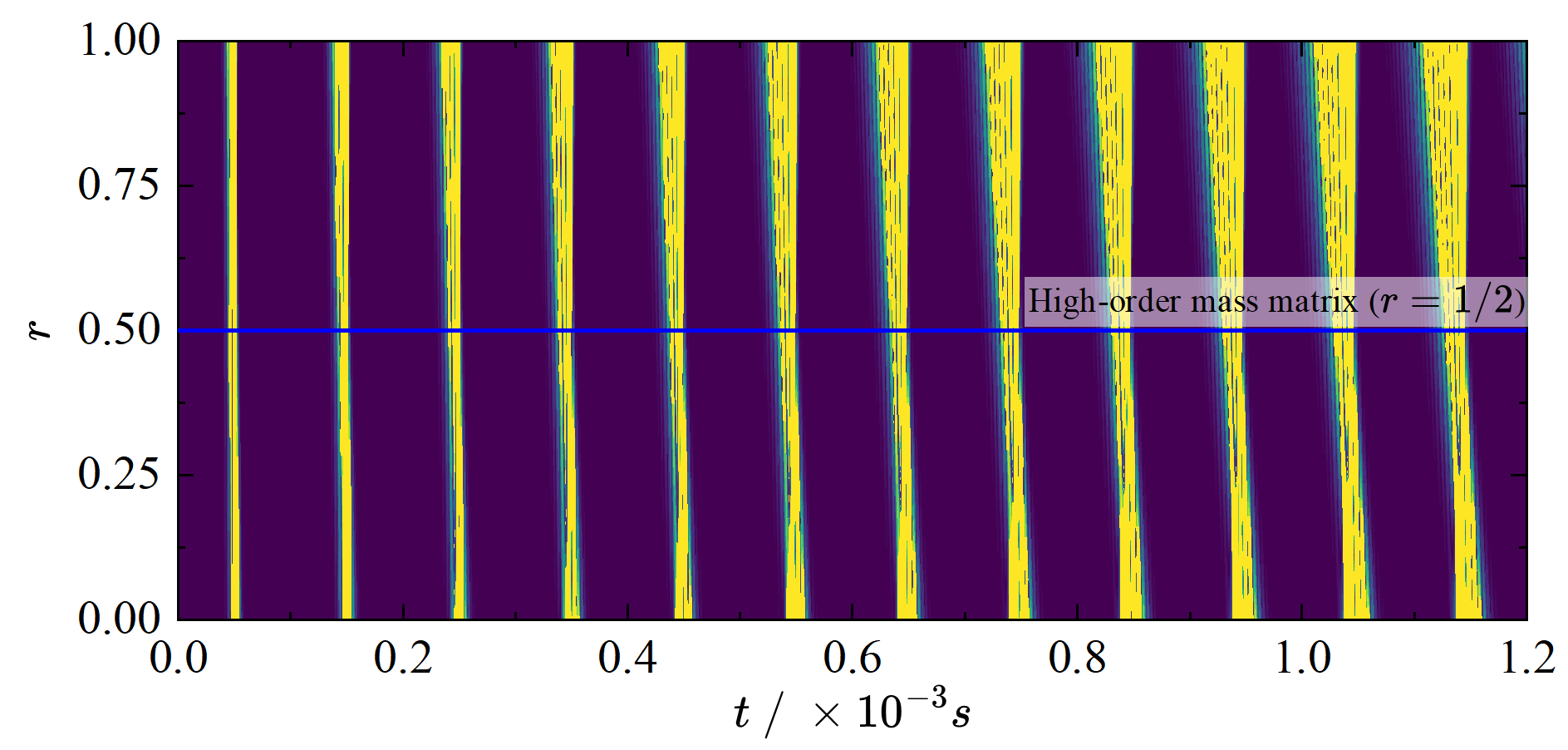}}
	\caption{The mid-point velocity of the bar predicted by the four-sub-step implicit method ($s=4$): (a-b) fourth-order accuracy with $\infrho=0.8$ and (c-d) fifth-order accuracy.}
    \label{fig:bar_s4_1}
\end{figure}

\cref{fig:bar_s4_1}(a-b) first investigates the fourth-order four-sub-step implicit method with $\infrho=0.8$. Three different CFL numbers, including the optimal value ($\approx0.9423$), are considered. It can be observed that all tested CFL values provide highly accurate predictions of the midpoint velocity, and the numerical solutions almost overlap with the reference solution. Therefore, unlike the lower-order cases, the advantage of the optimal CFL number is not clearly distinguishable. This is mainly because the fourth-order method already provides sufficiently small dispersion errors, such that the influence of CFLs on the final numerical response becomes relatively weak. Nevertheless, the error contour in \cref{fig:bar_s4_1}(b) indicates that the optimal CFL still corresponds to the region with improved dispersion behavior, which is consistent with the theoretical analysis.

For the fifth-order four-sub-step implicit method ($\varpi_1=0$), the numerical frequency relation  given by \cref{eq:omega_2node_4} reduces to
\begin{equation}\label{eq:omega_s4_1}
\overline{\omega}_\mathrm{s=4}=\omega+ \dfrac{(2r-1)\dx^{2}}{24c_0^{2}}\omega^3+\dfrac{20r^{2}-20r+1}{1920c_0^4}\dx^4\omega^5+ \bigg[ \dfrac{1400r^{3}-2100r^{2}+546r-3}{967680c_0^6}\dx^6+\dfrac{\varpi_3}{336}\dt^6 \bigg]\omega^7+\mathcal{O}(\omega^9). 
\end{equation}
Similarly, the elimination of the $\omega^3$ term still requires the high-order mass matrix, i.e., $r=1/2$. Otherwise, the frequency error remains dominated by the spatial discretization and stays at $\mathcal{O}(\omega^3)$ regardless of the increased temporal order. When $r=1/2$ is adopted, \cref{eq:omega_s4_1} becomes
\begin{equation}
\overline{\omega}_\mathrm{s=4}=\omega- \dfrac{1}{480c_0^4}\dx^4\omega^5+\mathcal{O}(\omega^7).
\end{equation}
It can be observed that the leading dispersion error is purely induced by the spatial discretization and no longer contains any adjustable temporal contribution. Therefore, the spatial and temporal dispersion errors cannot be compensated by modifying the CFL number, and no real optimal CFL number exists for the fifth-order four-sub-step method.

\cref{fig:bar_s4_1}(c-d) further examines the fifth-order four-sub-step implicit method by varying the mass weighting parameter $r$. According to the analysis, the high-order mass matrix with $r=1/2$ is required to eliminate the term $\mathcal{O}(\omega^3)$ and achieve higher-order accuracy. However, the numerical results show that the minimum error does not exactly occur at $r=1/2$. A possible explanation is related to the stability of the fifth-order method. Since the fifth-order four-sub-step implicit method is $A(\alpha)$-stable rather than unconditionally stable, a narrow frequency range may exist where the spectral radius slightly exceeds unity, resulting in weak numerical amplification of specific frequency components. Although this amplification is very small, it modifies the overall error distribution and may offset the dispersion improvement obtained by adopting the high-order mass matrix. Consequently, the theoretically optimal choice $r=1/2$ based on dispersion analysis does not necessarily produce the minimum response error in the complete numerical simulation.

These results demonstrate that, for wave propagation problems involving spatial discretization, the final numerical accuracy is determined by the combined effects of dispersion and stability. The dispersion-based optimal parameters provide valuable guidance for improving frequency accuracy; however, the spectral characteristics of the algorithm should also be considered when evaluating the practical performance of high-order methods.

\paragraph{The five- and six-sub-step implicit methods}
For the five- and six-sub-step implicit methods, the dispersion error is given by \cref{eq:omega_2node_56}. It is observed that the leading errors are independent of the temporal parameters, and are only determined by the spatial discretizations. This indicates that, after increasing the number of sub-steps to five or six, the temporal dispersion errors are sufficiently reduced, and the numerical accuracy is dominated by the spatial discretization.

Therefore, further increasing the temporal accuracy or optimizing the sub-step parameters cannot effectively improve the dispersion accuracy unless a higher-order spatial discretization scheme or an improved mass matrix is adopted. Since the attainable frequency accuracy of the five- and six-sub-step methods is mainly restricted by the spatial discretization, additional numerical validations for \cref{eq:omega_2node_56} are not performed herein. These results further demonstrate that, for spatially discretized problems, improving the temporal accuracy alone does not necessarily lead to higher overall accuracy when the spatial discretization becomes the dominant error source.


\subsubsection{Spatial discretization with three-node quadratic elements} 
The spatial discretization considered in the previous subsection is based on low-order linear finite elements, which limits the overall accuracy and prevents the advantages of high-order time integration algorithms from being fully demonstrated. Therefore, this subsection further improves the spatial discretization accuracy by employing 100 three-node quadratic finite elements. When three-node quadratic finite elements are adopted, the corresponding elemental mass and stiffness matrices are given, respectively, as
\begin{subequations}
	\begin{equation}\label{eq:2me}
		\mbf{m}^e=(1-r)\cdot\frac{\rho A\dx}{6}\begin{bmatrix}
			1 & 0 & 0 \\ 0 & 4 & 0 \\ 0 & 0 & 1
		\end{bmatrix}+r\cdot\dfrac{\rho A\dx}{30}\begin{bmatrix}
		4 & 2 & -1 \\ 2 & 16 & 2\\ -1 & 2 & 4
		\end{bmatrix}
	\end{equation}
	and 
	\begin{equation}
		\mbf{k}^e=\frac{EA}{3\dx}\begin{bmatrix}
			7 & -8 & 1 \\ -8 & 16 & -8\\ 1 & -8 & 7
		\end{bmatrix}
	\end{equation}
\end{subequations} 
where $\dx$ denotes the length of three-node finite elements. The parameter $r$ in \cref{eq:2me} controls the distribution of the elemental mass matrix between lumped and consistent forms. In this case, the numerical frequencies $\overline{\omega}$ for various sub-step implicit methods are directly given as 
\begin{subequations}
\begin{align}
\overline{\omega}_\mathrm{s=1}=&\omega-\left( \gamma_1^{2}-\gamma_1+ \dfrac{1}{3}  \right)\dt^{2}\omega^{3}+ \dfrac{(3r-1)\dx^4}{2880c_0^4}\omega^5+\mathcal{O}(\omega^7)\label{eq:omega_3node_1}\\
\overline{\omega}_\mathrm{s=2}=&\omega+ \dfrac{\vartheta_1}{6}  \dt^{2}\omega^{3}+ \left[ \dfrac{3r-1}{2880c_0^4}\dx^4- \dfrac{\vartheta_3}{20}\dt^4 \right]\omega^5+\mathcal{O}(\omega^7)\label{eq:omega_3node_2}\\
\overline{\omega}_\mathrm{s=3}=&\omega+ \left[ \dfrac{3r-1}{2880c_0^4}\dx^4+ \dfrac{\phi_2}{30}\dt^4 \right]\omega^5- \left[ \dfrac{63r^2+25}{2419200c_0^6}\dx^6+ \dfrac{\phi_4}{252}\dt^6 \right]\omega^7 +\mathcal{O}(\omega^9)\label{eq:omega_3node_3}\\
\overline{\omega}_\mathrm{s=4}=&\omega+ \left[ \dfrac{3r-1}{2880c_0^4}\dx^4- \dfrac{\varpi_1}{120}\dt^4 \right]\omega^5+ \left[\dfrac{\varpi_3}{336}\dt^6-\dfrac{63r^2+25}{2419200c_0^6}\dx^6\right]\omega^7+\mathcal{O}(\omega^9)\label{eq:omega_3node_4}\\
\overline{\omega}_\mathrm{s=5}=&\omega+\dfrac{3r-1}{2880c_0^4}\dx^4\omega^5-\left[\dfrac{63r^2+25}{2419200c_0^6}\dx^6+ \dfrac{\psi_2}{840}\dt^6\right]\omega^7+\mathcal{O}(\omega^9)\label{eq:omega_3node_5}\\
\overline{\omega}_\mathrm{s=6}=&\omega+\dfrac{3r-1}{2880c_0^4}\dx^4\omega^5+\left[\dfrac{\varphi_1}{5040}\dt^6-\dfrac{63r^2+25}{2419200c_0^6}\dx^6\right]\omega^7+\mathcal{O}(\omega^9)\label{eq:omega_3node_6}
\end{align}
\end{subequations}
where $\vartheta_{1,3}$, $\phi_{2,4}$, $\varpi_{1,3}$, $\psi_2$ and $\varphi_1$ are given by \cref{eq:s2_amp_phas,eq:s3_amp_phas,eq:s4_amp_phase,eq:s5_amp_phase,eq:s6_amp_phase}, respectively. Obviously, the quadratic finite element eliminates the third-order spatial dispersion error $\dx^3$ that exists in the linear finite elements. 

\paragraph{The single- and two-sub-step implicit methods}

For the single-sub-step implicit method, the leading error is still governed by the temporal contribution associated with the coefficient $\gamma_1^2-\gamma_1+1/3$. Although the spatial discretization error is reduced to $\mathcal{O}(\omega^5)$ by employing quadratic finite elements, the third-order temporal dispersion error cannot be eliminated because $\gamma_1^2-\gamma_1+1/3=0$ has no real solution. Therefore, the single-sub-step method cannot achieve further improvement in dispersion accuracy through the real parameter adjustment.

\begin{figure}[htbp]
	\centering 
        \includegraphics[scale=1.0]{bar_s1_leg}\\
    \subfigure[$\infrho=1-\sqrt{3}~\&~r=1$]{
		\includegraphics[scale=1.0]{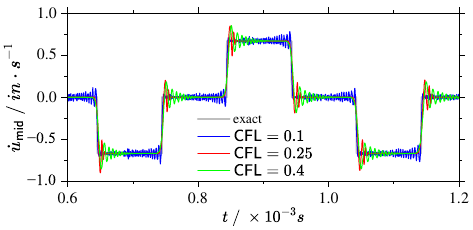}}
    \subfigure[$\infrho=1-\sqrt{3}~\&~r=1$]{
        \includegraphics[scale=1.0]{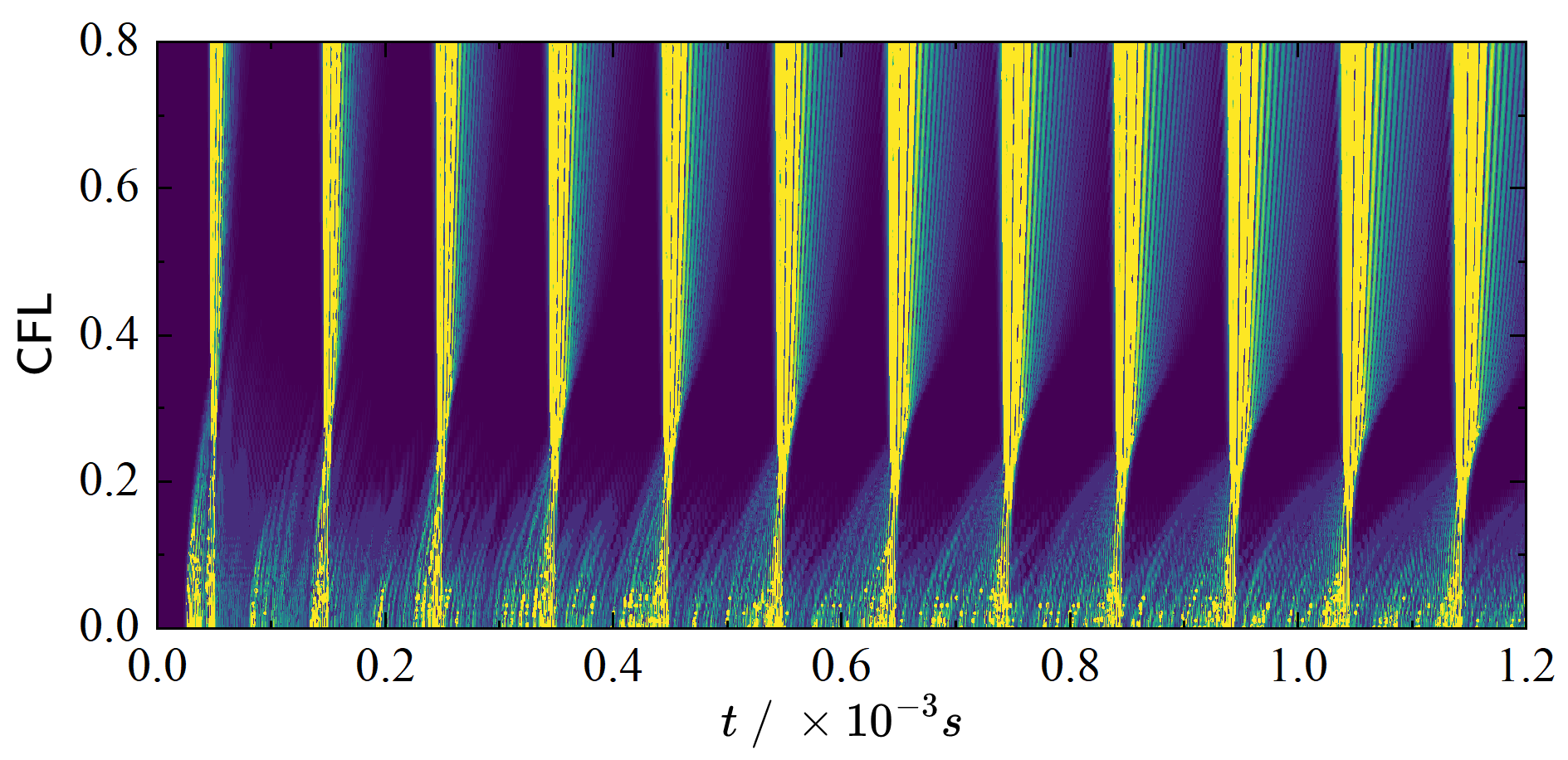}}
    \subfigure[$r=1$]{
		\includegraphics[scale=1.0]{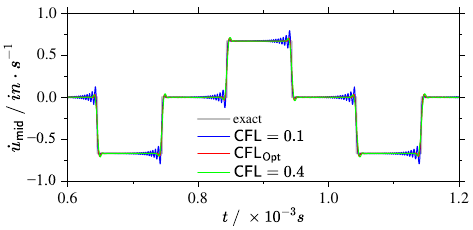}}
    \subfigure[$r=1$]{
        \includegraphics[scale=1.0]{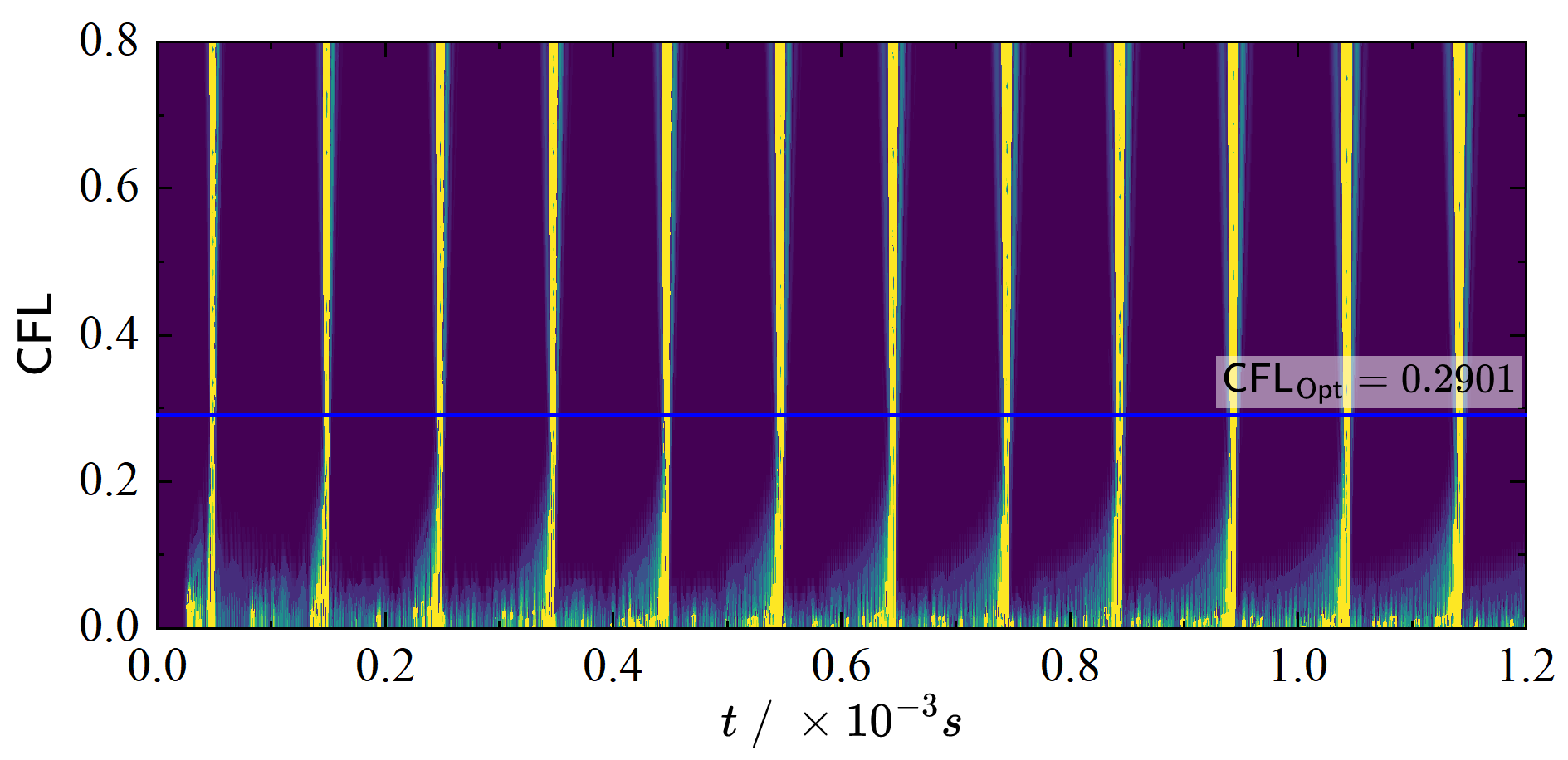}}
	\caption{The mid-point velocity of the bar predicted by the two-sub-step implicit method ($s=2$): (a-b) second-order accuracy and (c-d) third-order accuracy.}
    \label{fig:bar_s2_2}
\end{figure}

For the second-order member, the leading dispersion error is dominated by the $\omega^3$ term due to $\vartheta_1\neq0$. Therefore, no theoretical optimal CFL number exists for eliminating the leading frequency error. As shown in \cref{fig:bar_s2_2}(a-b), numerical responses obtained with different CFL numbers exhibit significant spurious oscillations, and the distribution of spurious oscillations around the wave front varies significantly with the CFL number. Interestingly, an empirical optimal region appears near $\mathsf{CFL}\approx0.25$, where the oscillations are minimized. This phenomenon does not originate from the cancellation of the leading dispersion error predicted by the frequency analysis, but rather from the modification of the propagation characteristics of high-frequency components. As the CFL number varies, the spurious oscillations gradually shift from the front of the wave to the rear of the wave, resulting in different error distributions. Therefore, the observed optimal CFL value should be regarded as an empirical choice for minimizing numerical oscillations rather than a theoretically derived dispersion-optimal parameter.

For the third-order member, the condition $\vartheta_1=0$ eliminates the third-order temporal dispersion term, and the leading frequency error is improved to $\mathcal{O}(\omega^5)$. Consequently, the spatial and temporal dispersion errors introduced by the quadratic finite element and the time integration scheme can be balanced through the optimal CFL number: 
\begin{equation}
\mathsf{CFL}_{\mathsf{Opt}}= \dfrac{1}{2}\sqrt[4]{ \dfrac{3r-1}{9+5\sqrt{3}} }. 
\end{equation}
For the conventional consistent mass matrix ($r=1$), the optimal CFL number is $\mathsf{CFL}_{\mathrm{Opt}}\approx0.2901$. As demonstrated in \cref{fig:bar_s2_2}(c-d), the numerical response corresponding to the optimal CFL number agrees well with the exact solution, and the wave-front oscillations are significantly reduced compared with non-optimal CFL selections. The error distribution in \cref{fig:bar_s2_2}(d) further confirms that the optimal CFL value effectively balances the spatial and temporal dispersion errors.

These results demonstrate that increasing the temporal accuracy alone does not necessarily guarantee improved dispersion performance unless the spatial discretization error is sufficiently reduced. With quadratic finite elements, the third-order two-sub-step method can fully exploit its higher temporal accuracy and achieve further dispersion improvement through CFL optimization. In contrast, the second-order member lacks such a theoretical optimization mechanism, and the observed optimal CFL behavior mainly reflects the redistribution of numerical high-frequency errors rather than the cancellation of the leading dispersion term.

\begin{figure}[htbp]
	\centering 
        \includegraphics[scale=1.0]{bar_s1_leg_down}\\
    \subfigure[$\infrho=0.8~\&~r=1$]{
		\includegraphics[scale=1.0]{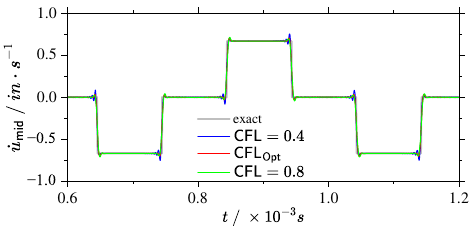}}
    \subfigure[$\infrho=0.8~\&~r=1$]{
        \includegraphics[scale=1.0]{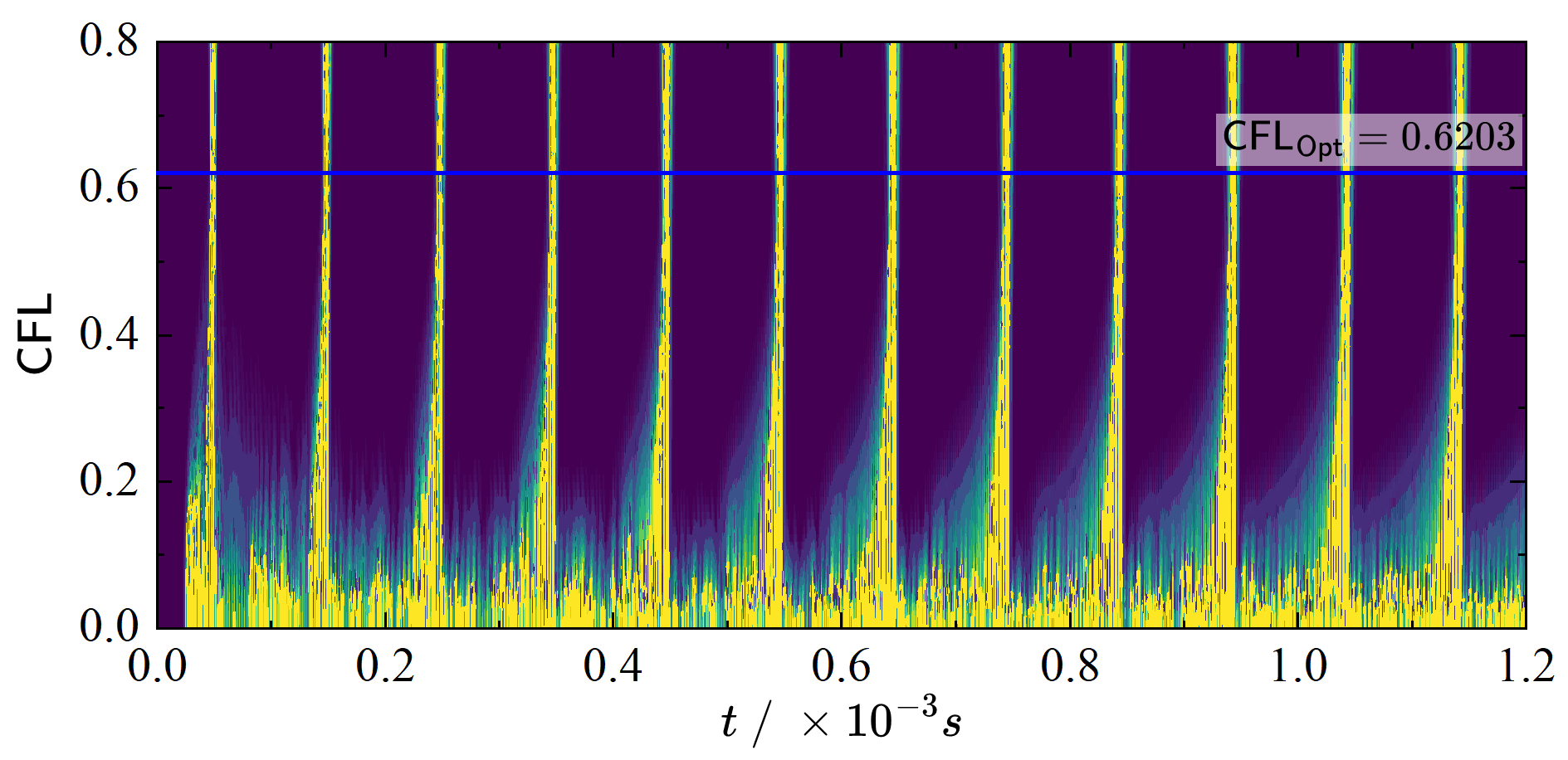}}
    \subfigure[$\gamma_1=1/2+\sqrt{3}\cos(\pi/18)/3~\&~r=0$]{
		\includegraphics[scale=1.0]{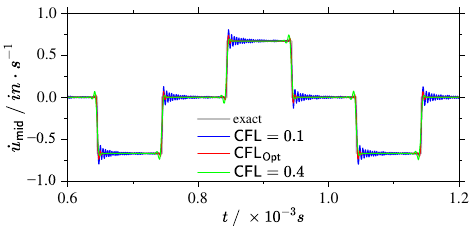}}
    \subfigure[$\gamma_1=1/2+\sqrt{3}\cos(\pi/18)/3~\&~r=0$]{
        \includegraphics[scale=1.0]{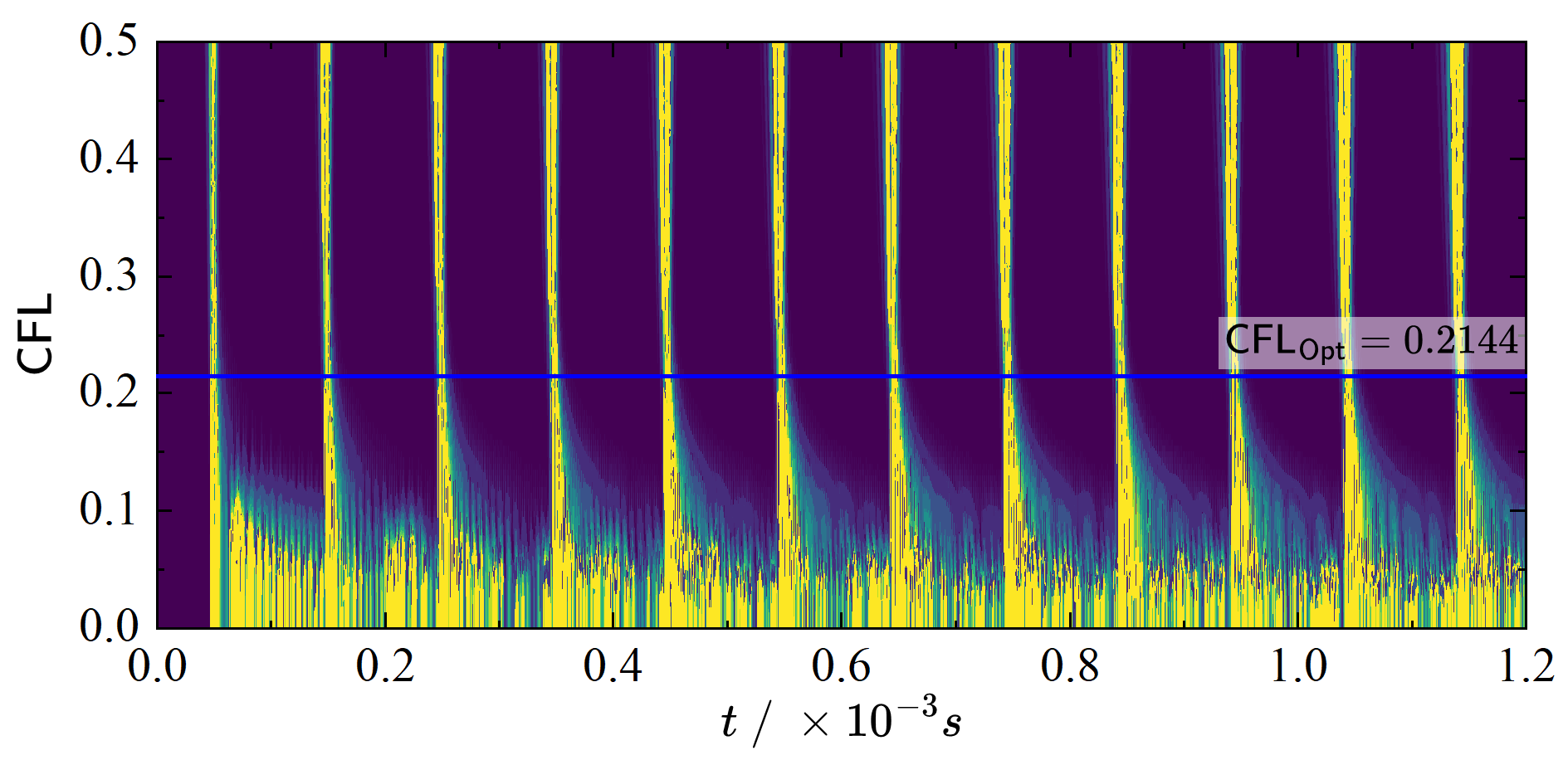}}
    \subfigure[$\phi_2=0~\&~\mathsf{CFL}=0.2$]{
		\includegraphics[scale=1.0]{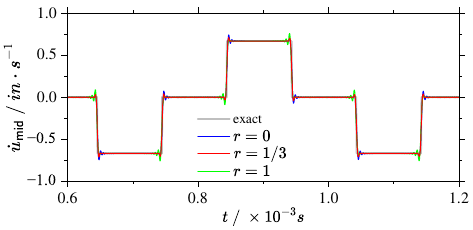}}
    \subfigure[$\phi_2=0~\&~\mathsf{CFL}=0.2$]{
        \includegraphics[scale=1.0]{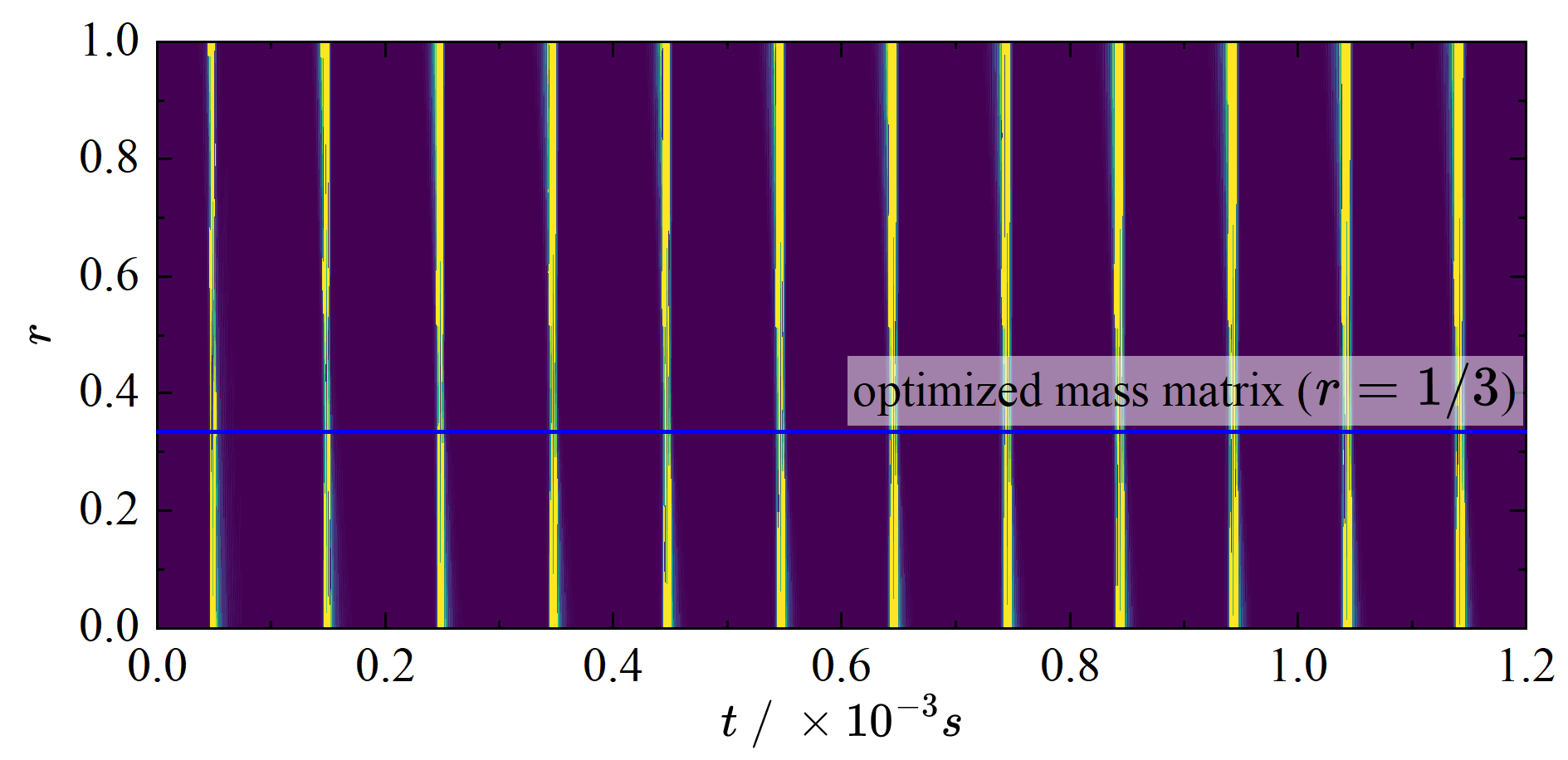}}
	\caption{The mid-point velocity of the bar predicted by the three-sub-step implicit method ($s=3$): (a-b) third-order accuracy, (c-d) fourth-order accuracy and (e-f) third-order accuracy with $\phi_2=0$.}
    \label{fig:bar_s3_2}
\end{figure}

\paragraph{The three-sub-step implicit method}
For the three-sub-step implicit method, the dispersion relation obtained with quadratic finite elements is given by \cref{eq:omega_3node_3}. Since the quadratic element eliminates the third-order spatial dispersion error, the leading frequency error is improved to $\mathcal{O}(\omega^5)$. 
\begin{itemize}
    \item When $\phi_2\neq0$, the spatial and temporal dispersion errors can be balanced by appropriately selecting the time step, leading to the following optimal CFL number:
\begin{equation}\label{eq:s3_cfl_opt}
\mathsf{CFL}_\mathsf{Opt}=\dfrac{1}{2} \sqrt[4]{\dfrac{1-3r}{6\phi_2}}=\dfrac{1}{2} \sqrt[4]{\dfrac{1-3r}{6(90\gamma_1^4-150\gamma_1^{3}+75\gamma_1^{2}-15\gamma_1+1)}}
\end{equation}
where the values of $\gamma_1$ and $r$ should ensure that $\mathsf{CFL}_\mathsf{Opt}$ is real and positive. As shown in \cref{fig:bar_s3_2}(a-b), where $r=1$ and $\infrho=0.8~(\gamma_1=0.34887)$ are adopted, the numerical response obtained using the optimal CFL number agrees well with the exact solution. Compared with other CFL selections, the optimized CFL value effectively reduces the dispersive oscillations around the wave front. The error distribution shown in \cref{fig:bar_s3_2}(b) further confirms that the dispersion error is minimized near the optimal CFL value, demonstrating the validity of the optimal CFL analysis. 
\item Furthermore, since \cref{eq:omega_3node_3} does not involve $\phi_1$, the fourth-order three-sub-step implicit method, which satisfies $\phi_1=0$, shares the same optimal CFL number given by \cref{eq:s3_cfl_opt} as the third-order member. As demonstrated in \cref{fig:bar_s3_2}(c-d) where $r=0$ and $\gamma_1=1/2+\sqrt{3}\cos(\pi/18)/3$ are adopted, the numerical solution obtained with the optimal CFL number exhibits improved agreement with the exact solution compared with non-optimal CFL selections. This confirms that increasing the temporal order from third to fourth order does not alter the optimal CFL condition when the leading dispersion term remains unchanged.
\item For the special case with $\phi_2=0$, corresponding to $\gamma_1=0.9756745886944403$, the fifth-order temporal dispersion contribution vanishes, and the dispersion error is dominated by the spatial discretization. According to the theoretical analysis, the optimized mass matrix with $r=1/3$ can eliminate the fifth-order spatial error and improve the leading frequency accuracy to $\mathcal{O}(\omega^7)$. The numerical comparisons in \cref{fig:bar_s3_2}(e-f) verify this conclusion. When different mass weighting parameters are considered, the optimized mass matrix ($r=1/3$) provides the most accurate prediction, while the lumped ($r=0$) and consistent ($r=1$) mass matrices exhibit relatively larger dispersion errors. The error distribution also shows that the use of the optimized mass matrix significantly suppresses the remaining numerical oscillations.
\end{itemize}

\paragraph{The four-sub-step implicit method}
For the fourth-order four-sub-step implicit method, the leading dispersion error in \cref{eq:omega_3node_4} is governed by the $\omega^5$ term containing both spatial and temporal contributions. Since $\varpi_1<0$ for the range of $\gamma_1$ given by \cref{eq:s4_uc}, an optimal CFL number exists when an appropriate mass weighting parameter satisfying $r<1/3$ is employed. By balancing the spatial and temporal dispersion errors, the $\omega^5$ term can be eliminated to derive
\begin{equation}
\mathsf{CFL}_\mathsf{Opt}=\sqrt[4]{\dfrac{3r-1}{24\varpi_1}}=\sqrt[4]{\dfrac{3r-1}{24(120\gamma_1^4 - 240\gamma_1^3 + 120\gamma_1^2 - 20\gamma_1 + 1)}}
\end{equation}
and thus the numerical frequency can be improved to $\mathcal{O}(\omega^7)$.

\begin{figure}[htbp]
	\centering 
        \includegraphics[scale=1.0]{bar_s1_leg_down}\\
    \subfigure[$\infrho=0.8~\&~r=0$]{
        \includegraphics[scale=1.0]{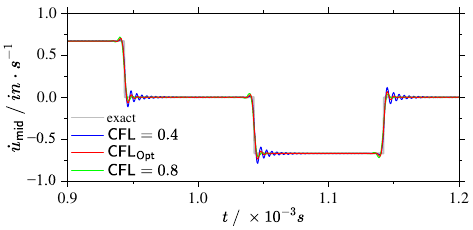}}
    \subfigure[$\infrho=0.8~\&~r=0$]{
        \includegraphics[scale=1.0]{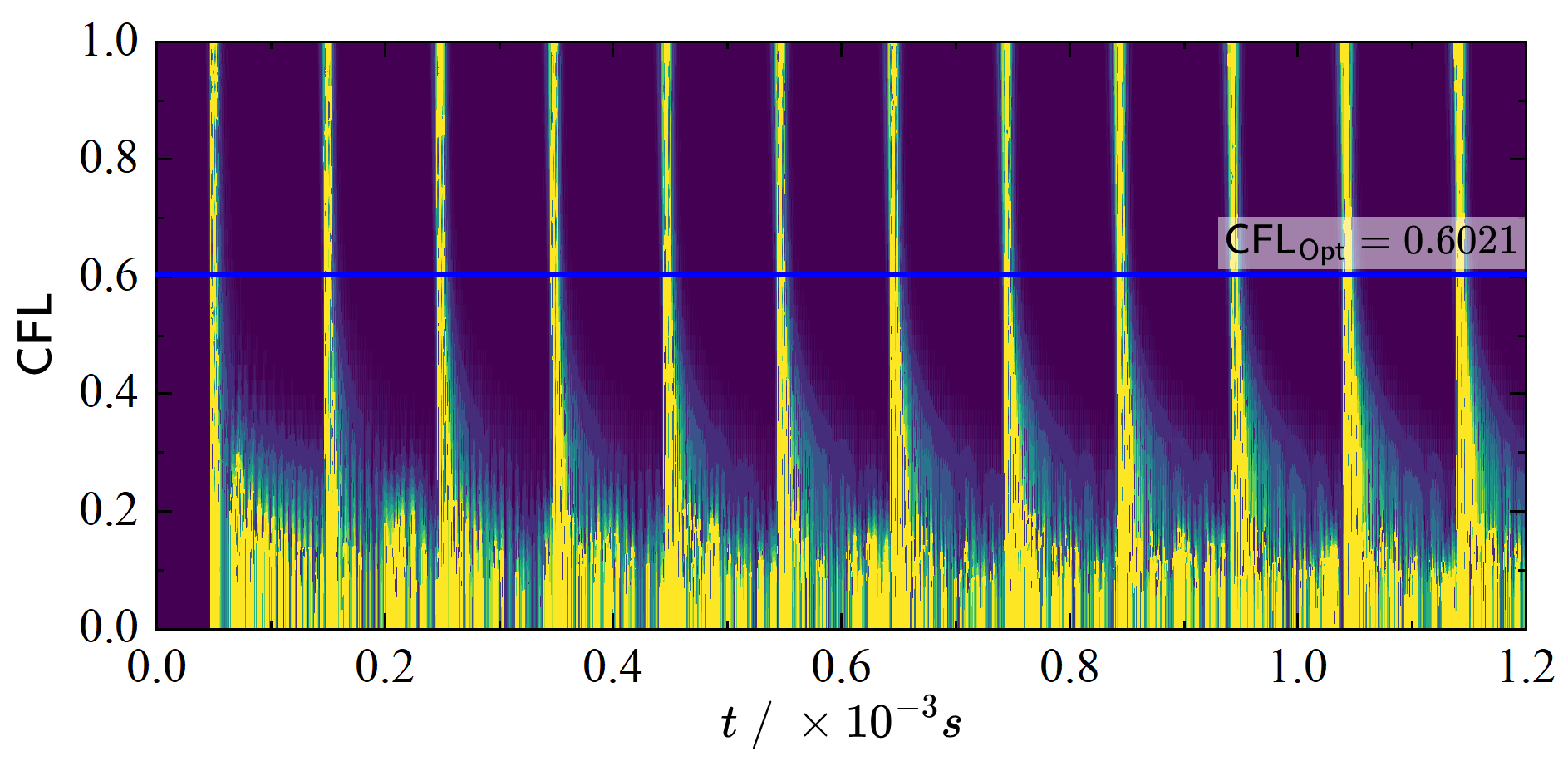}}
    \subfigure[$\gamma_1=1.3453664197803336~\&~r=1/3$]{
        \includegraphics[scale=1.0]{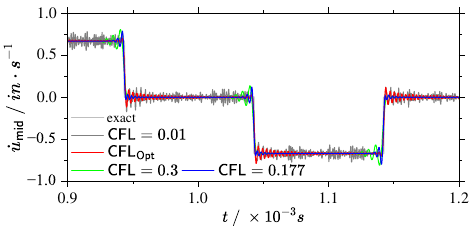}}
    \subfigure[$\gamma_1=1.3453664197803336~\&~r=1/3$]{
        \includegraphics[scale=1.0]{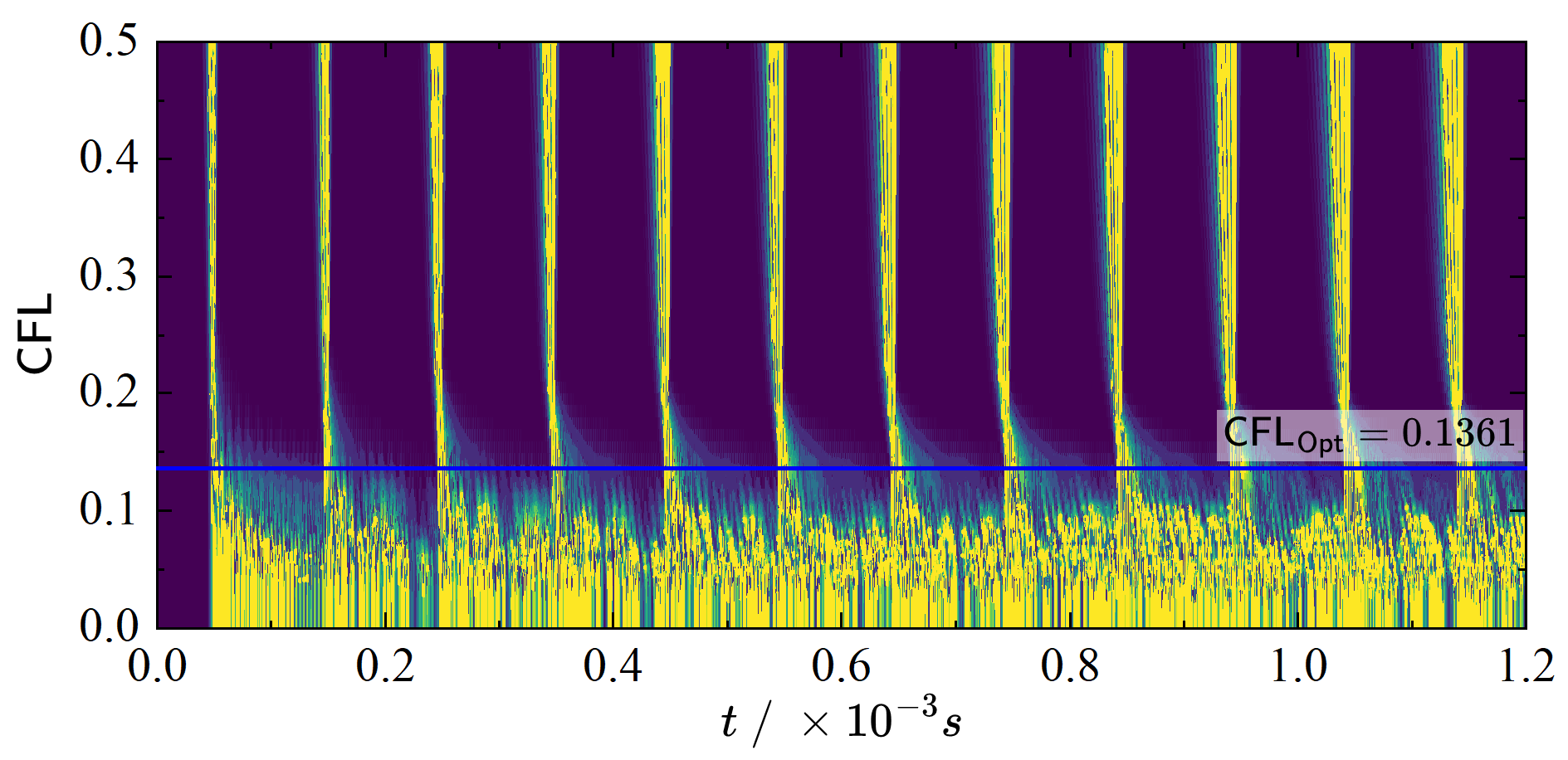}}
	\caption{The mid-point velocity of the bar predicted by the four-sub-step implicit method ($s=4$): (a-b) fourth-order accuracy and (c-d) fifth-order accuracy with $A(\alpha)$-stability.}
    \label{fig:bar_s4_2}
\end{figure}

As shown in \cref{fig:bar_s4_2}(a-b), the numerical responses obtained with different CFL numbers are compared under $\infrho=0.8$ and $r=0$. The optimal CFL number is computed as $\mathsf{CFL}_{\mathsf{Opt}}=0.6021$. The numerical solution corresponding to $\mathsf{CFL}_{\mathsf{Opt}}$ agrees better with the exact solution than those obtained with $\mathsf{CFL}=0.4$ and $\mathsf{CFL}=0.8$, especially near the wave fronts where dispersion errors are more significant. The corresponding error distribution in \cref{fig:bar_s4_2}(b) further confirms that the optimal CFL number effectively reduces the numerical error introduced by the spatial and temporal discretizations. 

For the $A(\alpha)$-stable fifth-order four-sub-step implicit method ($\varpi_1=0$), which requires $\gamma_1=1.3453664197803336$, the optimized mass matrix with $r=1/3$ is needed to eliminate the remaining spatial error. Consequently, the leading dispersion error in \cref{eq:omega_3node_4} is reduced to $\left[\varpi_3\dt^6/336-\dx^6/(75600c_0^6)\right]\omega^7$, and thus the optimal CFL number for the fifth-order member is derived as
\begin{equation}
\mathsf{CFL}_\mathsf{Opt}=\dfrac{1}{\sqrt[6]{225\varpi_3}}\approx0.136115999461561. 
\end{equation}
Therefore, the $A(\alpha)$-stable fifth-order four-sub-step implicit method requires not only the optimized mass matrix $r=1/3$ but also an optimal CFL number to achieve further improvement in dispersion accuracy to $\mathcal{O}(\omega^9)$. 

The numerical results shown in \cref{fig:bar_s4_2}(c-d) indicate that the smallest error is achieved at approximately $\mathsf{CFL}=0.177$, rather than exactly at the theoretical value. This discrepancy is mainly attributed to the $A(\alpha)$-stability of the fifth-order four-sub-step implicit method. Although the dispersion analysis predicts the optimal CFL number by balancing the phase errors, the $A(\alpha)$-stability property introduces a narrow frequency range in which the spectral radius may slightly exceed unity, resulting in weak numerical amplification. Such amplitude distortion modifies the overall numerical error and prevents the theoretically predicted CFL number from achieving the minimum total error in practical computations. Therefore, the deviation between the theoretical and numerical optimal CFL numbers does not indicate the failure of the dispersion analysis, but rather reflects the combined influence of phase accuracy and amplitude behavior. A similar phenomenon has also been observed previously in \cref{fig:bar_s4_1}(c-d). 

\paragraph{The five- and six-sub-step implicit methods}
For the five- and six-sub-step implicit methods, the dispersion relations are given by \cref{eq:omega_3node_5,eq:omega_3node_6}, respectively. It can be observed that the optimized mass matrix with $r=1/3$ is required to eliminate the leading $\omega^5$ spatial error. After this optimization, the leading dispersion terms of both methods are reduced to $\mathcal{O}(\omega^7)$, and the remaining spatial and temporal contributions can be further balanced by selecting an optimal CFL number: 
\begin{equation}
\mathsf{CFL}_\mathsf{Opt}=\begin{cases}
\dfrac{1}{\sqrt[6]{-90\psi_2}}=\dfrac{1}{\sqrt[6]{-90(4200\gamma_1^6-10920\gamma_1^5+8400\gamma_1^4-2800\gamma_1^{3}+455\gamma_1^{2}-35\gamma_1+1)}}, & s=5\\
\dfrac{1}{\sqrt[6]{15\varphi_1}}=\dfrac{1}{\sqrt[6]{15(5040\gamma_1^6 - 15120\gamma_1^5 + 12600\gamma_1^4 - 4200\gamma_1^3 + 630\gamma_1^2 - 42\gamma_1 + 1)}}, & s=6
\end{cases}
\end{equation}
where $\gamma_1$ should be selected to make $\mathsf{CFL}_\mathsf{Opt}$ real and positive. Therefore, both five- and six-sub-step implicit methods require the optimized mass matrix $r=1/3$ to fully exploit their high-order temporal accuracy. Further improvement of dispersion accuracy can be achieved through an optimal CFL number to eliminate the remaining leading frequency error.

\begin{figure}[htbp]
	\centering 
        \includegraphics[scale=1.0]{bar_s1_leg_down}\\
    \subfigure[$\infrho=0.8~\&~r=1/3$]{
        \includegraphics[scale=1.0]{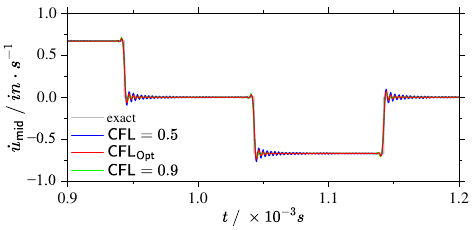}}
    \subfigure[$\infrho=0.8~\&~r=1/3$]{
        \includegraphics[scale=1.0]{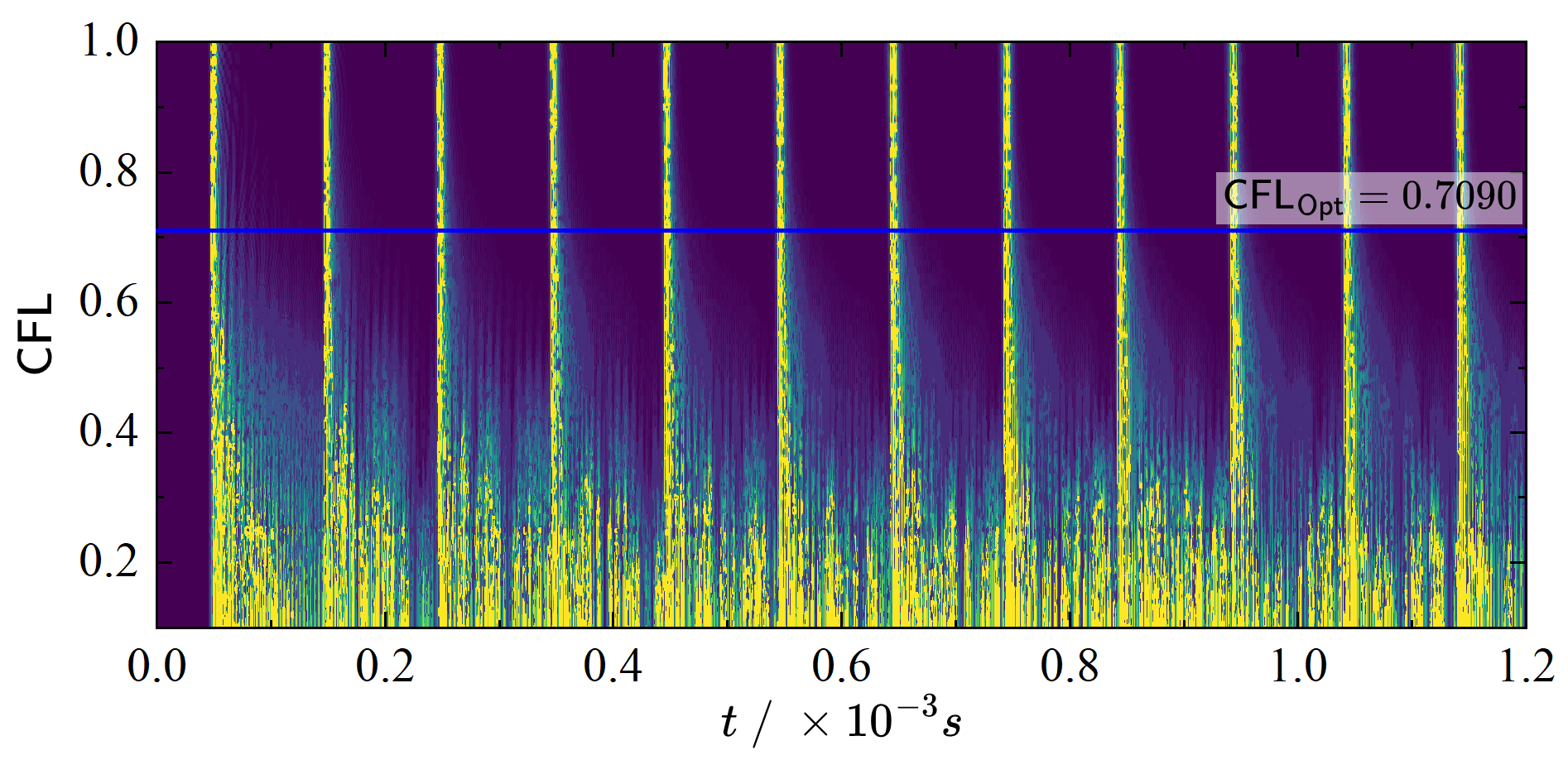}}
    \subfigure[$\infrho=0.0~\&~r=1/3$]{
        \includegraphics[scale=1.0]{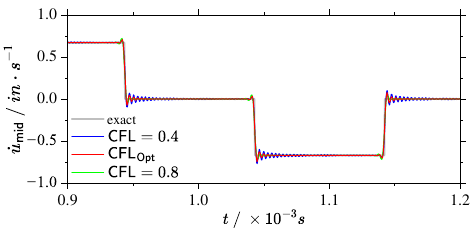}}
    \subfigure[$\infrho=0.0~\&~r=1/3$]{
        \includegraphics[scale=1.0]{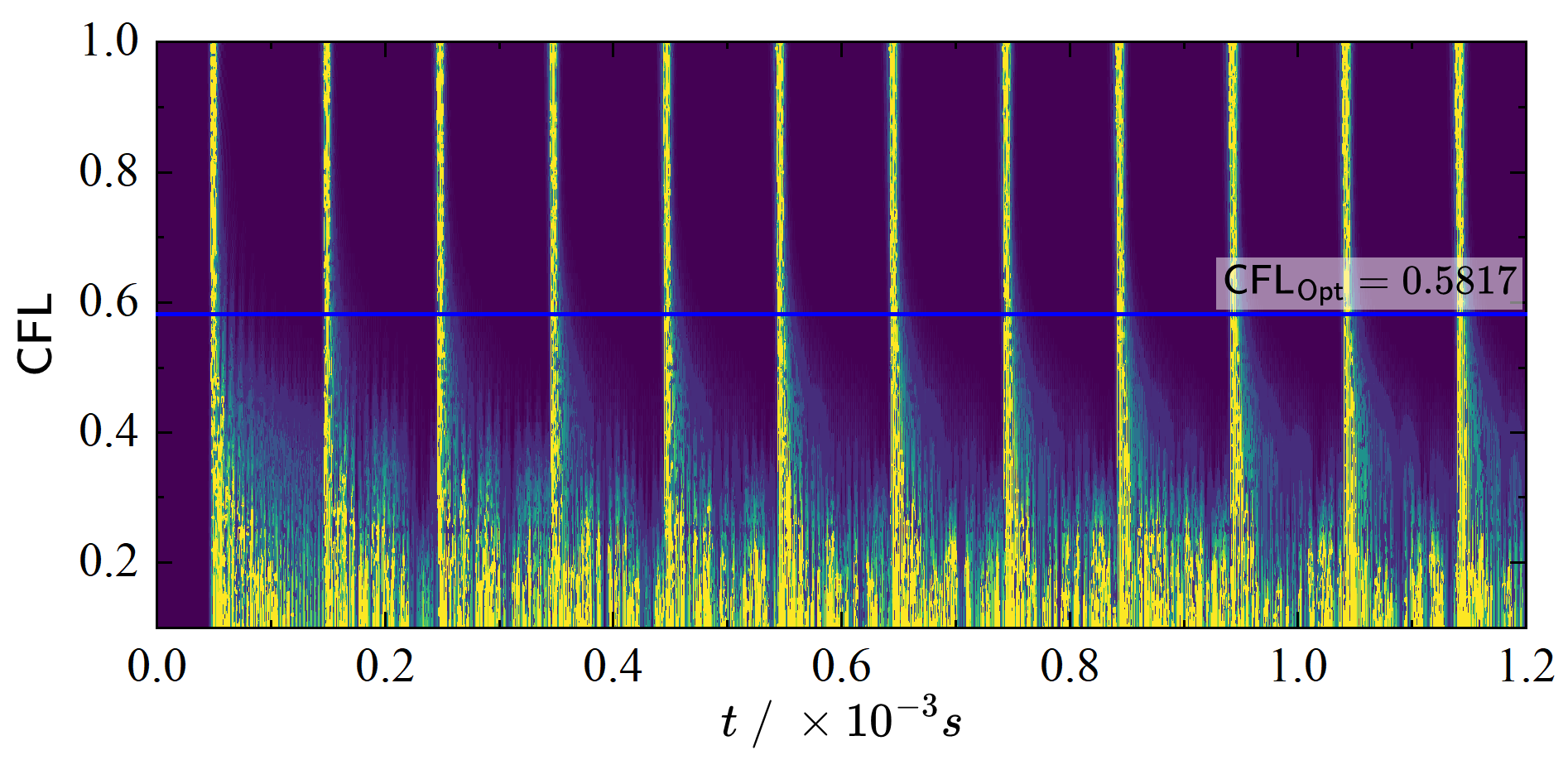}}
	\caption{The mid-point velocity of the bar predicted by the sixth-order six-sub-step implicit method ($s=6$).}
    \label{fig:bar_s6_2}
\end{figure}

The numerical performance of the six-sub-step implicit method is investigated in \cref{fig:bar_s6_2} as a representative example of the five- and six-sub-step methods. Since both methods share the same dispersion optimization strategy, i.e., employing the optimized mass matrix $r=1/3$ and selecting an optimal CFL number to balance the remaining $\omega^7$ terms, only the six-sub-step method is considered here for brevity. As shown in \cref{fig:bar_s6_2}(a-b), the numerical results obtained with $\infrho=0.8$ and $r=1/3$ are compared for different CFL numbers. The optimal CFL number is $\mathsf{CFL}_{\mathsf{Opt}}=0.7090$. The numerical responses obtained with the optimal CFL number show good agreement with the exact solution. Similarly, \cref{fig:bar_s6_2}(c-d) presents the results of the strongly dissipative method with $\infrho=0$. In this case, the optimal CFL number is $\mathsf{CFL}_{\mathsf{Opt}}=0.5817$. The numerical results indicate that the solution obtained with the optimal CFL number provides high accuracy and effectively suppresses spurious oscillations near wave fronts. 


\section{Conclusions}\label{sec:conclusions}

This study has developed a generalized directly self-starting $s$-sub-step implicit integration framework for first- and second-order transient problems. By relaxing the constraint that fixes the last sub-step at the end of each time interval, the sub-step locations are treated as independent design variables and can therefore be systematically exploited to improve the temporal accuracy. The proposed framework is uniformly represented in the RK form while directly operating on the original governing equations, and the identical effective matrices are retained over all sub-steps. The principal conclusions are summarized as follows:

\begin{itemize}

\item For a prescribed number of sub-steps $s$, two complementary algorithm configurations can be systematically constructed. The first achieves $s$th-order accuracy while retaining user-controllable high-frequency dissipation through the parameter $\infrho$. The second fully exploits the available sub-step locations to satisfy one additional order condition and consequently achieves $(s+1)$th-order accuracy without increasing the number of sub-steps. In the latter case, however, the sub-step sizes are completely determined and the high-frequency dissipation level becomes fixed. The proposed construction has been developed up to $s=6$, providing implicit methods with accuracy as high as seventh order.

\item The stability and spectral analyses demonstrate that the proposed framework provides favorable stability together with effective control of undesirable high-frequency responses. The $s$th-order members admit a full ranges of the user-specified parameter $\infrho\in[0,~1]$ while preserving the required unconditional stability. For some $(s+1)$th-order high-order members, the additional accuracy condition leads to $A(\alpha)$-stability rather than strict unconditional stability. Nevertheless, the corresponding stability angles are extremely close to $90^\circ$, and only very weak amplification is observed over narrow frequency intervals.

\item The analytical amplitude and phase errors reveal a clear parity-dependent superconvergence pattern. For the $s$th-order methods in undamped systems, odd numbers of sub-steps exhibit one-order phase superconvergence, whereas even numbers of sub-steps exhibit one-order amplitude superconvergence. Moreover, an appropriate optimization of $\gamma_1$ can eliminate additional leading error terms. In particular, the three- and five-sub-step methods can achieve phase accuracies three orders higher than their formal orders, whereas the four- and six-sub-step methods can achieve amplitude accuracies three orders higher than their formal orders. These results demonstrate that the sub-step locations can be utilized not only for controlling high-frequency dissipation or elevating the formal order, but also for selectively optimizing the spectral accuracy.

\item The dispersion analyses for the elastic bar demonstrate that the accuracy of spatially discretized problems is governed by the interaction between spatial and temporal discretization errors. With two-node linear finite elements, the spatial dispersion error may dominate and consequently mask the advantages of high-order temporal integrators. In particular, simply reducing the time step does not necessarily make the numerical solution approach the exact solution of the original PDE, because the temporal integrator may increasingly resolve the spurious high-frequency components introduced by the spatial discretization. Appropriate CFL numbers can instead balance the spatial and temporal dispersion errors and substantially improve the numerical accuracy. When three-node quadratic elements are employed, the reduced spatial error enables the high-order temporal accuracy to be more fully exploited. Further improvements can be obtained through optimized mass matrices and CFL numbers, although for some $A(\alpha)$-stable high-order members the weak amplification over a narrow frequency range may cause the practically optimal CFL number to deviate slightly from that predicted solely by the dispersion analysis.
\end{itemize}

Numerical examples confirm the theoretically predicted $s$th- and $(s+1)$th-order convergence rates and demonstrate the capability of the proposed methods to suppress spurious high-frequency components while accurately preserving low-frequency responses. Increasing the number of sub-steps systematically increases the attainable temporal accuracy, and therefore generally improves the numerical solution when the same time step is employed. This improvement is accompanied by an increased computational cost because more sub-step solutions are required within each time step. Hence, the number of sub-steps should be selected according to the required balance between accuracy and computational efficiency.

\backmatter

\bmhead{Acknowledgements}
This work was supported by National Natural Science Foundation of China (No.~12502066 and No.~12272105), Fundamental Research Funds for the Central Universities (No.~HIT.DZJJ.2026019), China Postdoctoral Science Foundation (No.~2024M764165), and Heilongjiang Postdoctoral Grant (No.~LBH-Z23153).

\bmhead{Conflict of interests}
The authors declare that they have no known competing financial interests or personal relationships that could have appeared to influence the work reported in this paper.

\bmhead{Data availability}
Data sharing not applicable to this article as no datasets were generated or analyzed during the current study.

\bmhead{Author Contributions}
\begin{description}
	\item[Jinze Li] Conceptualization, Methodology, Writing - Original Draft \& Review, Visualization, Formal analysis, Funding acquisition.
	\item[Yaokun Liu] Data Curation, Visualization, Writing - Review.
	\item[Kewei Chen] Resources, Supervision, Writing - Review.
	\item[Hua Li] Supervision, Writing - Review.
	\item[Kaiping Yu] Funding acquisition, Supervision, Writing - Review.
\end{description}

\begin{appendices}

\section{algorithmic parameters for more than three sub-steps}\label{app:a}
\subsection{Four sub-steps: $s=4$}\label{app:s4}

\begin{subequations}\label{eq:s_eq_4}
The four-sub-step implicit method is expressed by substituting $s=4$ into \cref{eq:but_alg} as 
\begin{equation}
\begin{NiceArray}{c|c}[cell-space-limits=3pt,columns-width=0.5cm]
\mbf{c} & \mbf{A}\\ \hline & \mbf{b}\T
\end{NiceArray}\quad=\quad\begin{NiceArray}{c|cccc}[cell-space-limits=3pt,columns-width=0.7cm]
\gamma_1 & \alpha_{11} \\ 
\gamma_2 & \alpha_{21} & \alpha_{22} \\ 
\gamma_3 & \alpha_{31} & \alpha_{32} & \alpha_{33} \\ 
\gamma_4 & \alpha_{41} & \alpha_{42} & \alpha_{43} & \alpha_{44} \\ \hline
& \beta_1 & \beta_2 & \beta_3 & \beta_4
\end{NiceArray}
\end{equation}
where the elements in $\mbf{A}$ and $\mbf{b}$ are computed as
\begin{align}
\alpha_{11}&=\alpha_{22}=\alpha_{33}=\alpha_{44}=\gamma_1 \\
\alpha_{21}&=\gamma_2-\gamma_1\\ 
\alpha_{31}&=\gamma_3-\alpha_{32}-\gamma_1 \\
\alpha_{41}&=\gamma_4-\alpha_{42}-\alpha_{43}-\gamma_1\\
\alpha_{32}&=\dfrac{(\gamma_1-\gamma_3)(\gamma_2-\gamma_3)(24\gamma_1 ^{3}-36\gamma_1 ^{2}+12\gamma_1-1)}{2(\gamma_2-\gamma_1)(2\gamma_2-1)(6\gamma_1 ^{2}-6\gamma_1+1)}  \\
\alpha_{42}&=\dfrac{6(\gamma_2-\gamma_3)^2(\gamma_3-\gamma_1)(24\gamma_1 ^{3}-36\gamma_1 ^{2}+12\gamma_1-1)\beta_3-(6\gamma_1 ^{2}-6\gamma_1+1) ^{2}(2\gamma_3-1)(2\gamma_2-1)}{12(\gamma_2-\gamma_1)(\gamma_2-\gamma_3)(2\gamma_2-1)(6\gamma_1 ^{2}-6\gamma_1+1)\beta_4}  \\
\alpha_{43}&=\dfrac{(1-2\gamma_2)(6\gamma_1 ^{2}-6\gamma_1+1)}{12(\gamma_2-\gamma_3)(\gamma_1-\gamma_3)\beta_4} \\ 
\beta_1&=\dfrac{3-4(\gamma_2+\gamma_3+\gamma_4)+6(\gamma_2\gamma_3+\gamma_2\gamma_4+\gamma_3\gamma_4)-12\gamma_2\gamma_3\gamma_4}{12(\gamma_1-\gamma_2)(\gamma_1-\gamma_3)(\gamma_1-\gamma_4)} \\
\beta_2&=\dfrac{3-4(\gamma_1+\gamma_3+\gamma_4)+6(\gamma_1\gamma_3+\gamma_1\gamma_4+\gamma_3\gamma_4)-12\gamma_1\gamma_3\gamma_4}{12(\gamma_2-\gamma_1)(\gamma_2-\gamma_3)(\gamma_2-\gamma_4)} \\
\beta_3&=\dfrac{3-4(\gamma_1+\gamma_2+\gamma_4)+6(\gamma_1\gamma_2+\gamma_1\gamma_4+\gamma_2\gamma_4)-12\gamma_1\gamma_2\gamma_4}{12(\gamma_3-\gamma_1)(\gamma_3-\gamma_2)(\gamma_3-\gamma_4)} \\
\beta_4&=\dfrac{3-4(\gamma_1+\gamma_2+\gamma_3)+6(\gamma_1\gamma_2+\gamma_1\gamma_3+\gamma_2\gamma_3)-12\gamma_1\gamma_2\gamma_3}{12(\gamma_4-\gamma_1)(\gamma_4-\gamma_2)(\gamma_4-\gamma_3)}
\end{align}
\end{subequations}
and the sub-step sizes $\gamma_j~(j=1,~\cdots,~4)$ are free parameters. 

\begin{figure}[htbp]
	\centering 
	\includegraphics[scale=1.0]{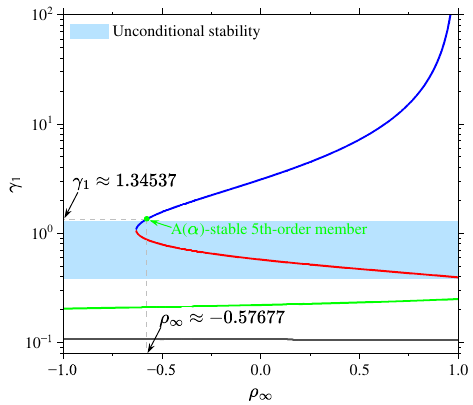}
	\caption{The variations of $\gamma_1$ with $\infrho\in[-1,~1]$ for the four-sub-step implicit method.}
	\label{fig:s4_r1_rho}
\end{figure}

For fourth-order accuracy, the first sub-step size $\gamma_1$ is related to the user-specified parameter $\infrho$ through
\begin{equation}
\dfrac{24\gamma_1^4-96\gamma_1 ^{3}+72\gamma_1 ^{2}-16\gamma_1+1}{24\gamma_1 ^{4}} = \left| \infrho \right| 
\end{equation}
whereas the remaining sub-step sizes $\gamma_j~(j=2,~3,~4)$ remain free parameters. For example, Burrage \cite{burrage_EfficientlyImplementableAlgebraically_1982} adopted $\gamma_2=3\gamma_1(2\gamma_1-1)^{2}/(12\gamma_1^{2}-6\gamma_1+1),~\gamma_3=1-\gamma_2$ and $\gamma_4=1-\gamma_1$ to achieve fourth-order accuracy and BN-stability within the four-sub-step framework \eqref{eq:s_eq_4} and Li et al. \cite{li_DirectlySelfstartingHigherorder_2022} used $\gamma_j=j\cdot\gamma_1~(j=2,~3)$ and $\gamma_4=1$ to splite the equal sub-step sizes for the first three sub-steps. The unconditional stability given in \cref{eq:Ey} requires $\gamma_1$ to satisfy
\begin{equation}\label{eq:s4_uc}
\gamma_1\in\left[\dfrac{3+\sqrt{3}}{12},~1.2805797612753055 \right],
\end{equation}
where $1.2805797612753055$ is an approximate root of $576\gamma_1^5 - 1224\gamma_1^4 + 768\gamma_1^3 - 204\gamma_1^2 + 24\gamma_1 - 1=0$. Accordingly, the first sub-step size $\gamma_1$ is determined by
\begin{equation}\label{eq:s4_r1_set}
\gamma_1\in\left\{\gamma_1~\bigg|~ \dfrac{3+\sqrt{3}}{12}\le\gamma_1\le1.2805797612753055,~24(1-\infrho)\gamma_1^4-96\gamma_1^{3}+72\gamma_1^{2}-16\gamma_1+1=0\right\}
\end{equation}
where $\infrho\in[1-2\cos(2\pi/9)/\cos(\pi/9),~1]$. \cref{eq:s4_r1_set} establishes the relationship between the first sub-step size $\gamma_1$ and the user-specified parameter $\infrho$. For a prescribed value of $\infrho$, the quartic equation $24(1-\infrho)\gamma_1^4-96\gamma_1^{3}+72\gamma_1^{2}-16\gamma_1+1=0$ yields four candidate solutions for $\gamma_1$. Since not all mathematical solutions necessarily satisfy the unconditional stability, the admissible solution is selected by additionally enforcing \cref{eq:s4_uc}. In this manner, a unique relationship between $\gamma_1$ and the user-specified parameter $\infrho$ can be established, as shown in \cref{fig:s4_r1_rho}.

For fifth-order accuracy, the sub-step sizes $\gamma_j~(j=1,~\cdots,~4)$ are determined by requiring the local truncation errors defined in \cref{eq:def_acc} to be of order $\mathcal{O}(\dt^6)$. The resulting conditions for achieving fifth-order accuracy with four sub-steps are
\begin{subequations}\label{eq:s4_r1234}
\begin{align}
\gamma_1&=\mathsf{RootOf}\left(120\gamma_1^4 - 240\gamma_1^3 + 120\gamma_1^2 - 20\gamma_1 + 1\right)\label{eq:s4_fifth_r1}\\
\gamma_2&= \dfrac{60\gamma_1^{3}-80\gamma_1^{2}+25\gamma_1-2}{5(24\gamma_1^{3}-36\gamma_1^{2}+12\gamma_1-1)}  \\
\gamma_3&=- \dfrac{120\gamma_1^5-360\gamma_1^4+342\gamma_1^{3}-126\gamma_1^{2}+19\gamma_1-1}{(6\gamma_1^{2}-6\gamma_1+1)(20\gamma_1^{2}-10\gamma_1+1)} \\
\gamma_4&=1-\gamma_1
\end{align}
\end{subequations}
where $\mathsf{RootOf}(f)$ denotes the set of all roots of $f=0$. \cref{eq:s4_fifth_r1} yields four candidate values of $\gamma_1$, none of which satisfies the unconditional stability given in \cref{eq:s4_uc}. Consequently, a strictly unconditionally stable fifth-order four-sub-step implicit method does not exist within the present framework \eqref{eq:but_alg}. Nevertheless, numerical investigations indicate that the stability region is maximized when
\begin{equation}
\gamma_1= \dfrac{1}{2}\left(1+\sqrt{\mu}+\sqrt{1-\mu+ \dfrac{1}{3\sqrt{\mu}} }\right) 
\end{equation}
where $\mu=1/3+2\cos\left(\arccos(\sqrt{10}/4)/3\right)/(3\sqrt{10})$. Therefore, this value, approximately $\gamma_1=1.3453664197803336$, is adopted in the present study as a practical choice for constructing the fifth-order four-sub-step implicit method. With this choice, the resulting fifth-order method is $A(\alpha)$-stable, where $\alpha$ is numerically estimated as $89.99173663^\circ$. This value is extremely close to the limiting case $\alpha=90^\circ$, corresponding to unconditional stability.

\subsection{Five sub-steps: $s=5$}
\begin{subequations}\label{eq:s_eq_5}
The five-sub-step implicit method is expressed by substituting $s=5$ into \cref{eq:but_alg} as 
\begin{equation}
\begin{NiceArray}{c|c}[cell-space-limits=3pt,columns-width=0.5cm]
\mbf{c} & \mbf{A}\\ \hline & \mbf{b}\T
\end{NiceArray}\quad=\quad\begin{NiceArray}{c|ccccc}[cell-space-limits=3pt,columns-width=0.7cm]
\gamma_1 & \alpha_{11} \\ 
\gamma_2 & \alpha_{21} & \alpha_{22} \\ 
\gamma_3 & \alpha_{31} & \alpha_{32} & \alpha_{33} \\ 
\gamma_4 & \alpha_{41} & \alpha_{42} & \alpha_{43} & \alpha_{44} \\ 
\gamma_5 & \alpha_{51} & \alpha_{52} & \alpha_{53} & \alpha_{54} &\alpha_{55} \\ \hline
& \beta_1 & \beta_2 & \beta_3 & \beta_4 & \beta_5
\end{NiceArray}
\end{equation}
where the elements in $\mbf{A}$ and $\mbf{b}$ are computed as
\begin{align}
\alpha_{11}&=\alpha_{22}=\alpha_{33}=\alpha_{44}=\alpha_{55}=\gamma_1 \\
\alpha_{21}&=\gamma_2-\gamma_1\\ 
\alpha_{31}&=\gamma_3-\alpha_{32}-\gamma_1 \\
\alpha_{41}&=\gamma_4-\alpha_{42}-\alpha_{43}-\gamma_1\\
\alpha_{51}&=\gamma_5-\alpha_{52}-\alpha_{53}-\alpha_{54}-\gamma_1\\
\alpha_{32}&=\dfrac{-120\gamma_1^4+240\gamma_1^{3}-120\gamma_1^{2}+20\gamma_1-1}{120\beta_5\alpha_{54}\alpha_{43}(\gamma_1-\gamma_2)}\\  
\alpha_{42}&= \dfrac{\left\{\begin{aligned}
&120\gamma_1^{3}\big[(\beta_5\alpha_{53}+\beta_4\alpha_{43})(\gamma_1-2)+\beta_5\alpha_{54}\alpha_{43}\big]+60(2\beta_5\alpha_{53}+2\beta_4\alpha_{43}-3\beta_5\alpha_{54}\alpha_{43})\gamma_1^{2}\\
&+20(3\beta_5\alpha_{54}\alpha_{43}-\beta_5\alpha_{53}-\beta_4\alpha_{43}-6\beta_5^{2}\alpha_{54}^{2}\alpha_{43}^{2})\gamma_1+120\beta_5^{2}\alpha_{54}^{2}\alpha_{43}^{2}\gamma_3-5\beta_5\alpha_{54}\alpha_{43}+\beta_5\alpha_{53}+\beta_4\alpha_{43}
\end{aligned}\right\}}{120\beta_5^{2}\alpha_{54}^{2}\alpha_{43}(\gamma_1-\gamma_2)} \\
\alpha_{52}&= \dfrac{\left\{ \begin{aligned}
&(\beta_5\beta_3\alpha_{54}-\beta_5\beta_4\alpha_{53}-\beta_4^{2}\alpha_{43})(120\gamma_1^4-240\gamma_1^{3}+120\gamma_1^{2}-20\gamma_1+1)-120\beta_5\beta_4\alpha_{54}\alpha_{43}\gamma_1^{3}\\
&-60\beta_5\alpha_{54}\alpha_{43}(2\beta_5\alpha_{54}-3\beta_4)\gamma_1^{2}+20\big[6\beta_5^{2}\alpha_{54}^{2}\alpha_{43}(1-\beta_5\alpha_{53}-\beta_5\alpha_{54})-3\beta_5\beta_4\alpha_{54}\alpha_{43}\big]\gamma_1\\
&+5\beta_5\alpha_{54}\alpha_{43}\big[24\beta_5^{2}\alpha_{54}(\alpha_{53}\gamma_3+\alpha_{54}\gamma_4)-4\beta_5\alpha_{54}+\beta_4\big]
\end{aligned}\right\}}{120\beta_5^{3}\alpha_{54}^{2}\alpha_{43}(\gamma_1-\gamma_2)} \\
\alpha_{43}&= \dfrac{\left\{ \begin{aligned}
&6(1+2\beta_1\gamma_5-2\gamma_2)\gamma_1^{2}+6\big[2\gamma_2+2\gamma_1(\gamma_1-\gamma_5)\beta_1-1\big]\gamma_1-12\beta_1\gamma_1^{3}+1-2\gamma_2\\
&-12\beta_5\alpha_{53}(\gamma_2-\gamma_3)(\gamma_1-\gamma_3)-12\beta_5\alpha_{54}(\gamma_2-\gamma_4)(\gamma_1-\gamma_4)
\end{aligned}\right\}}{12\beta_4(\gamma_1-\gamma_3)(\gamma_2-\gamma_3)} \\
\alpha_{53}&= \dfrac{\left\{ \begin{aligned}
&-120\beta_5^{2}\alpha_{54}^{2}(\gamma_1-\gamma_4)(\gamma_2-\gamma_4)-10\beta_5\alpha_{54}(2\gamma_2-1)(6\gamma_1^{2}-6\gamma_1+1)\\
&-\beta_4(120\gamma_1^{3}\gamma_2-60\gamma_1^{3}-180\gamma_1^{2}\gamma_2+80\gamma_1^{2}+60\gamma_1\gamma_2-25\gamma_1-5\gamma_2+2)
\end{aligned}\right\}}{120\beta_5^{2}\alpha_{54}(\gamma_1-\gamma_3)(\gamma_2-\gamma_3)} \\
\alpha_{54}&= \dfrac{\left\{ \begin{aligned}
&60\beta_1\gamma_1^{2}(\gamma_1-\gamma_4)(\gamma_1-\gamma_5)+60\beta_2\gamma_1\gamma_2(\gamma_2-\gamma_4)(\gamma_2-\gamma_5)+60\beta_3\gamma_1\gamma_3(\gamma_3-\gamma_4)(\gamma_3-\gamma_5)\\
&+10\gamma_1(1-\gamma_1)(6\gamma_2\gamma_3-3\gamma_2-3\gamma_3+2)-5\gamma_1(6\gamma_4\gamma_5-4\gamma_4-4\gamma_5+3)-10\gamma_2\gamma_3+5\gamma_2+5\gamma_3-3
\end{aligned}\right\} }{60\beta_5(\gamma_1-\gamma_4)(\gamma_2-\gamma_4)(\gamma_3-\gamma_4)} \\ 
\beta_1&=\dfrac{\left\{ \begin{aligned}
&60\gamma_2\gamma_3\gamma_4\gamma_5-30(\gamma_2\gamma_3\gamma_4+\gamma_2\gamma_3\gamma_5+\gamma_2\gamma_4\gamma_5+\gamma_3\gamma_4\gamma_5)\\
&+20(\gamma_2\gamma_3+\gamma_2\gamma_4+\gamma_2\gamma_5+\gamma_3\gamma_4+\gamma_3\gamma_5+\gamma_4\gamma_5)-15(\gamma_2+\gamma_3+\gamma_4+\gamma_5)+12
\end{aligned} \right\} }{60(\gamma_1-\gamma_2)(\gamma_1-\gamma_3)(\gamma_1-\gamma_4)(\gamma_1-\gamma_5)}\\
\beta_2&=\dfrac{\left\{ \begin{aligned}
&60\gamma_1\gamma_3\gamma_4\gamma_5-30(\gamma_1\gamma_3\gamma_4+\gamma_1\gamma_3\gamma_5+\gamma_1\gamma_4\gamma_5+\gamma_3\gamma_4\gamma_5)\\
&+20(\gamma_1\gamma_3+\gamma_1\gamma_4+\gamma_1\gamma_5+\gamma_3\gamma_4+\gamma_3\gamma_5+\gamma_4\gamma_5)-15(\gamma_1+\gamma_3+\gamma_4+\gamma_5)+12
\end{aligned} \right\} }{60(\gamma_2-\gamma_1)(\gamma_2-\gamma_3)(\gamma_2-\gamma_4)(\gamma_2-\gamma_5)}\\
\beta_3&=\dfrac{\left\{ \begin{aligned}
&60\gamma_1\gamma_2\gamma_4\gamma_5-30(\gamma_1\gamma_2\gamma_4+\gamma_1\gamma_2\gamma_5+\gamma_1\gamma_4\gamma_5+\gamma_2\gamma_4\gamma_5)\\
&+20(\gamma_1\gamma_2+\gamma_1\gamma_4+\gamma_1\gamma_5+\gamma_2\gamma_4+\gamma_2\gamma_5+\gamma_4\gamma_5)-15(\gamma_1+\gamma_2+\gamma_4+\gamma_5)+12
\end{aligned} \right\} }{60(\gamma_3-\gamma_1)(\gamma_3-\gamma_2)(\gamma_3-\gamma_4)(\gamma_3-\gamma_5)}\\
\beta_4&=\dfrac{\left\{ \begin{aligned}
&60\gamma_1\gamma_2\gamma_3\gamma_5-30(\gamma_1\gamma_2\gamma_3+\gamma_1\gamma_2\gamma_5+\gamma_1\gamma_3\gamma_5+\gamma_2\gamma_3\gamma_5)\\
&+20(\gamma_1\gamma_2+\gamma_1\gamma_3+\gamma_1\gamma_5+\gamma_2\gamma_3+\gamma_2\gamma_5+\gamma_3\gamma_5)-15(\gamma_1+\gamma_2+\gamma_3+\gamma_5)+12
\end{aligned} \right\} }{60(\gamma_4-\gamma_1)(\gamma_4-\gamma_2)(\gamma_4-\gamma_3)(\gamma_4-\gamma_5)}\\
\beta_5&=\dfrac{\left\{ \begin{aligned}
&60\gamma_1\gamma_2\gamma_3\gamma_4-30(\gamma_1\gamma_2\gamma_3+\gamma_1\gamma_2\gamma_4+\gamma_1\gamma_3\gamma_4+\gamma_2\gamma_3\gamma_4)\\
&+20(\gamma_1\gamma_2+\gamma_1\gamma_3+\gamma_1\gamma_4+\gamma_2\gamma_3+\gamma_2\gamma_4+\gamma_3\gamma_4)-15(\gamma_1+\gamma_2+\gamma_3+\gamma_4)+12
\end{aligned} \right\} }{60(\gamma_5-\gamma_1)(\gamma_5-\gamma_2)(\gamma_5-\gamma_3)(\gamma_5-\gamma_4)}
\end{align}
\end{subequations}
and the sub-step sizes $\gamma_j~(j=1,~\cdots,~5)$ are free parameters. 

\begin{figure}[htbp]
	\centering 
	\includegraphics[scale=1.0]{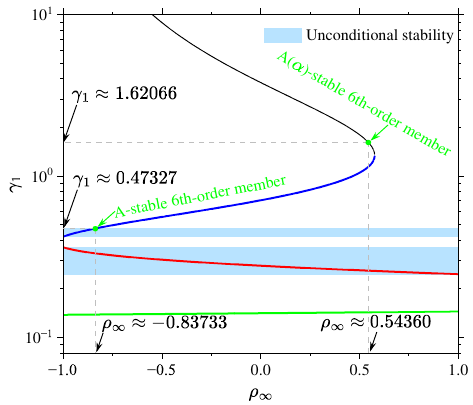}
	\caption{The variations of $\gamma_1$ with $\infrho\in[-1,~1]$ for the five-sub-step implicit method \eqref{eq:s_eq_5}.}
	\label{fig:s5_r1_rho}
\end{figure}

For fifth-order accuracy, the first sub-step size $\gamma_1$ is related to the user-specified parameter $\infrho$ through
\begin{equation}\label{eq:s5_D_infty}
\dfrac{120\gamma_1^5-600\gamma_1^4+600\gamma_1^3-200\gamma_1^2+25\gamma_1-1}{120\gamma_1 ^{5}} = \left| \infrho \right| 
\end{equation}
whereas the remaining sub-step sizes $\gamma_j~(j=2,~3,~4,~5)$ remain free parameters. For example, Li et al. \cite{li_DirectlySelfstartingHigherorder_2022} used $\gamma_j=j\cdot\gamma_1~(j=2,~3,~4)$ and $\gamma_5=1$ to splite the equal sub-step sizes for the first four sub-steps. The unconditional stability given in \cref{eq:Ey} requires $\gamma_1$ to satisfy
\begin{equation}\label{eq:s5_uc}
\gamma_1\in\left[0.2465051931428203,~ \dfrac{5+\sqrt{5}}{20}  \right]\cup\left[ \dfrac{15+\sqrt{105}}{60},~0.4732683912582953 \right],
\end{equation}
where $0.2465051931428203$ is an approximate root of $240\gamma_1^5 - 600\gamma_1^4 + 600\gamma_1^3 - 200\gamma_1^2 + 25\gamma_1 - 1=0$ and $0.4732683912582953$ is an appropriate root of $720\gamma_1^5 - 1800\gamma_1^4 + 1200\gamma_1^3 - 300\gamma_1^2 + 30\gamma_1 - 1=0$. Accordingly, the first sub-step size $\gamma_1$ is determined by
\begin{equation}\label{eq:s5_r1_set}
\gamma_1\in\left\{\gamma_1~\bigg|~ \text{Eq.}~\eqref{eq:s5_uc},~ 120  (1 + \infrho)  \gamma_1^5 - 600  \gamma_1^4 + 600  \gamma_1^3 - 200  \gamma_1^2 + 25  \gamma_1 - 1=0\right\}
\end{equation}
where $\infrho\in[-1,~1]$. \cref{eq:s5_r1_set} establishes the relationship between the first sub-step size $\gamma_1$ and the user-specified parameter $\infrho$. For a prescribed value of $\infrho$, the quintic equation $120  (1 + \infrho)  \gamma_1^5 - 600  \gamma_1^4 + 600  \gamma_1^3 - 200  \gamma_1^2 + 25  \gamma_1 - 1=0$ yields five candidate solutions for $\gamma_1$. Since not all mathematical solutions necessarily satisfy the unconditional stability, the admissible solution is selected by additionally enforcing \cref{eq:s5_uc}. In this manner, a unique relationship between $\gamma_1$ and the user-specified parameter $\infrho$ can be established, as shown in \cref{fig:s5_r1_rho}.

\cref{fig:s5_r1_rho} illustrates the variation of the first sub-step size $\gamma_1$ with respect to the user-specified parameter $\infrho$ for the five-sub-step implicit method. Similar to the lower-order cases, the first sub-step size $\gamma_1$ is determined by the high-frequency spectral radius through the dissipation control condition \eqref{eq:s5_D_infty}. However, the fifth-order five-sub-step implicit method \eqref{eq:s_eq_5} can maintain unconditional stability while allowing $\infrho$ to vary over the entire range of $[-1,~1]$. This indicates that the fifth-order five-sub-step implicit method preserves both high-order accuracy and controllable high-frequency dissipation. In particular, by prescribing different values of $\infrho$, users can directly adjust the asymptotic spectral radius without violating the stability requirement.

For sixth-order accuracy, the sub-step sizes $\gamma_j~(j=1,~\cdots,~5)$ are determined by requiring the local truncation errors defined in \cref{eq:def_acc} to be of order $\mathcal{O}(\dt^7)$. The resulting conditions for achieving sixth-order accuracy with five sub-steps are
\begin{subequations}\label{eq:s5_r12345}
\begin{align}
\gamma_1&=\mathsf{RootOf}\left(720\gamma_1^5 - 1800\gamma_1^4 + 1200\gamma_1^3 - 300\gamma_1^2 + 30\gamma_1 - 1 \right) \label{eq:s5_r1} \\
\gamma_2&= \dfrac{(6\gamma_1-1)(30\gamma_1^{3}-45\gamma_1^{2}+15\gamma_1-1)}{3(120\gamma_1^4-240\gamma_1^3+120\gamma_1^{2}-20\gamma_1+1)}  \\
\gamma_3&= \dfrac{60\gamma_1^{3}\gamma_2-40\gamma_1^{3}-80\gamma_1^{2}\gamma_2+50\gamma_1^{2}+25\gamma_1\gamma_2-14\gamma_1-2\gamma_2+1}{5\gamma_1(2\gamma_1-1)(12\gamma_1\gamma_2-6\gamma_1-12\gamma_2+5)-5\gamma_2+2}  \\
\gamma_4&= \dfrac{5\gamma_1(\gamma_1-1)(6\gamma_2\gamma_3-4\gamma_2-4\gamma_3+3)+5\gamma_2\gamma_3-3\gamma_2-3\gamma_3+2}{10\gamma_1(\gamma_1-1)(6\gamma_2\gamma_3-3\gamma_2-3\gamma_3+2)+10\gamma_2\gamma_3-5\gamma_2-5\gamma_3+3} \\
\gamma_5&=1-\gamma_1. 
\end{align}
\end{subequations}

It is also observed in \cref{fig:s5_r1_rho} that a specific value of $\gamma_1\approx0.4732683912582953$ corresponds to the $A$-stable sixth-order member. At this point, the associated high-frequency dissipation parameter is fixed as $\left|\infrho\right|\approx0.83733$, indicating that the enhancement of accuracy by one additional order eliminates the independent control of numerical dissipation. Therefore, the five-sub-step implicit method provides two possible algorithm configurations: for arbitrary appropriate $\infrho\in[-1,~1]$, it achieves fifth-order accuracy with controllable high-frequency dissipation; whereas for the specific choice $\gamma_1=0.4732683912582953$, it achieves sixth-order accuracy with a fixed high-frequency dissipation level.
\begin{remark}
    It is worth mentioning that \cref{eq:s5_r1} also yields another admissible solution, $\gamma_1\approx1.6206645011137183$, which leads to an $A(\alpha)$-stable five-sub-step implicit method with sixth-order accuracy and $\alpha=89.94077997^\circ$. Nevertheless, this member suffers from a significantly larger relative period error and $A(\alpha)$-stability. Consequently, it is not pursued further in the present work.
\end{remark}

\subsection{Six sub-steps: $s=6$}
\begin{subequations}\label{eq:s_eq_6}
The six-sub-step implicit method is expressed by substituting $s=6$ into \cref{eq:but_alg} as 
\begin{equation}
\begin{NiceArray}{c|c}[cell-space-limits=3pt,columns-width=0.5cm]
\mbf{c} & \mbf{A}\\ \hline & \mbf{b}\T
\end{NiceArray}\quad=\quad\begin{NiceArray}{c|cccccc}[cell-space-limits=3pt,columns-width=0.7cm]
\gamma_1 & \alpha_{11} \\ 
\gamma_2 & \alpha_{21} & \alpha_{22} \\ 
\gamma_3 & \alpha_{31} & \alpha_{32} & \alpha_{33} \\ 
\gamma_4 & \alpha_{41} & \alpha_{42} & \alpha_{43} & \alpha_{44} \\ 
\gamma_5 & \alpha_{51} & \alpha_{52} & \alpha_{53} & \alpha_{54} &\alpha_{55} \\ 
\gamma_6 & \alpha_{61} & \alpha_{62} & \alpha_{63} & \alpha_{64} &\alpha_{65} & \alpha_{66} \\ \hline
& \beta_1 & \beta_2 & \beta_3 & \beta_4 & \beta_5 & \beta_6
\end{NiceArray}
\end{equation}
where the elements in $\mbf{A}$ and $\mbf{b}$ are computed as
\begin{align}
    \alpha_{11}&=\alpha_{22}=\alpha_{33}=\alpha_{44}=\alpha_{55}=\alpha_{66}=\gamma_1 \\
\alpha_{21}&=\gamma_2-\gamma_1\\ 
\alpha_{31}&=\gamma_3-\alpha_{32}-\gamma_1 \\
\alpha_{41}&=\gamma_4-\alpha_{42}-\alpha_{43}-\gamma_1\\
\alpha_{51}&=\gamma_5-\alpha_{52}-\alpha_{53}-\alpha_{54}-\gamma_1\\
\alpha_{61}&=\gamma_6-\alpha_{62}-\alpha_{63}-\alpha_{64}-\alpha_{65}-\gamma_1\\
\alpha_{32}&=\dfrac{-720\gamma_1^5+1800\gamma_1^4-1200\gamma_1^{3}+300\gamma_1^{2}-30\gamma_1+1}{720\beta_6\alpha_{65}\alpha_{54}\alpha_{43}(\gamma_2-\gamma_1)} \\ 
\alpha_{42}&=\dfrac{\left\{\begin{aligned}
&720\big[(\alpha_{65}\alpha_{53}+\alpha_{64}\alpha_{43})\beta_6+\alpha_{54}\alpha_{43}\beta_5\big]\gamma_1^5\\
&+360\big[(2\alpha_{65}\alpha_{54}\alpha_{43}-5\alpha_{64}\alpha_{43}-5\alpha_{65}\alpha_{53})\beta_6-5\alpha_{54}\alpha_{43}\beta_5\big]\gamma_1^4\\
&-240\big[(6\alpha_{65}\alpha_{54}\alpha_{43}-5\alpha_{64}\alpha_{43}-5\alpha_{65}\alpha_{53})\beta_6-5\alpha_{54}\alpha_{43}\beta_5\big]\gamma_1^3\\
&+60\big[(12\alpha_{65}\alpha_{54}\alpha_{43}-5\alpha_{64}\alpha_{43}-5\alpha_{65}\alpha_{53})\beta_6-5\alpha_{54}\alpha_{43}\beta_5\big]\gamma_1^2\\
&+30\big[24\alpha_{65}^{2}\alpha_{54}^{2}\alpha_{43}^{2}\beta_6^{2}- (4\alpha_{65}\alpha_{54}\alpha_{43}-\alpha_{64}\alpha_{43}-\alpha_{65}\alpha_{53})\beta_6+\alpha_{54}\alpha_{43}\beta_5 \big]\gamma_1\\
&-720\alpha_{65}^{2}\alpha_{54}^{2}\alpha_{43}^{2}\beta_6^{2}\gamma_3+(6\alpha_{65}\alpha_{54}\alpha_{43}-\alpha_{64}\alpha_{43}-\alpha_{65}\alpha_{53})\beta_6-\alpha_{54}\alpha_{43}\beta_5 
\end{aligned}\right\}}{720\beta_6^{2}\alpha_{65}^{2}\alpha_{54}^{2}\alpha_{43}(\gamma_2-\gamma_1)}\\
\alpha_{52}&=\dfrac{\left\{\begin{aligned}
&720\big[\alpha_{54}\alpha_{43}\beta_6(\alpha_{65}\beta_4-2\alpha_{64}\beta_5)-(\alpha_{64}^2\alpha_{43}+\alpha_{65}\alpha_{64}\alpha_{53}-\alpha_{65}\alpha_{63}\alpha_{54})\beta_6^2-\alpha_{54}^2\alpha_{43}\beta_5^2 \big]\gamma_1^5\\
&+360\big[5\alpha_{54}^2\alpha_{43}\beta_5^2-\alpha_{54}\alpha_{43}\beta_6(2\alpha_{65}\alpha_{54}\beta_5-10\alpha_{64}\beta_5+5\alpha_{65}\beta_4)\\
&\qquad-(\alpha_{64}\alpha_{43}(2\alpha_{65}\alpha_{54}-5\alpha_{64})-5\alpha_{65}(\alpha_{64}\alpha_{53}-\alpha_{63}\alpha_{54}))\beta_6^2 \big]\gamma_1^4\\
&+240\big[5\alpha_{65}\alpha_{54}\alpha_{43}\beta_4\beta_6+2\alpha_{54}\alpha_{43}\beta_5\beta_6(3\alpha_{65}\alpha_{54}-5\alpha_{64})-5\alpha_{54}^2\alpha_{43}\beta_5^2\\
&\qquad-(3\alpha_{65}\alpha_{54}\alpha_{43}(\alpha_{65}\alpha_{54}-2\alpha_{64})+5\alpha_{65}(\alpha_{64}\alpha_{53}-\alpha_{63}\alpha_{54})+5\alpha_{64}^2\alpha_{43})\beta_6^2 \big]\gamma_1^3\\
&+60\big[5\alpha_{54}^2\alpha_{43}\beta_5^2-\alpha_{54}\alpha_{43}\beta_6(12\alpha_{65}\alpha_{54}\beta_5-10\alpha_{64}\beta_5+5\alpha_{65}\beta_4)\\
&\qquad+(6\alpha_{65}\alpha_{54}\alpha_{43}(3\alpha_{65}\alpha_{54}-2\alpha_{64})+5\alpha_{65}(\alpha_{64}\alpha_{53}-\alpha_{63}\alpha_{54})+5\alpha_{64}^2\alpha_{43})\beta_6^2 \big]\gamma_1^2\\
&+30\big[24\alpha_{65}^3\alpha_{54}^2\alpha_{43}\beta_6^3(\alpha_{53}+\alpha_{54})+\alpha_{54}\alpha_{43}\beta_6(4\alpha_{65}\alpha_{54}\beta_5-2\alpha_{64}\beta_5+\alpha_{65}\beta_4)-\alpha_{54}^2\alpha_{43}\beta_5^2\\
&\qquad-(4\alpha_{65}\alpha_{54}\alpha_{43}(3\alpha_{65}\alpha_{54}-\alpha_{64})+\alpha_{65}(\alpha_{64}\alpha_{53}-\alpha_{63}\alpha_{54})+\alpha_{64}^2\alpha_{43})\beta_6^2 \big]\gamma_1\\
&-720\alpha_{65}^3\alpha_{54}^2\alpha_{43}\beta_6^3(\alpha_{53}\gamma_3+\alpha_{54}\gamma_4)-\alpha_{54}\alpha_{43}\beta_6(6\alpha_{65}\alpha_{54}\beta_5-2\alpha_{64}\beta_5+\alpha_{65}\beta_4)+\alpha_{54}^2\alpha_{43}\beta_5^2\\
&+(6\alpha_{65}\alpha_{54}\alpha_{43}(5\alpha_{65}\alpha_{54}-\alpha_{64})+\alpha_{65}(\alpha_{64}\alpha_{53}-\alpha_{63}\alpha_{54})+\alpha_{64}^2\alpha_{43})\beta_6^2
\end{aligned}\right\}}{720\beta_6^3\alpha_{65}^3\alpha_{54}^2\alpha_{43}(\gamma_2-\gamma_1)}\\
\alpha_{62}&=\dfrac{\left\{ \begin{aligned}
&-1440\beta_6\beta_5\beta_4\alpha_{65}\alpha_{54}\alpha_{43}\gamma_1^5-720\beta_6\alpha_{65}\alpha_{54}\alpha_{43}\big[\beta_6(\beta_4\alpha_{65}-\beta_5\alpha_{64})-\beta_5^{2}\alpha_{54}-5\beta_4\beta_5\big]\gamma_1^4\\
&+240\beta_6\alpha_{65}\alpha_{54}\alpha_{43}\big[3\beta_6\beta_5(\alpha_{65}\alpha_{54}-2\alpha_{64})-6\beta_5^{2}\alpha_{54}+6\beta_6\beta_4\alpha_{65}-10\beta_5\beta_4\big]\gamma_1^3\\
&-120\beta_6\alpha_{65}\alpha_{54}\alpha_{43}\big[3\beta_6\alpha_{65}\alpha_{54}(3\beta_5-2\beta_6\alpha_{65})-6\beta_5^{2}\alpha_{54}-6\beta_6\beta_5\alpha_{64}+6\beta_6\beta_4\alpha_{65}-5\beta_5\beta_4\big]\gamma_1^{2}\\
&-60\beta_6\alpha_{65}\alpha_{54}\alpha_{43}\big[6\beta_6\alpha_{65}\alpha_{54}(2\beta_6\alpha_{65}-\beta_5)-12\beta_6^3\alpha_{65}^{2}\alpha_{54}(\alpha_{63}+\alpha_{64}+\alpha_{65})\\
&\qquad+2\beta_5^{2}\alpha_{54}+2\beta_6(\beta_5\alpha_{64}-\beta_4\alpha_{65})+\beta_5\beta_4\big]\gamma_1\\
&+2\beta_6\alpha_{65}\alpha_{54}\alpha_{43}\big[60\beta_6^{2}\alpha_{65}^{2}\alpha_{54}-360\beta_6^{3}\alpha_{65}^{2}\alpha_{54}(\alpha_{63}\gamma_3+\alpha_{64}\gamma_4+\alpha_{65}\gamma_5)\\
&\qquad+3\beta_6(\beta_5\alpha_{64}-\beta_4\alpha_{65}-5\beta_5\alpha_{65}\alpha_{54})+\beta_5(3\beta_5\alpha_{54}+\beta_4)\big]\\
&+\big[\beta_5^{2}\alpha_{54}\alpha_{43}(\beta_5\alpha_{54}+2\beta_6\alpha_{64})+\beta_6^{2}(\beta_5\alpha_{64}-\beta_4\alpha_{65})(\alpha_{64}\alpha_{43}+\alpha_{65}\alpha_{53})\\
&\qquad-\beta_6^{2}\alpha_{65}\alpha_{54}(\beta_5\alpha_{63}-\beta_3\alpha_{65})\big]\times\big[720\gamma_1^5-1800\gamma_1^4+1200\gamma_1^{3}-300\gamma_1^{2}+30\gamma_1-1\big]
\end{aligned}\right\}}{720\beta_6^4\alpha_{65}^3\alpha_{54}^{2}\alpha_{43}(\gamma_2-\gamma_1)} \\ 
\alpha_{43}&=\dfrac{\left\{ \begin{aligned}
&(6\gamma_1^{2}-6\gamma_1+1)(1-2\gamma_2)-12\beta_5\big[\alpha_{53}(\gamma_1-\gamma_3)(\gamma_2-\gamma_3)+\alpha_{54}(\gamma_1-\gamma_4)(\gamma_2-\gamma_4)\big]\\
&-12\beta_6\big[\alpha_{63}(\gamma_1-\gamma_3)(\gamma_2-\gamma_3)+\alpha_{64}(\gamma_1-\gamma_4)(\gamma_2-\gamma_4)+\alpha_{65}(\gamma_1-\gamma_5)(\gamma_2-\gamma_5)\big]
\end{aligned}\right\}}{12\beta_4(\gamma_1-\gamma_3)(\gamma_2-\gamma_3)} \\
\alpha_{53}&=\dfrac{\left\{ \begin{aligned}
&\beta_4\big[60\gamma_1^3(1-2\gamma_2)+20\gamma_1^2(9\gamma_2-4)-5\gamma_1(12\gamma_2-5)+5\gamma_2-2\big]\\
&-120[(\beta_6\alpha_{64}+\beta_5\alpha_{54})^2-\beta_6\beta_4\alpha_{65}\alpha_{54}](\gamma_1-\gamma_4)(\gamma_2-\gamma_4)\\
&-120\beta_6\alpha_{63}(\beta_6\alpha_{64}+\beta_5\alpha_{54})(\gamma_1-\gamma_3)(\gamma_2-\gamma_3)-120\beta_6\alpha_{65}(\beta_6\alpha_{64}+\beta_5\alpha_{54})(\gamma_1-\gamma_5)(\gamma_2-\gamma_5)\\
&-10(6\gamma_1^2-6\gamma_1+1)(2\gamma_2-1)(\beta_6\alpha_{64}+\beta_5\alpha_{54})
\end{aligned}\right\}}{120(\beta_5^{2}\alpha_{54}+\beta_6\beta_5\alpha_{64}-\beta_6\beta_4\alpha_{65})(\gamma_1-\gamma_3)(\gamma_2-\gamma_3)} \\
\alpha_{54}&=\dfrac{\left\{\begin{aligned}
&60\beta_6\alpha_{64}(\gamma_1-\gamma_4)(\gamma_2-\gamma_4)(\gamma_3-\gamma_4)+60\beta_6\alpha_{65}(\gamma_1-\gamma_5)(\gamma_2-\gamma_5)(\gamma_3-\gamma_5)\\
&+10\gamma_1(6\gamma_2\gamma_3-3\gamma_2-3\gamma_3+2)(\gamma_1-1)+10\gamma_2\gamma_3-5\gamma_2-5\gamma_3+3
\end{aligned}\right\}}{60\beta_5(\gamma_4-\gamma_3)(\gamma_4-\gamma_2)(\gamma_4-\gamma_1)}\\
\alpha_{63}&=\dfrac{\left\{\begin{aligned}
&180(\beta_6\beta_5\alpha_{64}+\beta_5^2\alpha_{54}-\beta_6\beta_4\alpha_{65})(8\beta_1\gamma_2-2\beta_1\gamma_1-2\gamma_2+1)\gamma_1^4\\
&+60\big[(\beta_6\beta_5\alpha_{64}+\beta_5^2\alpha_{54}-\beta_6\beta_4\alpha_{65})(6\beta_1\gamma_1^2-24\beta_1\gamma_1\gamma_2+12\gamma_2-5)\\
&\qquad-3\beta_6\beta_5\alpha_{65}\alpha_{54}(2\beta_1\gamma_1-6\beta_1\gamma_2+2\gamma_2-1)\big]\gamma_1^3\\
&+15\big[4\beta_6\beta_5\alpha_{65}\alpha_{54}(6\beta_1\gamma_1^2-18\beta_1\gamma_1\gamma_2+9\gamma_2-4)-12\beta_6^2\alpha_{65}^2\alpha_{54}(2\beta_1\gamma_1-4\beta_1\gamma_2+2\gamma_2-1)\\
&\qquad-3\beta_5^2\alpha_{54}(8\gamma_2-3)-3\beta_6\beta_5\alpha_{64}(8\gamma_2-3)+3\beta_6\beta_4\alpha_{65}(8\gamma_2-3)\big]\gamma_1^2\\
&-3\big[120\beta_6^3\alpha_{65}^2\alpha_{64}\alpha_{54}(\gamma_2-\gamma_4)+120\beta_6^3\alpha_{65}^3\alpha_{54}(\gamma_2-\gamma_5)+5\beta_6\beta_5\alpha_{65}\alpha_{54}(12\gamma_2-5)\\
&\qquad-60\beta_6^2\alpha_{65}^2\alpha_{54}(2\beta_1\gamma_1^2-4\beta_1\gamma_1\gamma_2+2\gamma_2-1)-(\beta_6\beta_5\alpha_{64}+\beta_5^2\alpha_{54}-\beta_6\beta_4\alpha_{65})(20\gamma_2-7)\big]\gamma_1\\
&+360\beta_6^3\alpha_{65}^2\alpha_{64}\alpha_{54}\gamma_4(\gamma_2-\gamma_4)+360\beta_6^3\alpha_{65}^3\alpha_{54}\gamma_5(\gamma_2-\gamma_5)-30\beta_6^2\alpha_{65}^2\alpha_{54}(2\gamma_2-1)\\
&+3\beta_6\beta_5\alpha_{65}\alpha_{54}(5\gamma_2-2)-(\beta_6\beta_5\alpha_{64}+\beta_5^2\alpha_{54}-\beta_6\beta_4\alpha_{65})(3\gamma_2-1)
\end{aligned}\right\}}{360\beta_6^3\alpha_{65}^2\alpha_{54}(\gamma_1-\gamma_3)(\gamma_2-\gamma_3)}\\
\alpha_{64}&=\dfrac{\left\{\begin{aligned}
&120\beta_6^2\alpha_{65}^2(\gamma_3-\gamma_5)(\gamma_2-\gamma_5)(\gamma_1-\gamma_5)-5\beta_5\gamma_3(5\gamma_1+\gamma_2)+\beta_5(14\gamma_1+2\gamma_2+2\gamma_3-1)\\
&-2\beta_6\alpha_{65}\big[10\gamma_1(6\gamma_2\gamma_3-3\gamma_2-3\gamma_3+2)(1-\gamma_1)-10\gamma_2\gamma_3+5\gamma_2+5\gamma_3-3\big]\\
&+120\beta_5\beta_4\gamma_1^3(\gamma_3-\gamma_4)(\gamma_2-\gamma_4)+120\beta_5^2\gamma_1^3(\gamma_3-\gamma_5)(\gamma_2-\gamma_5)+120\beta_6\beta_5\gamma_1^3(\gamma_3-\gamma_6)(\gamma_2-\gamma_6)\\
&+120\beta_5\beta_1\gamma_1^3(\gamma_3-\gamma_1)(\gamma_2-\gamma_1)-20\beta_5\gamma_1^2\gamma_2(9\gamma_3-4)+20\beta_5\gamma_1\gamma_3(4\gamma_1+3\gamma_2)-25\beta_5\gamma_1(2\gamma_1+\gamma_2)
\end{aligned}\right\}}{120\beta_6^2\alpha_{65}(\gamma_4-\gamma_3)(\gamma_4-\gamma_2)(\gamma_4-\gamma_1)}\\
\alpha_{65}&=\dfrac{\left\{\begin{aligned}
&5\gamma_2\gamma_1(6\gamma_3\gamma_4-2\gamma_3-2\gamma_4+1)-(10\gamma_3\gamma_4-5\gamma_3-5\gamma_4+3)(\gamma_1+\gamma_2)+5\gamma_3\gamma_4-3\gamma_3-3\gamma_4+2\\
&-60\beta_6\gamma_1(\gamma_1-\gamma_6)(\gamma_2-\gamma_6)(\gamma_3-\gamma_6)(\gamma_4-\gamma_6)-60\beta_5\gamma_1(\gamma_1-\gamma_5)(\gamma_2-\gamma_5)(\gamma_3-\gamma_5)(\gamma_4-\gamma_5)
\end{aligned}\right\}}{60\beta_6(\gamma_5-\gamma_4)(\gamma_5-\gamma_3)(\gamma_5-\gamma_2)(\gamma_5-\gamma_1)}\\
\beta_1&=\dfrac{\left\{ \begin{aligned}
&60\gamma_2\gamma_3\gamma_4\gamma_5\gamma_6-30(\gamma_2\gamma_3\gamma_4\gamma_5+\gamma_2\gamma_3\gamma_4\gamma_6+\gamma_2\gamma_3\gamma_5\gamma_6+\gamma_2\gamma_4\gamma_5\gamma_6+\gamma_3\gamma_4\gamma_5\gamma_6)\\
&+20(\gamma_2\gamma_3\gamma_4+\gamma_2\gamma_3\gamma_5+\gamma_2\gamma_3\gamma_6+\gamma_2\gamma_4\gamma_5+\gamma_2\gamma_4\gamma_6+\gamma_2\gamma_5\gamma_6+\gamma_3\gamma_4\gamma_5+\gamma_3\gamma_4\gamma_6+\gamma_3\gamma_5\gamma_6+\gamma_4\gamma_5\gamma_6)\\
&-15(\gamma_2\gamma_3+\gamma_2\gamma_4+\gamma_2\gamma_5+\gamma_2\gamma_6+\gamma_3\gamma_4+\gamma_3\gamma_5+\gamma_3\gamma_6+\gamma_4\gamma_5+\gamma_4\gamma_6+\gamma_5\gamma_6)\\
&+12(\gamma_2+\gamma_3+\gamma_4+\gamma_5+\gamma_6)-10
\end{aligned} \right\}}{60(\gamma_2-\gamma_1)(\gamma_3-\gamma_1)(\gamma_4-\gamma_1)(\gamma_5-\gamma_1)(\gamma_6-\gamma_1)}  \\ 
\beta_2&=\dfrac{\left\{ \begin{aligned}
&60\gamma_1\gamma_3\gamma_4\gamma_5\gamma_6-30(\gamma_1\gamma_3\gamma_4\gamma_5+\gamma_1\gamma_3\gamma_4\gamma_6+\gamma_1\gamma_3\gamma_5\gamma_6+\gamma_1\gamma_4\gamma_5\gamma_6+\gamma_3\gamma_4\gamma_5\gamma_6)\\
&+20(\gamma_1\gamma_3\gamma_4+\gamma_1\gamma_3\gamma_5+\gamma_1\gamma_3\gamma_6+\gamma_1\gamma_4\gamma_5+\gamma_1\gamma_4\gamma_6+\gamma_1\gamma_5\gamma_6+\gamma_3\gamma_4\gamma_5+\gamma_3\gamma_4\gamma_6+\gamma_3\gamma_5\gamma_6+\gamma_4\gamma_5\gamma_6)\\
&-15(\gamma_1\gamma_3+\gamma_1\gamma_4+\gamma_1\gamma_5+\gamma_1\gamma_6+\gamma_3\gamma_4+\gamma_3\gamma_5+\gamma_3\gamma_6+\gamma_4\gamma_5+\gamma_4\gamma_6+\gamma_5\gamma_6)\\
&+12(\gamma_1+\gamma_3+\gamma_4+\gamma_5+\gamma_6)-10
\end{aligned} \right\}}{60(\gamma_1-\gamma_2)(\gamma_3-\gamma_2)(\gamma_4-\gamma_2)(\gamma_5-\gamma_2)(\gamma_6-\gamma_2)}  \\   
\beta_3&=\dfrac{\left\{ \begin{aligned}
&60\gamma_1\gamma_2\gamma_4\gamma_5\gamma_6-30(\gamma_1\gamma_2\gamma_4\gamma_5+\gamma_1\gamma_2\gamma_4\gamma_6+\gamma_1\gamma_2\gamma_5\gamma_6+\gamma_1\gamma_4\gamma_5\gamma_6+\gamma_2\gamma_4\gamma_5\gamma_6)\\
&+20(\gamma_1\gamma_2\gamma_4+\gamma_1\gamma_2\gamma_5+\gamma_1\gamma_2\gamma_6+\gamma_1\gamma_4\gamma_5+\gamma_1\gamma_4\gamma_6+\gamma_1\gamma_5\gamma_6+\gamma_2\gamma_4\gamma_5+\gamma_2\gamma_4\gamma_6+\gamma_2\gamma_5\gamma_6+\gamma_4\gamma_5\gamma_6)\\
&-15(\gamma_1\gamma_2+\gamma_1\gamma_4+\gamma_1\gamma_5+\gamma_1\gamma_6+\gamma_2\gamma_4+\gamma_2\gamma_5+\gamma_2\gamma_6+\gamma_4\gamma_5+\gamma_4\gamma_6+\gamma_5\gamma_6)\\
&+12(\gamma_1+\gamma_2+\gamma_4+\gamma_5+\gamma_6)-10
\end{aligned} \right\}}{60(\gamma_1-\gamma_3)(\gamma_2-\gamma_3)(\gamma_4-\gamma_3)(\gamma_5-\gamma_3)(\gamma_6-\gamma_3)}  \\  
\beta_4&=\dfrac{\left\{ \begin{aligned}
&60\gamma_1\gamma_2\gamma_3\gamma_5\gamma_6-30(\gamma_1\gamma_2\gamma_3\gamma_5+\gamma_1\gamma_2\gamma_3\gamma_6+\gamma_1\gamma_2\gamma_5\gamma_6+\gamma_1\gamma_3\gamma_5\gamma_6+\gamma_2\gamma_3\gamma_5\gamma_6)\\
&+20(\gamma_1\gamma_2\gamma_3+\gamma_1\gamma_2\gamma_5+\gamma_1\gamma_2\gamma_6+\gamma_1\gamma_3\gamma_5+\gamma_1\gamma_3\gamma_6+\gamma_1\gamma_5\gamma_6+\gamma_2\gamma_3\gamma_5+\gamma_2\gamma_3\gamma_6+\gamma_2\gamma_5\gamma_6+\gamma_3\gamma_5\gamma_6)\\
&-15(\gamma_1\gamma_2+\gamma_1\gamma_3+\gamma_1\gamma_5+\gamma_1\gamma_6+\gamma_2\gamma_3+\gamma_2\gamma_5+\gamma_2\gamma_6+\gamma_3\gamma_5+\gamma_3\gamma_6+\gamma_5\gamma_6)\\
&+12(\gamma_1+\gamma_2+\gamma_3+\gamma_5+\gamma_6)-10
\end{aligned} \right\}}{60(\gamma_1-\gamma_4)(\gamma_2-\gamma_4)(\gamma_3-\gamma_4)(\gamma_5-\gamma_4)(\gamma_6-\gamma_4)}  \\ 
\beta_5&=\dfrac{\left\{ \begin{aligned}
&60\gamma_1\gamma_2\gamma_3\gamma_4\gamma_6-30(\gamma_1\gamma_2\gamma_3\gamma_4+\gamma_1\gamma_2\gamma_3\gamma_6+\gamma_1\gamma_2\gamma_4\gamma_6+\gamma_1\gamma_3\gamma_4\gamma_6+\gamma_2\gamma_3\gamma_4\gamma_6)\\
&+20(\gamma_1\gamma_2\gamma_3+\gamma_1\gamma_2\gamma_4+\gamma_1\gamma_2\gamma_6+\gamma_1\gamma_3\gamma_4+\gamma_1\gamma_3\gamma_6+\gamma_1\gamma_4\gamma_6+\gamma_2\gamma_3\gamma_4+\gamma_2\gamma_3\gamma_6+\gamma_2\gamma_4\gamma_6+\gamma_3\gamma_4\gamma_6)\\
&-15(\gamma_1\gamma_2+\gamma_1\gamma_3+\gamma_1\gamma_4+\gamma_1\gamma_6+\gamma_2\gamma_3+\gamma_2\gamma_4+\gamma_2\gamma_6+\gamma_3\gamma_4+\gamma_3\gamma_6+\gamma_4\gamma_6)\\
&+12(\gamma_1+\gamma_2+\gamma_3+\gamma_4+\gamma_6)-10
\end{aligned} \right\}}{60(\gamma_1-\gamma_5)(\gamma_2-\gamma_5)(\gamma_3-\gamma_5)(\gamma_4-\gamma_5)(\gamma_6-\gamma_5)}  \\ 
\beta_6&=\dfrac{\left\{ \begin{aligned}
&60\gamma_1\gamma_2\gamma_3\gamma_4\gamma_5-30(\gamma_1\gamma_2\gamma_3\gamma_4+\gamma_1\gamma_2\gamma_3\gamma_5+\gamma_1\gamma_2\gamma_4\gamma_5+\gamma_1\gamma_3\gamma_4\gamma_5+\gamma_2\gamma_3\gamma_4\gamma_5)\\
&+20(\gamma_1\gamma_2\gamma_3+\gamma_1\gamma_2\gamma_4+\gamma_1\gamma_2\gamma_5+\gamma_1\gamma_3\gamma_4+\gamma_1\gamma_3\gamma_5+\gamma_1\gamma_4\gamma_5+\gamma_2\gamma_3\gamma_4+\gamma_2\gamma_3\gamma_5+\gamma_2\gamma_4\gamma_5+\gamma_3\gamma_4\gamma_5)\\
&-15(\gamma_1\gamma_2+\gamma_1\gamma_3+\gamma_1\gamma_4+\gamma_1\gamma_5+\gamma_2\gamma_3+\gamma_2\gamma_4+\gamma_2\gamma_5+\gamma_3\gamma_4+\gamma_3\gamma_5+\gamma_4\gamma_5)\\
&+12(\gamma_1+\gamma_2+\gamma_3+\gamma_4+\gamma_5)-10
\end{aligned} \right\} }{60(\gamma_1-\gamma_6)(\gamma_2-\gamma_6)(\gamma_3-\gamma_6)(\gamma_4-\gamma_6)(\gamma_5-\gamma_6)}
\end{align}
\end{subequations}
and the sub-step sizes $\gamma_j~(j=1,~\cdots,~6)$ are free parameters. 

\begin{figure}[htbp]
	\centering 
	\includegraphics[scale=1.0]{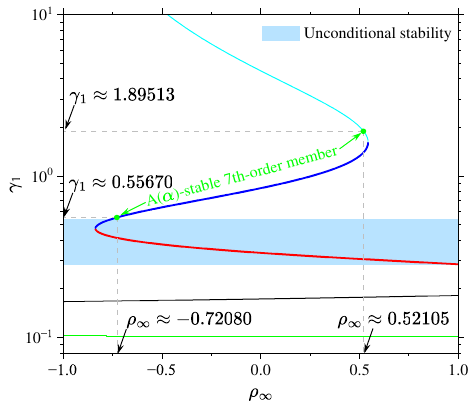}
	\caption{The variations of $\gamma_1$ with $\infrho\in[-1,~1]$ for the six-sub-step implicit method \eqref{eq:s_eq_6}.}
	\label{fig:s6_r1_rho}
\end{figure}

For sixth-order accuracy, the first sub-step size $\gamma_1$ is related to the user-specified parameter $\infrho$ through
\begin{equation}\label{eq:s6_D_infty}
\dfrac{720\gamma_1^{6}-4320\gamma_1^5+5400\gamma_1^4-2400\gamma_1^3+450\gamma_1^2-36\gamma_1+1}{720\gamma_1^{6}} = \left| \infrho \right|
\end{equation}
whereas the remaining sub-step sizes $\gamma_j~(j=2,~\cdots,~6)$ remain free parameters. For example, Li et al. \cite{li_DirectlySelfstartingHigherorder_2022} used $\gamma_j=j\cdot\gamma_1~(j=2,~\cdots,~5)$ and $\gamma_6=1$ to splite the equal sub-step sizes for the first five sub-steps. The unconditional stability given in \cref{eq:Ey} requires $\gamma_1$ to satisfy
\begin{equation}\label{eq:s6_uc}
\gamma_1\in\left[0.2840646380117983,~0.5409068780733081\right],
\end{equation}
where $0.2840646380117983$ is an approximate root of $1440\gamma_1^6 - 4320\gamma_1^5 + 5400\gamma_1^4 - 2400\gamma_1^3 + 450\gamma_1^2 - 36\gamma_1 + 1=0$ and $0.5409068780733081$ is an appropriate root of $34560\gamma_1^7 - 106560\gamma_1^6 + 97920\gamma_1^5 - 39600\gamma_1^4 + 8160\gamma_1^3 - 888\gamma_1^2 + 48\gamma_1 - 1=0$. Accordingly, the first sub-step size $\gamma_1$ is determined by
\begin{equation}\label{eq:s6_r1_set}
\gamma_1\in\left\{\gamma_1~\bigg|~ \text{Eq.}~\eqref{eq:s6_uc},~ 720(1 + \infrho)\gamma_1^6 - 4320\gamma_1^5 + 5400\gamma_1^4 - 2400\gamma_1^3 + 450\gamma_1^2 - 36\gamma_1+1=0\right\}
\end{equation}
where $\infrho\in[-0.83733,~1]$. \cref{eq:s6_r1_set} establishes the relationship between the first sub-step size $\gamma_1$ and the user-specified parameter $\infrho$. For a prescribed value of $\infrho$, the sextic equation $720(1 + \infrho)\gamma_1^6 - 4320\gamma_1^5 + 5400\gamma_1^4 - 2400\gamma_1^3 + 450\gamma_1^2 - 36\gamma_1+1=0$ yields six candidate solutions for $\gamma_1$. Since not all mathematical solutions necessarily satisfy the unconditional stability, the admissible solution is selected by additionally enforcing \cref{eq:s6_uc}. In this manner, a unique relationship between $\gamma_1$ and the user-specified parameter $\infrho$ can be established, as shown in \cref{fig:s6_r1_rho}.

\cref{fig:s6_r1_rho} illustrates the variation of the first sub-step size $\gamma_1$ with respect to the user-specified parameter $\infrho$ for the six-sub-step implicit method defined by \cref{eq:s_eq_6}. For the sixth-order formulation, the admissible values of $\gamma_1$ are obtained by combining the high-frequency dissipation with unconditional stability. The stable solution branch passes through the shaded unconditionally stable region, thereby allowing the sixth-order six-sub-step method to retain controllable high-frequency dissipation over the admissible range of $\infrho\in[0,~1]$.


For seventh-order accuracy, the sub-step sizes $\gamma_j~(j=1,~\cdots,~6)$ are determined by requiring the local truncation errors defined in \cref{eq:def_acc} to be of order $\mathcal{O}(\dt^8)$. The resulting conditions for achieving seventh-order accuracy with six sub-steps are
\begin{subequations}\label{eq:s6_r123456}
\begin{align}
\gamma_1&=\mathsf{RootOf}\left( 5040\gamma_1^6 - 15120\gamma_1^5 + 12600\gamma_1^4 - 4200\gamma_1^3 + 630\gamma_1^2 - 42\gamma_1 + 1 \right)\label{eq:s6_r1}  \\
\gamma_2&= \dfrac{2520\gamma_1^4(\gamma_1-2)+84\gamma_1^{2}(35\gamma_1-8)+63\gamma_1-2}{7\big[360\gamma_1^4(2\gamma_1-5)+300\gamma_1^{2}(4\gamma_1-1)+30\gamma_1-1\big]}  \\
\gamma_3&= \dfrac{420\gamma_1^4(3\gamma_2-2)-420\gamma_1^{3}(5\gamma_2-3)+63\gamma_1^{2}(15\gamma_2-8)-7\gamma_1(21\gamma_2-10)+7\gamma_2-3}{7\big[180\gamma_1^4(2\gamma_2-1)-60\gamma_1^{3}(12\gamma_2-5)+45\gamma_1^{2}(8\gamma_2-3)-3\gamma_1(20\gamma_2-7)+3\gamma_2-1\big]} 
  \\
\gamma_4&= \dfrac{35\gamma_1\gamma_2\gamma_3(2\gamma_1-1)(6\gamma_1-5)-7(4\gamma_1-1)(\gamma_2+\gamma_3)(10\gamma_1^{2}-10\gamma_1+1)+42\gamma_1^{2}(5\gamma_1-6)-14\gamma_2\gamma_3+63\gamma_1-4}{7\big[60\gamma_1\gamma_2\gamma_3(2\gamma_1-1)(\gamma_1-1)-(\gamma_2+\gamma_3)(60\gamma_1^{3}-80\gamma_1^{2}+25\gamma_1-2)+10\gamma_1^{2}(4\gamma_1-5)-5\gamma_2\gamma_3+14\gamma_1-1\big]}  \\
\gamma_5&= \dfrac{7\gamma_1(\gamma_1-1)\big[10\gamma_2\gamma_3(3\gamma_4-2)+15\gamma_4-5(\gamma_2+\gamma_3)(4\gamma_4-3)-12\big]-7(\gamma_2+\gamma_3)(3\gamma_4-2)+7\gamma_2\gamma_3(5\gamma_4-3)+14\gamma_4-10}{7\big[5\gamma_1(\gamma_1-1)\big(6\gamma_2\gamma_3(2\gamma_4-1)+4\gamma_4-2(\gamma_2+\gamma_3)(3\gamma_4-2)-3\big)-(\gamma_2+\gamma_3)(5\gamma_4-3)+5\gamma_2\gamma_3(2\gamma_4-1)+3\gamma_4-2\big]}  \\
\gamma_6&=1-\gamma_1
\end{align}
\end{subequations}
where $\mathsf{RootOf}(f)$ denotes the set of all roots of $f=0$. \cref{eq:s6_r1} yields six candidate values of $\gamma_1$, none of which satisfies the unconditional stability given in \cref{eq:s6_uc}. Consequently, a strictly unconditionally stable seventh-order six-sub-step implicit method does not exist within the present framework \eqref{eq:but_alg}.  Nevertheless, numerical investigations indicate that the stability region is maximized and the relative period error is minimized when
\begin{equation}
\gamma_1=0.5566999420873480.
\end{equation}
Therefore, it is adopted in the present study as a practical choice for constructing the seventh-order six-sub-step implicit method. With this choice, the resulting seventh-order method is $A(\alpha)$-stable, where $\alpha$ is numerically estimated as $89.99999995^\circ$. This value is extremely close to the limiting case $\alpha=90^\circ$, corresponding to unconditional stability. 
\begin{remark}
    It is worth mentioning that \cref{eq:s6_r1} also yields another admissible solution, $\gamma_1\approx1.8951305922313202$, which leads to an $A(\alpha)$-stable six-sub-step implicit method with seventh-order accuracy and $\alpha=89.85071945^\circ$. Nevertheless, this member suffers from a significantly larger relative period error and less favorable stability properties. Consequently, it is not pursued further in the present work.
\end{remark}

A notable feature of \cref{fig:s6_r1_rho} is the existence of two distinct seventh-order candidates obtained by imposing one additional order condition. The first candidate corresponds to $\gamma_1\approx0.55670$ and $\infrho\approx-0.72080$. This point lies outside the strict unconditionally stable region but remains very close to its boundary. The resulting seventh-order method is therefore $A(\alpha)$-stable rather than strictly $A$-stable. Since its stability angle is close to $90^\circ$, this member retains stability characteristics that are close to unconditional stability and is adopted as the preferred seventh-order configuration. \cref{eq:s6_r1} also yields another candidate, $\gamma_1\approx1.89513$, associated with $\infrho\approx0.52105$. This solution is located substantially outside the unconditionally stable interval and consequently exhibits less favorable stability. As discussed previously, it also produces a larger relative period error than the member with $\gamma_1\approx0.55670$. Therefore, although both parameter sets satisfy the seventh-order conditions, only the former is emphasized in the present study.

\end{appendices}


\bibliography{Reference}

\end{document}